**The Edge Geometry of Regular N-gons (Part I for N ≤ 25)**

**G. H. Hughes**

In [H5] (*First Families of Regular Polygons and their Mutations*) we defined the First Family members of a regular N-gon to be 'conforming' to the star polygons of N and hence they survive in the complement of the singularity set W under the outer-billiards map $\tau$. Since all the dynamics of $\tau$ are based on this complement, these tiles play an important role in determining the topology of W. Because this topology is very complex we suggested that it might be fruitful to begin with a study the 'edge-geometry' of these regular polygons to yield some insight into the global topology.

Every regular N-gon has a locally invariant region which will include at least the S[1] and S[2] 'tiles 'in the First Family of N. Typically these regions will be bounded by the 'major' resonances of N so they would be expected to include least 1/3 or 1/4 of all the S[k] in the First Family of N. The tiles bound together in an invariant region should share a similar geometry and dynamics. For N = 60 shown below the inner invariant region extends out to S[25].

In [H2] (*Outer-billiards, digital filters and kicked Hamiltonians*) we describe two additional mappings that can be used to generate portions of this edge geometry of $\tau$. These mappings are the digital-filter map *Df* of Chou and Lin [Cl] and the dual-center map *Dc* of E. Goetz [Go]. In the Appendix of [H5] we show how these two mappings can be used to trace the evolution of W.

In [H5] we show that the evolution of the $\tau$-web can be reduced to a simple 'shear and rotation' and in [H2] we show how the *Df* map can also be reduced to a shear and rotation – but on a toral space. The *Dc* map is *defined* to be a shear and rotation in the complex plane Therefore it is not surprising that these two mappings have singularity sets that are locally conjugate to $\tau$. This supports our premise that this geometry is inherent in the N-gon. The Df and Dc maps offer a significant computational advantage, but their dynamics are quite different

The $\tau$- singularity set W (a.k.a. the 'web') can be generated by iterating the extended edges of the N-gon as shown below. If these extended edges are truncated they form classical star polygons and the cyan First Family S[k] tiles are defined to be conforming to these nested star polygons. For N = 14 the maximal S[k] is S[5] – also known as D. Here it is congruent to N but for N-odd it is a 2N-gon with edge length identical to N. This D tile is always globally maximum among convex tiles that can evolve in W (with any extended edge length), and rings of these D tiles guarantee that this resulting 'generalized star polygon' is invariant, so it is sufficient to iterate just the star polygon edges to define what we call the web W. If full extended edges are iterated, concentric rings of these D tiles will guarantee global stability (and introduce no additional scaling or geometry). See Appendix II of [H5].

**Figure 1** The web development for the regular tetradecagon known as N = 14.

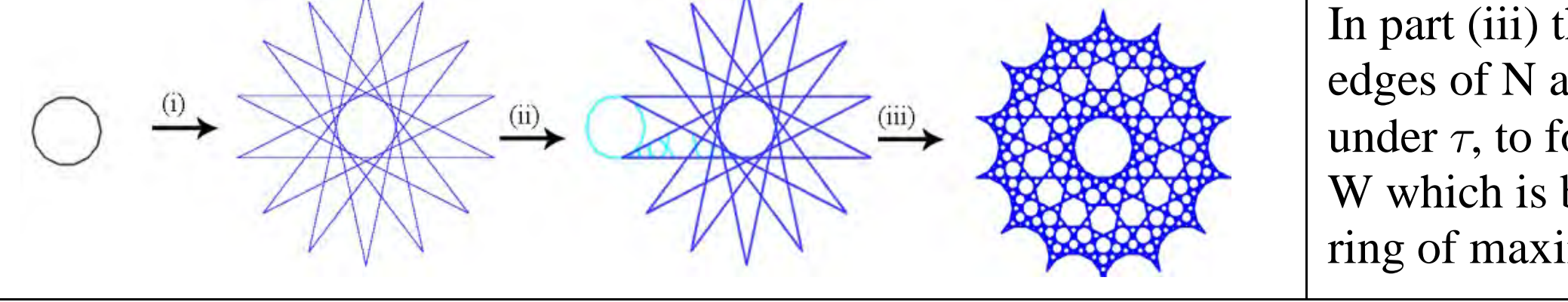

In part (iii) the truncated edges of N are iterated under $\tau$, to form the 'web' W which is bounded by a ring of maximal 'D' tiles.

These ‘generalized star polygons’ in (iii) share the same scaling and dihedral symmetry as N so they are inherent in N and their geometry is determined by the matching cyclotomic field $\mathbb{Q}_N$. Each S[k] in the First Family of N defines a star[k] point and a scale[k]. The number of ‘primitive’ scales (gcd(k,N) = 1) matches the ‘algebraic complexity’ of N, namely $\phi(N)/2$ where $\phi$ is the Euler totient function. This is the rank of the maximal real subfield $\mathbb{Q}_N^+$ of $\mathbb{Q}_N$.

Based on a 1949 result of C.L.Siegel communicated to S. Chowla [Ch], the primitive scales form a unit basis for $\mathbb{Q}_N^+$, so this is what we call the ‘scaling field’ of N. The traditional generator of $\mathbb{Q}_N^+$ is $\lambda_N = 2\cos(2\pi/N)$ which is $\zeta + \zeta^{-1}$ where $\zeta = \exp(2\pi i/N)$. We will typically use primitive scales as alternate generators of $\mathbb{Q}_N^+$ because the resulting scaling polynomials will be more meaningful than the generic polynomials in $\lambda_N$.

Our ultimate goal is to understand the topology of W, but currently the only cases where the topology of W is understood are the ‘linear’ or ‘quadratic’ cases of N = 3,4,5,6,8,10 and 12. The N = 14 case described here has $\phi(N)/2 = 3$ so N = 14 and the matching N = 7 are classified as ‘cubic’ and there are two non-trivial primitive scales so W is probably multi-fractal.

As noted by J. Moser in 1978, this $\tau$-web can be regarded as a discontinuous version of the phase space for a Hamiltonian system, so all three maps will have geometry related to the classical 1969 Standard Map of Chirikov. This connection is outlined in [A3] of [H5]

At this time the $\tau$-representations are the most meaningful since there is a well-developed theory of dynamics for both regular and non-regular N-gons. But for computational purposes the Dc maps are more efficient and we will sometimes use these alternative maps to generate W.

Our primary concern here is the geometry local to N and this will always include the S[1] and S[2] tiles of N, but this geometry is typically shared by adjacent tiles in an invariant region local to N. In this way the edge geometry interacts with the overall geometry of N. In Appendix II of [H5] we describe the global evolution of W based on concentric rings of D tiles and show that the initial Ring 0 guarantees that the canonical ‘inner star’ region is invariant. But invariance exists at all scales and typically this inner region is sub-divided into secondary regions as shown below for N = 60. We conjecture that the boundaries are driven by edge-sharing of adjacent tiles.

When N is twice-odd this edge-sharing occurs early and is easy to predict but for N odd or twice-even, neighboring S[k] cannot share edges so these edge sharing ‘hugs’ are harder to predict. For N = 60 below the black invariant region extends all the way to S[26] which by the GFFT has a step-4 right-side family that includes neighboring S[13] and S[17] tiles which share an edge. This shared edge in the early web creates a dynamical barrier between magenta and black points.

**Figure 2** The invariant regions for N = 60 where edge sharing of tiles defines boundaries.

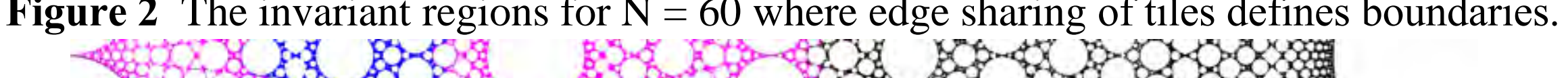

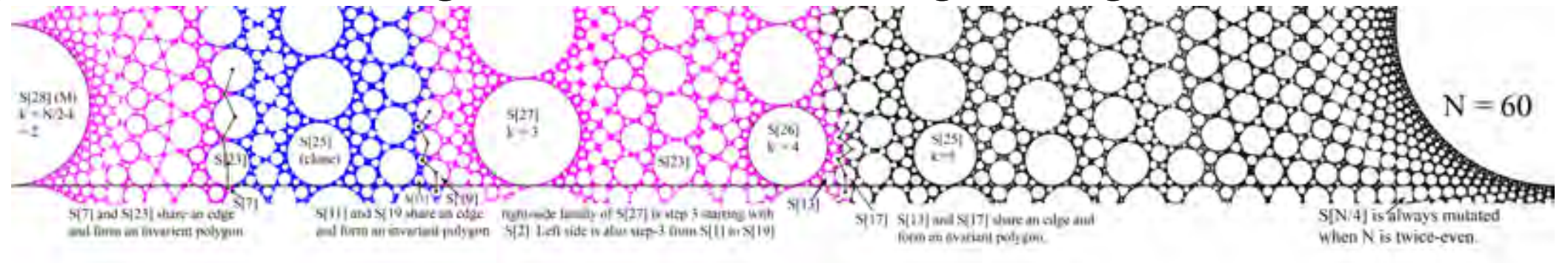

The concentric rings of D tiles described in Appendix II of [H5] are linked by shared vertices but they share S[1] tiles with their neighbors so they are also linked by shared edges and it appears that this edge sharing of D or N with S[1] can be regarded as a 'template' for all known instances of edge sharing. Below we will define two tiles to be 'dual' if they have the form S[k] and S[k′] so one is the 'derivative' of the other. We will see that dual tiles have a close algebraic relationship which makes them capable of edge sharing. For N even, every S[k] is self-dual since $k' = N/2-k$ has $k'' = k$. Since the D tiles always evolve as step-1, they are dual to the S[1]. Therefore they have a stronger link than just the edge-sharing of the secondary resonances.

When N is twice odd the matching 'M' tiles can link up in this stronger fashion as shown on the left below with N = 14. Here the S[5] of N plays the role of a local N = 7 and the tiles in the chain share vertices as well as an edge with their S[1] tile (which is an S[2] of N). This 'peck' and 'hug' is a strong version of invariance similar to the rings of D tiles. The first isolated edge sharing here is the S[Floor[N/4]] = S[3] with neighbor S[4] just like N = 26 with S[6] and S[7].

**Figure 3** Invariant regions for N = 14 and N = 26

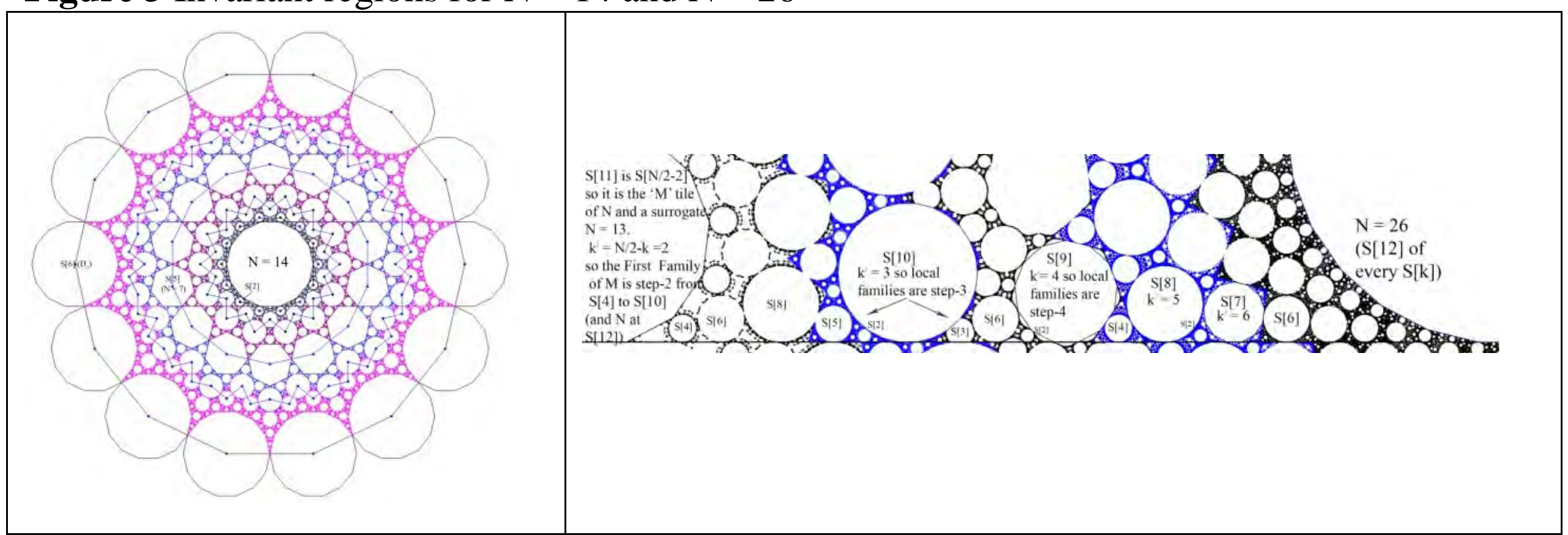

To return to the relatively cold and hug-free world of N = 60, it is up to N = 26 with $k' = N/2-k = 4$ to foster a step-4 family that must be compatible with S[25] which evolves earlier. So the right-side family of S[26] must include S[25], and counting down mod-4 it will be begin with a clone of S[5] as show here. But by serendipity this family includes an S[13] which is a perfect match with S[17] which has $k' = 13$. (The original S[17] also has $k' = 13$ but there are no S[13] neighbors to form a bond.) On careful inspection the magenta edge of S[13] will be a symmetric fit inside the edge of S[17] since $sS[13]/sS[17] = \mathrm{Tan}[13\pi/60]^2 \approx .6557$.

**Figure 4** By the GFFT, the S[26] tile must have a step-4 right-side family since $k' = 4$

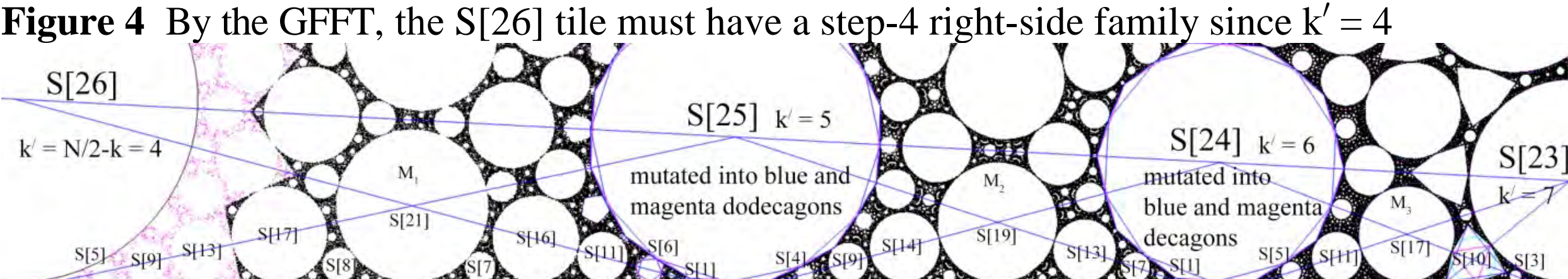

Based on these few examples it seems unlikely that we will find a simple heuristic that can be used to predict the bounds of these invariant regions, but below we present some suggestions.

**Recommendations for determining the primary $\tau$-invariant region of a regular n-gon**

**Definition**; The last surviving S[k] tile of N in this invariant region we call the 'shepherd' $S_N$

(i) When N is twice-odd $S_N$ = S[Floor[N/4]] and it will be the S[N/2-k] tile in the family of S[Floor[N/4]+1] and hence will have an edge which is a subset of an edge of S[Floor[N/4]+1]. (See Appendix II for details of N = 202 along with samples of orbital 'Projections'.)

(ii) When N is odd or twice-even the S[k] in the First Family of N are always even so the geometry is different from the 'egalitarian' twice-odd case. Here there are many cases where S[N/3] is a good approximation for $S_N$ but there are also cases where N/3 is far off as witnessed by N = 60 with $S_N$ = 25 instead of 20 or N = 243 with $S_N$ = 96 instead of 81. See N = 300 below.

**Figure 5**. For large k values the S[k] can only interact through their local step-k families. As noted above the twice-odd family is the only one where the early S[k] can share an edge with their neighbor, so the typical case shown below will involve sharing within these local step-k families. The GFFT can be used to identify the primary tiles in these families but they are just clones of the actual S[k] in the early web, so their location will depend on the local star points of the 'parent'. Often the local 'M' tiles act as level-2 parents. They are the only tiles in both right and left-side families of neighboring S[k], so their role is pivotal.

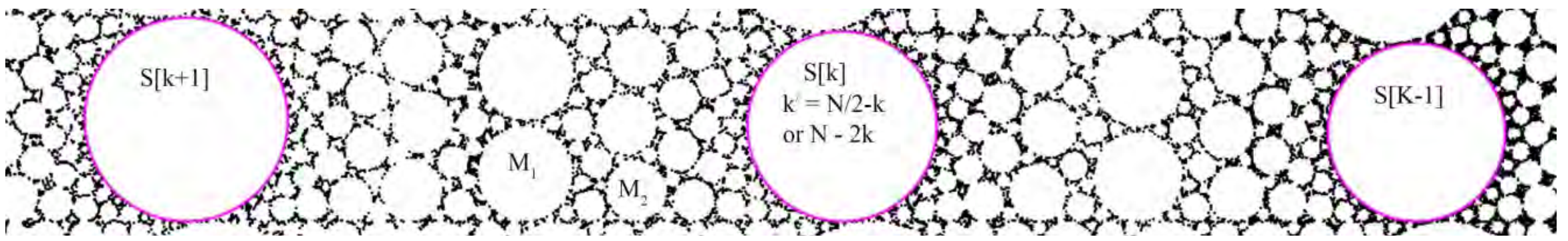

At the lowest level we would know whether any of these S[k] in the First Family of N have a local step-k family that begins with an S[k′]. If this is true then S[k′] will share an edge with S[k] and create a potential barrier for the dynamics. For N = 26 in Figure 3 this occurs with S[6] which has k′ = 17-12 = 5 and it shares an edge with a DS[5] to create the initial boundary.

N = 36 below has an example of a level-2 transition where S[15] fosters a step-3 family which includes an S[11] with k′ = 7 . Meanwhile S[14] has a step-4 left side web which counts down from 15 and this includes an S[7] which lives in both black and magenta worlds.

**Figure 6**. N = 36 has its first transition at S[15] but it is the local 'M' at S[11] that supports the edge sharing.

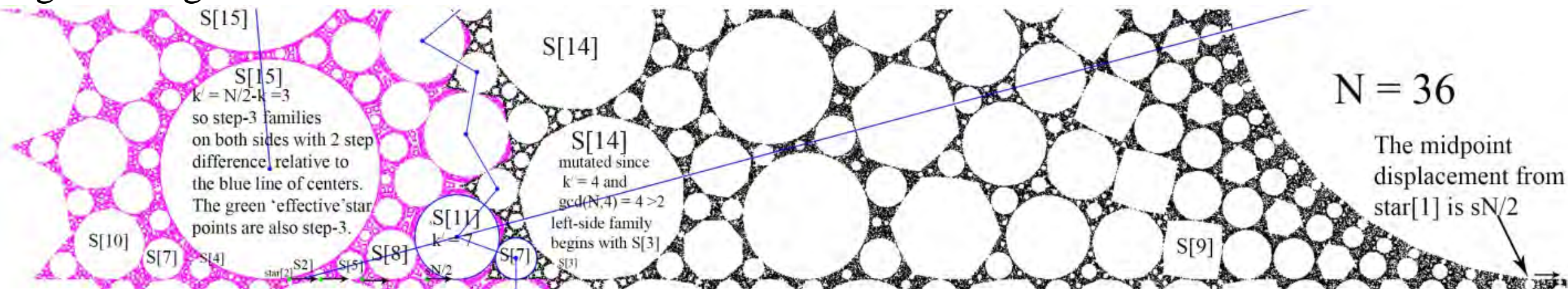

As noted earlier when N increases the S[k] tiles are torn apart but in the early web there will still be ‘triads’ like the one shown below for N = 200 where two neighboring S[k] are related through an intermediary. Here S[67] has $k' = 33$ and its step-33 right-side family which begins with S[66] will have an S[66-33] as the local ’S[1]’. This will guarantee that S[67] and a clone of S[33] will share an edge. The GFFT says that the midpoint of S[33] will be star[33] of S[67] + {sN/2,0} as shown below in green. All S[k] have this same canonical displacement but in general the ‘effective’ star points of a tile must be consistent with the right-side neighbor which evolves first.

**Figure 7** For N = 200 , Floor[N/3] is an accurate predictor for the ‘shepherd’ $S_N$ at S[66].

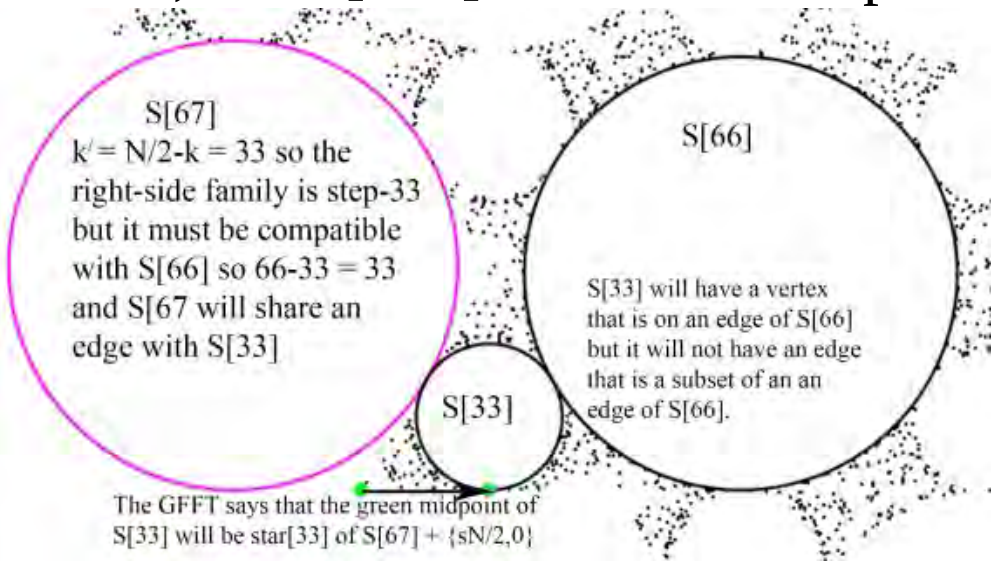

In the First Family Theorem of [H5] we showed that each star[k] of N defines an S[k] tile. To so this it was necessary to match the (right-side) star points of the proposed S[k] with the (left-side) star points of N. This meant that star[k] of N must be matched with star[N/2-k] of S[k]. Therefore the step-$k'$ nature of the S[k] has its origin in the need to assure conformity of S[k] to the star polygon constraints of N and it is independent of any mapping. This implies that there may be a close connection between any S[k] and the matching S[$k'$].

**Definition**: For an arbitrary S[k] of N define the *dual* of S[k] to be S[$k'$].

For N even, this is a reciprocal relationship since N/2-(N/2-k) = k and it begins with N and S[1], since N is congruent to D at S[N/2-1]. Lemma 3.2 of [H5] shows how any regular N-gon can be regarded as the S[1] tile of a next-generation D, so N and S[1] have a natural duality and this will apply to any S[k] and S[$k'$] pair when S[$k'$] is the primary tile in the local step-$k'$ family of S[k]. This is the edge sharing described above. The fundamental reason for this duality is the fact that N has reflective symmetry to match that of tan and cotan, namely $\tan[\pi/2 - \theta] = \text{cotan}[\theta]$.

Side of S[k] = sS[k] = SideFromRadius[rS[k],N] = $2r[S[k]\cdot\text{Tan}[\pi/N]$ . This simplifies to $2\text{Cot}[\pi k'/N]\cdot\text{Tan}[\pi k/N]^2$ so $sS[k']/sS[k] = \text{Tan}[\pi k'/N]^2$ which acts as a scale factor. In the special case of S[1] and N at step-1, $\text{Tan}[\pi/N]^2$ is GenerationScale[N] = hS[1]/hN] = hS[1]. For N-odd, D has the same side as N, so replace N with D and the same relationship holds with S[k] and S[$k'$] of D. (When working inside N odd we refer to these as DS[k] to avoid confusion and the connection is that every S[k] of N is congruent to DS[2k] and N/2-k of D matches N-2k of N.)

**Figure 8** For N = 103 below S[34] has $k' = N-2k = 35$ which is still be true if it is regarded as a DS[68]. Its left-side step-35 family counts down from S[35] which is DS[70] so it contains a DS[35] = DS[$k'$]. S[41] on the left has $k' = 21$ and S[42] = DS[84] counts down mod-21 to 21.

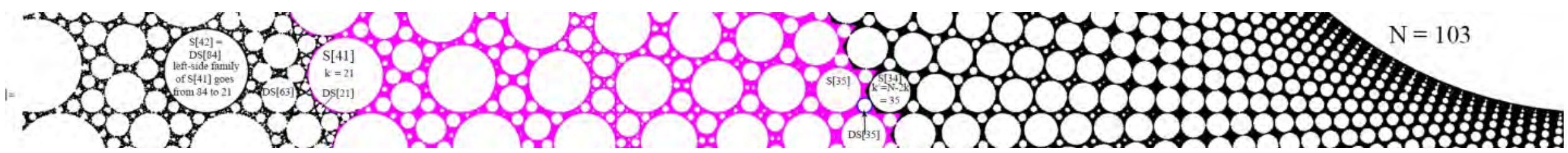

**Figure 9** For N = 300 the initial invariant region extends out to S[128] which is farther than the predicted S[100]. The plot below extends only out to S[142] at about x = -10, but the D tile at S[149] is at about x = -20. This initial black region was generated by a single point, namely star[1] of N displaced outwards by $\epsilon$ which is $\{10^{-8}, 10^{-8}\}$. This is about ¼ of our default 35-place accuracy. This took about 500 million iterations because typically points with long orbits will cover in a very uneven manner, But using multiple initial points could disguise transitions. The four transitions here are all very similar and we will just show detail of S[136]-S[139] .

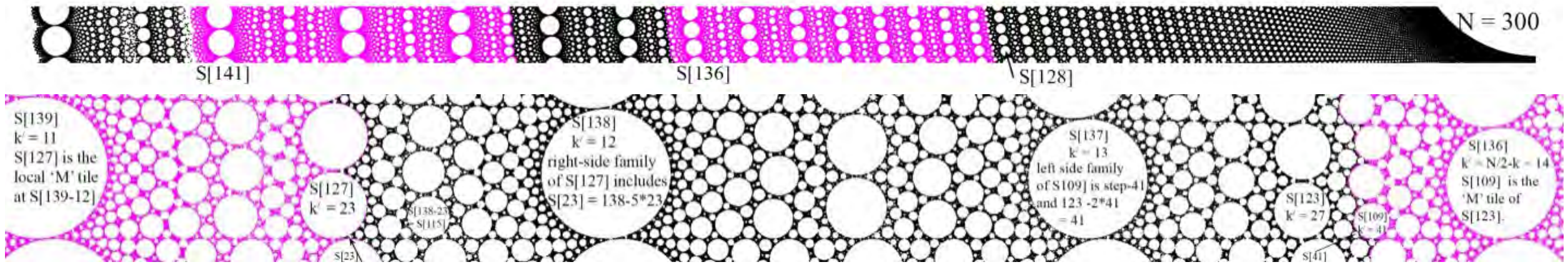

The S[136] transition is what we call level-3 because it entails two intermediate 'M' tiles – namely S[123] and S[109]. The left-side transition above and the earlier S[128] transition were both level 2 since they required just a single M tile. It is easy to construct an M tile like S[123]. It is the unique tile in the local families of both S[136] and S[137]. We give the steps below under 'Nuts and Bolts'. S[123] has k′ = N/2-k = 27 so its right-side family includes an S[109]. This is the target tile which is on the magenta-black border. S[109] has k′ = 41 and indeed it supports a small S[41] because 123-2(41) = 41. This tile will share an edge with S[109]. For S[138] on the left the local M tile is S[127] with k′ =23 and it supports an S[23] because 138-5(23) = 23.

**Nuts and Bolts**: The Generalized First Family Theorem (GFFT) of [H5] is a generalization of the fact that for N even, the S[k] can only be conforming to the extended edges of N if the right-side star point that matches star[k] of N is star[N/2-k] of S[k]. We define this to be an 'effective' star point of S[k] and the GFFT shows that by the symmetry or 'duality' of S[k] and N, any multiple of star[N/2-k] could also generate a tile conforming to S[k] acting as N. So S[k] can be 'parents'.

Using N = 300 as an example where S[137] has k′ = 13 the 'effective' star points are found by counting backwards mod-13 from the earlier S[136] so star[123] of S[137] will be effective and generate an S[123] and this will continue down to S[6]. The first step is to find <u>all</u> the star points of S[137] and that is easy.

Step1**:MidS137 = {cS[137],-1}; StarS137 = Table[MidS137 + {hS]137]*Tan[k*Pi/300], 0}, {k,1,149}]** Step 2: **MidS123 = {StarS137[123][[1]],-1} + {sN/2,0}, cS123 = MidS123 + {0,hS[123]}** Step 3: **S123 = TranslationTransform[cS123-cS[123]] [S[123]]**

**Note**: Most of the images shown here of invariant regions were generated by a single orbit, but we believe the conjecture by Richard Schwartz that for regular polygons, periodic points will have 'full measure', so any 'random' point that we choose as generator will likely be periodic with just a long period. But if multiple test points yield similar results it is possibility that these invariant regions might have some underlying dynamical structure that these 'test' orbits could help to uncover. In the pages to follow and in Part 2 we will occasionally track orbits and look for symmetry and structure using the 'Projections' that Richard Schwartz found useful for N = 8 in [S3]. See Appendix II for more about Projections. We close with a wish that in our very imperfect world, you and your dual will eventually find yourselves in proximity for a hug.

Returning to the web of N = 60 in Figure 2, the only First Family tiles that are in their correct position are the four ‘parents’ S[23], S[24], S[25] and S[26] . All the other ‘satellite’ S[k] are displaced copies of actual S[k] and the GFFT allows us to calculate their location (adjusted for mutations).These satellite families include central ‘penultimate’ tiles which we call $M_k$. Because these tiles span two families their steps will be additive, so $M_1$ will be step 4+5 and $M_2$ will be step 5+6 and $M_3$ will be step 6+7. Some of these tiles foster their own families but recursion is very limited.

**Figure 10** (same as Figure 4)

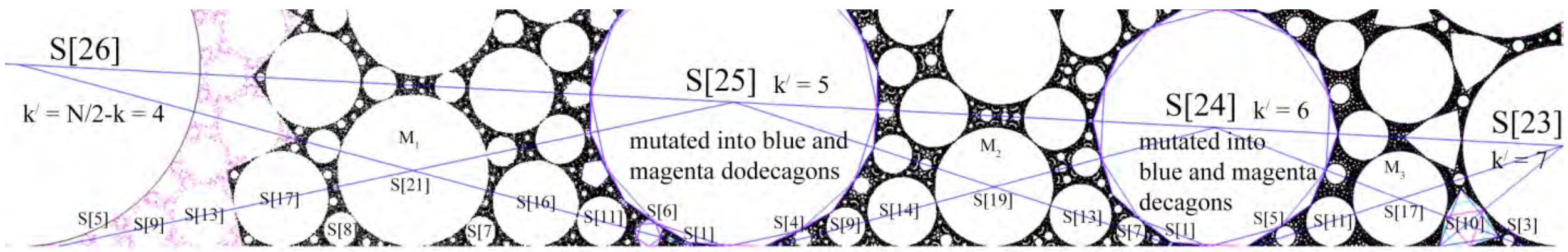

As noted above the step-4 family of S[26] includes S[13] and S[17] which share an edge. This means that one edge of S[13] is a subset of an edge of S[17]. This implies that S[13] and S[17[ linked together can create a dynamical barrier under rotation about N.

It is not surprising that a secondary tile like S[6] shown below can inherit the same k′ = N/2-6 = 24 step sequence as the original S[6], because for N twice-even S[25] is an N-gon that can play the role of a step-5 surrogate N. This means that both S[6]s can share a similar mutation. Since N/gcd(60,24) = 5, they will both be the weave to two pentagons with ‘base’ spanning gcd(60,24) star points, but this can occur in multiple ways and the Mutation Conjecture of [H5] predicts the minimal surviving star point will be the (absolute) minimal of N/2-1-jk′. This correctly predicts that S[25] here will have star[1] as minimal and also predicts star[5] for S[6] of N = 60 .

**Figure 4** Detail of S[25] showing the canonical two step difference between the left and right sides for N twice-even.

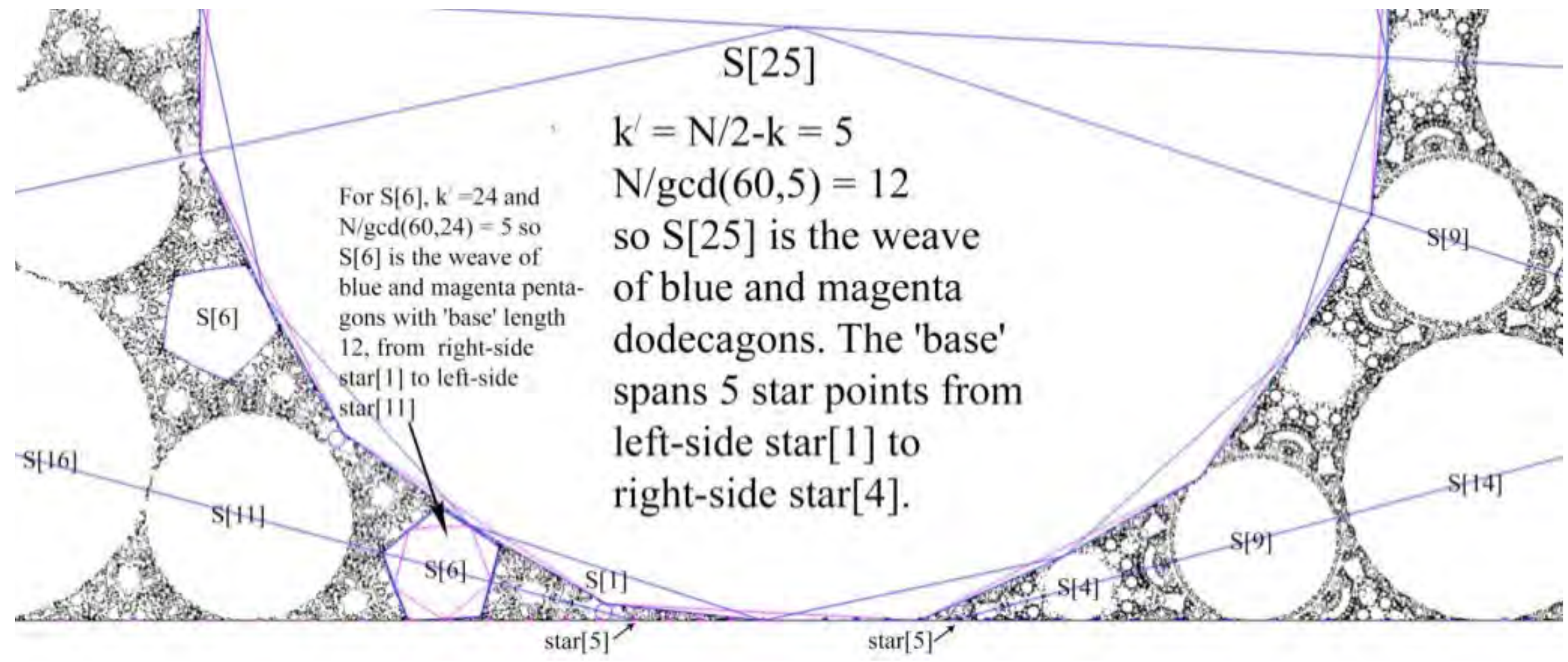

**Figure 5** Detail of S[24] and S[23] also showing a mutation in S[10]

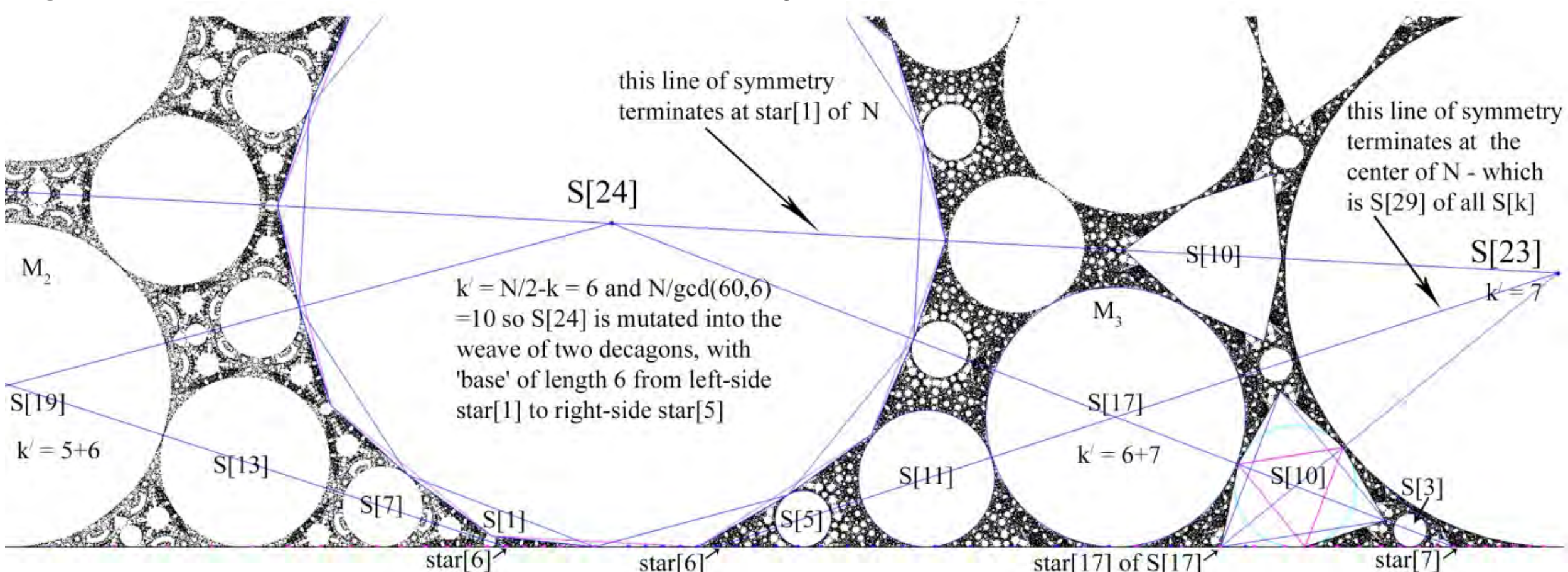

In the pages to follow it will be clear that most displaced S[k] will have displaced mutations so S[6] (and S[10]) are not exceptions. For S[24], $k' = 6$ and N/gcd(60,6) = 10 so the mutation is two decagons with base star[1] to star[5]. For S[10], $k' = 20$ so N/gcd(60,20) = 3 and the base will be 20 steps. The original S[10] goes from star[9] to star[11] as predicted but here it is star[3] to star[17]. This star[17] point also defines a step-13 family of $M_3$ consisting of just S[10] and S[23] and this implies that S[10] can play the part of an M tile with step 13 +7 to match $k'$. Indeed S[10] has its own short step-20 family consisting of S[3] and S[23]. Following the blue lines of symmetry all of this geometry references back to N which is the S[29] tile of every S[k].

The S[k] of N = 60 define what we call the first generation of N. To define the next generation we can use S[2] as a surrogate N or D. If S[2] has an S[2][2] that survives the web W, this can be used to continue a chain converging to the 'foot' of S[2] (which is locally congruent to star[1] of N). For the 8k+2 family such chains appear to exist, but in general all we know is that for N-even, S[2] will have a step-4 family of DS[k] where the matching S[1] of N is DS[N/2-2]. Here for N = 60, counting backwards mod-4 from DS[28] yields the DS[k] shown below. This 8k+4 family is the only one where S[2] is mutated since N/gcd(N/2-2,N) = N/4. So S[2] below is the weave of two regular 15-gons and the 8k+4 Conjecture says that for the mod-16 subfamily 12+16j the DS[4] will have an alternative 'parent' which we call Px. For N =12, Px = S[3] of N.

**Figure 6** The geometry local to S[2] for N = 60

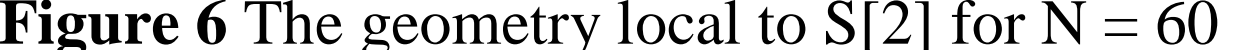

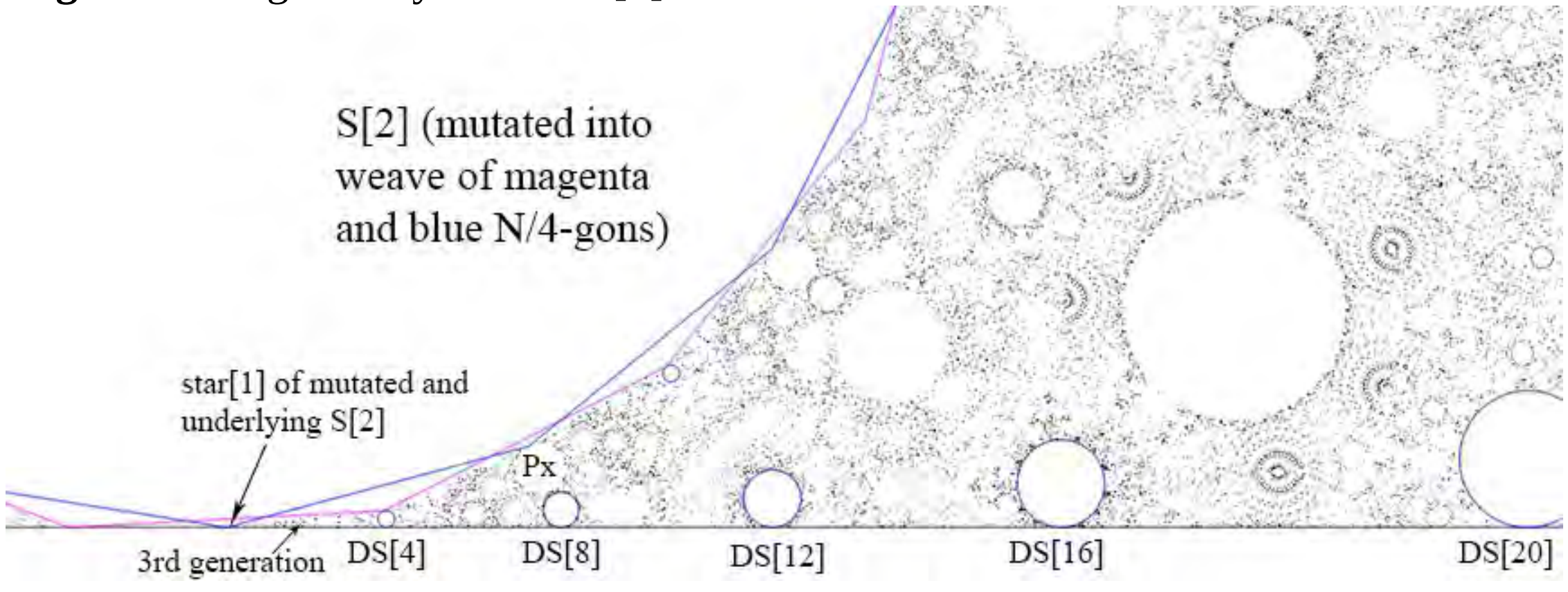

Below is a table illustrating what we call the '8-fold-way' for regular N-gons.

**Table 1**: A classification of web geometry on the edges of a regular N-gon - based on the Rule of 4 for N even (top) and the Rule of 8 for N odd (bottom). The examples show the local webs of the S[2] and S[1] tiles of N. The DS[k] are the 'next-generation' tiles of S[2] predicted by the Edge Conjecture so they arise in the early magenta web. The limiting web is black.

| 8k family | 8k + 2 family | 8k + 4 family | 8k + 6 family |
|---|---|---|---|
| DS[2] N = 16 S[2] DS[6] S[1] | DS[3] N = 18 S[2] DS[7] S[1] | DS[4] N = 20 S[2] DS[8] S[1] | DS[1] DS[5] N = 22 S[2] DS[9] S[1] |
| **8k + 1 family** | **8k + 3 family** | **8k+ 5 family** | **8k+7 family** |
| DS[5] N = 17 S[2] DS[13] S[1] | DS[7] N = 19 S[2] DS[15] S[1] | DS[1] DS[9] N = 21 S[2] DS[17] S[1] | DS[3] DS[11] N = 23 S[2] DS[19] S[1] |

**Organization of the three sections of this paper**

**Section 1** The singularity set of the outer-billiards map

(i) Definition of the outer-billiards map and the primitive domains (atoms)

(ii) Definition of the singularity set W (also known as the 'web')

(iii) Default web based on the star polygons of N (iv) The three maps

**Section 2** The local evolution of the singularity set W

(i) The web W is partitioned by the star points of N so it can be regarded as the disjoint union of the local webs of the S[k] tiles.

(ii) When N is even each S[k] is formed in a step $k' = N/2-k$ fashion and when N is odd these indices are doubled to $k' = N-2k$. The 'effective' star points of the S[k] will be also be step-$k'$. and the Generalized First Family Theorem of [H5] defines the right-side 'families' of the S[k].

(iii) The S[1] and S[2] tiles will have retrograde (ccw) steps 2 and 3 for N even and twice this for N odd. This simplifies the analysis of the joint S[1]-S[2] web and shows that S[2] has an effective step-4 web for N even and step-8 web for N odd.

(iv) The Edge Conjecture predicts that the 'next-generation' DS[k] tiles S[2] will exist at least in a mod-4 fashion for N even and a mod-8 fashion for N odd. These are called the Rule of 4 and and Rule of 8 and together they define 8 distinct classes of dynamics for N of the form 8k + j.

(v) The 8k+2 Conjecture predicts that if N has this form then S[1] and S[2] can be regarded as M[1] and D[1] of an infinite family of M[k] and D[k] converging to star[1] of N or D. This sequence will have temporal scaling N/2+1 as determined by difference equations for the periods.

**Section 3** Catalog of edge geometry for $N \leq 25$

**Appendix I** 'Deep-Field ' maps for = 19 and N = 200

**Appendix II** The invariant regions of N = 202 and analysis of orbits using 'projections'.

## Section 1. The singularity set of the outer-billiards map

**Definition 1.1** (The outer-billiards map $\tau$) Suppose that P is a convex n-gon in Euclidean space with origin internal to P. If p is a point external to P that does not lie on a blue 'trailing edge' of P, then the (clockwise) outer-billiards image of p is a 'central' reflection about the nearest clockwise vertex of P, so $\tau(p) = 2c-p$ where c is the nearest clockwise vertex of P as shown on the left below.

**Figure 1.1** The geometry of $\tau$ on the left and the geometry of the primitive domains $T_k$ on the right. These domains are partitioned by P and the blue trailing edges of P.

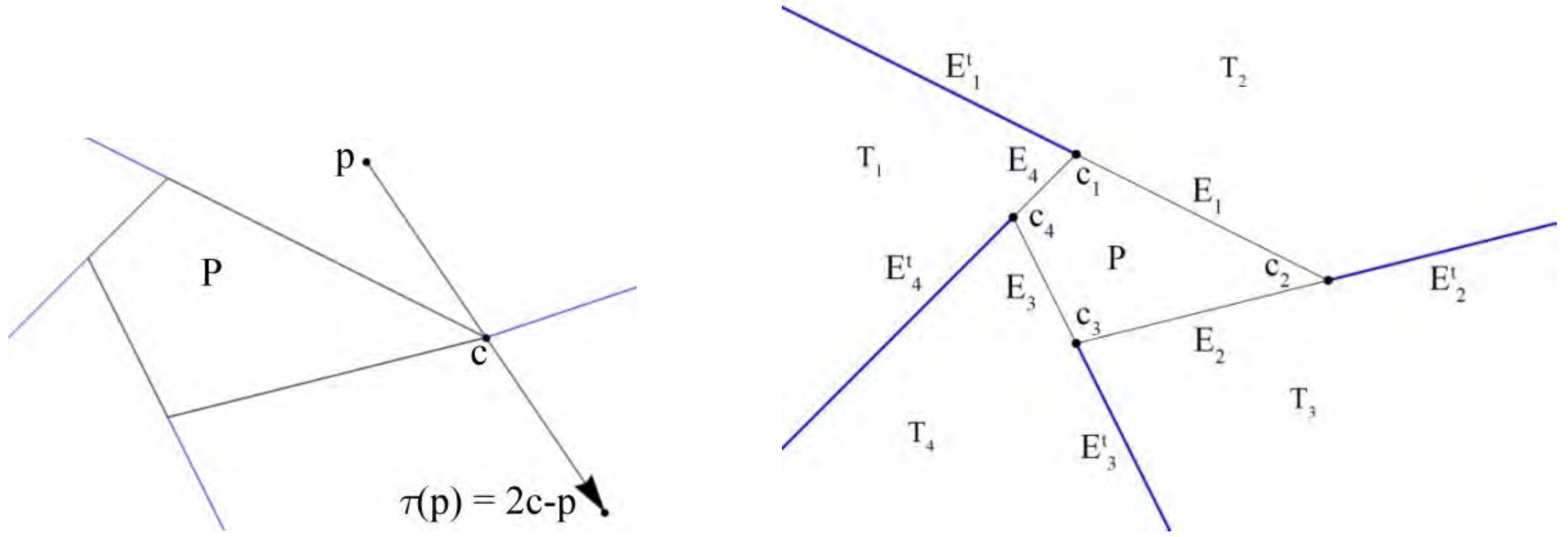

Since $\tau$ is not defined on P or on the trailing edges, the level-0 web is defined to be $W_0 = E \cup E^t$. where E is the set of edges of N and the $E^t$ are the extended 'trailing' edges of N. (For a counterclockwise $\tau$ these would be extended forward edges.) $W_0$ is called the level-0 exceptional (or singular) set of $\tau$. Since $W_0$ is connected, the complement of $W_0$ external to P consists of n disjoint open (convex) sets which are known as level-0 tiles or 'atoms'. Using these primitive $T_k$ tiles, the mapping $\tau$ can be defined as $\tau(p) = \tau_k(p)$ if $p \in T_k$ where $\tau_k(p) = 2c_k-p$. Therefore the domain of $\tau_k$ is $T_k$, which we write as $Dom(\tau_k) = T_k$.

It follows that $Dom(\tau)$ is the union of the primitive tiles, $\cup T_k$. Since the $T_k$ are called the 'level-0' tiles, after iteration k of the web algorithm (defined below) any new tiles which arise will be called level-k tiles.

By definition $Dom(\tau^2)$ is $Dom(\tau) - \tau^{-1}(W_0)$. The union of $W_0$ and $\tau^{-1}(W_0)$ is called the level-1 web, $W_1$. In general the level k (forward) web is defined to be:

$$\text{(i)}\quad W_k^f = \bigcup_{j=0}^{j=k} \tau^{-j}(W_0^f) \quad \text{where } W_0^f = E \cup E^t$$

The level k inverse web is defined in a similar fashion using $\tau$ and the extended forward edges:

$$\text{(ii)}\quad W_k^i = \bigcup_{j=0}^{j=k} \tau^{j}(W_0^i) \quad \text{where } W_0^i = E \cup E^f$$

At each iteration these webs are distinct and there are computational advantages to using both webs so we will define the level k web $W_k$ to be the union of $W_k^f$ and $W_k^i$.

It is not necessary to implement $\tau^{-1}$ explicitly because this map can be obtained from τ by reversing the orientation of the generating polygon. As long as the origin of the coordinate system is inside the polygon, the orientation can be reversed by taking a reflection about the vertical axis – which we call Tr.

**Figure 1.2** The mapping Tr and its inverse

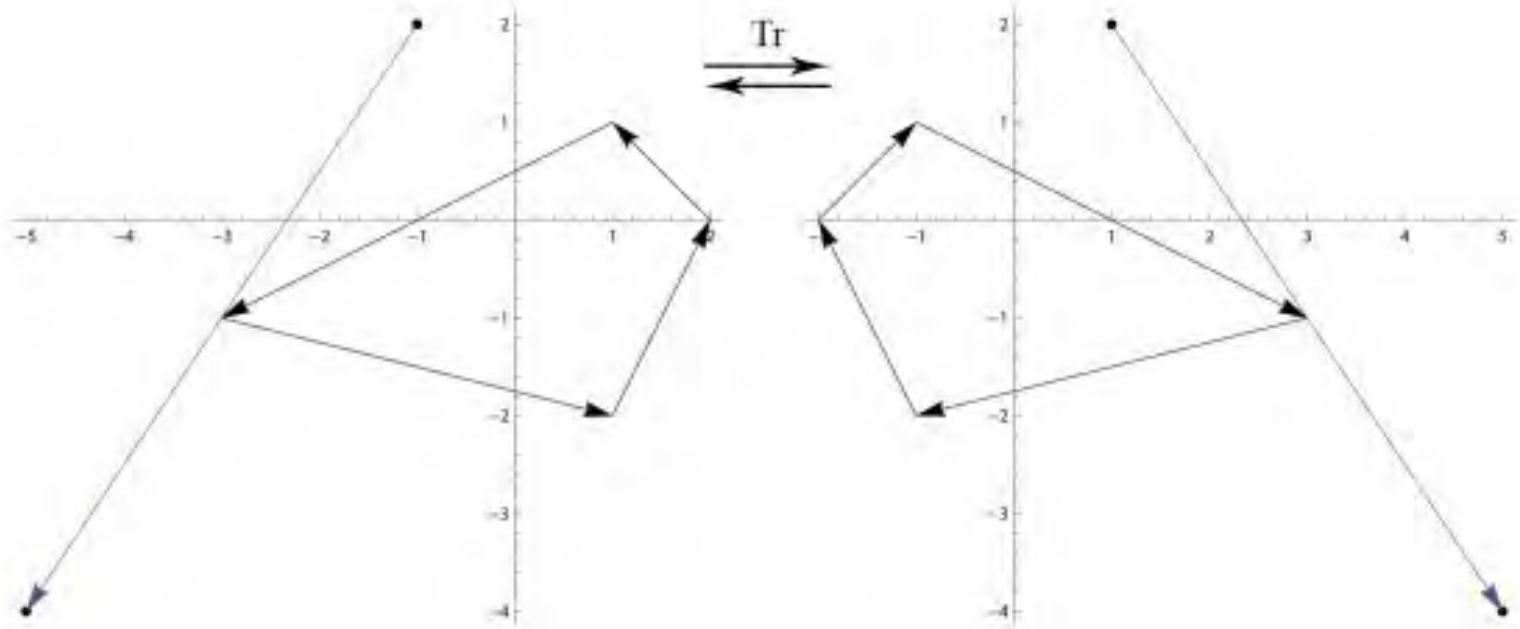

Our arbitrary choice is to implement τ clockwise, but for a counterclockwise matrix P, it is just a matter of reflecting P to get Tr[P], then using the clockwise τ, and reflecting the result back as shown above. For regular polygons centered at the origin, P and Tr[P] are identical (except for orientation) so the procedure for generating the web is very simple: implement τ for a clockwise P and recursively apply τ to the forward edges to obtain the inverse web $W_k^i$. Sometimes this is sufficient, but if a true web $W_k$ is needed, apply Tr to the resulting inverse web to get $W_k^f$. This works because there is no loss of generality in assuming that the original matrix was counterclockwise. Since $W_k^i$ and $W_k^f$ are reflections, any analysis can be done with either web.

**Figure 1.3** The two level-5 webs in the vicinity of the First Family tiles S[1] and S[2] on the edge of N = 11

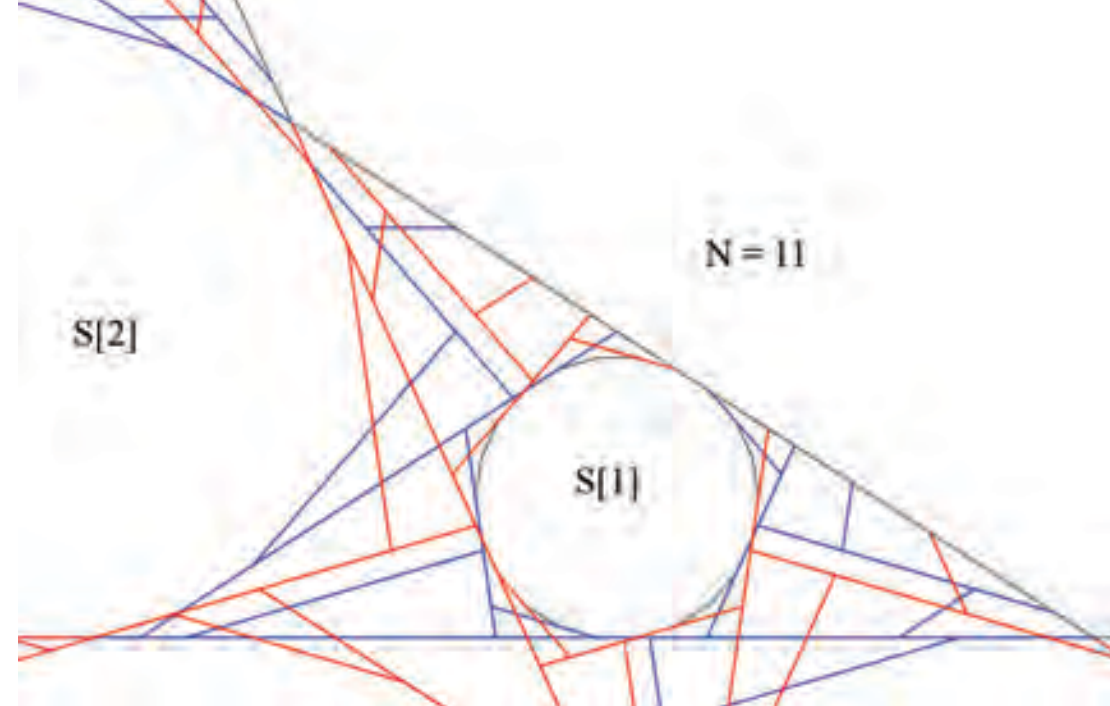

The web $W_5^i$ is in blue and the forward web $W_5^f$ is in magenta. These webs are just reflections of each other about the center of N, so using both webs is computationally very efficient. The tiles S[1] and S[2] get their names from the fact that their centers have $\tau$-orbits around N that skip 1 or 2 vertices on each iteration.

**The star polygon web W for a regular polygon**

As described in [H5], every regular N-gon defines a sequence of nested star-polygons formed by extending the edges until they meet. The 'maximal' star polygon in this process is simply $W_0 = W_0^f \cup W_0^i$ and this will be our traditional choice for initial $W_0$ because is easy to show that the resulting web W will be $\tau$- invariant. This maximal star polygon will always be bounded by a ring of 'D' tiles which are maximal among all regular tiles that can be formed in the web W.

By our remarks above, $\tau^{-1}$ is $\tau$ applied to a horizontal reflection of N, so $W_k$ and $W_k^i$ are also related by a simple reflection and it is our convention to first generate $W_k^i$ by mapping the 'forward' extended edges under $\tau$ and if desired a reflection gives $W_k$ also. In the limit W and $W^i$ must be identical but at every iteration they differ, so it is efficient to utilize both for analysis.

**Example 1.4** The star polygon webs of N = 7 and N = 14, The forward and trailing edges are shown in blue and magenta. Here we generate the level-k (inverse) webs $W_k^i$ by iterating the blue forward edges under $\tau^k$ for k = 0,1,2,3 and 50. The magenta trailing edges are shown for reference. At every stage these images could be enhanced by taking the union with the horizontal reflection. In the limit it would not matter.

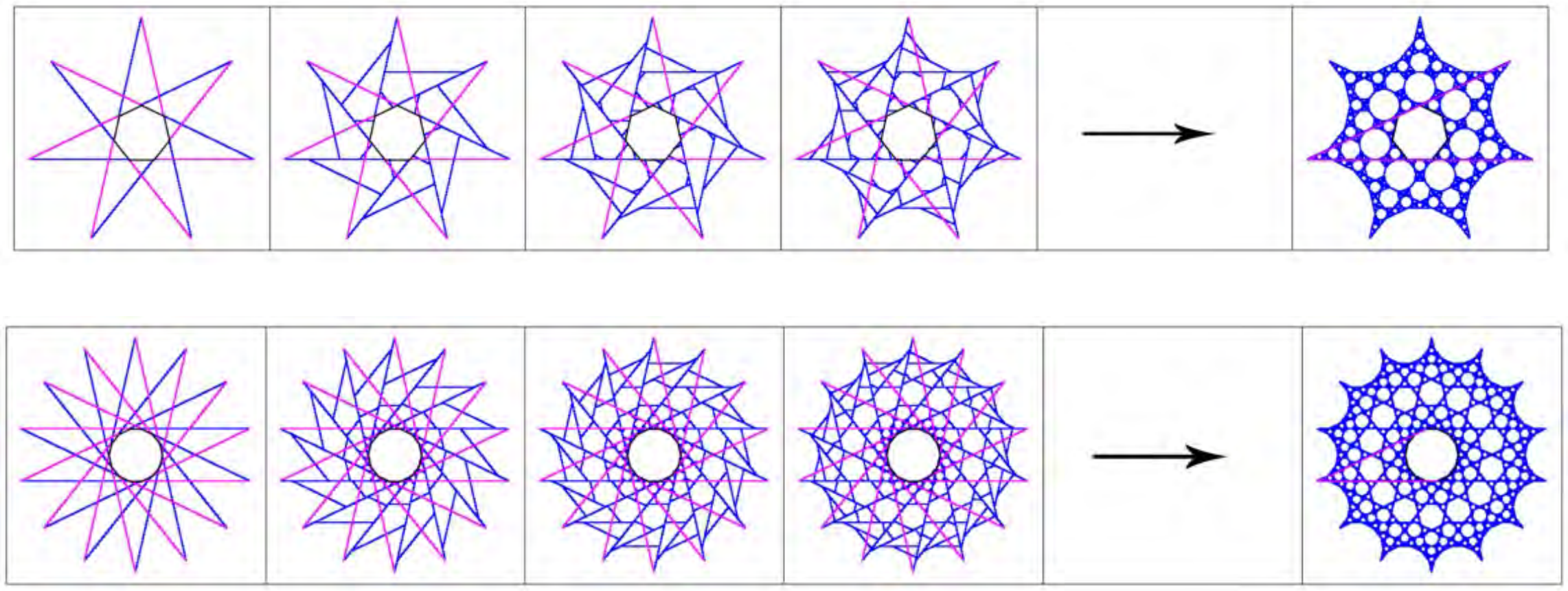

We call these 'generalized star polygons'. They retain the dihedral symmetry group $\mathcal{D}_N$ of N. Because the region bounded by D tiles is invariant under $\tau$, our default 'region of interest' will be the regions outlined above. By Lemma 4.1 of [H5] the orbital 'step-sequence' of D is maximal at step <N/2> which is N/2-1 for N even and Floor[N/2] for N odd. The points in these star regions cannot have steps that exceed that of D. The Twice-Odd Lemma of [H5] says that the web of N = 7 can be faithfully embedded in the web of N = 14. Except for scaling, they both share the magenta 'darts' outlined on the right. This will typically be our region of interest.

Implementing singularity sets by iterating all the extended edges of N $\tau$ is very inefficient. Because of rotational symmetry it should be sufficient to iterate a subset of the edges and this is where the Df and Dc maps are more efficient. They are based on just one or two extended edges. So for N = 14 shown below, the 14 primitive regions reduce to just three for the Digital Filter map and two for the Dual Center map.

**Figure 1.5** The primitive partitions (atoms) for the outer-billiards map, digital-filter map and dual-center map for N = 14 – followed by the singularity sets at level-200.

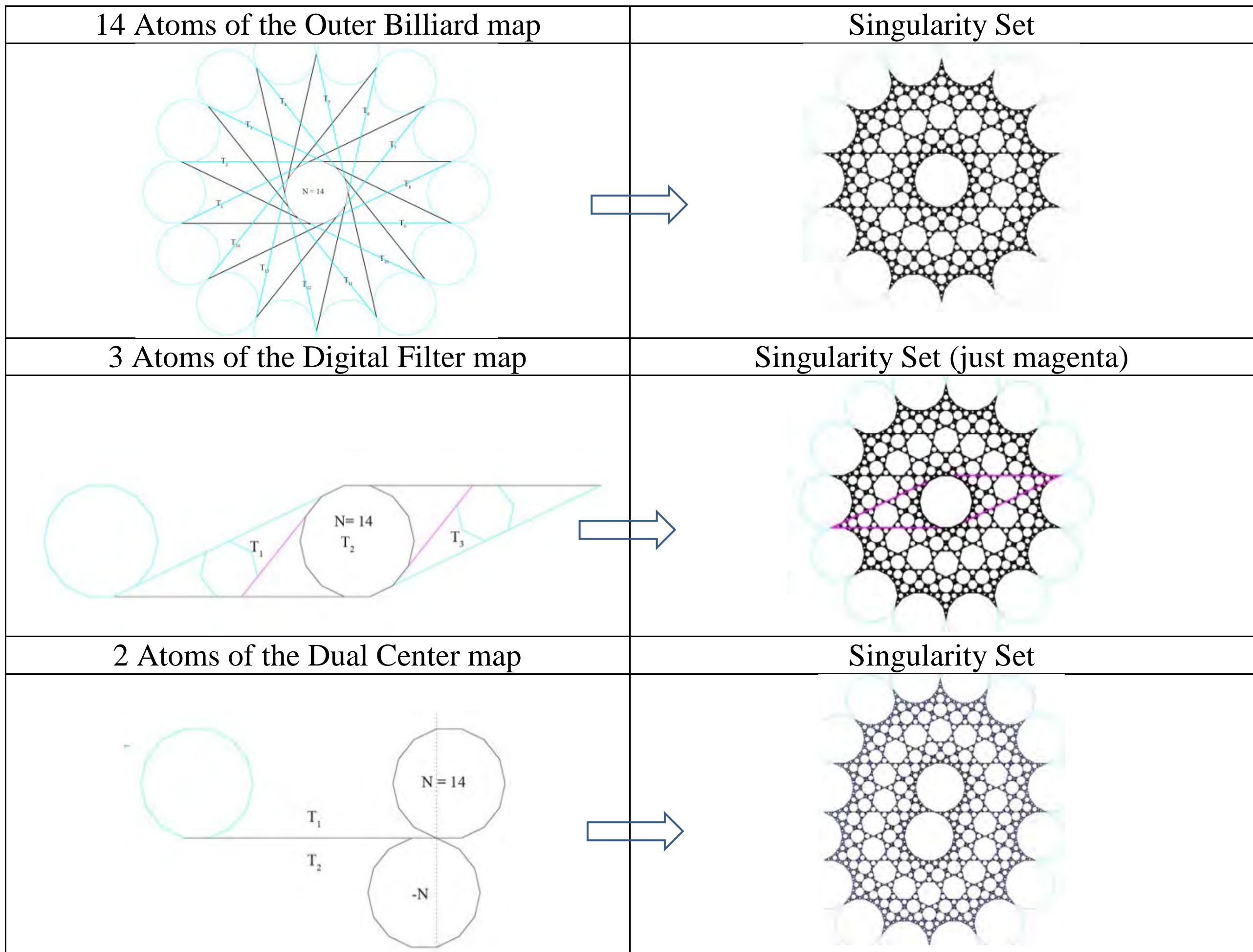

**Definition 1.2 ( The three maps)**

(i) (See Definition 1.1) The outer-billiards map $\tau(p) = 2c\text{-}p$ is not defined on N or on the trailing edges of N. The web W shown here is $W_{200}$ which results from 200 iterations of the 14 black forward edges. This ‘star’ region is invariant under $\tau$ and serves as a template for the global web.

(ii) The Df map is a toral map $[-1,1)^2 \rightarrow [-1,1)^2$ defined by Df[{x,y}] = {y,f(-x+ wy)} where w =2Cos[2π/N] yields an elliptical rotation and f(z) = Mod[z +1, 2]-1 models a sawtooth register. Df is rectified above to get a pure rotation on the magenta rhombus. The web shown here results from 200 iterations of a single forward edge of N (even only) and it is locally identical to W.

(iii) The dual-center map is globally defined as $Dc(z) = e^{-2\pi i/N}(z - \text{Sign}[\text{Im}[z]])$ where Sign[x] is 1,0 or -1 depending on whether x is positive, 0, or negative. The web shown here is 200 iterations of the single forward edge shown above. This web is identical to W on the upper (or lower) half-planes and both N and –N are parts of the web.

In terms of time-efficiency and computational efficiency there is a significant increase when using the digital-filter map or dual-center map. This is a major factor for large N. Below is timing for 20 iterations of the 50000 points in H =Table[{x,0},{x, .5, 1,.00001}] on a modest computer for N = 14:
(i) **Timing[Nest[Tau, H[[k]], 20] , {k,1,50000} ]]** = 120.4 seconds
(ii) **Timing[Nest[Df, H[[k]], 20] , {k,1,50000} ]]** = 2.9 seconds (plus 3.4 seconds to rectify)
(iii) **Timing[Nest[Dc, H[[k]], 20] , {k,1,50000} ]]** = 3.5 seconds.

This timing is a little misleading because by convention one iteration of the $\tau$-web iterates H on all 14 edges of N and these edges interact to produce a web which is accelerated relative to Df or Dc. The acceleration factor is variable, but for N = 14 it is about 20 times for a large-scale web, but this decreases for small-scale webs where the inter-edge interaction is less significant.

The biggest issue with Df or Dc is that the dynamics are very different from $\tau$ so they do not reproduce the global $\tau$-web faithfully. The horizontal axis defined by the extended 'base' edge of N is a line of discontinuity of all three maps but to preserve the rotational symmetry relative to this horizontal axis it will be necessary to rotate the web cw one step. This is still very fast.

Df and Dc are comparable in efficiency, but Df is only defined for N even and the points must be rectified, so we will typically use the dual-center map here, along with $\tau$. Since all exact calculations take place inside the cyclotomic field of N, there are sometimes advantages to using a complex-valued implementation of the web.

The key to the simplicity of the dual center map is the extra level of symmetry which results from shifting the origin to a vertex of N. Now the webs of N and –N can be generated together in a very efficient fashion and these two interact locally in the same fashion as the web development of W described in [H5]. This means that Dc First Families will be faithful to the First Family Theorem as illustrated below. This top-bottom juxtaposition of families actually occurs in the global $\tau$-web so the Dc web is locally compatible with the $\tau$ web. Since N = 14 and N = 7 share the same web, it would be possible to study them together, but typically we will study N = 7 at the origin where the dynamics are different. This is true for $\tau$ also.

**Figure 1.6** The First Family of N = 14 and N = 7 as generted by Dc

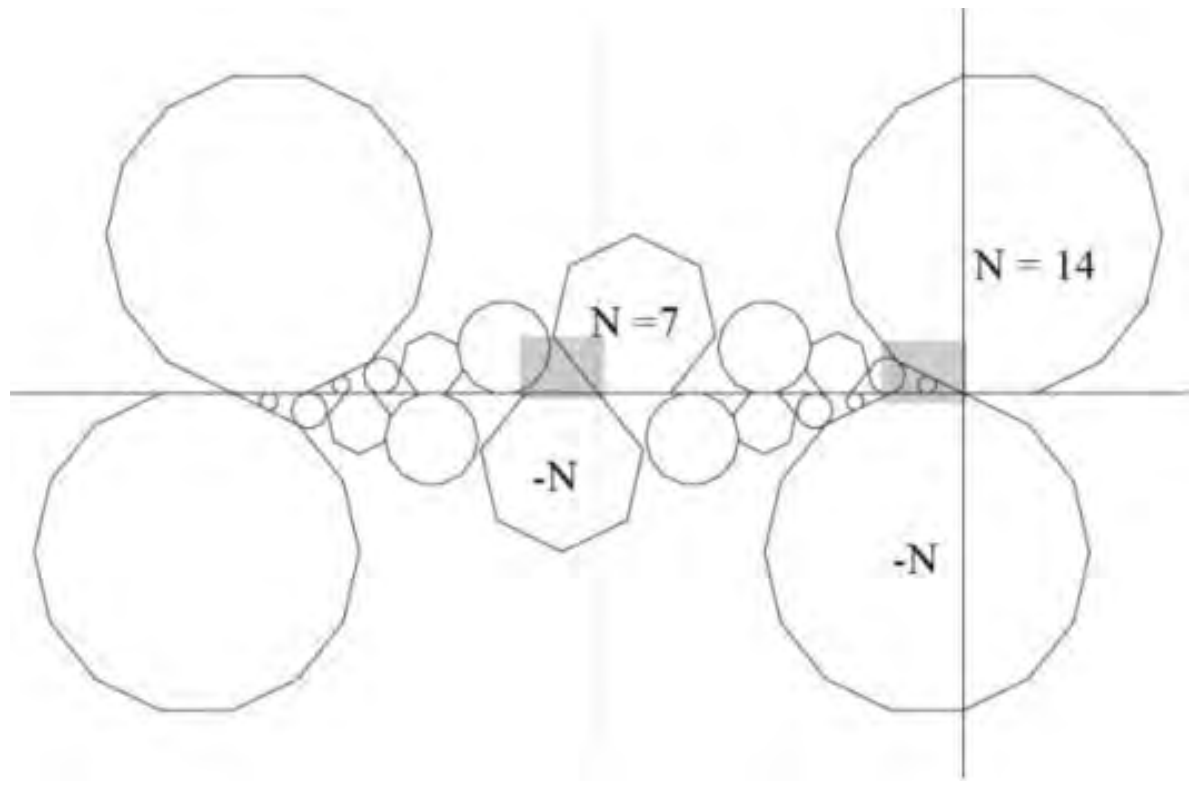

**Figure 1.7** – Iterations 0,1,2,4,7 and 8 of the web of the dual-center map Dc showing the construction of N = 7 and –N from the interval [-1,1]

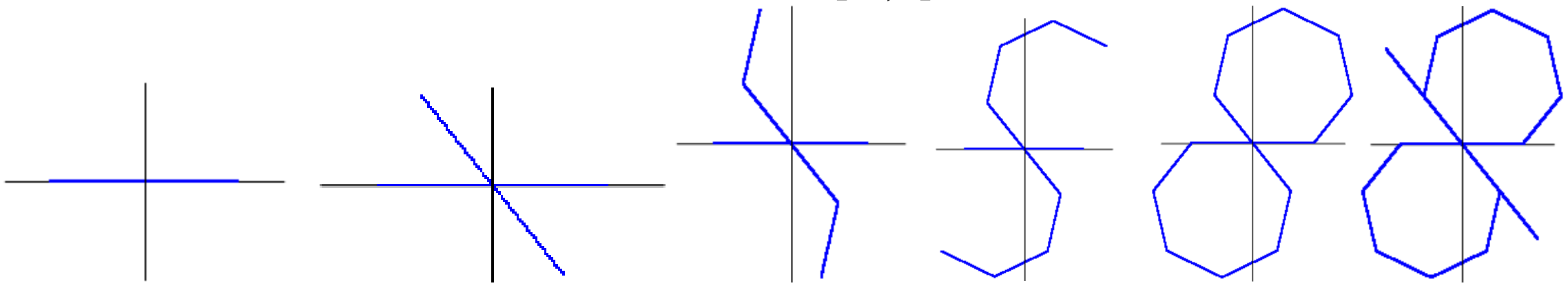

In Mathematica, the map is: **Dc[z_]:=Exp[-I*w]*(z-Sign[Im[z]])** where w = 2Pi/N to 35 decimal places. This extends z by 1,-1 or 0 depending on whether z is above, below or on the real axis. Then Dc rotates the result about the origin by $2\pi/N$. The following code generates the $8^{th}$ iteration on the right above :

**H = Table[x,{x, -1, 1, .001}]; Web = Table[NestList[Dc, H[[k]], 8],{k, 1, Length[H]}];**
**RealWeb ={Re[#],Im[#]}&/@Web; Graphics[Point[RealWeb], Axes->True]**

Under exact arithmetic, these extensions in level 8 would not exist and the web would be periodic after 7 iterations. Using Dc with an approximate rotation angle w will generate such extensions because when the two intervals [-1,0] and [0,1] map back to the x-axis there will be a small systematic error in the Sign function and most of the [-1,0] points will think they are negative while the [0,1] points will tend to be positive. This is easy to see using the Tally function for the 2000 points in the $7^{th}$ iteration of Dc.

**S[z_]:=Sign[Im[z]]; SF7 = S/@Nest[Dc,H,7]; Tally[SDc7] = {{-,944},{,972},{0,85}}**

This tally does not disclose which points are which, but a separate tally makes it clear that the majority of the [-1,0] points think they are negative and the bias is the opposite for the [0,1] points. This explains the extensions above. On careful examination there are some correct blue points which map back to the edges of N or –N.

It is necessary use an approximate w for extended calculations because exact evaluation of the Sign function may not be feasible. Therefore these extensions will occur in a more-or-less random fashion and there is no reason to iterate any points in the interval [-1,1]. Since both N and –N are known, the solution is use an initial interval of the form [1, x] or [-x,-1]. For the local geometry on the edges of N it is sufficient to set x = 2. Because of the reflective edge symmetry of N, it would be sufficient to iterate just [-1.5,-1] . The resulting web will also have +/- symmetry, so it is efficient to combine the resulting $W_k$ with $–W_k$ as shown by the following example from [H5].

Below is an example from N = 14 where we iterate 1,000 points in the interval H = {-2,-1} at a depth of 5,000. (The interval {-1,1} will generate N and –N in a period N orbit.) Here we crop these 5 million web points and their negatives and reflections to the desired region. (Less than 1 minute to generate and 1 minute to crop on a modest computer.)

**Figure 1.8** (The edge geometry of N = 14) F[z_]:=Exp[-I*w]*(z-Sign[Im[z]]); w =N[ 2*Pi/14, 35];(35 decimal places); H=Table[x,{x, -2, -1, .001}]; Web = Table[NestList[F, H[[k]], 5000],{k, 1, Length[H]}]; RealWeb = {Re[#], Im[#]}&/@Web; WebPoints = Crop[Union[RealWeb, -RealWeb, Reflection[RealWeb]]]; (about 450,000 points)

**Graphics[{AbsolutePointSize[1.0], Point[WebPoints]}}]**

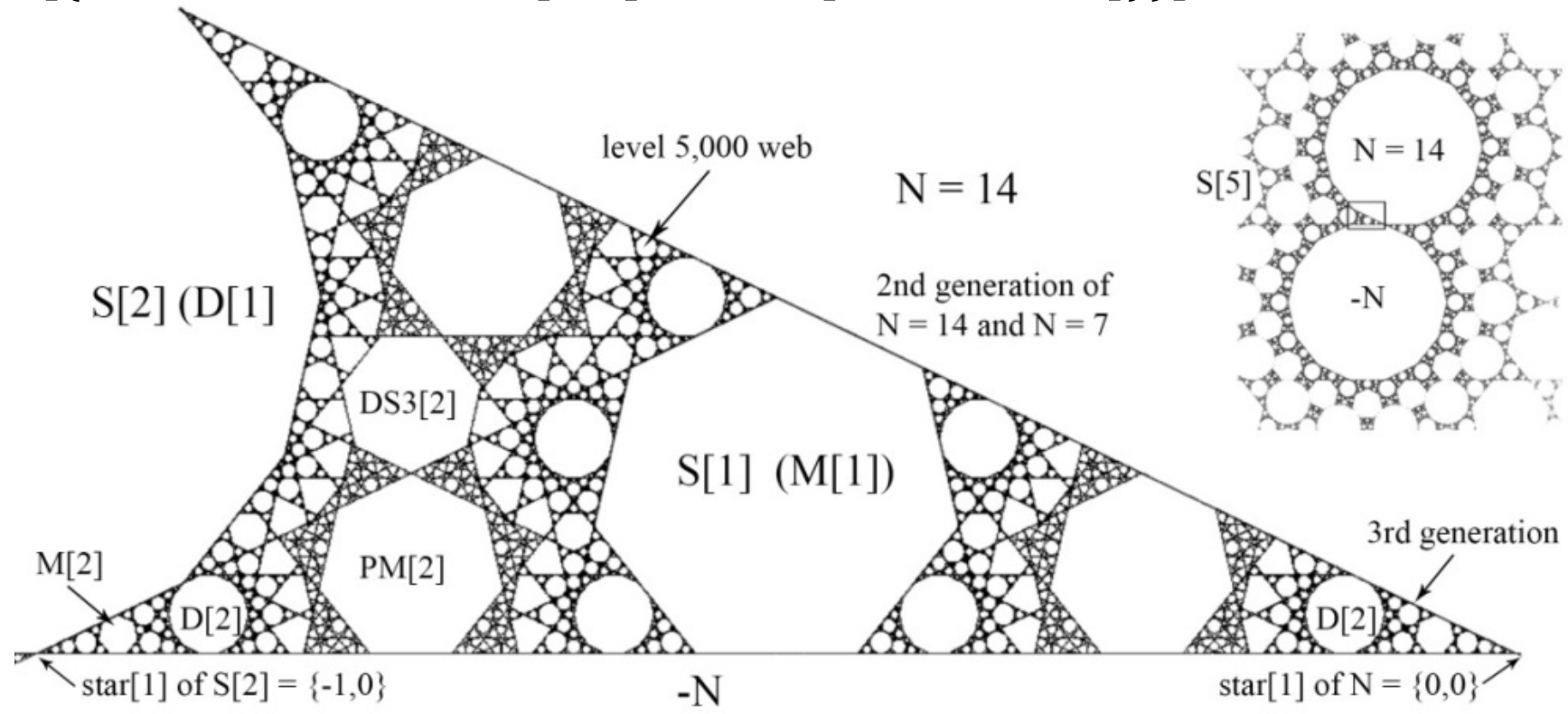

To generate the full star polygon web for the dual-center map as shown in the insert above, we can use the blue interval shown on the left below. This interval includes [-1,1] so that N and –N will be part of the web. Shown here are levels 0,1, 2, 3, 8, 25 50 and 200.

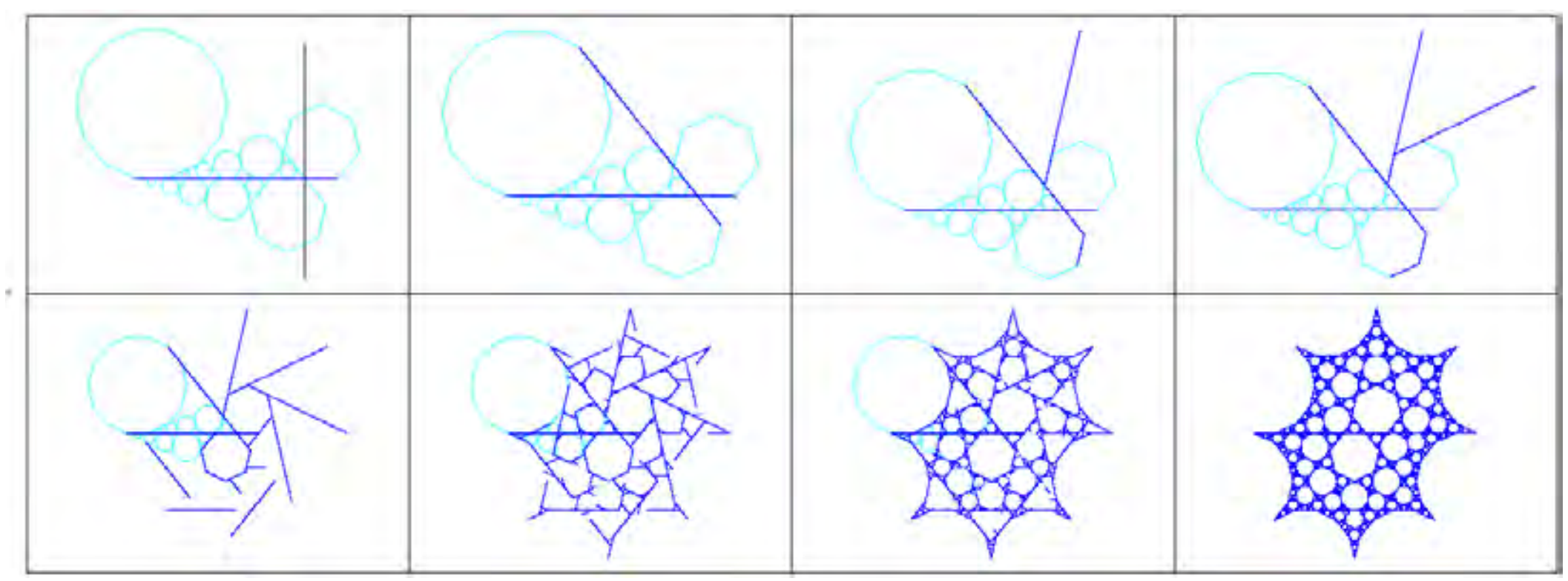

The regions above and below the real axis are identical to the outer-billiards web, Since the $\tau$-web has $2\pi$ symmetry about the center of N, a local rotation of this Dc web will yield correct geometry above and below the horizontal axis. This web algorithm is very fast and there is a bonus of plus-minus symmetry as in the Df map. This maps the lower web to the upper as shown below in magenta for level-8. In addition there is reflective symmetry in green. The reflective symmetry here is relative to the center of N so it only applies to the upper half plane. These symmetries yield an accelerated web with amazing efficiency and make it possible in two lines of code to generate webs with greater detail than ever before.

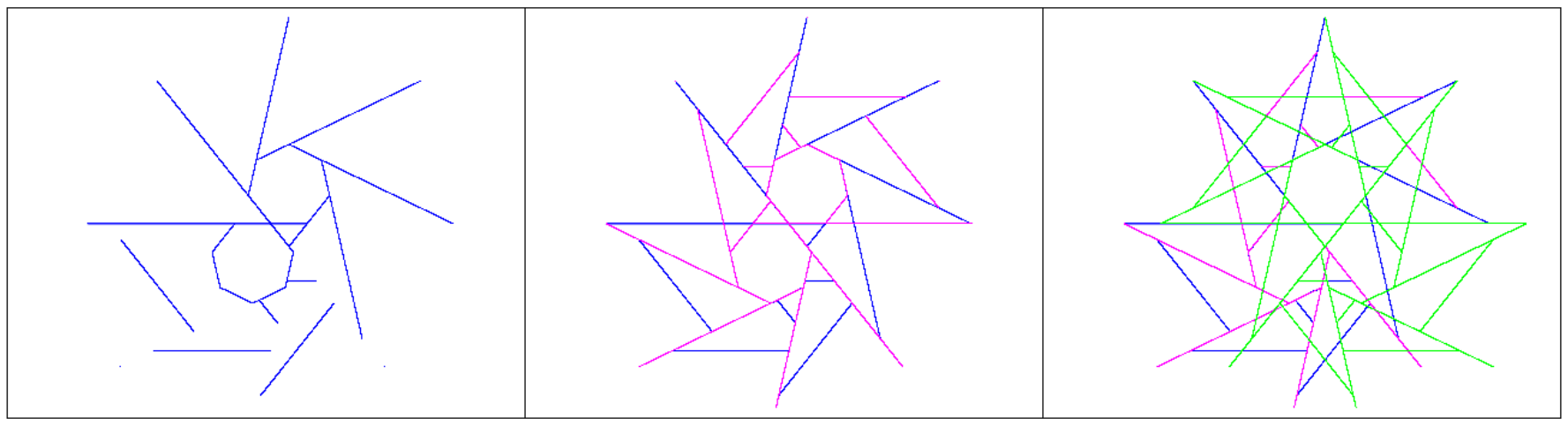

**Note**: This web algorithm is so efficient that processing the graphical data can be a major issue. For Mathematica the raw data is typically points in 35-decimal place 'postscript' format, so the data files are large. This size can be an issue when converting to raster form for display. This is usually done with Photoshop but the data points need to be in high-speed memory for processing and this can be a major issue for billions of points. Since these images are generated inside Mathematica, which is itself a very powerful image processor, an alternative is to save the raw data but use Mathematica to generate on-screen images which can be preserved by software or by a simple screen capture on a 4K monitor with resolution of 3840 by 2160. Some graphics cards allow the user to double this resolution to 8K and we will occasionally do this. On a typical monitor the screen resolution is a relatively coarse 72 dpi but inside Photoshop a 4K image captured with Screen Save will be 23Mb and 33 by 50 inches, which will still have good quality at 200 dpi (and smaller size) for printing. See the Appendix for more on postscript files.

**Section 2. The evolution of the outer-billiards web W**

In Section 4 of [H5] we described the evolution of the early web, and we will summarize the results here. For every regular N-gon, the star points partition the extended edges so the web W can be regarded as the disjoint union of local webs determined by these intervals. We show that in the early web, each star[k] to star[k+1] interval defines an S[k] tile of the FFT, so W can be regarded as the *disjoin union of the local S[k] webs*.

For edge geometry our interest is primarily the local webs of S[1] and S[2] but the case of N = 60 of the Introduction shows that this geometry is typically shared by multiple S[k] in an invariant region local to N. For N = 22 this invariant region extends to S[5] and the graphic below shows the intricate fashion in which these five local webs interact.

**Figure 2.1** The local webs of S[1] through S[5] for N = 22

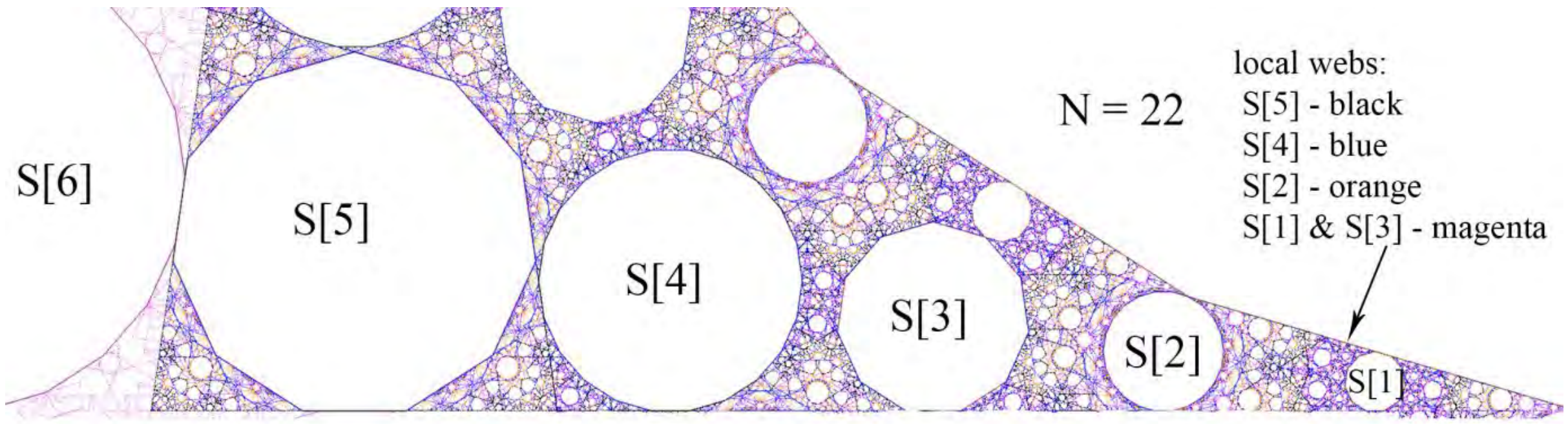

We show in [H5] that each of the star[k] intervals evolve by an iteration of 'shear and rotate' where the shear is of constant magnitude sN and the rotation angle for each S[k] interval is the 'star[k] angle' $k'\phi$ with $\phi = 2\pi/N$ and $k' = N/2-k$ for N even and $k' = N-2k$ for N odd. We will use N = 22 to illustrate the even case. By rotational symmetry it is sufficient to analyze the evolution of W in a single domain (atom) of $\tau$ or $\tau^{-1,}$ and this is done below for N = 22.

The magenta trailing edges of N are lines of discontinuity of $\tau$ so they partition each domain of $\tau^{-1}$. For N even there are N/2-1 star points which partition the horizontal forward extended edge L as shown below. The outermost region will be unbounded. For N = 22 we will see how this 10th region defines the S[10] (D) tile of the First Family. This tile will be congruent to N but it will evolve in a 'retrograde' fashion so if the center was shifted to D, its extended edges would generate N. By symmetry D also generates a left-side copy of N and this defines endless rings. Our default region is the first ring bounded by 22 Ds. In this region S[9] is in a central position and the Twice-Odd Lemma of [H5] says that its 'in-situ' web geometry will match N = 11.

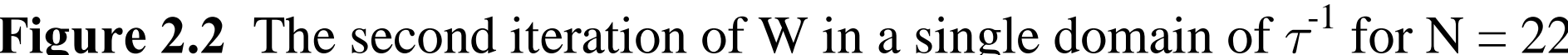

**Figure 2.2** The second iteration of W in a single domain of $\tau^{-1}$ for N = 22

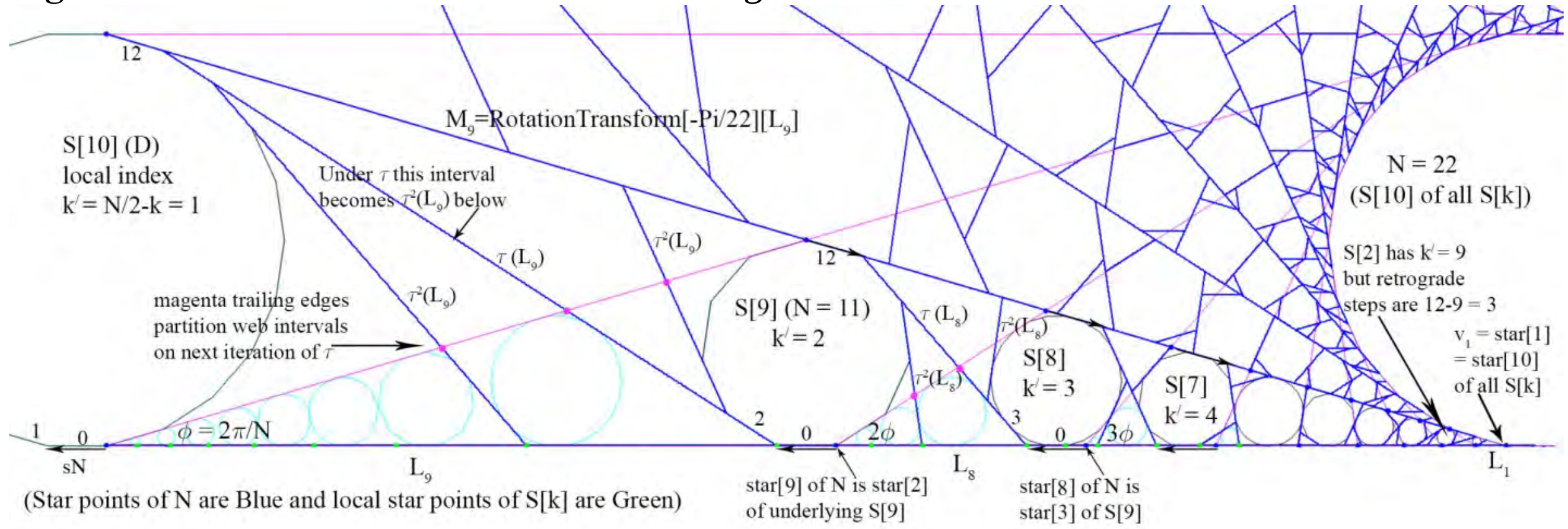

As indicated above, the evolution local to each star[k] is a simple outwards 'shear' of magnitude sN (as $\tau$ changes from one target vertex to the next) and then a variable rotation that would align each magenta trailing edge with the horizontal 'base' edge L. These rotations are the star angles $k'\phi$ where $k' = N/2-k$ for N even. So D with $k' = 1$ evolves in the same fashion as N itself. The star[1] point of D is first extended outwards by the shear and then rotated by $\phi$. This occurs recursively because the shear applies to all points on the horizontal line L as well as its image $\tau$(L). The rotation is a constant $\phi$ in the magenta region defined by D, so D will be maximal among all S[k] (and among all possible regular tiles in the web W).

For S[9] of N = 22, the star[1] point is star[9] of N and now the rotation is $2\phi$ . This is also applied recursively so S[9] evolves from the interval star[9]-star[10] in a step-2 fashion – which is why S[9] (and all odd S[k]) will be N/2-gons.This analysis applies only to the interior of the region between consecutive blue edges. On the opposite side of these blue forward edges the shears must be reversed as show by the arrows above. Two consecutive forward edges will be $\phi$ apart. Therefore for D the top edge will be N/2+1 = 12, and that remains true for all S[k].

When N is twice-odd this offset between top and bottom will be even, so tiles like S[9] will be formed in a redundant fashion from either cycle. Therefore all odd S[k] will be N/2 gons with skips 2,4,6,… For the even S[k] the rotation angle will be odd and hence not synchronized top to bottom, so the even S[k] will be regular N-gons formed from both cycles. (When N is twice odd, these $k' = N/2-k$ steps in the web evolution of S[k] will cause 'mutations' in the S[k] when gcd($k'$,N) >1, because the web cycles will be shortened. When N is twice even the N/2-1 offset between top and bottom cycles is odd so the cycles are no longer redundant and the mutation condition is a more forgiving gcd($k'$,N) >2. Therefore for N = 12, the S[4] tile is not mutated.)

When N is odd the (relative) shears are unchanged from the even case and the star angles are compatible since they are of the form $(\phi/2)(N-2k) = \phi(N/2-k)$. Since the S[k] are 2N-gons their local indices are $k' = 2(N/2-k) = N-2k$ and D will again have index 1 with rotation angle $\phi/2$. Therefore D will be a 2N-gon with the same side as N. Since $k' = N-2k$ must be odd, the primary web cycle be odd and it will be shortened iff gcd(N-2k,2N) > 1. The top cycle will also be odd since it is based on edge N+2 of D, so these two cycles will be synchronized mod 2 as in the twice-odd case, but now both cycles are odd relative to the base edge, so there will be no gender-based mutations in the S[k], but of course D will have the full spectrum of gender changes in the S[k] tiles – which we call the DS[k].

**Local webs of the S[k]**

Our convention for numbering the S[k] is right to left but the star angles increase left to right so for N even the $k'$ steps of S[k] are N/2-k and this is doubled to N-2k for N odd. This means that the critical S[1] and S[2] tiles will have webs with large step sizes, but conceptually there is no issue with regarding large cw steps as small ccw steps as shown here with N = 22.

**Figure 2.3** The S[k] tiles of N = 22 showing the level-1 local webs with edge steps 2 to 10

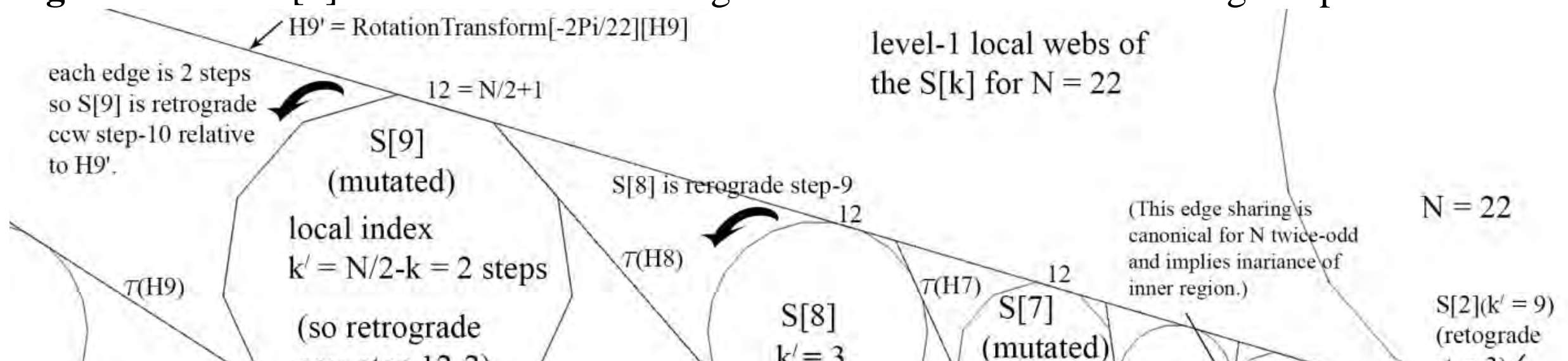

Every S[k] defines two local webs. Our clockwise convention for N means that the 'traditional' local web of an S[k] would be defined by iteration of the left-side interval between S[k] and S[k+1]. These are the intervals H9, H8, etc shown above. But iterating an interval like H8 also defines what we call the right-side web of S[9]. Both webs are generated by the same algorithm so the only distinction is the choice of initial interval.

Since cw and ccw rotations have edge steps differing by N/2+1 for N even (and double this for N odd), the retrograde right-side web of an S[k] will be step k+1 for N even and 2k+2 for N odd. This means that the critical S[1] and S[2] tiles will have relatively small step sizes when their web evolution is regarded as ccw instead of the traditional cw.

This is purely a conceptual issue but it makes sense that at least the 'next-generation' vertex based tiles of any N-gon will tend to have opposite web polarity, so it is natural to view the S[2] of N = 22 as evolving ccw with 3 steps rather than cw with 9 steps. Likewise the S[k] of S[2] (which we will call the DS[k]) should have a cw web orientation. What is very special about S[1] and S[2] for N even or odd, is that their webs evolve together and complement each other.

**Figure 2.4** The combined web evolution of S[1] and S[2] for N = 22 begins with an initial (magenta and blue) interval that spans both tiles. This will be our default initial interval 0.

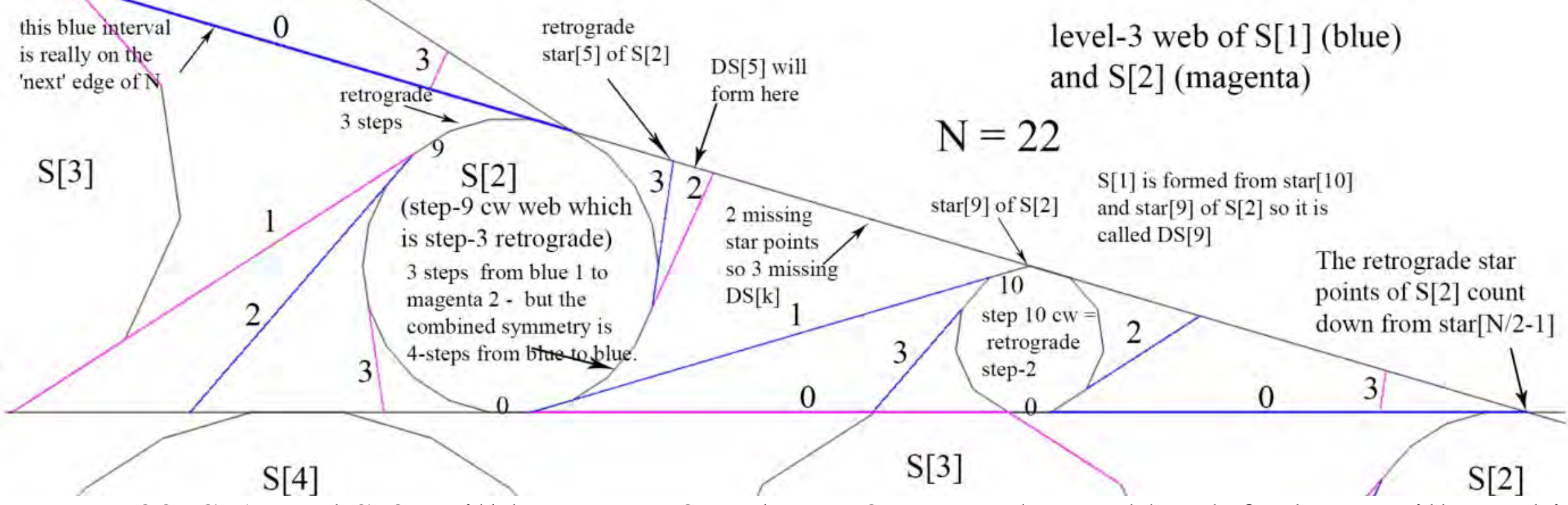

For N = 22, S[1] and S[2] will have step-2 and step-3 ccw webs, and by default we will combine their webs together as shown here. There will be a 1-step offset between these webs so the result is what we call a hybrid step-4 web for S[2]. This extra step is perfect for generating the shared next generation S[k] tile of S[2], which we call the DS[k].

For N-even, the GFFT of [H5] explains the consequences of a step-k′ web for any S[k] and this can be applied to S[2] with its step-4 hybrid web. In particular, every effective star point of S[2] should support a DS[k] tile in step-4 family of S[2]. The reference for these families is N itself which must be the DS[N/2-1] tile of S[2], so S[1] is DS[N/2-2] and it can be used as reference to count backwards mod-4. Therefore for N = 22, the ‘surviving’ DS[k] are just DS[5] and DS[1].. This called the **Rule of Four** . It does not exclude the possibility that other ‘volunteer’ DS[k] tiles may evolve in the web – just that these predicted DS[k] will evolve early and naturally. For N = 40 below, S[1] is DS[18] and the early survivors are DS[14], DS[10], DS[6] and DS[2].

**Figure 2.5** The web of N = 40, showing the symmetry local to S[2]

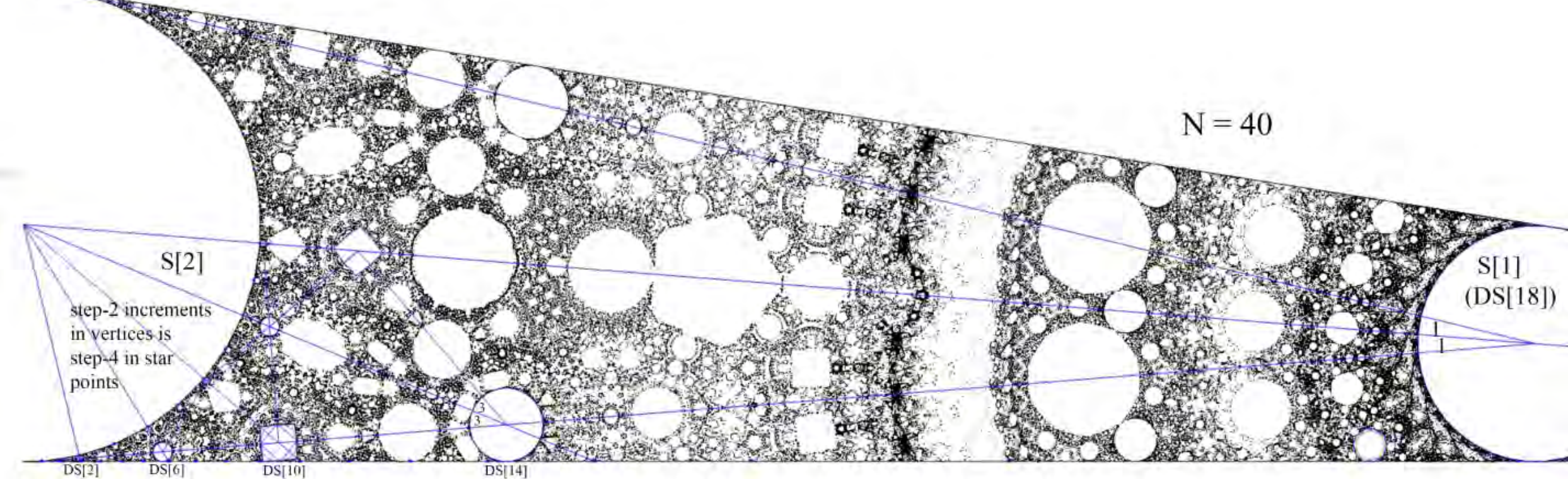

The primary lines of symmetry are determined N and locally they are based on S[1] which has k′ = N/2-k = 18 and retrograde step k + 1 = 2. At each effective star point of S[2], the web splits and DS[14] will have k′ = 6 . This determines the local symmetry of DS[14] as shown by the blue lines. The S[10] tile of N will be mutated since k′ = N/2-k = 10 and N/gcd(N,10) = 4. Based on S[1], these DS[k] seem to inherit these same k′ = N/2-k steps relative to S[2], so the DS[10] of N = 40 will share the same mutation as S[10] of N – namely the weave of two squares. Of course the S[4] and S[5] of N will also be mutated into pairs of octagons and pentagons respectively.

In general the web local to S[1] is not strongly linked to S[2] so S[1] will have its own lines of symmetry. It is easy to find these lines of symmetry for any N-gon independent of the web, and they are very uniform so DS[5] and DS[1] of N = 22 will have k′ = 6 and 10. This implies that these webs have a thread that links N = 40 back to the first non-trivial case of N = 8 where S[1] is DS[2]. The DS[k] centers provide the links like a road map.

**Figure 2.10** The step-1 symmetry lines for N = 22 in magenta

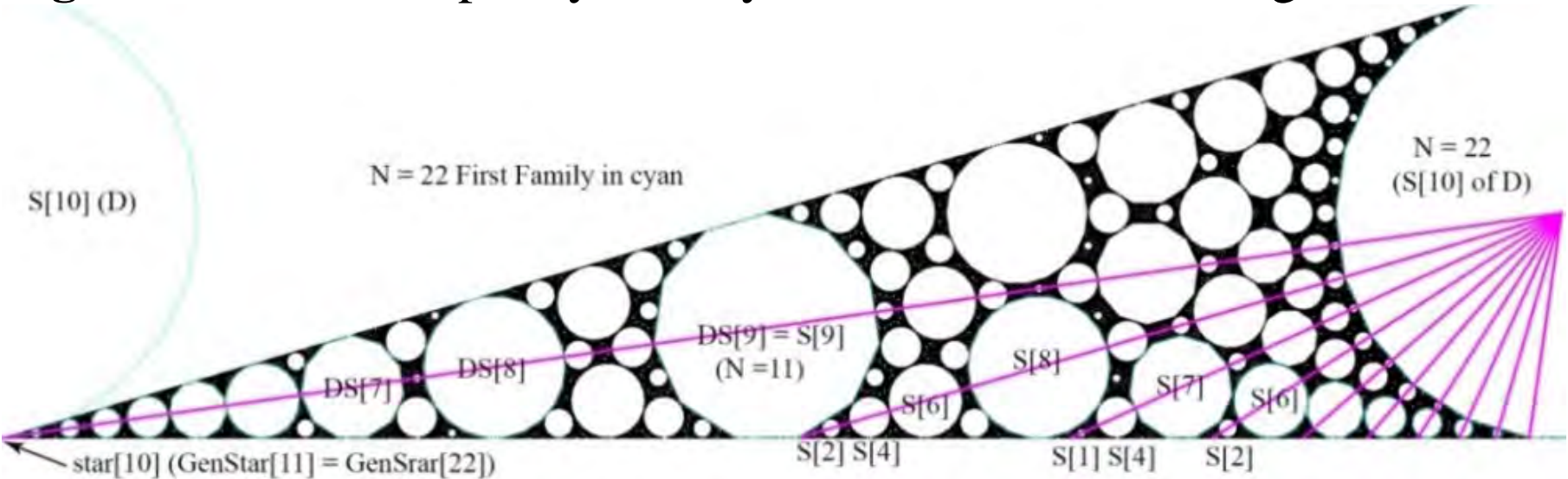

These lines joining the origin to star points showcase the members of the step-k families as discussed in the Generalized First Family Theorem of [H5] and the Introduction. At this time there is no version of this theorem for the tiles in the second generation because the DS[k] are relatively isolated, but there is still a lot to be learned about the shared symmetry. The most promising self-similarity occurs in the 8k+2 family as discussed below.

The **N-odd** case is very similar but now S[1] and S[2] are 2N-gons, so the web steps are doubled from N/2-k to N -2k and the retrograde steps are also doubled at k′= 2k+2. This means that S[1] will be retrograde step 4 and S[2] will be step-6, but the combined web of S[1] and S[2] will be step-8 as shown below. The DS[k] are now based on star[2] of S[2] and they will also be step-8 to match the effective star points. This is what we call the **Rule of 8**.

**Figure 2.6** The level-3 right-side webs of S[1] (blue) and S[2] (magenta) for N = 23.

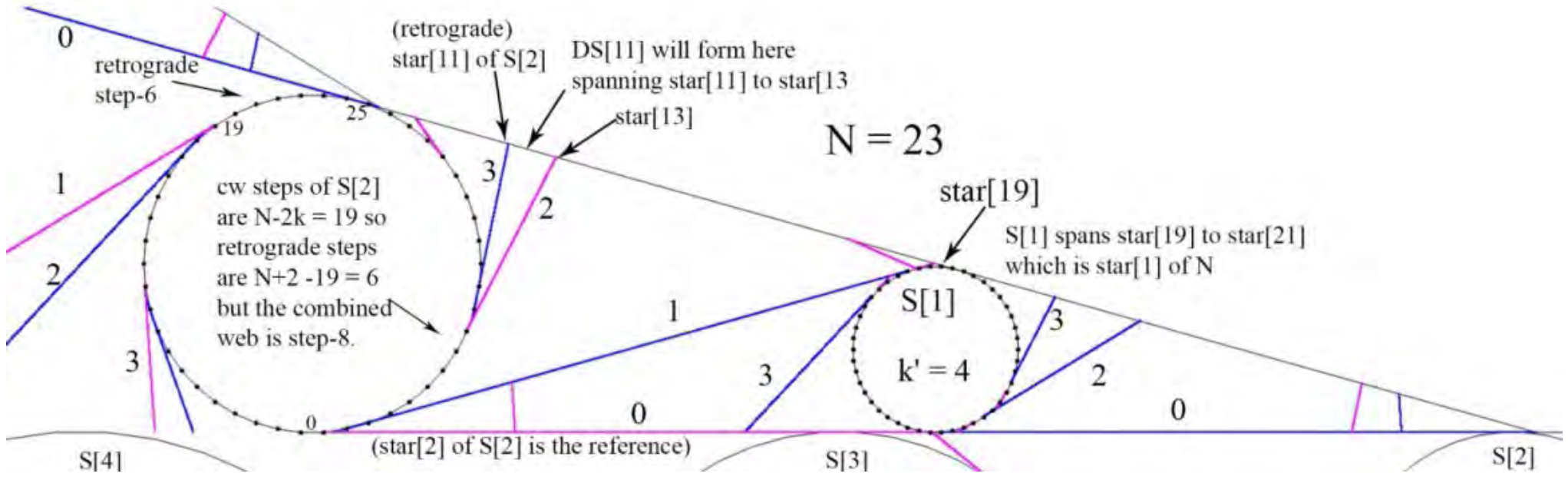

Once again the reference for the DS[k] is S[1] which is now DS[N-4]. Since this is odd all the predicted mod-8 DS[k] will have odd k. Overall this web development is very similar to N-even, but the GFFT will have to be revised to account for the increased displacement of the DS[k] since the web is now based on star[2] of S[2] instead of star[1]. This occurs because S[2] is 2N-gon and it implies that the DS[k] will be larger than the S[k] tiles in the First Family of S[2]. For N = 17 below, the increase in height of DS[5] over S[5] is proportional to the difference between displacements x and y. This difference is the distance from star[1] to star[2] of S[2] and it implies that hDS[]/hS[k[ must be 1/scale[2] of 2N. This is what we call an a D-M relationship in the normal First Family of 2N, since it matches hD[1]/hM[1]. Therefore all the DS[k] for N-odd will have this D-M relationship with the (virtual) S[k] in the First Family of S[2]. We call this the 'star[2] family' of S[2].

**Figure 2.7** The geometry of the DS[k] tiles for S[2] of N = 17

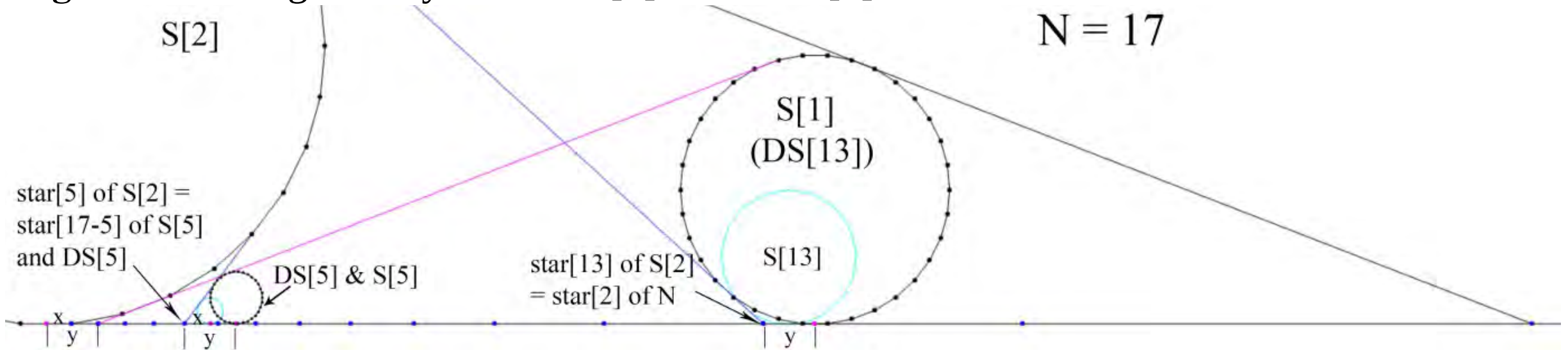

It is easy to construct these DS[k] for an arbitrary k and arbitrary odd N. Here are the steps

(i) hDS[k] = hS[k]/scale[2] = hS[2]*Tan[Pi/2N]*Tan[k*Pi/2N]/(Tan[Pi/2N]/Tan[Pi/N])
(ii) Midpoint = StarS2[[k]]+{hDS[k]*Tan[(N-k)*Pi/2N],0}; cDS[k] = Midpoint +{0,hDS[k]};
(iii) rDS[k] = RadiusFromHeight[hDS[k],2N] = hDS[k]/Cos[Pi/2N];
(iv) DS[k] = RotateCorner[cDS[k] + {rDS[k],0}, 2N, cDS[k]] (rotate the vertex at 3:00 for N twice-even or use a star[1] point of DS[k] otherwise)

**Figure 2.8**. N = 43 showing that exactly 4 DS[k] with odd k can fit between existing pairs so the natural fit is k odd (or k even) only.

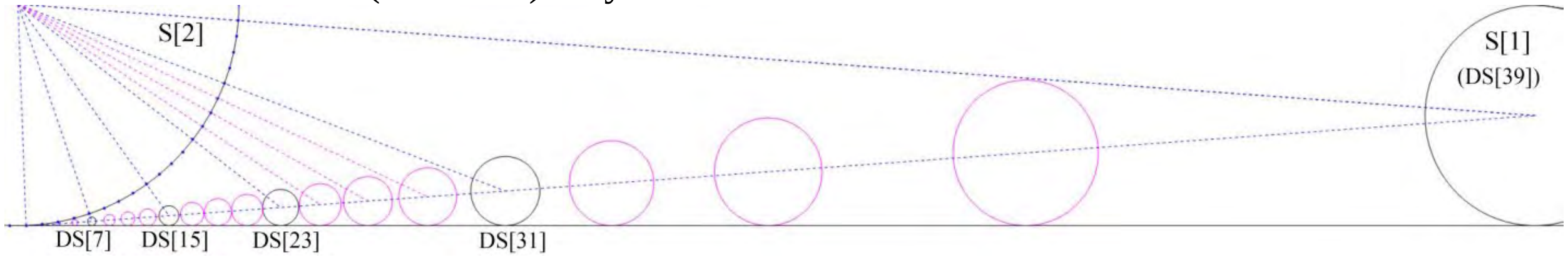

Below is a symmetry plot for N = 39 where the step-8 symmetry of S[2] is partitioned into 4 vertex-based steps corresponding to the step-8 effective star points.

**Figure 2.9** The web symmetry of S[2] local to S[2] for N = 39 in the 8k +7 family

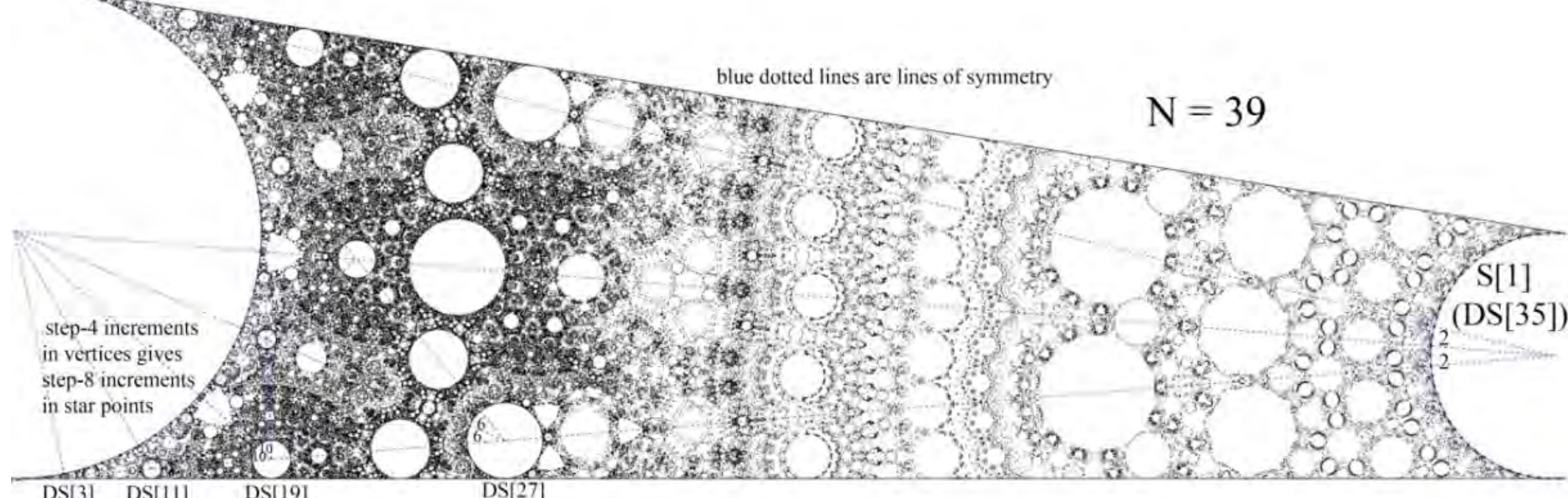

Even though N is odd, the DS[k] will evolve relative to S[2] which is a 2N-gon. Therefore the web steps of the DS[k] will be $k' = 2N/2 - k = N - k$ just like the even case. This will imply that the 2nd generations for N-even and N-odd are similar although the relative steps here are doubled. Compare this with the symmetry lines of the First Family where the 'parent' N-gon is step-1.

**The Edge Conjecture and Recursion**

Below we summarize these results about the evolution local to the S[2] tile of N.

**Definition 2.1** (The DS[k] tiles of S[2]): For any regular N-gon, the tiles on the edges of N which are strongly conforming to S[2] and exist in the limiting web W, will be called the 'DS[k] tiles of S[2]'. For N even these DS[k] will be strongly conforming to star[1] of S[2] and hence identical to the S[k] tiles in the (ideal) First Family of S[2]. For N odd the DS[k] will be strongly conforming to star[2] of S[2] and hence span the region from star[k] to star[k+2]. In general they will not be part of the First Family of S[2].They will be called the 'star[2] family' of S[2].

**Edge Conjecture.** For an arbitrary regular N-gon, every potential DS[k] tile on the edges of N which satisfies the Rule of 4 for N-even or the Rule of 8 for N-odd, will be among those that exist in the web W. The Rule of 4 says that counting down from S[1] at DS[N/2-2], the DS[k] will exist mod 4. The Rule of 8 says that counting down from S[1] at DS[N-4], the DS[k] will exist mod 8. (Other DS[k] may exist as secondary or 'volunteer' tiles.)

The natural question for any existence conjecture is whether it can be applied recursively to the DS[k] and potentially generate sequences of tiles. Theoretically any DS[k] can support an extended family and when N is even the DS[1] and DS[2] tiles of any S[2] can serve as 'next-generation' S[1] and S[2] to define sequences of tiles converging to the star[1] point of S[2].

**Definition 2.1 (**Ideal generations at GenStar[N]**)** For any regular N-gon with N > 4 , define D[0] = D and M[0] = M (the penultimate tile of D). Then for any natural number k > 0, define M[k] and D[k] to be DS[1] and DS[2] of D[k-1], so M[k] will be the penultimate tile of D[k]. The (ideal) *kth generation* of N is defined to be the (ideal) First Family of D[k-1]. Therefore M[k] and D[k] will be 'matriarch' and 'patriarch' of the next generation – which is generation k+1. By Lemma 3.1 of [H5], hM[k]/hM[k-1] = hD[k+1]/hD[k] will be GenScale[N] (or GenScale[N/2] if N is twice-odd).

**Figure 2.10** Ideal generations at GenStar[N] for N = 9 and N = 18

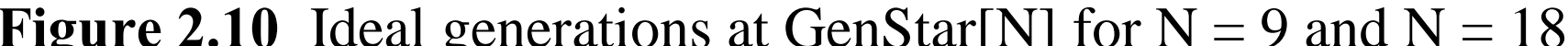

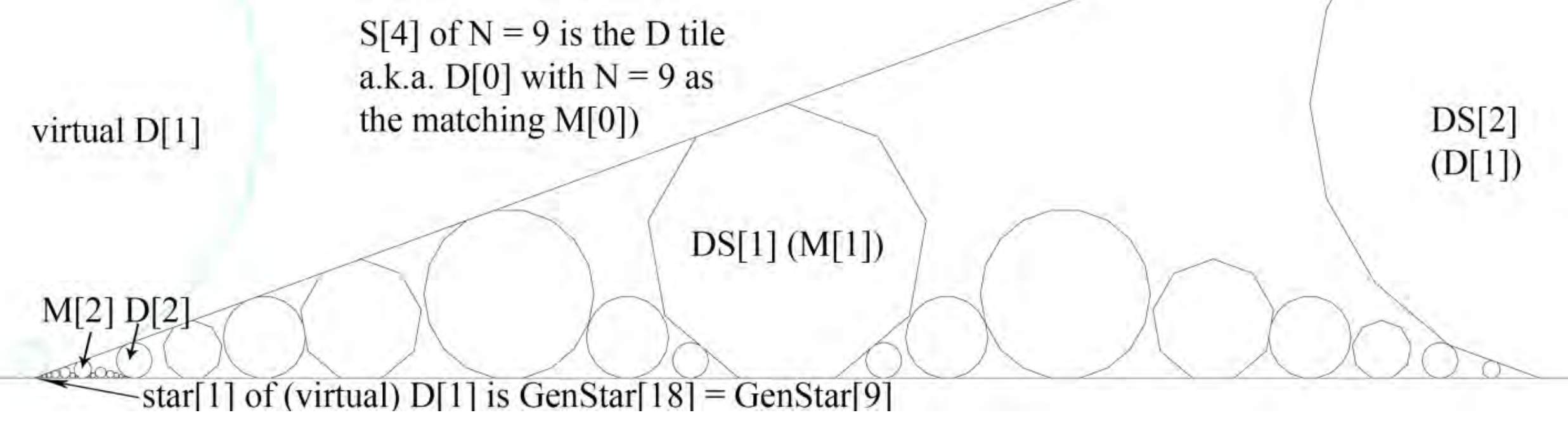

(The Mathematica code for these families is **FFM1 = TranslationTransform[cM[1]][First Family*GenScale]** and **FFM2 = TranslationTransform[cM[2]][FirstFamily*GenScale^2]**.)

In [H3](2013) we conjectured that if N is in the 8k+2 family, it will support self-similar families of S[1][k] and S[2][k] tiles for all k. Here we reproduce the 8k+2 Conjecture of [H5].

**The 8k+2 Conjecture**

Suppose N is a regular polygon in the 8k+2 family. Since N is even, it is interchangeable with D and Definition 3.3 describes a well-defined (ideal) sequence of M[k] and D[k] tiles converging to star[1] of N which is also known as GenStar[N] or GenStar[N/2]. By reflective symmetry there will be an equivalent sequence converging to star[1] of S[2] which is the DS[2] of N. It is our convention to study the sequence there because S[2] can serve as a 'surrogate' N or D.

In this sequence N itself can be regarded as D[0] and the matching M[0] will be the S[N/2-2] penultimate tile of N. M[1] and D[1] will be the DS[1] and DS[2] of N which are simply the S[1] and S[2] of N since N is even. Then in the ideal sequence of Definition 2.1, for any positive integer k, hM[k+1]/hM[k] = hD[k+1]/hD[k] = GenScale[N/2] since N is twice-odd.

(i) We conjecture that these M[k] and D[k] tiles will always exist and converge to GenStar[N] with geometric scaling GenScale[N/2] and temporal scaling given by (ii) below

(ii) If $D_k$ is the $\tau$-period of D[k], it will satisfy the second-order difference equation

$D_k = nD_{k-1} + (n+1)D_{k-2}$ where n = N/2 and $D_1 = n$ and $D_2 = n^2$

The solution to this equation is $D_k = -\dfrac{n\left((-1)^k-(1+n)^k\right)}{2+n}$.

Therefore the periods of the D[k] for N = 2n will satisfy this closed-form equation.

(iii) The ratio of the periods will be $D_k/D_{k-1} = \dfrac{(-1)^k-(1+n)^k}{(-1)^{1+k}-(1+n)^{-1+k}}$

This sequence of ratios will clearly approach n + 1 = N/2 + 1 in the limit. It will begin with n and exceed n+1 on the second iteration and then alternate low-high relative to the limit.

(iv) Since 'most' M[k] form on the edges of the D[k-1] they will have the same limiting ratios and also satisfy the same basic difference equation as the D[k], but the new initial conditions will be $M_1$ =period of M[1] = n and $M_2$ =period of M[2] = n(3(n-1)/2 + 2 ), (This later period can be explained by noting that as shown below, 'most' D[k] arise in groups of 2 with 3 M[k] in each group plus the 2 M[k] of the single D[k]). The solution is

$$M_k = \frac{n\left((-1)^{1+k}+(-1)^k n+3(1+n)^k\right)}{2(2+n)}$$ (period of the M[k] for n = 2N)

(v) Each of the predicted DS[k] in the Edge Conjecture will also survive in each generation.

**Note**:Part (v) implies that the entire region local to S[2] should share the geometric scaling of the M[k] and D[k] tiles, but the remaining DS[k] may have very different temporal scaling. Working back from D[2] toward star[1] of S[2] we know both the geometric and temporal scaling so the local fractal dimension should be Log[N/2+1]/Log[1/GenScale[N/2] for N = 8k+2. This begins with 1.2411 for N = 10 and approaches ½ with N. This does not imply that the subsequent generations will be self-similar and evidence says that typically they are not.

N = 26 has algebraic complexity 6 and it is the first 8k+2 case which is neither quadratic nor cubic in complexity. We will use it to help explain the difference equations for the evolution of the D[k] and M[k]. The Edge Conjecture predicts that S[1] will be DS[11] and there will also be DS[7]s and DS[3]s as shown below. The solid blue center lines track the local web symmetry and this also applies to the First Family as shown by the S[2] line of symmetry.

**Figure 2.11** The geometry of N = 26 showing 'clusters of D[2] tiles anchored by S[3]s (not shown). The blue lines are lines of symmetry

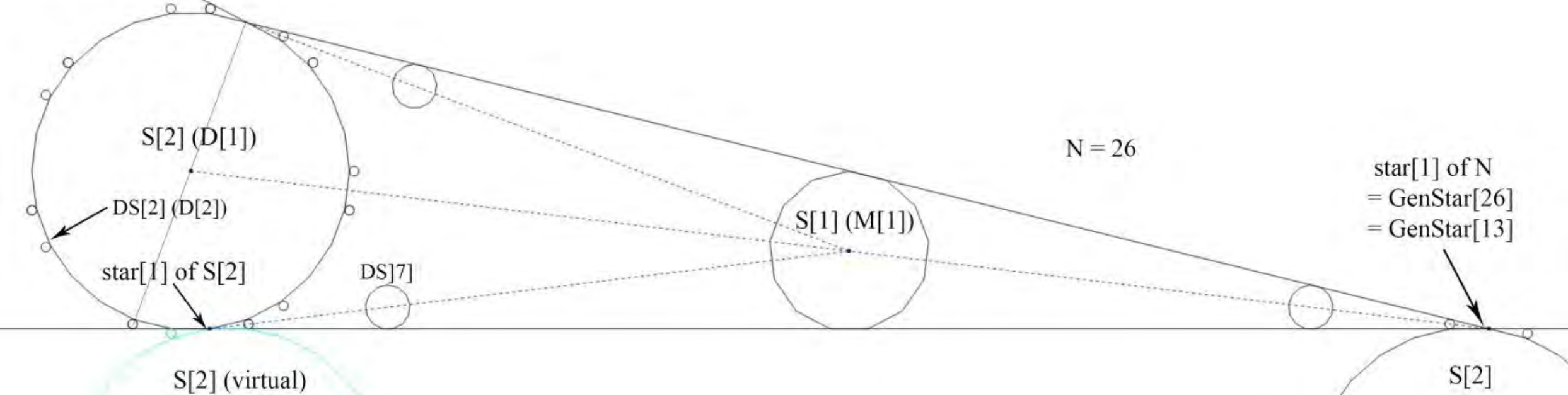

Because the combined web of S[1] and S[2] is step-4 there will be k (from 8k+2) of these S[3] 'clusters' on either side of this line of symmetry. Each cluster has 2 D[2]s and there a lone D[2] on the line of symmetry for a total of 4k+1 D[2]s. This is the 'n = N/2' in the difference equation. Since the period of any S[k] is N/gcd(k,N), S[2] has period n also. Therefore the period of the D[2]s is $n^2$. Every member of the 8k+2 family will have D[1] at period n and D[2] at period $n^2$. These are the initial conditions of the p[k] equation. (Note that the S[2] local to N is in a different group of S[2]s so its small D[2] is not counted here. However when D[1] is replaced with D[2] to count the next generation D[3]s, these 'outliers' must be counted.)

**Figure 2.12** At D[1] there will be 13 magenta D[3] around each D[2] - each with their outliers. So there will be 14 outlier D[3]s which are just visible here as magenta dots.

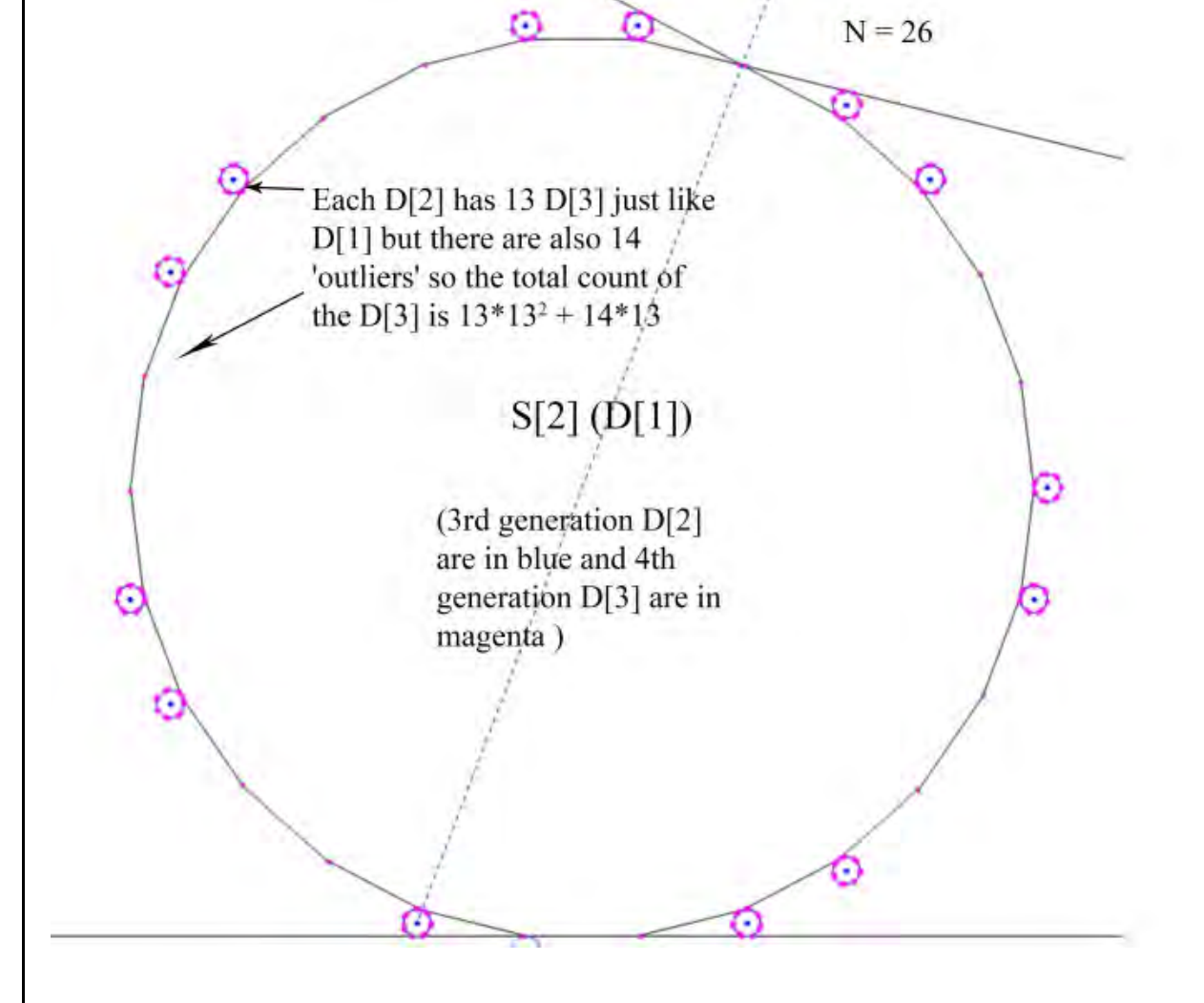

The periods of the first three D[k] are:
D[1] at 13, D[2] at $13^2$ and D[3] at $13^3 + 14\cdot 13$
so the proposed difference equation for the periods $D_k$ of the D[k] is

$$D_k = nD_{k-1} + (n+1)D_{k-2}$$

where n = N/2 and $D_1 = n$ and $D_2 = n^2$
Mathematica's solution is:

**RSolve**[{D[k]= D[k-1]n + D[k-2](n+1), D[1]=n, D[2]=n^2}, D[k], k]

$$D[k] \to -\frac{n\left((-1)^k - (1+n)^k\right)}{2+n}$$

as in the 8k+2 Conjecture

Therefore: $D[n,k] = -\dfrac{n\left((-1)^k - (1+n)^k\right)}{2+n}$ gives the $\tau$-period of each D[k] for N = 2n

**Table**[D[13,k], {k,1,8}] = {13, 169, 2379, 33293, 466115, 6525597, 91358371, 1279017181}

For dynamics it is necessary to adopt a global perspective as shown below for N = 10. Here it is clear that D[1] and D[2] have periods of 5 and $5^2$ and then D[3] must have period $155 = 5{\cdot}5^2 + 6{\cdot}5$. This applies equally to the blue or magenta regions, but our default region is blue.

**Figure 2.13** N = 10 showing the canonical 'decomposition' of the orbit of S[2] for N-even.

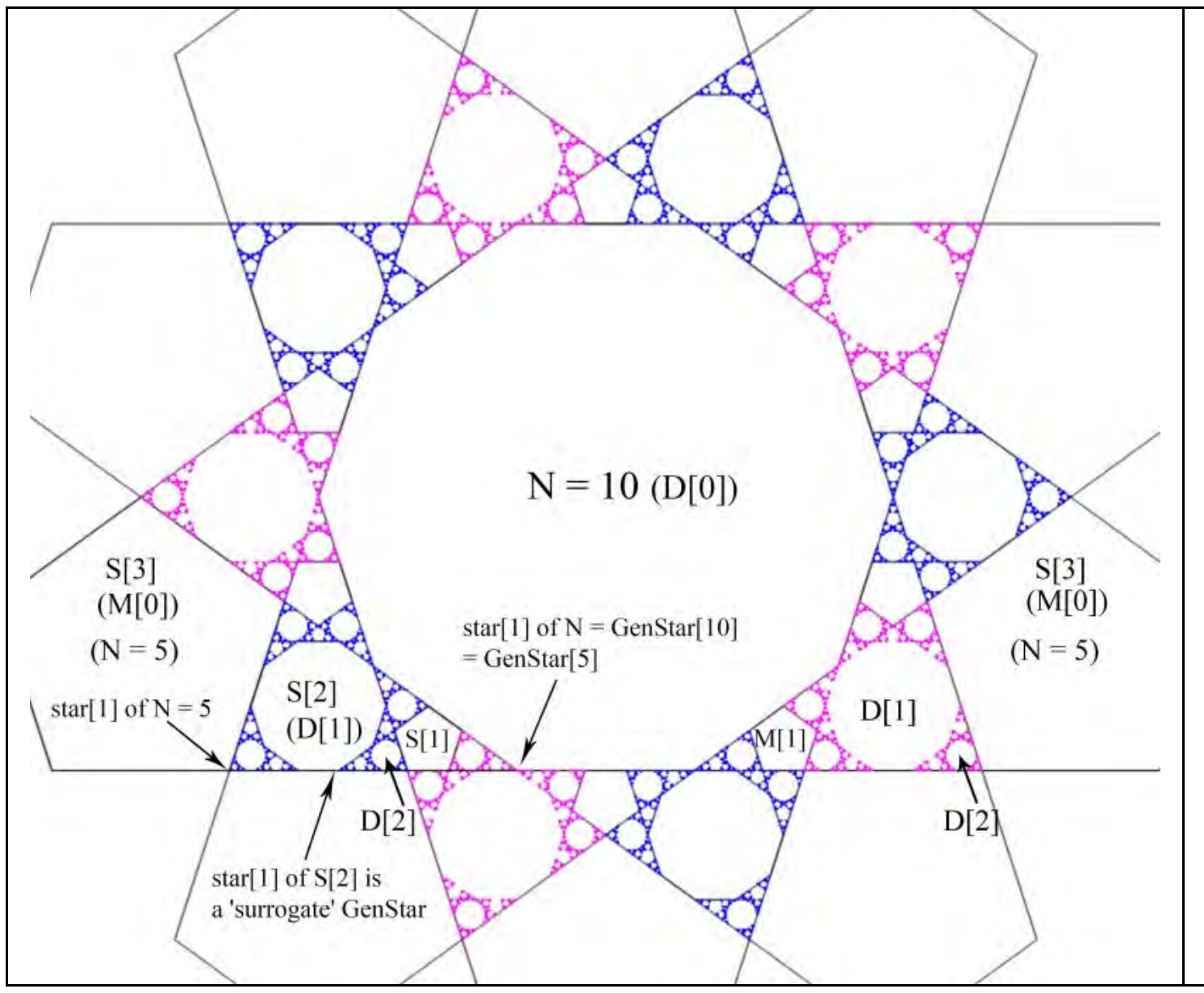

Any twice-odd pair like N = 26 and N = 13 will share the same web, but the dynamics of the embedded N = 13 may be very different from the dynamics local to N. Here with N = 10 there is only one non-trivial scale so N = 5 and N = 10 share the same geometric scaling and since this geometry is self-similar, they also share the same temporal scaling. This implies that N = 10 and N = 5 must share the same difference equations for the D[k] and M[k] (but possibly with distinct initial conditions). The resulting four equations are given in Section 5 of [H5] and will be reproduced below.

This invariance occurs at all scales and these regions are typically 'nested' as shown below.

**Figure 2.14** For N = 26, the 13 'even' and 'odd' S[2] still map to other but only locally. The blue and magenta regions below are the orbits of points close to S[2] and S[1] respectively. On close inspection it is clear that the 'even' S[2] will have magenta D[k] and M[k].

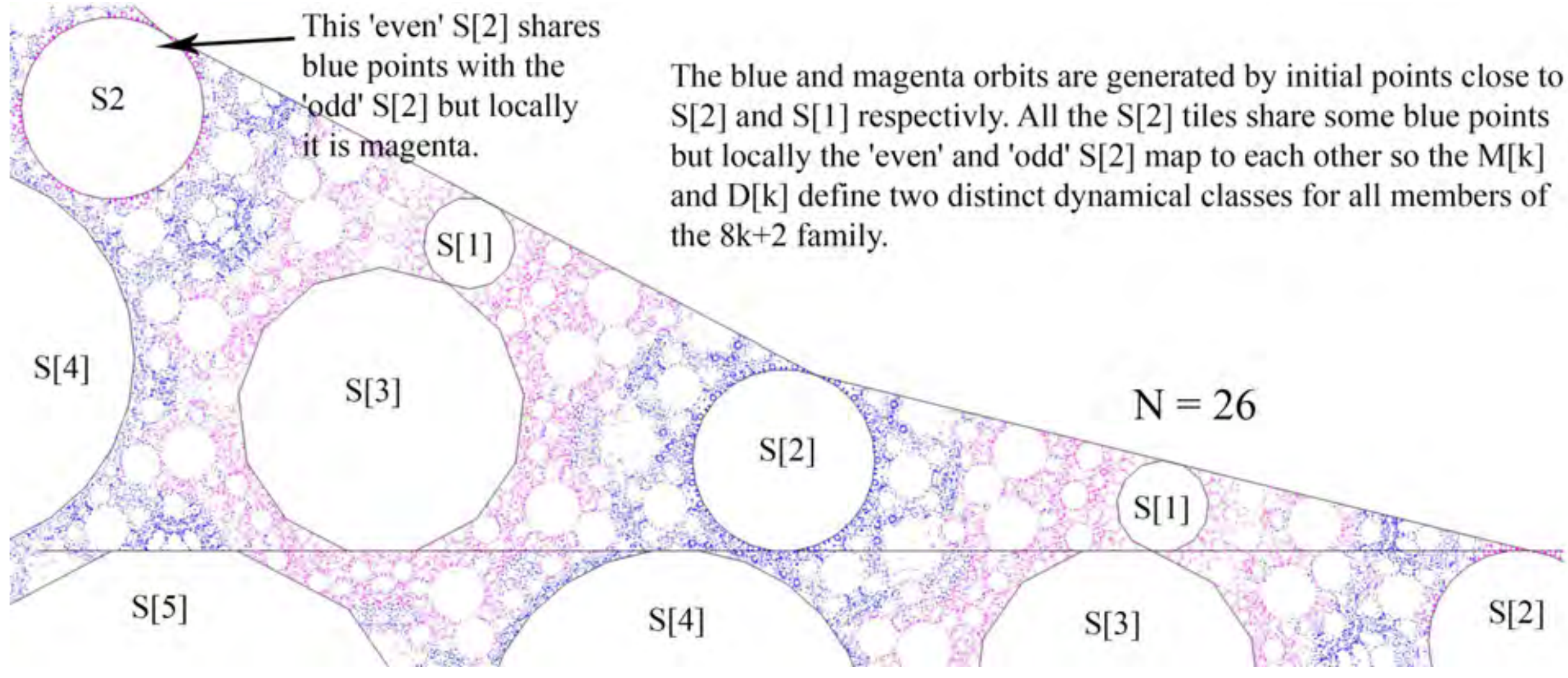

**Towers of S[k] Tiles**

As noted earlier every N-gon has an invariant local region that includes the S[1] and S[2] tiles. All the tiles in this region would be expected to share similar dynamics and neighboring tiles like S[3] and S[4] can have a significant effect on the geometry. Because of rotational symmetry, rotated copies of S[k] will share this same geometry. Here we look at these rotations and the resulting vertical 'towers' of S[k] which also share star points. Once again the odd cases have shared indices that are double the even case.

**Figure 2.15** (N even) S[k] and S[k+2] will share star[k+1] as shown here for N = 20.

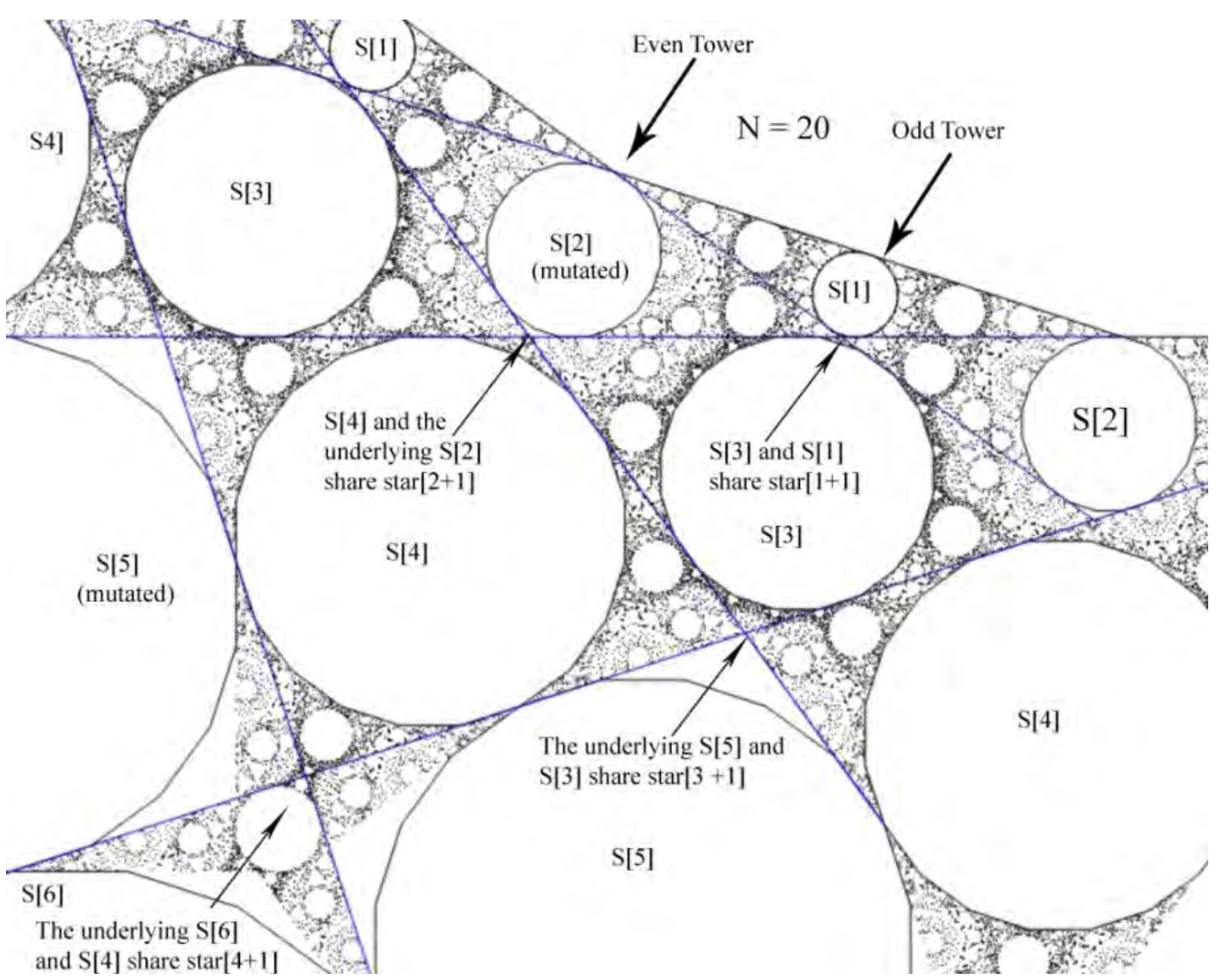

The twice-odd case is the same when the mutated odd Sk] are replaced with underlying S[k].

**Figure 2.16** (N odd) S[k] and S[k+2] will share star[2k+2] as shown for N = 17

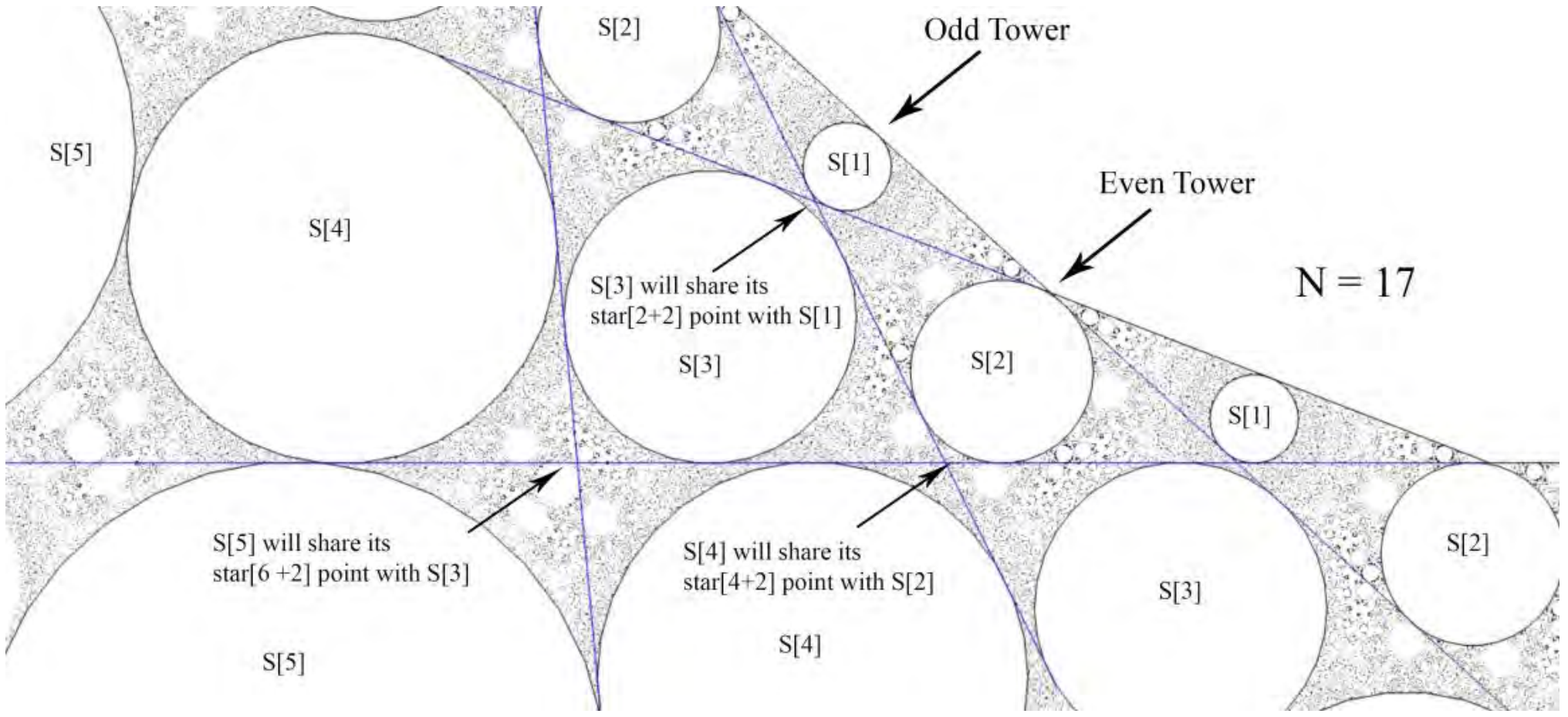

**Summary of Edge Geometry**

**Figure 2.17** Classification of Regular N-gons

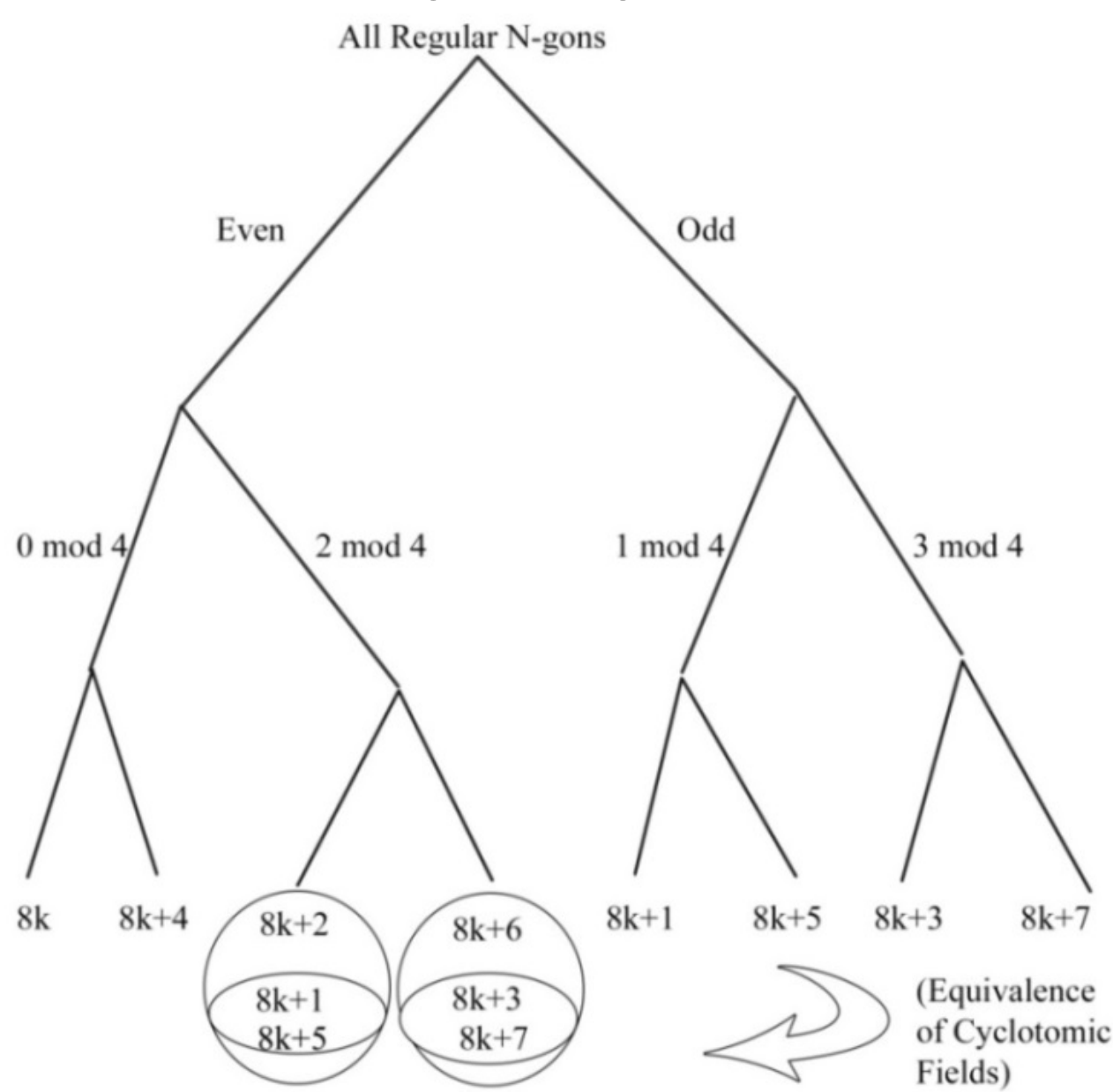

These embeddings are of theoretical interest but at least for the purposes of edge geometry, we will regard the odd cases as relatively independent of the matching twice-odd case. This seems prudent until we fully understand how they are related. Their geometric relationship goes back to the First Family Theorem of [H5] where we show that in terms of the S[k], the N-odd case is a simple doubling of the $s_k$ trigonometric parameters for N-even. This allows for the unification of First Families of N-odd and 2N but locally they evolve in a very different fashion.

The mutations in the S[k] and DS[k] follow directly from the $k'$ web steps and here we expand the original **Mutation Conjecture** of [H5] to include the DS[k] tiles in the second generation. As explained above, when N is even, the S[k] of N (or D) and the S[k] of S[2] (which we call DS[k]) appear to have the same $k' = N/2\text{-}k$ web steps, so mutations in the DS[k] should be matched by equivalent mutations in the S[k] of N (or D acting as N for N twice-odd). So for N = 30, D and S[2] should share mutations in the predicted DS[5] and DS[9] (and a volunteer DS[3]).

For N twice-odd the DS[k] will have odd k values so they will be N/2-gons which can be regarded as the 'underlying' S[k] or DS[k]. Therefore mutations should occur whenever $\gcd(N/2,k') > 1$ and the mutation will be the weave of two $(N/2)/\gcd(N/2,k')$ – gons with base spanning $\gcd(N/2,k')$ star points of DS[k]. For N twice-even the underlying S[k] and DS[k] are all N-gons, so the S[k] of N and the DS[k] will be mutated when $\gcd(k',N) > 2$ and the mutation will be the weave of two $N/(\gcd(k',N)$ – gons with base spanning $\gcd(k',N)$ star points of S[k] or DS[k]. In all N-even cases, the initial base star point should be the minimum of $N/2\text{-}1\text{-}jk'$ with respect to the true underlying N-gon S[k] or DS[k]. For N-odd, S[2] is twice-odd so the web steps are $(2N/2)\text{-}k = N\text{-}k$ , but the DS[k] are now oversized 2N-gons so the DS[k] mutations are relative to these underlying 2N-gons with mutation criteria $\gcd(2N,k') >2$ . The initial base star point for mutations will now be the minimum of $N\text{-}2\text{-}jk'$ . See N = 15.

## Section 3. Catalog of singularity set edge geometry for regular N-gons with N $\leq$ 25

Below are examples of the 8 possible edge cases based on the Rule of 4 for N even and the Rule of 8 for N odd. This is the '8 fold-way' for regular polygons.

**Table 3.1** The 8 dynamical classes of edge geometry for N even (top) and N odd (bottom)

| 8k family (8,16,24,32,40,48) | 8k + 2 family (10,18,26,34,42,50) | 8k + 4 family (4,12,20,28,36,44) | 8k + 6 family (6,14,22,30,38,46) |
|---|---|---|---|
| DS[2] N = 16 DS[6] S[2] S[1] | DS[3] N = 18 DS[7] S[2] S[1] | DS[4] N = 20 DS[8] S[2] S[1] | DS[1] DS[5] N = 22 DS[9] S[2] S[1] |
| **8k + 1 family (9,17,25,33,41,49)** | **8k + 3 family (3,11,19,27,35,43)** | **8k+ 5 family (5,13,21,29,37,45)** | **8k+7 family (7,15,23,31,29,47)** |
| DS[5] N = 17 DS[13] S[2] S[1] | DS[7] N = 19 DS[15] S[2] S[1] | DS[1] DS[9] N = 21 DS[17] S[2] S[1] | DS[3] DS[11] N = 23 DS[19] S[2] S[1] |

As illustrated above, these even and odd cases are linked because any odd N-gon can be regarded as embedded in the 2N-case. Therefore the 8k+1 and 8k+5 families can be embedded in 8k+2 families while the 8k+3 and 8k+7 families can be embedded in 8k+6. This embedding of an odd N-gon in 2N is algebraically sound (with the proper transformation), but the dynamics are different and the web local to N may be very different from the web local to D or 2N.

The **8k+2 Conjecture** described above says that these N-gons will have an edge geometry driven by sequences of self-similar D[k] and M[k] tiles with known geometric and temporal scaling, The **8k+4 Conjecture** says that the geometry is driven by the mutation of S[2]. The **8k+1 Conjecture** says that these families will have a volunteer DS[2] to go along with the predicted DS[5]. The **8k+7 Conjecture** says that the predicted DS[3] will generate dual DS[1]s at S[N-3] with step-2 webs which support at least S[1] tiles. The **Twice-even S[1] Conjecture** says that since S[1] has a step-2 web it can support 'step-2'tiles called Skx which are D tiles relative to S[k] (like N-odd).These Skx include S[2] and S3x is an S[2] tile of S[3] for N$\geq$ 12.

Every N-gon has a local web which is invariant and this web would be expected to contain at least 1/4 of the S[k], so there is a link between edge geometry and the large scale geometry. Both are driven by the cyclotomic field and the corresponding scaling field $S_N$ with complexity φ(N)/2. Hopefully the examples below may shed some light on the issue of 'nature' (algebraic complexity) vs. 'nurture' (web and edge complexity under $\tau$ ). As N increases there appears to be a surprising amount of diversity within the 'algebraic families' shown below.

**Table 3.2** Algebraic Complexity of regular N-gons for N $\leq$ 50

| $\phi$(N)/2 | 1 | 2 | 3 | 4 | 5 | 6 | 8 | 9 | 10 | 11 | 12 | 14 | 15 | 18 | 20 | 21 | 23 |
|---|---|---|---|---|---|---|---|---|---|---|---|---|---|---|---|---|---|
| N | 3<br>4<br>6 | 5<br>8<br>10<br>12 | 7<br>9<br>14<br>18 | 15<br>16<br>20<br>24<br>30 | 11<br>22 | 13<br>21<br>26<br>28<br>36<br>42 | 17<br>32<br>34<br>40<br>48 | 19<br>27<br>38 | 25<br>33<br>44<br>50 | 23<br>46 | 35<br>39<br>45 | 29 | 31 | 37 | 41 | 43<br>49 | 47 |

**•N = 3**

| | |
|---|---|
| 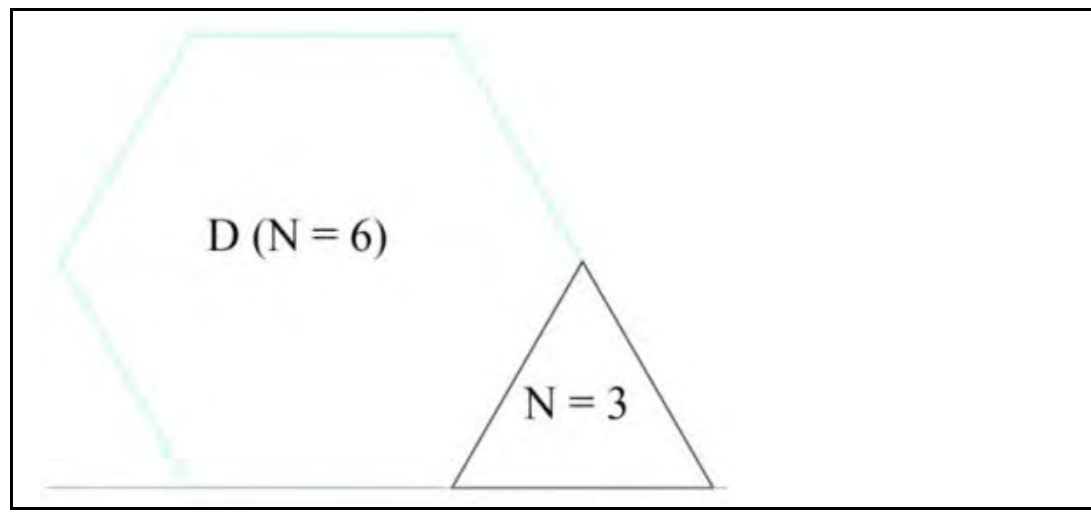  | The scaling fields $S_3 = S_6$ are generated by x = GenScale[3] = Tan[π/3]·Tan[π/6] = 1. By the Twice-Odd Lemma there is a natural equivalence of the webs for N = 3 and N = 6. |

**•N = 4**

| | |
|---|---|
| 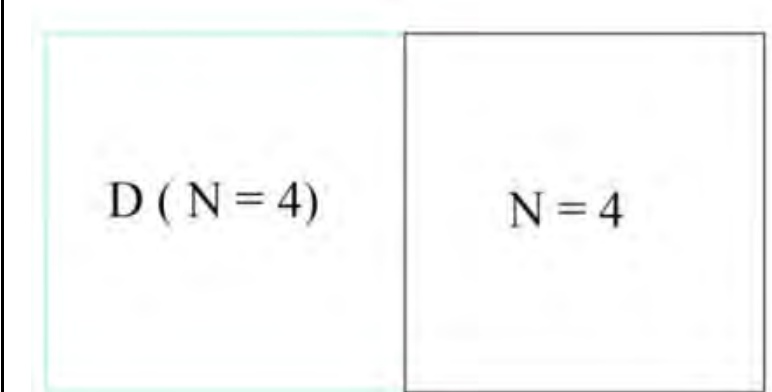  | The scaling field $S_4$ is generated by x = GenScale[4] = Tan[π/4]$^2$ = 1 |

**• N = 5**

N = 5 is the first 'quadratic' polygon and the charter member of the 8k+5 family. The only non-trivial scale is scale[2] = GenScale[5] = Tan[π/5]·Tan[π/10] = $\sqrt{5} - 2$. The Twice-odd Lemma implies that N = 5 and N = 10 share the same web so the five D tiles shown below are surrogate N = 10s where the 8k+2 Conjecture predicts that there will be well-defined sequences of M[k] and D[k|] tiles converging to GenStar[5], which is star[1] of D. The 8k+2 Conjecture predicts that this sequence will have temporal scaling 6. Assuming that this joint web is self-similar this same temporal scaling should apply to the sequence of M[k] and D[k] converging to star[1] of N and here we repeat the 'renormalizing' argument from [T] and [H3] to support that hypothesis.

**Figure 5.1** N = 5 and N = 10 share the same web geometry but different dynamics. Here we derive the temporal scaling at star[1] of N = 5 using the self-similar blue 'darts' below.

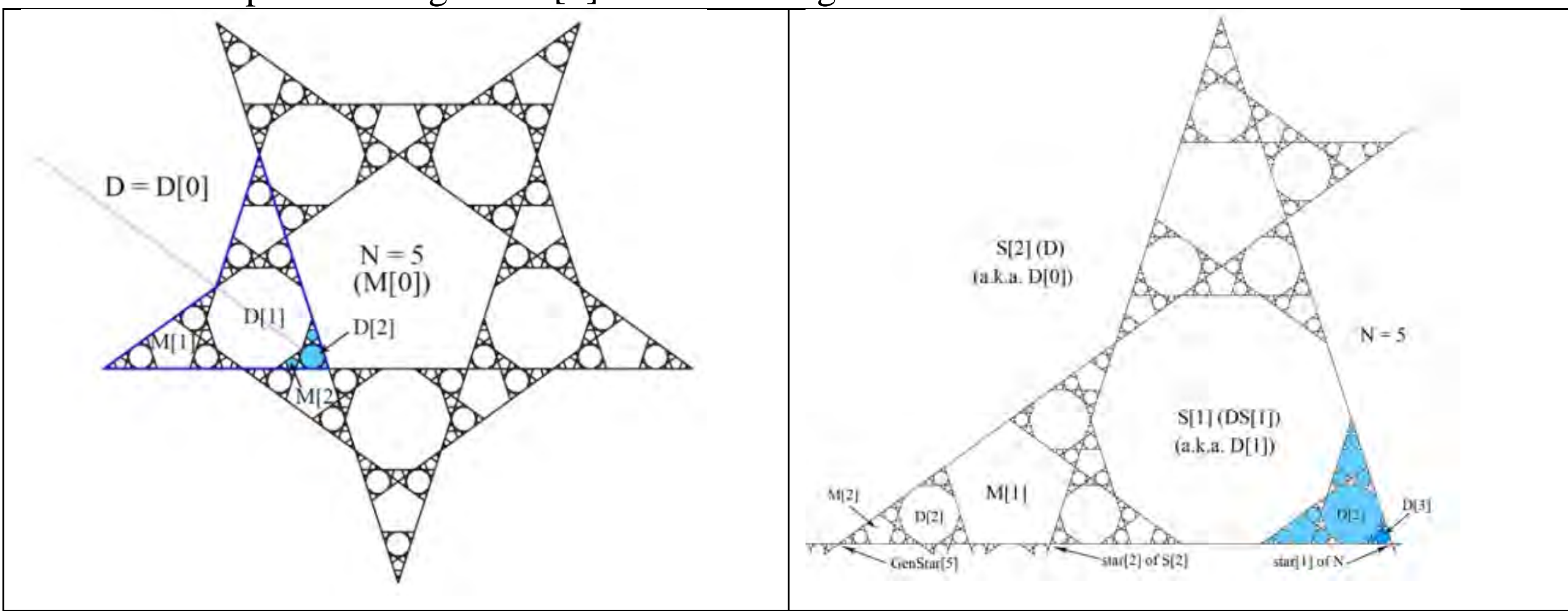

The blue 'darts' shown here are anchored by D[k] tiles so they scale geometrically by GenScale[5] and we want to show that their temporal scaling is 6. The light-blue and dark-blue darts are the 2nd and 3rd 'generations'. Note that 5 copies of the 1st generation dart will almost tile the full generation '0' but tiling the 1st generation dart with light blue darts will take at least 7 copies and it will take at least 41 dark-blue darts to cover the previous dart. Since each dart is anchored by a D[k] tile, the limiting tiling is equivalent to determining the growth of the D[k]. In each dart. D[k] is surrounded by 2 M[k] and in the limit each M[k] accounts for 3 next-generation D[k+1], so the D[k] should have a temporal growth of 6.

In [H3] we derived the following difference equations to describe this growth of decagons and pentagons: $d_n = 3d_{n-1} + 2p_{n-1}$ ; $p_n = 6d_{n-1} + 2p_{n-1}$. This yields the table below for one blue dart.

| Generation | decagons - $d_n$ | pentagons - $p_n$ |
|---|---|---|
| 1 | $d_1 = 1$ (D[1]) | $p_1 = 2$ (M[1])s |
| 2 | $d_2 = 3d_1 + 2p_1 = 7$ | $p_2 = 6d_1 + 2p_1 = 10$ |
| 3 | $d_3 = 3d_2 + 2p_2 = 41$ | $p_3 = 6d_2 + 2p_2 = 62$ |
| n | $d_n = 3d_{n-1} + 2p_{n-1}$ | $p_n = 6d_{n-1} + 2p_{n-1}$ |

Eliminating $p_n$ gives the second-order difference equation: $d_n = 5d_{n-1} + 6d_{n-2}$ . This equation must also describe the growth of the D[k] for N =10 as in the 8k+2 Conjecture, and indeed it is identical to $D_k = nD_{k-1} + (n+1)D_{k-2}$ where $n = N/2$. It is a simple linear second-order homogeneous difference equation which can be solved by assuming that $D_k = C_1 z_1^k + C_2 z_2^k$ where $z_1$ and $z_2$ depend on n and n+1 and $C_1$ and $C_2$ depend on the initial conditions.

For N = 10 (and all 8k+2 cases) the initial conditions are simply $D_1 = N/2$ and $D_2 = (N/2)^2$ because these are the periods of D[1] and D[2]. For N = 5 the initial conditions for the table above with one dart are $D_1 = 1$ and $D_2 = 7$ and the full star region will have $D_1 = 5$ and $D_2 = 5(7)$. For this default star region, Mathematic gives the following (simplified) solutions:

$$D[n,k] = \frac{n\left((-1)^k + (-1)^k n + 3(1+n)^k + n(1+n)^k\right)}{(1+n)(2+n)}$$

For the M[k] pentagons the initial conditions in this same region are 2N and $2N^2$ to yield

$$M[n,k] = -\frac{2n\left((-1)^k - (n+1)^k\right)}{2+n}$$

Therefore in both cases $z_1 = -1$ and $z_2 = 1+n$ and this explains why in both cases the ratio must approach 1+n in an alternating fashion. The first few periods are given below

**Table**[D[5,k], {k,1,9}] ={5, 35, 205,1235, 7405, 44435, 266605, 1599635, 9597805}
**Table**[M[5,k], {k,1,9}] ={10, 50, 310, 1850, 11110, 66650, 399910, 2399450, 14396710}

The first few $D_k$ ratios are : {35/5, 205/35, 1235/205, 7405/1235..} = {7, 5.857, 6.024, 5.9959}

In Section 4 of [H3] (2013) we gave this same difference equation for N = 5, but these results and others were obtained years earlier and communicated to Richard Schwartz at Brown. He passed them on to his colleague S.Tabachnikov who had already studied this case and in the following years we exchanged e-mails as their interests evolved to other areas. We did not realize that N = 5 was just a special case of a family of 8k+2 polygons based on N = 10. For us personally it was a simple case of 'bias' in favor of the beauty and elegance of primes like N = 5. What can N = 10 offer that is not already known for N = 5 ? Now there are likely an infinite number of 'generalized' N = 10s but we know almost nothing about the matching N/2-gons. This seems to be a very gender-specific world, especially when dynamics are involved.

•N = 6

<table>
<tr>
<td>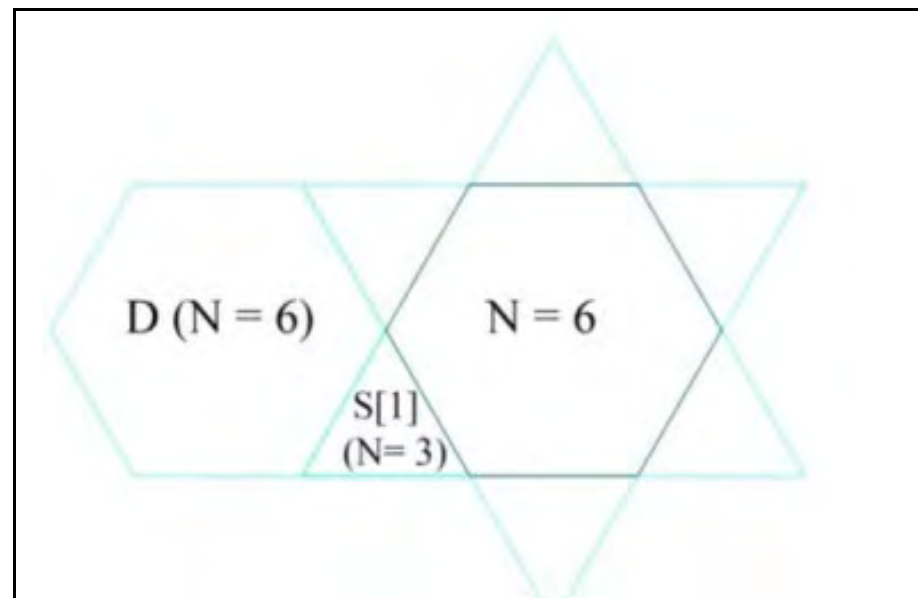
</td>
<td>The scaling fields $S_6 = S_3$ are generated by x = GenScale[3] = $\text{Tan}[\pi/3]\cdot\text{Tan}[\pi/6] = 1$. If the origin is shifted to the center of S[1], this web will be identical to the web for N = 3.</td>
</tr>
</table>

•N = 7

N = 7 and N = 9 are 'cubic' polygons so they have a second non-trivial primitive scale along with the minimal scale[<N/2>] which we call GenScale[N] (or just GenScale when N is understood). Recall that <N/2> is Floor[N/2] for N odd and N/2-1 for N even, so for N = 7, GenScale = scale[3] = $\text{Tan}[\pi/7]/\text{Tan}[3\pi/7] = \text{Tan}[\pi/7]\cdot\text{Tan}[\pi/14] \approx 0.109914$. For both N = 7 and N = 9, the competing primitive scale is scale[2] which is also hS[1]/hS[2]. These two scales are noncommensurate and their interaction is the major source of complex dynamics beyond the 'single-scale' fractal structure of N = 5, 10, 8 and 12.

In the twice-odd family, N and 2N share the same cyclotomic field and the same web, but their local geometry may be very different and we would like to know how their geometries are related. For large N values the separation between N and D can make it difficult to find relationships, but N = 7 and N = 14 form a very compact joint family where there are natural family connections between N = 7 and D. Note the edge sharing of DS[3] and DS[4] leading to the canonical invariant regions. The overview shown below shows these invarient regions as defined by N and D. DS[4] (S[2]) is the only tile that is in the First Families of both D and N and it plays a critical role on the boundary between the regions. Here we will concentrate on the inner region local to N. This is the only case where D is the S[3] of N.

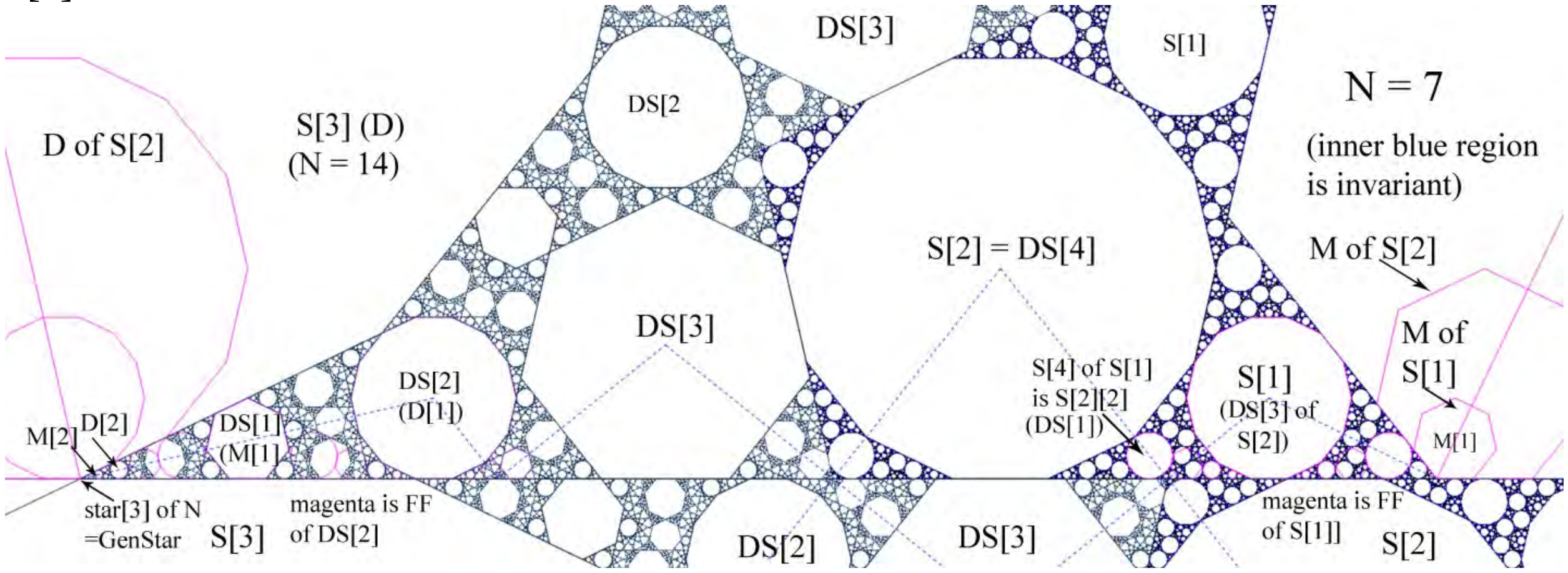

N = 7 can be regarded as the charter member of the 8k+7 family. This is a very interesting family and the 8k+7 Conjecture makes a number of predictions. First it says that the DS[3] predicted by the Edge Conjecture will evolve along with dual DS[1]s on the edges of S[2] so together they span 4 edges of S[2]. These DS[1]s will be S[2]s of S[2] so they can play the role of D[2]s and foster 'next-generation' families. Because N is odd the normal S[k] in the First Family of S[2] are replaced by the DS[k] 'star[2]' families which are displaced relative to the S[k] and hence scaled up by 1/scale[2] of 2N. So DS[1] is a 'D' tile relative to the S[1] of S[2] and is a perfect match for an S[2] of S[2]. The S[2] tile will act as a local D tile and here the DS[3] tile of S[2] is identical to the S[1] of N. We apologize for the possible confusion of this D with the D of N.

The 8k+7 Conjecture also says that the DS[1]s will be S[N-3] tiles of DS[3], so here they are S[4]s of S[1].The evolution of these DS[1]s will be an important issue throughout the 8k+7 family. Since S[2] is a 2N-gon, the DS[k] will be expected to evolve with $k' = 2N/2\text{-}k = N\text{-}k$. Since Mod[2N,N-1] is 2 for N odd, DS[1] will be expected to have a step-2 limiting web.
For DS[3],Mod[2N,N-3] is typically 6, but N = 7 is an exception since DS[3] is S[1] with an expected $k' = 2k + 2 = 4$ which is consistent with the fact that S[1] is DS[N-4] for N odd. In the graphic below note that the DS[1] tiles are step-4 relative to S[1].

**Figure 7.3** The local web showing lines of convergence to star[1] of N and star[2] of S[2]

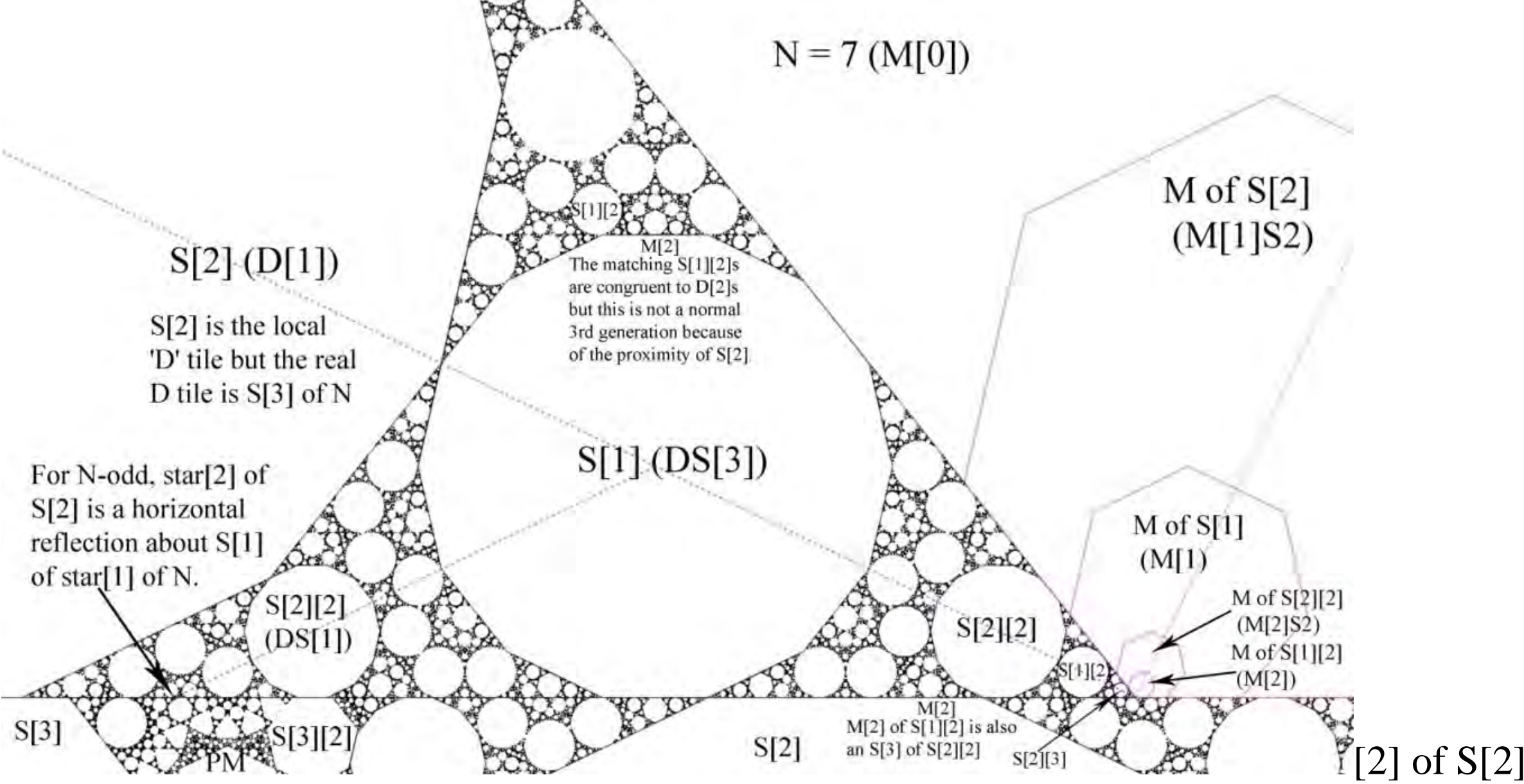

Here our interest is primarily the convergence at star[1] of N which matches the convergence at star[2] of S[2]. As shown in Figure 7.2 earlier, there is also a convergent sequence at star[3] of S[2]. We noted above that the DS[1] volunteers can serve as next-generation S[2][2]s which we sometimes call D[2]s. But there are no matching M[2]s to provide a framework for a 3rd generation. By contrast S[1] has both D[2] and M[2] tiles. This occurs because S[1] is congruent to an S[2] tile of D so its edges can support an 'ideal' 3rd generation scaled by $\text{GenScale}^2$ relative to the First Family of N as in Definition 2.3 However these familes on the edges of S[1] are only a modified version of a 'normal' 3rd generation.

This mutated family on the edges of S[1] does support further generations but in an irregular gashion using both GenScale and scale[2]. The S[2]-S[1] sequence has the same issues and it is no surprise that the geometry at star[1] of N is replicated on the edges of S[1].

The S[2][k] sequence ends temporarily with the tiny S[2][3] shown above, because there is no matching S[1][3]. This S[1][3] should occur as an edge tile of the (virtual) M[2] shown in magenta. This breakdown seems to result from the fact that the M[2] tiles have two distince dynamical roles – as M tiles of the S[1][k] and S[3] tiles of the S[2][k]. This is illustrated above with the M[2] adjacent to the right-side S[2][2].

Geometrically these two roles of M[2] are compatible since hM[2]/hN = $\text{GenScale}^2$ and this is identical to an S[3] of S[2][2] which has hS[3]/hN = $(\text{GenScale}^2/\text{scale}[2])\cdot\text{scale}[2]$. But the scaling occurs sequentially and dynamically these roles are very different , so the web evolution appears to have a form of 'sensitive dependence on initial conditions' where the progeny of M[2] feel conflicting influence from the local webs of S[2] and S[1] with their different scaling. Therefore the S[2]-S[1] sequence breaks down at S[1][3], but it eventually continues with an S[2][5] playing the role of the original S[2]. Below is an enlargement of this star[1] region.

**Figure 7.4** The star[1] region showing the break-down of the S[2][k]-S[1][k] at S[1][3].

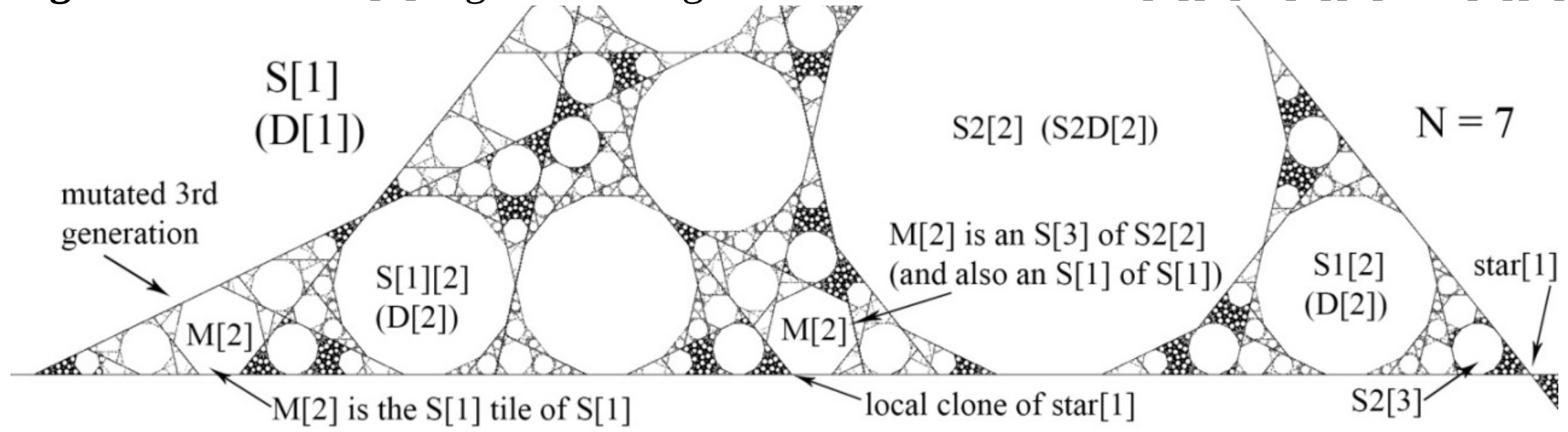

It is clear that the location where S[1][3] should form has a complex geometry of 'dark matter'. This same issue occurs with the two real M[2]s above, because the geometry at star[2] of S[1] is just a reflection of the geometry at star[1] of N. Therefore the $3^{rd}$ generations on the edges of S[1] are mutated because of this conflict between scale[2] (represented by S[2]) and scale[3] (represented by S[1]). At D these $3^{rd}$ generation families evolve 'normally' on the edges of S[2] of D acting a D[1], because there is no conflict with scale[2[. Families self-similar to the First Family can arise on any scale so it is possible for S[2] to support such families as shown below with the S[2][3] tiles, but they will be 'scale[2]' families which are not always compatible with the GenScale families at D or S[1].

**Figure 7.5** Two levels of detail of the star[1] region

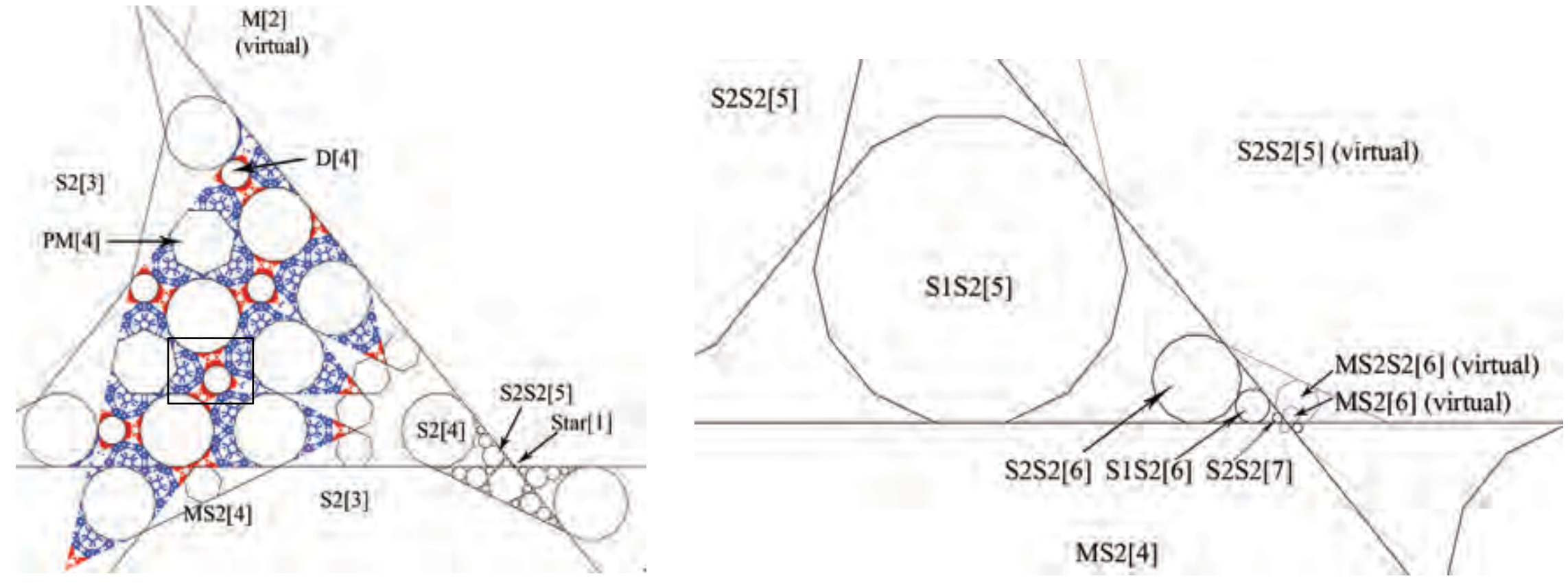

The left side above shows that the break-down after S[2][3] (a.k.a. S2[3]) is localized and there is an S2[4] at star[1] that can be used to continue the sequence. This S2[4] is a D tile in a 5th generation scale[2] family on the edges of S2[3]. The matriarch of this family is an MS2[4] which is an S2-scaled M[4] so hMS2[4] = GenScale$^4$/scale[2].This will be the scaling for each new restart but the next cycle based on S2S2[5] will be only be partially self-similar because MS2[4] shares a vertex with star[1], so there will be no 'subterranean' family as on the left. The continuation depends on the fact that MS2[4] supports a whole arc of S[2] tiles which are what we call S2S2[5]s. As shown above, the virtual S2S2[5] adjacent to star[1] has an S[2] which is the needed S2S2[6] and the virtual M tile of this S2S2[6] will generate S1S2[6] as the (real) S[2]. This virtual M tile has a matching MS2[6] which generates the tiny S2S2[7] . But this MS[2][6] has conflicting dynamical roles like M[2], so it has no S[1] to serve as the S1S2[7].

**Figure 7.6** The second break-down occurs on the edges of MS2[6] but the virtual S1S2[7] inside MS2[4] will generate an S2S2[8] with matriarch MS2S2[8] and scale GenScale$^8$/scale[2]$^2$.

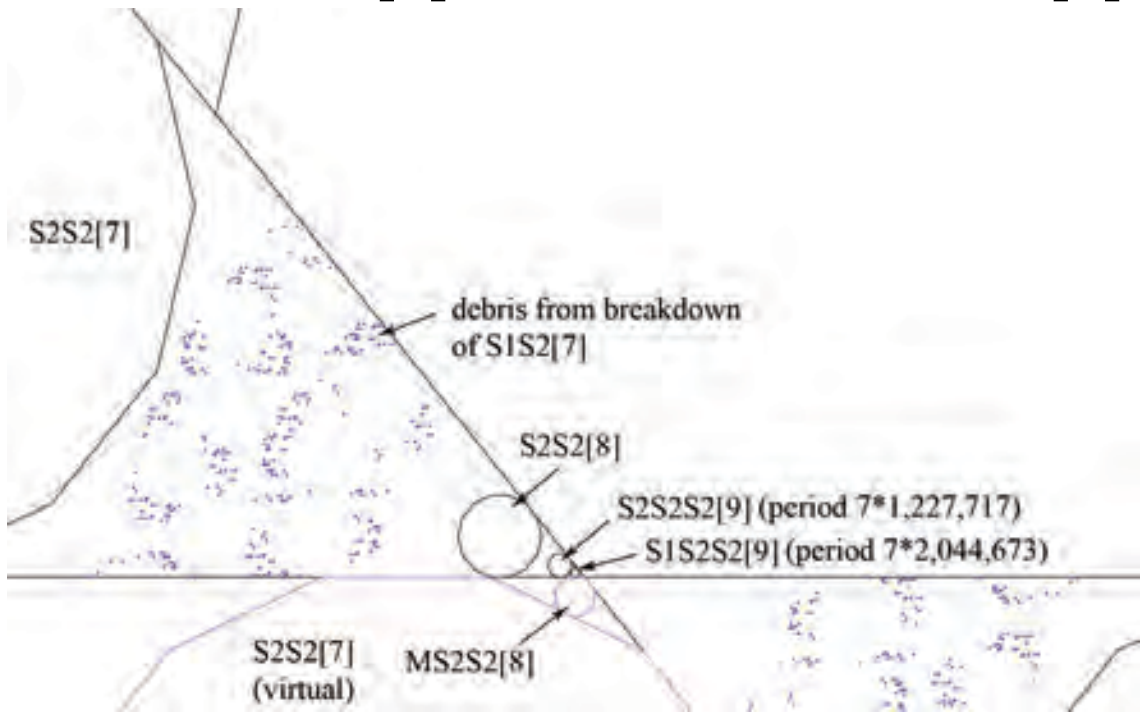

We cannot prove that this sequence continues, but the centers of these S[1][k] and S[2][k] are known and the periods of the existing tiles have been tracked up to the 10th generation where S[1][10] and S[2][10] have $\tau$-periods 509698714 and 63001498 with ratios 1253.8158 and 1253.861 with respect to the 9th generations. This is strong evidence for their mutual survival.

As expected the geometry of these 'chaotic' regions is a mixture of scales. The first such region shown below is presided over by a virtual MS2[4] with a central S1[4] (D[4]) embedded in a 'sea' of larger S2[4]s. This region also has PM[4] tiles which are weakly conforming to S[2][3]. They are modified 'M' tiles with divided allegiance in a manner similar to M[2].

**Figure 7.7** Enlargement of the region outlined in Figure 7.5

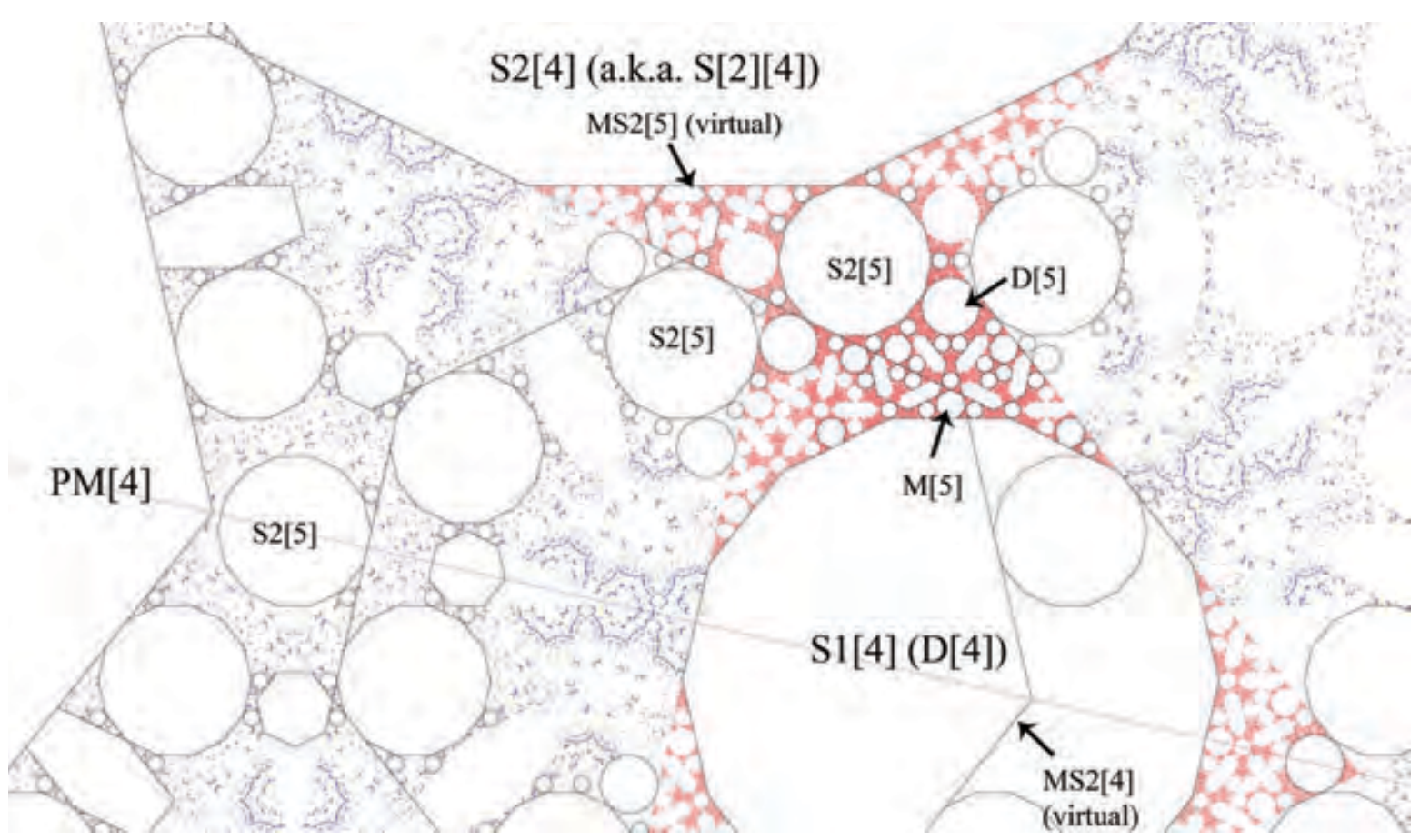

•N = 8

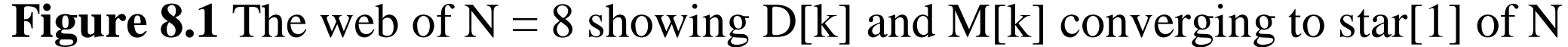

**Figure 8.1** The web of N = 8 showing D[k] and M[k] converging to star[1] of N

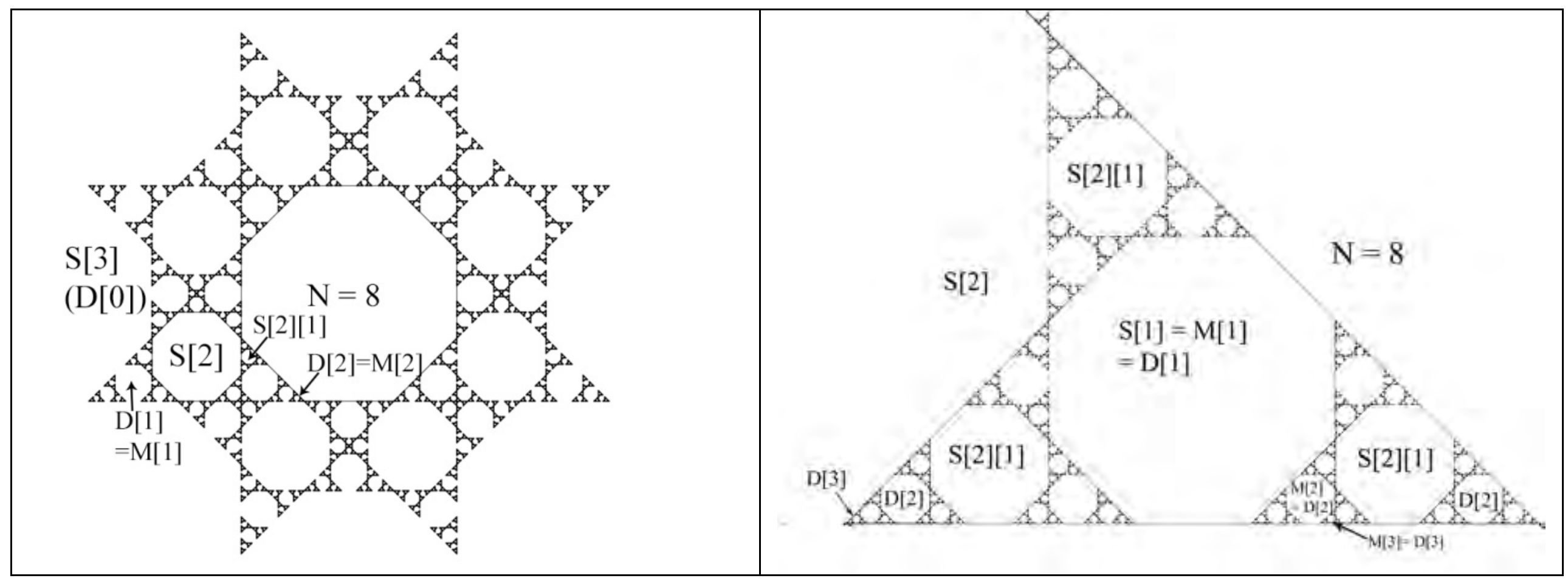

N = 8 has quadratic complexity like N = 5, 10 and 12, so GenScale[N] will be the only non-trivial primitive scale and the web W should be self-similar. N = 8 is the charter member of the 8k family and the Edge Conjecture says that there should be a DS[2] which can possibly act as a D[1] and foster generations converging to star[1] of S[2] and star[1] of N. The problem is that this rarely occurs in the 8k family and N = 8 and N = 16 are the only known cases where D[1] actually fosters next-generations. The main problem is that a 'single' D[1] will usually need a supporting M[1], to form a joint web that mimics the original S[1]-S[2] geometry.

Since N = 8 is twice even there is no distinction between an M[k] tile and a D[k] tile so either one can be used to describe the convergence at star[1] of N or star[1] of S[2]. Our convention is to scale by 'generations' based on GenScale[N] which here is scale[3] = $\mathrm{Tan}^2[\pi/8]$. Any scaling could be used as long it matches the temporal scaling but scale[2] = Tan[π/8] is not primitive.

Unlike the twice-odd case the S[1] of N is now an N-gon, so it can be used as reference for scaling. This S[1] of N is what we call M[1] and the M[k] scale by GenScale[N] = $\mathrm{Tan}[\pi/N]^2$. However M[1] is an edge tile of N and any M[k] sequence converging to star[1] of N will appear 'jagged' relative to the strict vertex scaling of D[k], but here for N = 8 converting to D[k] just means skipping one intermediate vertex, so the vertex path from the D[1] to D[2] goes through S[2][1].

On the left above the primery 'dart' would be anchored by S[2], but M[1] or D[1] can also be used as reference as on the right Here it is clear that D[1] to D[2] has a temporal scale of 9 and this continues for D[2] to D[3] in a self-similar fashion, so the D[k] (and M[k]) scale by 9, and this same limiting scaling will clearly apply to the entire default star region. Therefore this web W will have similarity dimension Log[9]/Log[1/GenScale[8]] ≈ 1.2465. Below we will compare this with N = 5 and 12.

•N = 9

N = 7 and N = 9 are the 'boundary' cases of the 8k+7 and 8k+1 families so it may not be prudent to make generalizations about these families based on these two cases. Concerning N = 9, the Edge Conjecture makes no predictions before DS[5], but the 8k+1 Conjecture says that volunteer DS[2]s will exist with potential for extended family structure. We will see that this is true for N = 9, but there is little evidence to show that this conjecture is the result of the embedding of 8k+1 in the 'well-behaved' 8k+2 family. As explained in Figure 17.3, one possible explanation for the existence of DS[2] is that star[5] is always 'effective' and this yields a local geometry at star[2] of S[2] which can support a DS[2] spanning star[2] to star[4].

**Figure 9.1** The edge geometry of N = 9 showing clusters of next-generation S1, S2 and S3 tiles anchored by blue S[3][2]s. These make up the second iteration of sequences of S[1][k], S[2][k] and S[3][k] converging to star[2] of S[2] or star[1] of N.

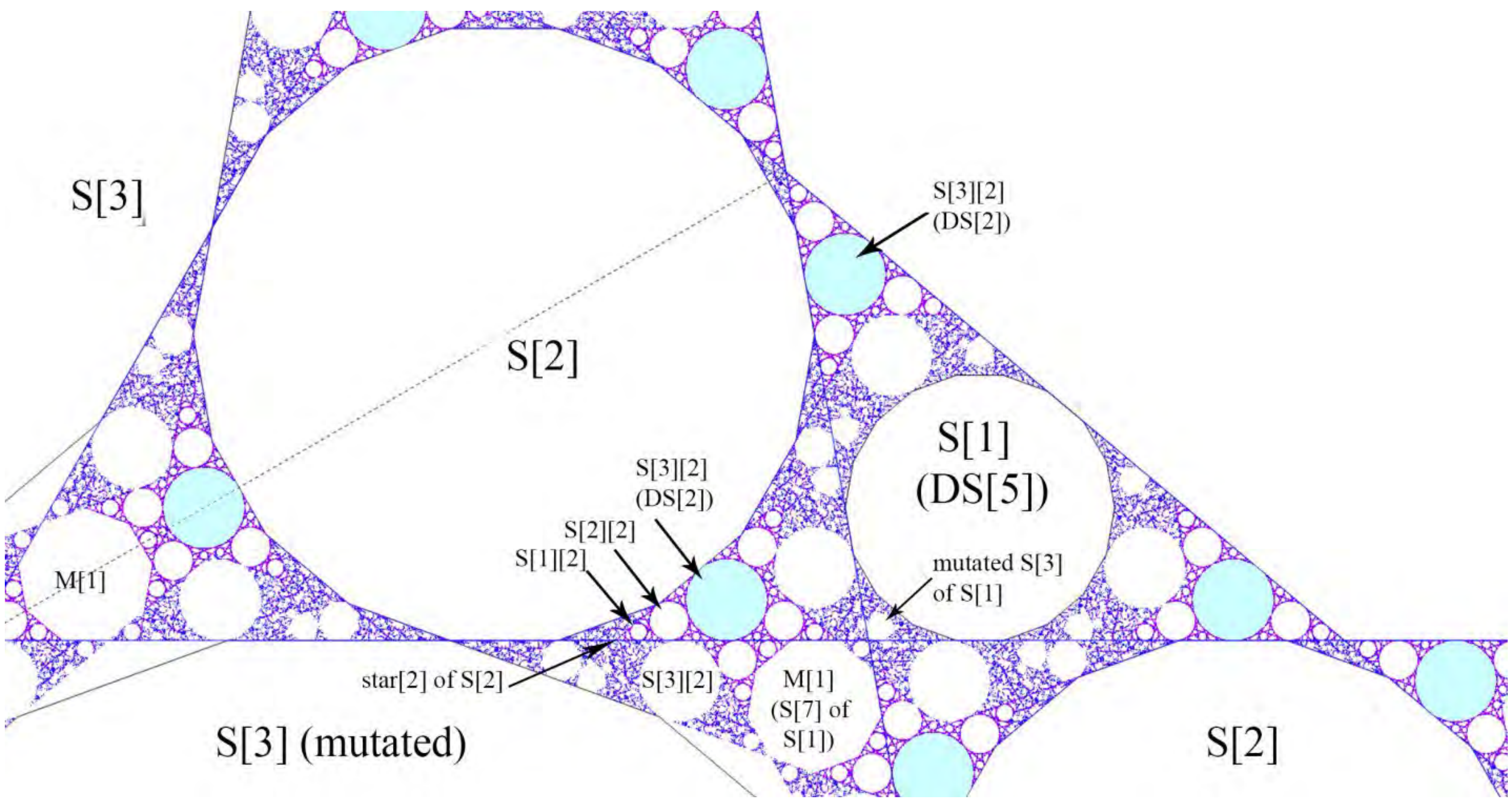

It is clear that the webs of N = 7 and N = 9 have much in common because any cubic N-gon will have two competing primitive scales with the potential for local single-scale behavior where one of the primitive scales is suppressed. N = 7 and N = 9 both have sequences of tiles converging to star[1] of N or the equivalent star[2] of S[2]. For N = 7 these sequences involve both scale[2] and scale[3] (GenScale[7]) in what appears to be a multi-fractal convergence with an unknown spectrum of temporal scaling. For N = 9 the star[1] convergence appears to be much simpler with a limiting temporal scaling of 20. This may be true because this convergence appears to take place inside 'islands' of single-scale self-similar dynamics which characterize N = 9 and 18, and to a lesser degree N = 7 and 14.

For N = 9 these DS[2] are also S[3][2]s as shown above, so they are the foundation for the clusters of S[1][2] and S[2][2] tiles . These DS[2] tiles will exist for other members of these families but only for N = 9 will they also be S[3][2]s and members of the extended family of S[1] which includes the S[k] tiles of M[1] – which is the S[7] tile of S[1]. The mutation of S[3] allows it to support M[1]s at its vertices, and this seems to a stabilizing influence since these M[1] tiles serve an important role. The next level of this sequence is shown in the enlargement below. The local geometry of S[1][2] shows signs of cubic instability, but these issues are localized and do not affect the convergence. This is true for N = 18 also.

**Figure 9.2** – Detail of the star[2] region of S[2] showing sequences of S1[k] ,S2[k] and S[3][k] tiles converging to star[2] of S[2].

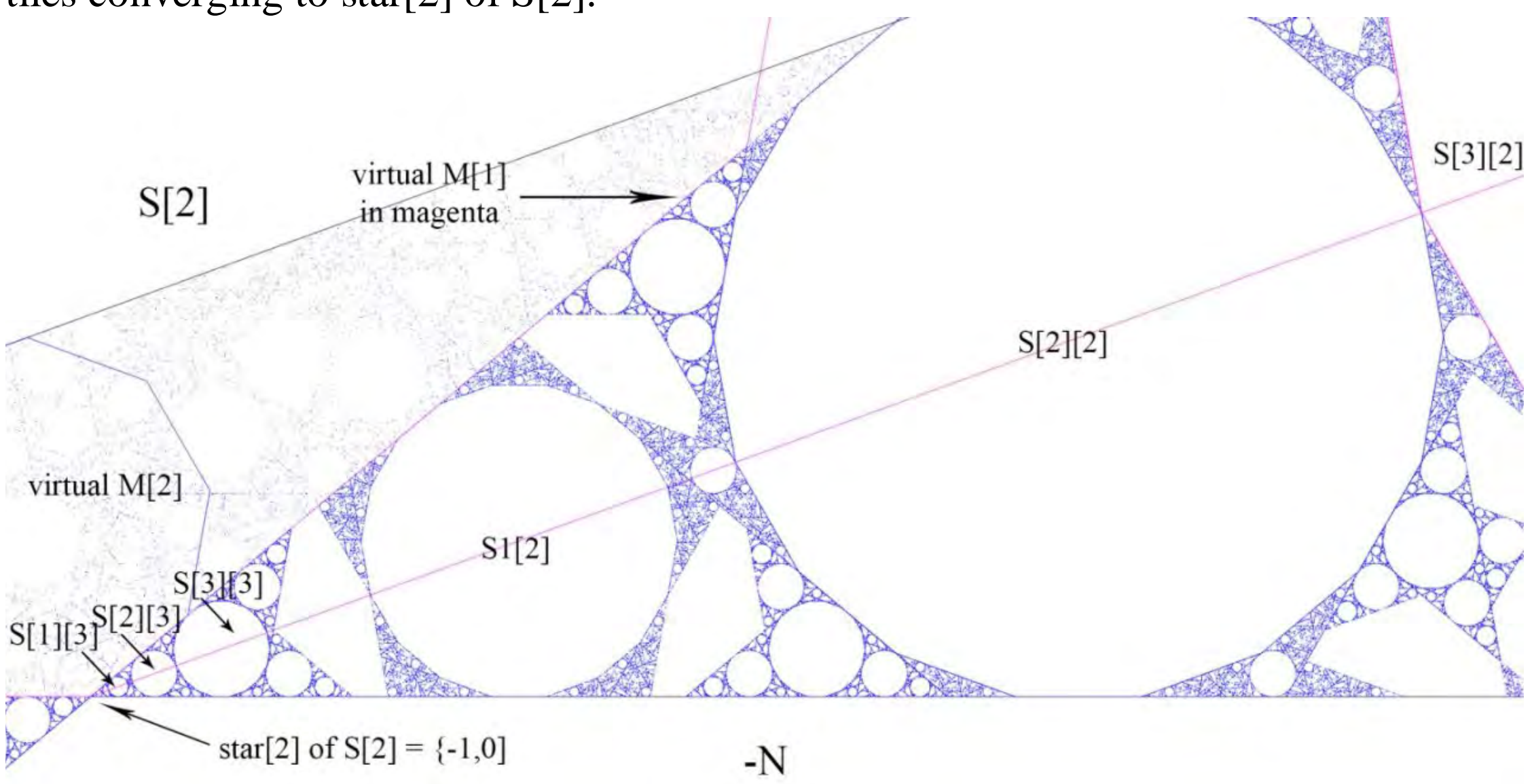

(This is a Dc plot so star[2] of S[2] is the {-1,0} vertex of –N.) It is an easy matter to track the tiles in this sequence to estimate their temporal scaling using the $\tau$-periods of their centers.

**Table 9.1** – The S[1][k] periods converging to star[2] of S[2]

| Tile | S[1][1] | S[1][2] | S[1][3] | S[1][4] | S[1][5] | S[1][6] |
|---|---|---|---|---|---|---|
| Period | 9 | 207 | 4005 | 79263 | 156374 | 30863799 |
| Ratio | | 23 (exact) | 19.3478 | 19.791 | 19.7285 | 19.7372 |

We conjecture that the temporal scaling of the S[1][k] is 20 and Figure 9.1 above can be used to give a plausible explanation for this value. The count for S[1][2]s on either side of the line of symmetry is 10 with one shared between sides. This yields a mean of 4 S[1][2]s per cluster and provides a likely temporal scaling of 20. The geometric scaling is hS[1][k+1]/hS[1][k]. Under the assumption of self-similarity this is simply hS[1][1]/hN since the sequence begins with N. By the First Family Theorem hS[1]/hN = $\text{Tan[Pi/9]}^2 \approx .132474$ so the local Hausdorff fractal dimension should be $-\text{Log[20]/Log[Tan[Pi/9]}^2\text{ ]} \approx 1.48203$,

This analysis was based on just the S[1] tiles but the surviving sequences appear to include S[2] and S[3] so they would share the same scaling. This applies to N = 18 also where the self-similarity begins with the 2nd generation at star[1] of S[2] and includes the same triple, but now the temporal scaling should be 9+1 and the geometric scaling is GenScale[9] = scale[4] of N = 18 = Tan[Pi/9]/Tan[4*Pi/9] to yield a local fractal dimension ≈ .838493.

**N = 10**

For N = 10 the 8k+2 Conjecture makes a number of assertions about the web and dynamics local to S[1] and S[2]. The Edge Conjecture makes no DS[k] predictions beyond S[1] at DS[3]. Since N = 10 has quadratic complexity, the web it shares with N = 5 has just one non-trivial scale and should be fractal in nature. As indicated earlier, N = 5 and N = 8 are the only non-trivial regular cases where the dynamics and singularity sets have been studied in detail. In [T] (1995) S. Tabachnikov derived the fractal dimension of W for N = 5 using 'normalization' methods and symbolic dynamics and in [S2] (2006) R. Schwartz used similar methods for N = 8. In [BC] (2011) Bedaride and Cassaigne reproduced Tabachnikov's results in the context of 'language' analysis and showed that N = 5 and N= 10 had equivalent sequences. (The Twice-Odd Lemma implies equivalent webs but not equivalent dynamics.) In N = 5 earlier we gave an independent analysis of the temporal scaling based on difference equations and here we will give the implications of the 8k+2 Conjecture for the matching cases of N = 10 and N = 5. As noted earlier, S[2] will have period N/2 so it will define two disjoint invariant regions and it is sufficient to track to just the blue regions shown here.

**Example 5.1** The temporal scaling of any 8k+2 N-gon (with N = 5 as the lone N/2 case)

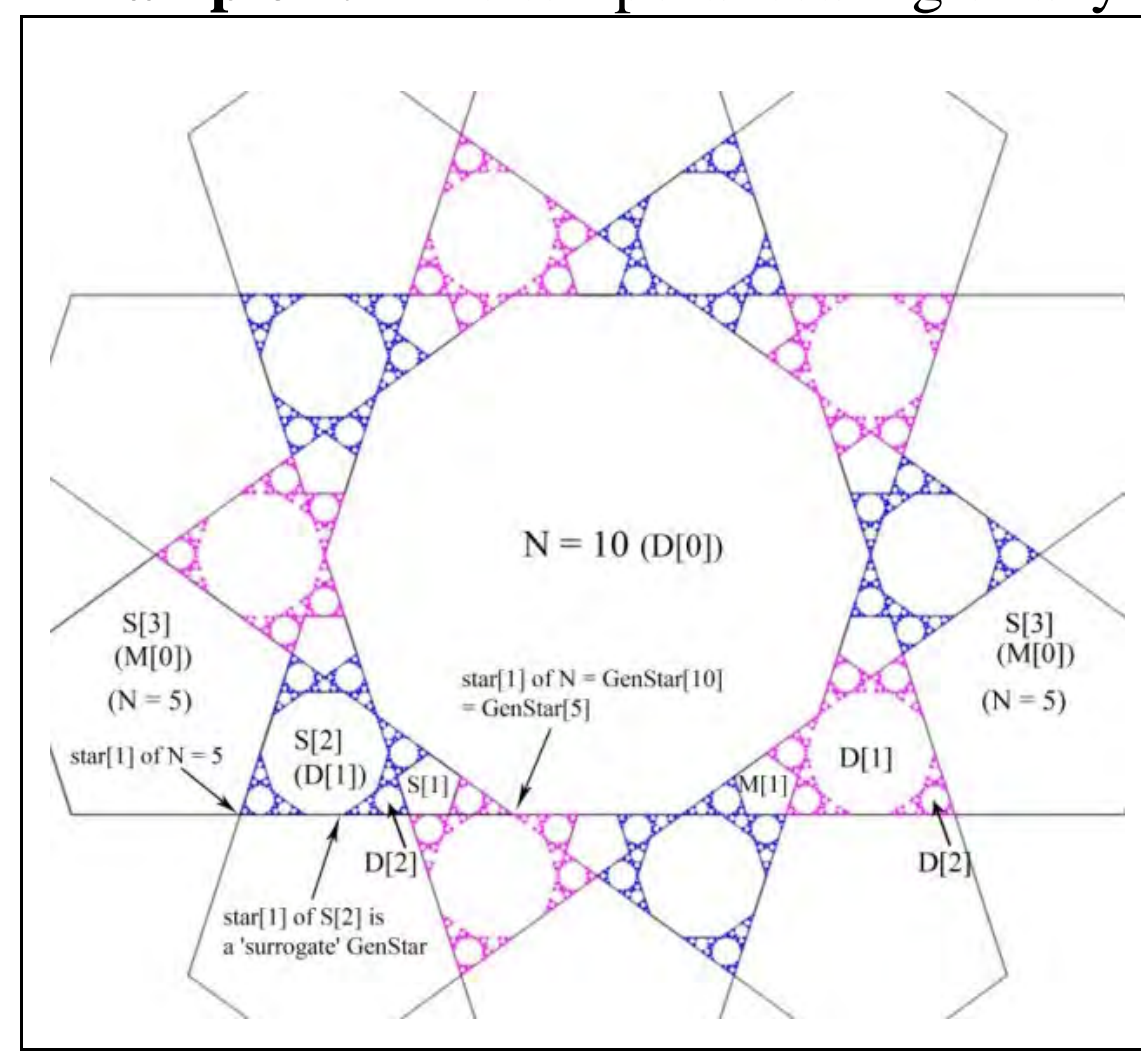

(i) The D[k] of N will have period $D_k$ where
$D_k = nD_{k-1} + (n+1)D_{k-2}$ where n = N/2 and $D_1 = n$ and $D_2 = n^2$

(ii) The M[k] of N will have period $M_k$ where
$M_k = nM_{k-1} + (n+1)M_{k-2}$ where n = N/2 and $M_1 = n$ and $M_2 = n(3(n-1)/2 + 2)$

(iii) The D[k] of N = 5 will have period $D_k$ where
$D_k = nD_{k-1} + (n+1)D_{k-2}$ where n = N and $D_1 = n$ and $D_2 = n(n+2)$

(iv) The M[k] of N = 5 will have period $M_k$ where
$M_k = nM_{k-1} + (n+1)M_{k-2}$ where n = N and $M_1 = 2n$ and $M_2 = 2n^2$ (these will be $2D_k$ for N = 10)

The solutions are:

| | |
|---|---|
| $(i)\ n = N/2: D[n,k] = -\dfrac{n\left((-1)^k - (1+n)^k\right)}{2+n}$ | $(ii)\ n = N/2, M[n,k] = \dfrac{n\left((-1)^{1+k} + (-1)^k n + 3(1+n)^k\right)}{2(2+n)}$ |

The first few periods are: $D_k$:{5, 25, 155, 925, 5555, 33325, 199955, 1199725, 7198355} and $M_k$ {5, 40, 230, 1390, 8330, 49990, 299930, 1799590,10797530}

| | |
|---|---|
| $(iii)\ n = 5, D[n,k] = \dfrac{n\left((-1)^k + (-1)^k n + 3(1+n)^k + n(1+n)^k\right)}{(1+n)(2+n)}$ | $(iv)\ n = 5, M[n,k] = -\dfrac{2n\left((-1)^k - (n+1)^k\right)}{2+n}$ |

The first few periods are: $D_k$ {5, 35, 205,1235, 7405, 44435, 266605, 1599635, 9597805} and $M_k$:{10, 50, 310, 1850, 11110, 66650, 399910, 2399450, 14396710} = $2D_k$ for N = 10.

• **N = 11**

This is the lone ‘quintic’ polygon (along with the matching N = 22). N = 11 is the first non-trivial member of the 8k+3 family. The Edge Conjecture predicts that for N odd, S[1] will be the DS[N-4] tile of S[2], so here S[1] should be the DS[7] of S[2].

N = 11 is the only case where S[2] and S[3] share a vertex and there is also a close relationship between S[1] and S[3]. They are in the same odd ‘tower’ so they share star[4] points. This point will also be star[N-4] of S[2]. Therefore S[1] will be the DS[N-4] in the star-2 family of S[2] as outlined in the Edge Conjecture. Note that this star[4] point is also the star[2] of N = 11, so it is like the floor wax that is also a shampoo and a dessert topping.

**Figure 11.1** When N is odd and greater than 5, S[3] will exist and share star[4] points with S[1]. For N = 11 this implies a strong connection between the geometry of S[1], S[2] and S[3].

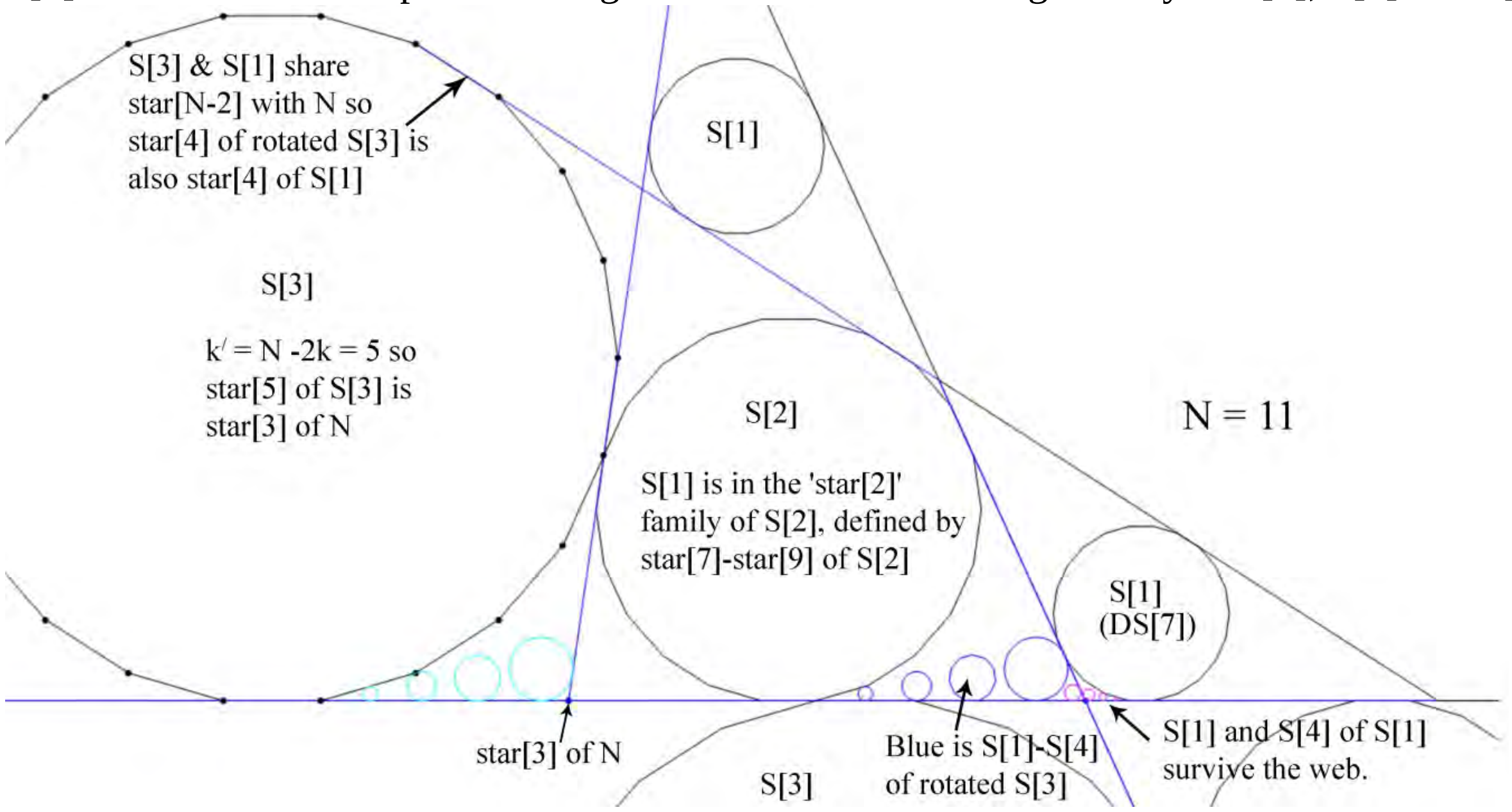

By contrast, for the twice-even family, S[1] and S[3] will share star[2] points and since S[1] is mutated in the twice-odd case, S[1] and S[3] will share star[1] points. (See N = 18 and N = 20). The shared star[4] point here implies that the matching S[4] tiles will have the same height relationships as S[1] and S[3], namely hS[1]/hS[3] = scale[3] = Tan[Pi/11]/Tan[3*Pi/11]. This is a scale[3] version of the scale[2] relationship between S[1] and S[2] and it would be possible to define ‘scale[3]’ families of S[1] in the same fashion as the ‘scale[2]’ families of S[2], but it seems that such families seldom exist. None of the S[1] through S[4] tiles of S[3] shown above actually survive the web, but S[4] and S[2] of S[1] both survive.

This gives the impression that the rotated S[3] does not play an important role in the geometry of this region, but that is not the case. The webs of S[1] are S[3] are closely linked and this grows stronger as N increases. It would be impossible to explain the geometry here without recourse to S[3]. For example the volunteer $D_1$ tile shown below shares a midpoint with the virtual S[1] of S[3]. Since it is (weakly) conforming to star[2] of S[2] its center must lie on the common line of symmetry joining star[2] with cS[1]. Therefore $D_1$ lives in both worlds and its characteristic polynomial can serve as 'DNA' to relate it to S[1], S[2] and S[3]. In addition one of the S[4] tiles which seems out of place is actually tied to S[3] as shown by the link to the center of S[3]. All the geometry local to this errant S[4] tile can be tied to the geometry of S[3]. Note that this local geometry is not shared by the remaining 4 DS[4]s. This is because these four are based on the virtual S[8] of S[1] shown here in magenta.

**Figure 11.2** The geometry local to S[2] shows the influence of S[3] and S[1]

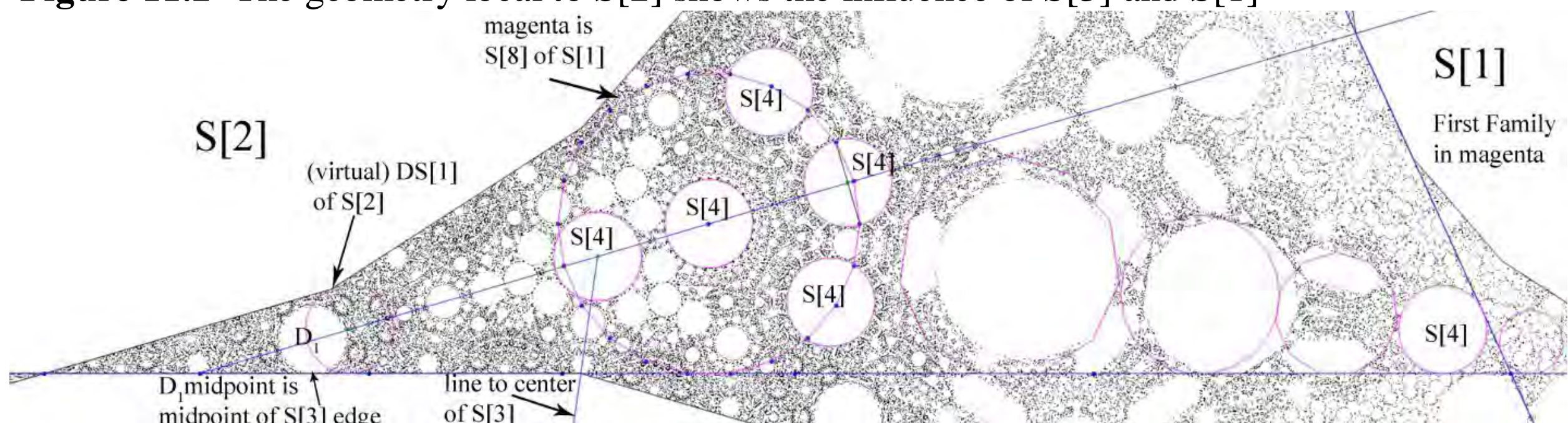

These are not approximate alignments and the S[8] case is among the strangest seen yet. The vertices of S[8] only define the centers of the three 'satellite' S[4]s in an indirect fashion. Just one alignment is shown here but they are all the same so these three satellites are rotations by 6Pi/22 about the center of S[8]. As indicated above the left-most S[4] defined by S[3] seems to be largely independent of the S[8] family and the green center displacements to the virtual DS[1] and to the outermost S[4] are only approximately equal. This was a mystery for a long time until S[8] was discovered as the source of the remaining 4 S[4]s. (These distances must all lie in $S_{11}$.)

Note: The 'normal' S[k] of S[2] typically do not exist but S[6] of S[2] would be congruent to S[4] of S[3] because the 'parents' have a scale[3]/scale[2] relationship. By the First Family Theorem, hS[4] (of S[3]) is hS[3]$\cdot s_1 \cdot s_4$ (where $s_k$ = Tan[kPi/N]). For the S[6] of S[2], hS[6] = hS[2]$\cdot s_1 \cdot s_6$. These are equivalent because hS[2]/hS[3] = scale[3]/scale[2] = $s_2/s_3$ (of N = 11) = $s_4/s_6$ (of N = 22). These two are only identical here because S[2] and S[3] share star[1] points.

Below are the three polynomials relating $hD_1$ with hS[1], hS[2] and hS[3] respectively in the scaling filed $S_{11}$ generated by x = GenScale[11] = scale[5] = Tan[Pi/11]/Tan[5Pi/11] ≈ .0422171

**AlgebraicNumberPolynomial[ToNumberField[hD1/ h__, GenScale],x] =**

$$-\frac{1}{8}+\frac{17x}{4}+\frac{17x^2}{2}+\frac{3x^3}{4}-\frac{3x^4}{8} \;;\quad \frac{1}{8}-\frac{9x}{4}+x^2+\frac{x^3}{4}-\frac{x^4}{8} \;;\quad \frac{1}{2}-11x-10x^2-x^3+\frac{x^4}{2}$$

This 'DNA' evidence points toward S[2] and S[3] as the closest relatives to $D_1$.

Mathematica can also find the relative horizontal or vertical displacements of $cD_1$ and cS[3].
**AlgebraicNumberPolynomial[ToNumberField[cD1[[2]]/ cS[3][[2]], GenScale],x]** gives
$2-12x-\frac{29x^2}{4}+\frac{x^4}{4}$. Internally $cD_1$[[2]] is stored in complex 'cyclotomic' form as

$$\frac{2\left(2+(-1)^{1/11}+(-1)^{4/11}+(-1)^{7/11}+2(-1)^{8/11}\right)}{1+2(-1)^{2/11}+(-1)^{3/11}+2(-1)^{4/11}+(-1)^{5/11}+2(-1)^{6/11}+(-1)^{8/11}}$$

**Figure 11.2** The edge geometry of N = 11 – the charter member of the 8k+3 family

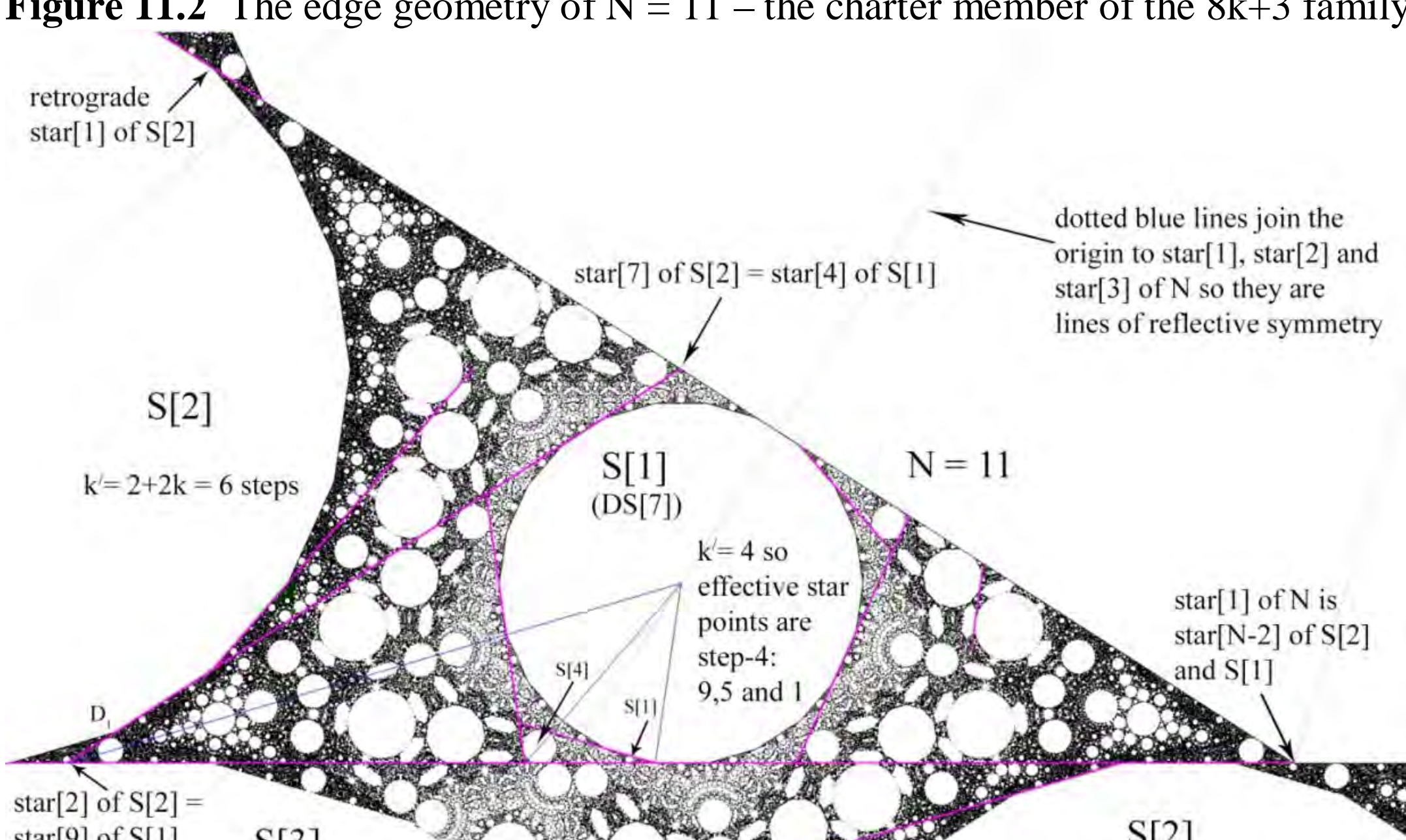

Both S[1] and S[2] have 'effective' star points which form early in the web and these star points depend on the step-sequences which are (retrograde) step-$k'$ where $k' = 2+2k$. So S[2] has a step-6 ccw web but when combined with the step-4 S[1] web it will have an extra 2 steps, so every predicted DS[k] will span 3 star points from star[k] to star[k+2]. Therefore S[1] is formed from the triplet star[9]-star[7] and there are no other predicted DS[k]. Of these 3 star points of S[2], only star[7] is called 'effective' because it occurs in the early web of S[2].

For N odd, all S[k] share their penultimate star points so star[N-2] of S[1] matches star[N-2] of S[2] and by reflective symmetry, star[N-2] of S[2] is equivalent to star[2] of S[2]. The effective star points of S[1] are always step-4 and there is no added contribution from S[2]. Since these effective star points count down from star[N/2-2] at star[2] of S[2], they are star[9], star[5] and star[1] of S[1]. The effective star points of S[2] are step-8 and count down from star[N-4] which is star[7] here. This point is also star[4] of S[1] which is not effective, but the neighboring star[5] is effective and forms early in the web. Together they define the S[4] tile of S[1].

Back in the First Family, each star[k] point corresponds to a (strongly) conforming S[k], and here we might hope that the effective star points of S[1] foster conforming tiles in a manner similar to the DS[k] of S[2]. This is not true in general, but here S[1] does have S[1] and S[4] tiles corresponding to star[1] and star[5]. There is also a volunteer at star[9] of S[1] which we call $D_1$ because it is weakly conforming to star[2] of D[1].

On the S[2] side of the 8k+3 family, there is a tendency for $D_k$ volunteers to exist between predicted DS[k]. This $D_1$ is a boundary case between DS[7] and S[2]. N = 19 has a conforming $D_1$ between DS[7] and DS[15] as well as a $D_2$ between DS[7] and S[2]. The one known exception is N = 35 which had no volunteer between DS[7] and DS[15]. This largest $D_1$ can often serve as a surrogate 'M' tile between S[1] and S[2].

**Figure 11.3** Since S[1] is congruent to S[2] of D, the 'second generation' tiles at D (N = 22) can be imported as virtual tiles of S[1]. The actual tiles that exist in the web of S[1] are black and the imported tiles that exist in the web at D are magenta. The transformation is **T[x]** = **TranslationTransform[cS[1]-cDS[2]] [x]**

This graphic provides some perspective of how the tile structure at S[1] differs from D[1]. This tile structure is what drives the web – or conversely. All these imported tiles are theoretically compatible with S[1] but only M[2] actually exists in both worlds. The magenta tiles such as S[5] and M[1] are actually in the First Family of S[1] but they do not exist in the web here on the edge of N = 11. Likewise the S[4] tile of S[1] exists here and not in the web at D[1].

Note how awkwardly M[1] of D is embedded in S[2]. This shows the lack of compatibility between S[1] and S[2] and helps to explain why this web geometry is so complex. We have found virtually no relationship with the First Family of M[1] and the actual tiles at S[2]. On the positive side N will have the same edge length as D – even though they have different genders.

**Figure 11.4** The shared geometry of S[1] and S[2]

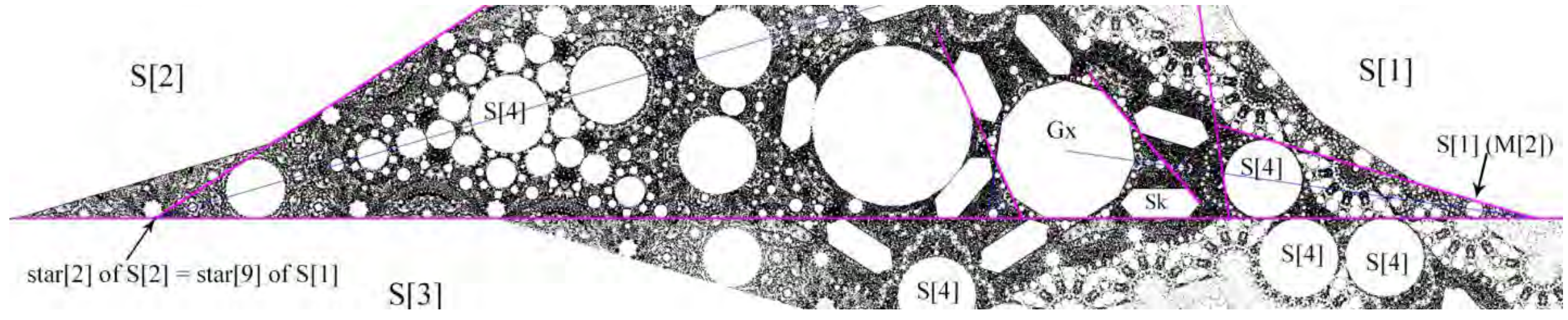

By convention we will choose to study the left side web of S[1] because the clockwise early web defines initial star points on this side of S[1] - and these 'effective' star points are step-4 as noted above. It is our contention that the overall dynamics and geometry of this region are driven by these effective star points - as in the generalized First Family Theorem. Of course N = 11 is special in many ways - and the true nature of the 8k+3 family and the efficacy of these step-4 star points will be better illustrated by larger N values – such as N = 19, 27, 35 and 43. In general these step-4 families will contain very few remnants of normal First Families.

**Figure 11.5** Detail of the web local to S[4] and the weakly conforming ‘volunteer’ Gx. In the First Family of S[1] there are also S[1]s as shown here. These would be M[2]s relative to D.

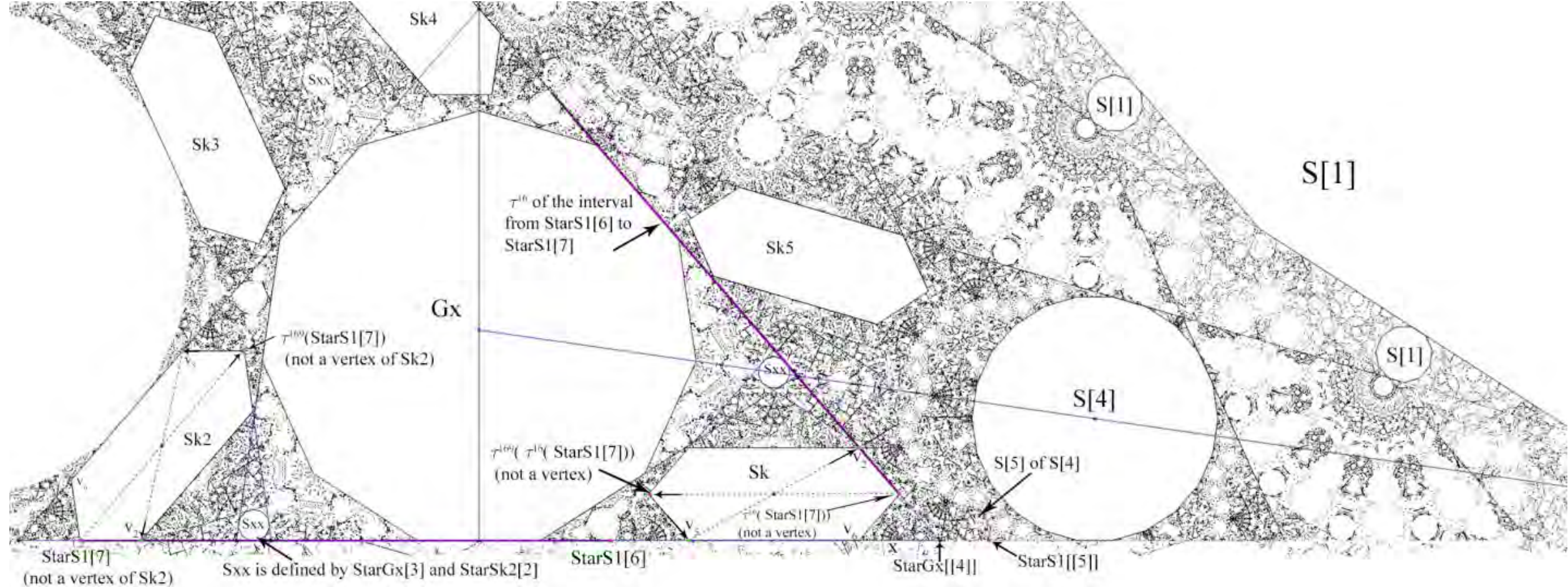

In the level-16 web of S[1], the interval from StarS1[6] to StarS1[7] defines a shared edge of Gx and Sk as shown above. The slope of this magenta interval must match an edge of N so it is known. Therefore this interval can be used to define StarGx[4] and also an edge of Sk. The lower end point of this interval is $\tau^{16}$(StarS1[7]) and it is necessary to find this point exactly. This is an issue which we have addressed in the Appendix of [H5]. The major star points such as StarS1[7] will typically not have extended orbits because they will map to trailing edges of N within a few iterations. But these star points typically have one-sided limiting orbits which can be calculated exactly using surrogate initial points and we will do this here for both Gx and Sk.

It appears above that StarS1[7] is vertex $v_6$ of Sk2 but there is a small horizontal offset involved and the same is true for $\tau^{16}$(StarS1[7]) and Sk. The Sk and Sk2 tiles are related by a simple rotation about the center of Gx, so we can work with either Sk or Sk2 once Gx is known.

**Calculations for Gx** (based on N = 11 at the origin and with hN = 1)

(i) Any exact point on an extended edge of Gx will define the matching star point, so to find StarGx[4] all that is needed is $p_1 = \tau^{16}(p_0)$ where $p_0$ = StarS1[7]. This is a simple calculation that Mathematica will do with exact arithmetic – but only if the correct vertex points in the orbit are known. Here it will likely generate an error because these star points will either have no image at all under $\tau$ or will soon map to a trailing edge of N. Therefore we will use a surrogate neighbor such as pn = $p_0$ + {0, .000001}. This will put pn inside Sk2 – which is an advantage because all points in a tile such as Sk2 must map together.

Tiles like Sk2 and Gx are prominent in the web, so it is easy to probe them with test points and find their periods. Every point inside Sk2 has period 338 – except the center which has period 169. These are called ‘period-doubling’ tiles . Here we only need to track 16 iterations of pn for Gx, but we will later track further iterations of pn to get the vertices of Sk – as shown by the dotted blue arrows above.

(ii) Ind = IND[pn,16] will generate the indices of the first 16 vertices of N in the orbit of pn. These are {8,11,3,7,10,3,6,9,11,1,2,4,6,8,9,10} which means that the vertices visited will be $c_8$,$c_{11}$, etc –where by convention $c_1$ is the 'top' vertex of N. (In general we prefer to use step-sequences of orbits which are the first differences of these indices – because all the Sk will have the same step sequences - but not indices.) We will assume that $p_0$ will have this same 'corner sequence' and if this is false the error will be large and easily detected.

(iii) The matching routine PIM[pn,8] will take these 16 indices and perform 8 iterations of the 'return map' $\tau^2$ to calculate the actual orbit of pn, so $\tau^{16}$(pn] will be the last element of PIM. Here we set P = PIM[$p_0$,8] and P[[8]] = $\tau^{16}(p_0) = p_1$ will be exact.

(iv) This point $p_1$ is on an extended edge of Gx with known slope and this defines an exact StarGx[4]. Since Gx is conforming to S[1], StarGx[5] = StarS1[[1]] so we can use the Two Star Lemma to find hGx and its location:
d = StarS1[1][[1]]-StarGx[4][[1]]; hGx = d/(Tan[4π/11]-Tan[π/11]) ≈ .013760647991661
AlgebraicNumberPolynomial[ToNumberField[hGx/hN,GenScale],x]

$$-\frac{5}{8}+\frac{59x}{4}+9x^2+\frac{x^3}{4}-\frac{3x^4}{8}$$

Midpoint Gx = StarS1[1]-{hGx*Tan[5*Pi/11],0)};cGx = MidGx + {0,hGx};
rGx = RadiusFromHeight[hGx,11];Gx = RotateCorner[cGx + {0, rG}, 11, cGx};
StarGx = Table[Midpoint – {hGx*Tan[k*Pi/11], 0}, {k,1,5}]

**Calculations for Sk**

(i) Now that cGx is known we can work with either Sk or Sk2 because they are related by RotationTransform[8π/11,cGx]. This implies that $\tau^{16}(p_0)$ is simply a rotation by 8π/11 of StarS1[7] – but this calculation depends on knowing cGx.

As indicated earlier, it is possible to use the same pn point as a surrogate and simply extend the orbit to 169 iterations. So IND[pn,169] will generate the indices, but our software for PIM is based on the return map $\tau^2$ and 169 is odd. One solution is to replace $p_0$ with $p_2 = \tau(p_0)$ (and np with $\tau$(np). This is easy since the first index is known to be 8, so $\tau(np) = 2c_8 - np$ and IND[$\tau$(np),170] will yield the corner sequence for $\tau$(np) which should also be the corner sequence for $p_2 = \tau(p_0)$. Therefore P = PIM[$p_2$,85] will give the exact $\tau^2$ orbit for $p_2$ which will end with px = P[[85]]] = $\tau^{169}(p_0)$

(ii) This px point defines the top edge of Sk2 which intersects the edge of Gx at what we call vertex 5 or $v_5$. This vertex can then be mapped to the opposite vertex $v_2$ using the same indices because all point in a tile such as Sk2 must map together with the same indices – even the center which is the unique point where $\tau^{169}$(cSk) = cSk so it has period 169 instead of 338. Therefore $v_2 = \tau^{169}(v_5)$ and it can be found by setting P1 = PIM[$\tau(v_5)$,85] and then $v_2$ = P1[[85]].

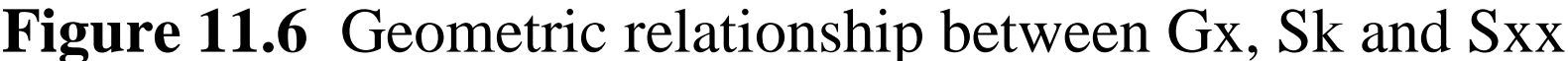

**Figure 11.6** Geometric relationship between Gx, Sk and Sxx

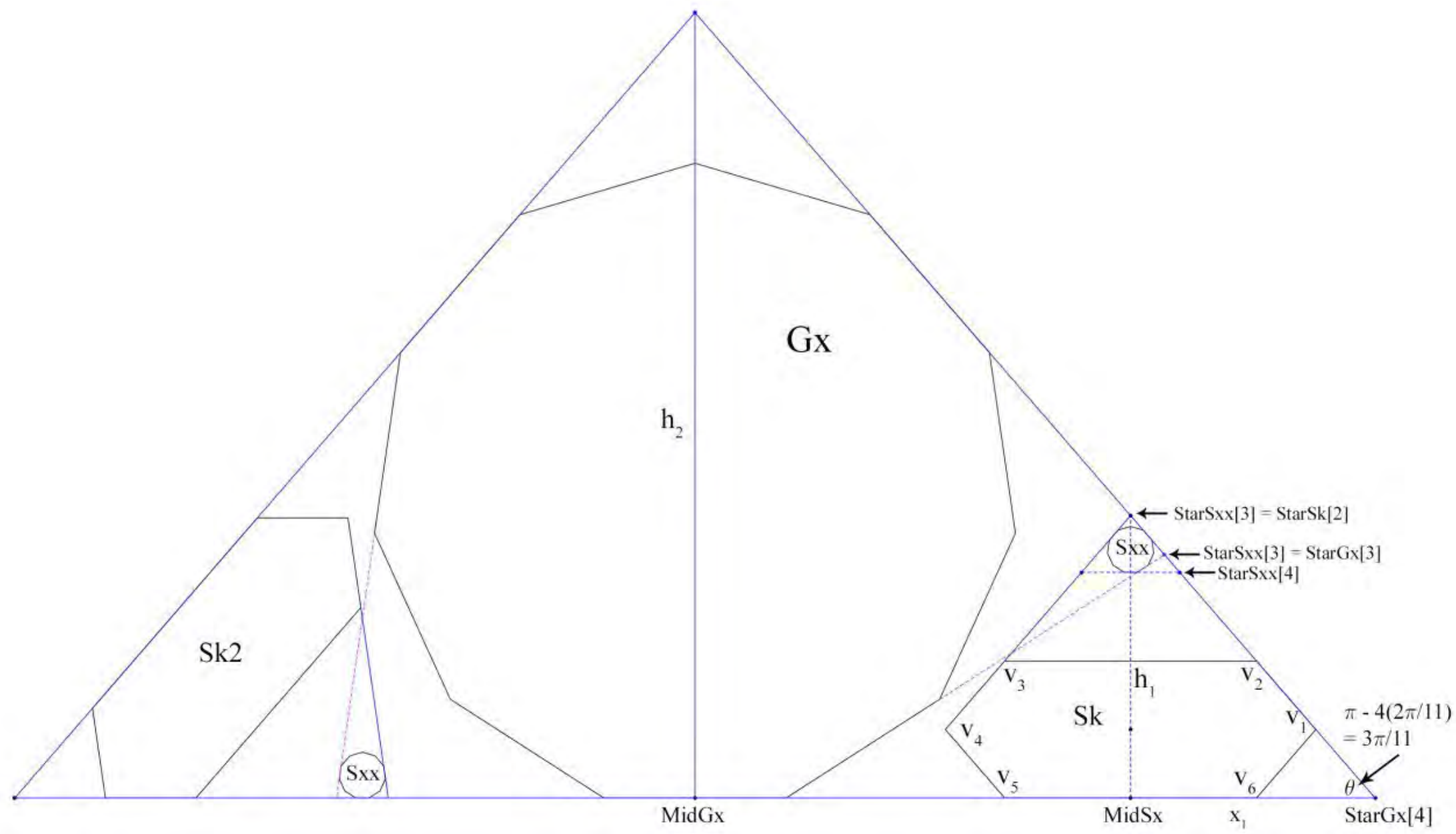

(iii) Any interior point can be used to find the center of a period-doubling tile because all points map to a reflection about the center under half of the period. The calculations above yield two candidate lines that must pass through the center. Using their intersection for the center avoids a distance calculation which might be an issue with the scaling field $S_{11}$. Now map cSk2 ,$v_5$ and $v_2$ to Sk using RotationTransform[8π/11, cGx] which is exact relative to $S_{11}$.

(iv) In Sk2, $v_2$ defines $v_6$ and the 'height' which we call $h_1$. Since $v_5$ is known the height defines $v_3$ and the horizontal center line defines $v_1$ and $v_4$ because the slopes of the edges are known. This defines Sk and the remaining Sk2, Sk3, Sk4 and Sk5 are either rotations or reflections about the center line joining Gx with StarS1[1].

(v) Sk clearly has a close geometric relationship with Gx and we will try to make this precise. Later we will show that these two tiles generate an 'offspring' called Sxx in a manner similar to the Sx tile of S[5] and Px at D.

There are two nested isosceles triangles here with heights $h_1$ and $h_2$. And Sk is embedded in the smaller triangle. The 'star[k] angles' of any odd N-gon have the form π - k$\phi$ where $\phi$ is 2π/11, so the angle θ at StarGx[4] is π -8π/11. Therefore Tan[3π/11] = $h_1/x_1$ = $h_2/x_2$ where $x_2$ is the distance from MidGx to StarGx[4] - which is hGx*Tan[4π/11. Therefore $h_1/h_2 = x_1/x_2$ and the embedded Sxx has a self-similar relationship with Gx. This is turn relates Gx and S[4] because the S[5] tile of S[4] is congruent to Sxx.

(vi) The last step for Sk is to show that the ratio each edge length with sN is in scaling field $S_{11}$. If $s_1$ and $s_2$ are the long and short sides of Sk their polynomials are shown below.

**AlgebraicNumberPolynomial[ToNumberField[$s_1$/sN],GenScale[11]],x] =**

$\frac{41}{8} - 118x - \frac{279x^2}{4} - x^3 + \frac{21x^4}{8}$. The matching polynomial for the short side $s_2$ is $\frac{3}{2} - \frac{277x}{8} - \frac{141x^2}{8} + \frac{x^3}{8} + \frac{5x^4}{8}$

**Figure 11.7** Detail of the web local to Sk. There are 'islands' here which are similar to the invariant regions which arise in the web of Mx at D. See N = 22.

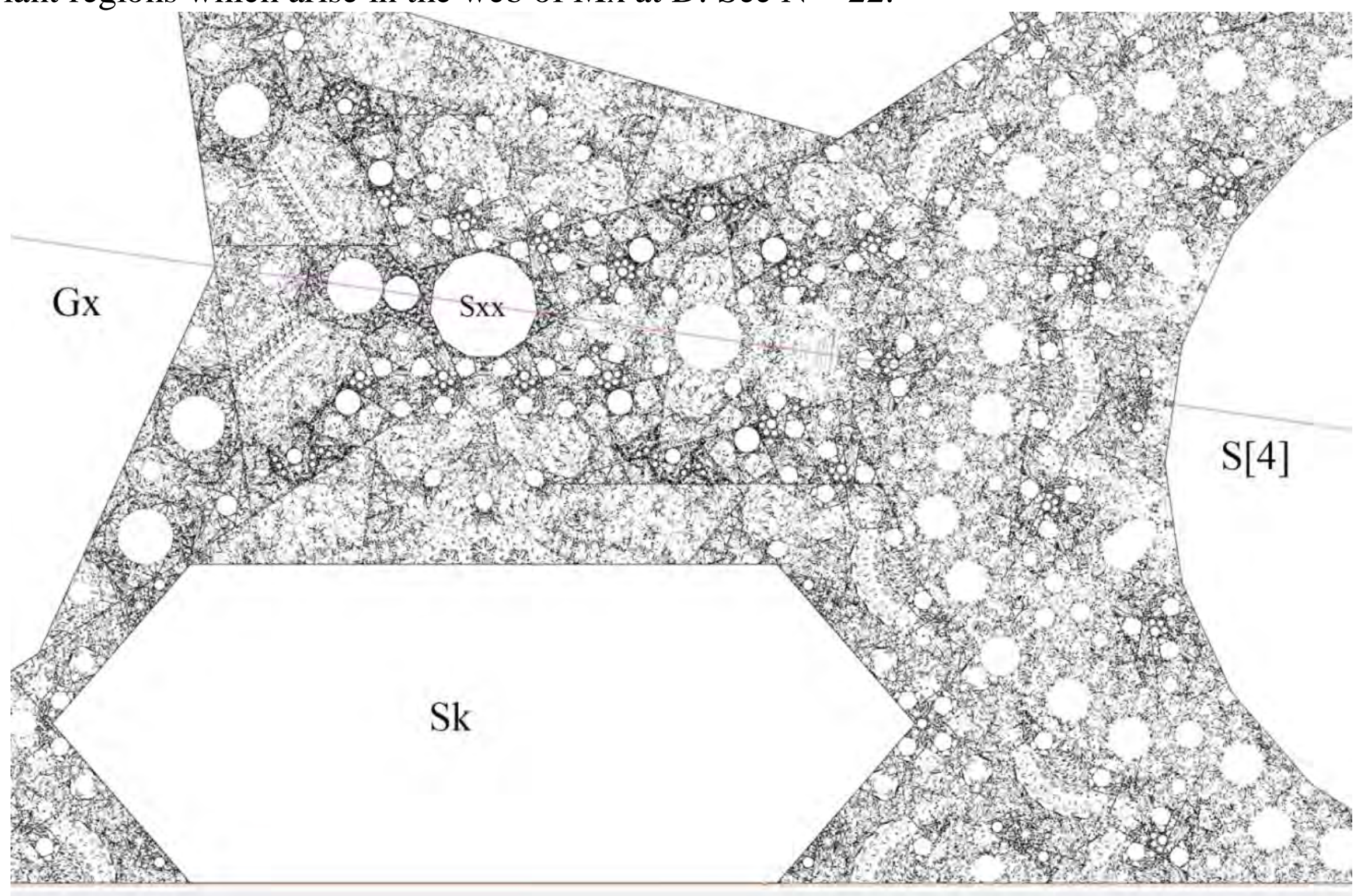

**Calculations for Sxx**

(i) We will find the parameters of the Sxx on left side of Gx. This is a simple application of the Two Star Lemma using StarGx[3] and what we call star[2] of Sk2. The local indices of Sxx are both 3, so d = StarGx[3][[1]]-StarSk2[2]; hSxx = d/(2Tan[3π/11])]

**AlgebraicNumberPolynomial[ToNumberField[hSxx/hN,GenScale[11]],x]** yields

$p[x] = -\frac{1}{8}+\frac{11x}{4}+\frac{11x^2}{2}+\frac{5x^3}{4}-\frac{3x^4}{8}$ Here hN = 1 so hSxx = p[GenScale[11]] ≈ 0.000992499045

(ii) To construction Sxx:

MidSxx = {StarGx[[3]][[1]]- hSxx*Tan[3Pi/11],0}; cSxx = {MidSxx[[1]],hSxx}

rSxx = RadiusFromHeight[hSxx,11]; Sxx[[1]] = cSxx + {0,rSxx};

Sxx = RotateVertex[Sxx[[1]],11, cSxx];

This polynomial for hSxx is relatively simple because it is also the polynomial for S[5] of S[4] – and this shows that Sxx is a true 3rd generation tile. But as with $D_1$ earlier the positional complexity of Sxx is much more difficult because it depends on star points of both Gx and Sk. For example the vertical coordinate of vertex 1 of Sxx is relatively simple since it is just the height plus radius, but the horizontal coordinate is another issue. The simplified form is shown below:

Sxx[[1]][[1]] = $2(-\frac{3}{2}-(-2\cos\left[\frac{\pi}{22}\right]\left(48064\cos\left[\frac{\pi}{22}\right]-198157\cos\left[\frac{3\pi}{22}\right]+410268\cos\left[\frac{5\pi}{22}\right]+382041\sin\left[\frac{\pi}{11}\right]-526796\sin\left[\frac{2\pi}{11}\right]\right)$

$$\left(-65410266379512+130820532376983\cos\left[\frac{\pi}{11}\right]-130820531468146\cos\left[\frac{2\pi}{11}\right]+130820529676436\sin\left[\frac{\pi}{22}\right]-130820529922657\sin\left[\frac{3\pi}{22}\right]+130820530531082\sin\left[\frac{5\pi}{22}\right]\right))/$$

$$\left(4\left(101\cos\left[\frac{\pi}{22}\right]-260\cos\left[\frac{3\pi}{22}\right]+276\cos\left[\frac{5\pi}{22}\right]+94\sin\left[\frac{\pi}{11}\right]-209\sin\left[\frac{2\pi}{11}\right]\right)^2+\left(33312-66408\cos\left[\frac{\pi}{11}\right]+65974\cos\left[\frac{2\pi}{11}\right]-66422\sin\left[\frac{\pi}{22}\right]+66076\sin\left[\frac{3\pi}{22}\right]-65840\sin\left[\frac{5\pi}{22}\right]\right)^2\right))\tan\left[\frac{\pi}{11}\right]$$

The expression for Sxx[[6]][[1]] (which is star[1] of Sxx) is much worse and would run for more than 30 pages of normal print .The first few terms of that expression are shown below:

**Figure 11.8** The first few terms in the trigonometric expression of the horizontal coordinate of vertex 6 of Sxx. For expressions like this that run for multiple pages, Mathematics asks if you want more (or less) . At each stage there are 'unresolved' numbered terms that are eventually evaluated. Note that the whole expression is in grey as a warning that it is partial.

But because calculations within the scaling field $S_{11}$ are so efficient, Mathematica only takes a few seconds to find the polynomial for these coordinates relative to GenScale[11] and sN.
**AlgebraicNumberPolynomial[ToNumberField[Sxx[[6]][[1]]/sN,GenScale[11]],x]** yields

p[x] = $-\frac{5}{8}-\frac{27x}{2}-\frac{33x^2}{4}-\frac{x^3}{2}+\frac{3x^4}{8}$ so Sxx[[6]][[1]] = sN* p[Tan[Pi/11*Tan[Pi/22] =

$$2\tan[\pi/11]\left[-\frac{5}{8}-\frac{27}{2}\tan\left[\frac{\pi}{22}\right]\tan\left[\frac{\pi}{11}\right]-\frac{33}{4}\tan\left[\frac{\pi}{22}\right]^2\tan\left[\frac{\pi}{11}\right]^2-\frac{1}{2}\tan\left[\frac{\pi}{22}\right]^3\tan\left[\frac{\pi}{11}\right]^3+\frac{3}{8}\tan\left[\frac{\pi}{22}\right]^4\tan\left[\frac{\pi}{11}\right]^4\right]$$

≈ - 0.7103831120784580 (with hN =1). This is identical to the lengthy trigonometric expression in Figure 11.8 above but Mathematica has a very hard time simplifying such an expression unless it is told to do so in the context of $S_{11}$.

Here the vertical coordinate of vertex 6 is -1, but in general the vertical coordinates must be determined relative to hN not sN. Once again Mathematic has no problem doing this.
**AlgebraicNumberPolynomial[ToNumberField[Sxx[[1]][[2]]/hN,GenScale],x]** yields

p[x] = $-\frac{13}{8}+\frac{59x}{4}+\frac{5x^2}{2}-\frac{7x^3}{4}+\frac{x^4}{8}$

**Figure 11.9** Detail of the web local to Sxx

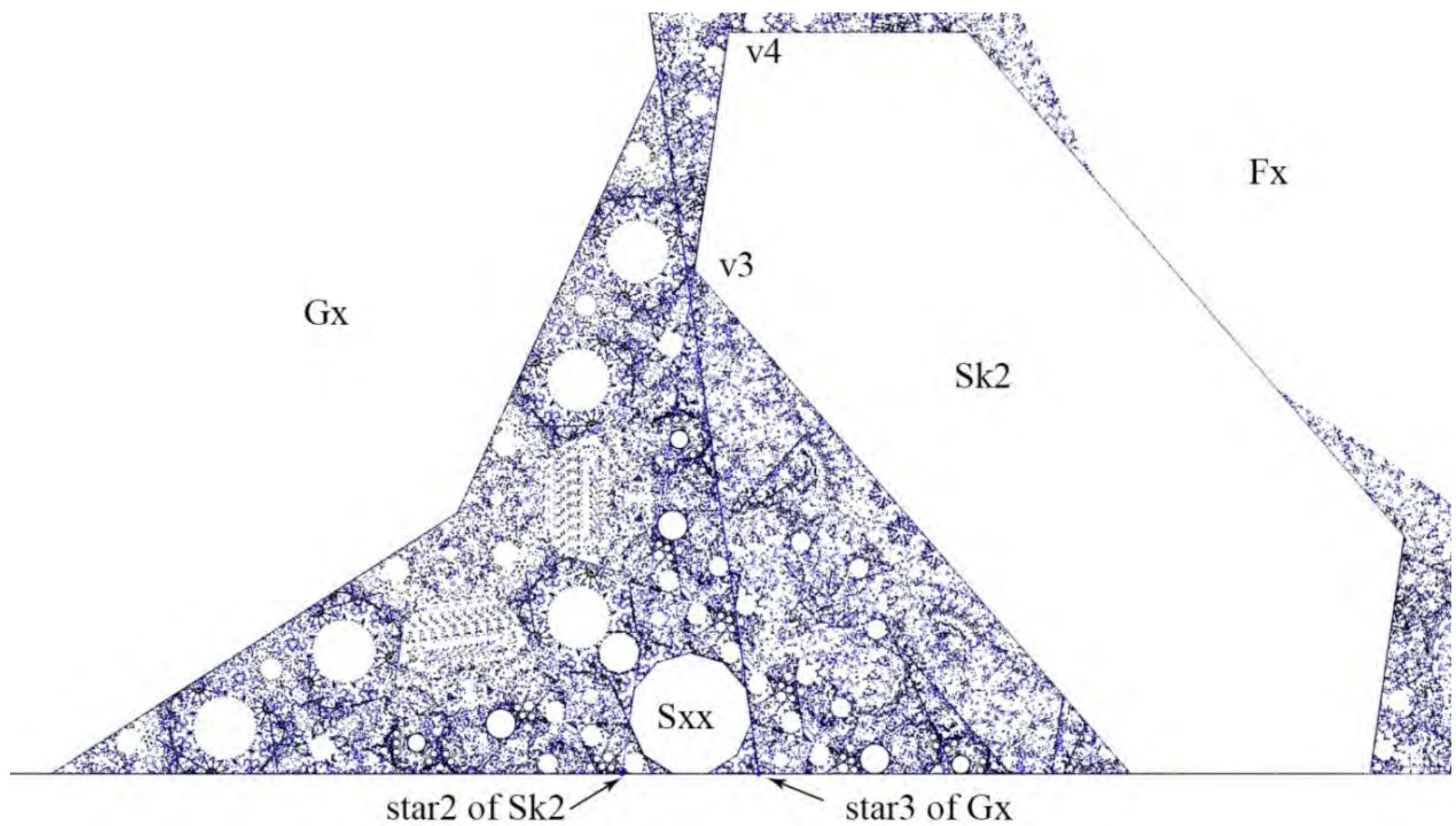

There is little doubt that Sxx will have some form of local extended family. Earlier we showed that geometrically Sxx is closely related to Gx, and based on this plot it seems that there is also a close relationship with Sk2 because Sxx appears to be the first tile in a sequence of tiles converging to the star point shared by Sxx and Sk2. The second tile in that sequence can be seen above. By convention the S[k] are scaled relative to D. For all N, S[1] is special because hS[1]/hS[k] is always scale[k], so here hS[1]/hD = hS[1]/hS[5]= scale[5] = GenScale[11]. But hS[1]/hS[2] = scale[2] of N = 11 – which is an awkward relationship so they share no First Family tiles.

| hS[1] /hD | hS[1]/hS[2] | hGx/hN | cS[k][[2]]/hN | hSxx/hN |
| --- | --- | --- | --- | --- |
| x | $\frac{1}{2}-x-\frac{x^2}{2}$ | $-\frac{5}{8}+\frac{59x}{4}+9x^2+\frac{x^3}{4}-\frac{3x^4}{8}$ | $\frac{3}{4}-\frac{161x}{4}-\frac{107x^2}{4}-\frac{3x^3}{4}+x^4$ | $-\frac{1}{8}+\frac{11x}{4}+\frac{11x^2}{2}+\frac{5x^3}{4}-\frac{3x^4}{8}$ |

**•N = 12**

N = 12 is a quadratic polygon and the first non-trivial member of the 8k+4 family. This is a very interesting family geometrically and algebraically. Since the mutation criteria for S[k] in the twice-even case is gcd(N/2-k,N) > 2, this family is the only one where S[2] is mutated. This is a simple star[1] to opposite-side star[3] mutation, so S[2] will always be the equilateral Riffle or 'weave' of two N/4-gons. This is compatible with the Edge Conjecture which predicts the survival of a DS[4]. S[3] is also mutated, but as always the M tile escapes mutation since it is S[N/2-2]. Since hS[1]/hS[2] = GenScale[12] = $\mathrm{Tan}[\pi/12]^2$, the M[k]will play a major role here.

**Fig. 12.1** The First Family

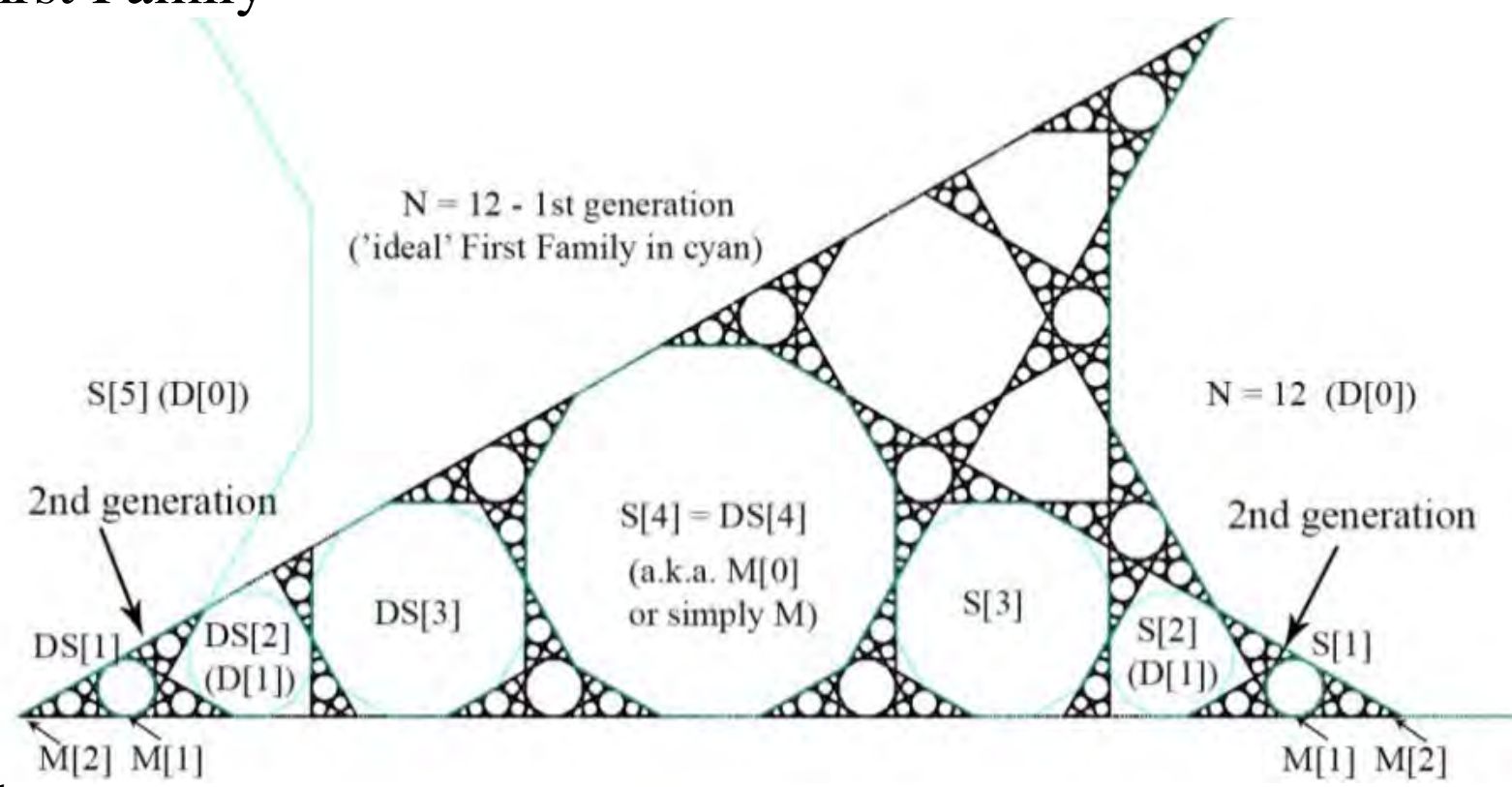

Below is the 2nd generation where the predicted DS[4] is S[1], also known as M[1]. The geometry local to M[1] looks like an unfinished web but this is consistent with the fact that S[2] is composed of two regular triangles which have trivial local webs. As N increases in the 8k+ 4 family, this local N/4-web of S[2] will be more complex, but apparently it will never yield 'normal' families for the DS[4]s at the magenta extended vertices. However the blue vertices will be also be vertices of the underlying S[2] and have a more promising geometry. It appears that for all members of this 8k+4 family, the blue star[1] of S[2] will be shared by an alternative 'parent' Px of DS[4]. For N = 12 this Px is the large S[3] tile of N shown here.

**Fig 12.2** The 2nd generation

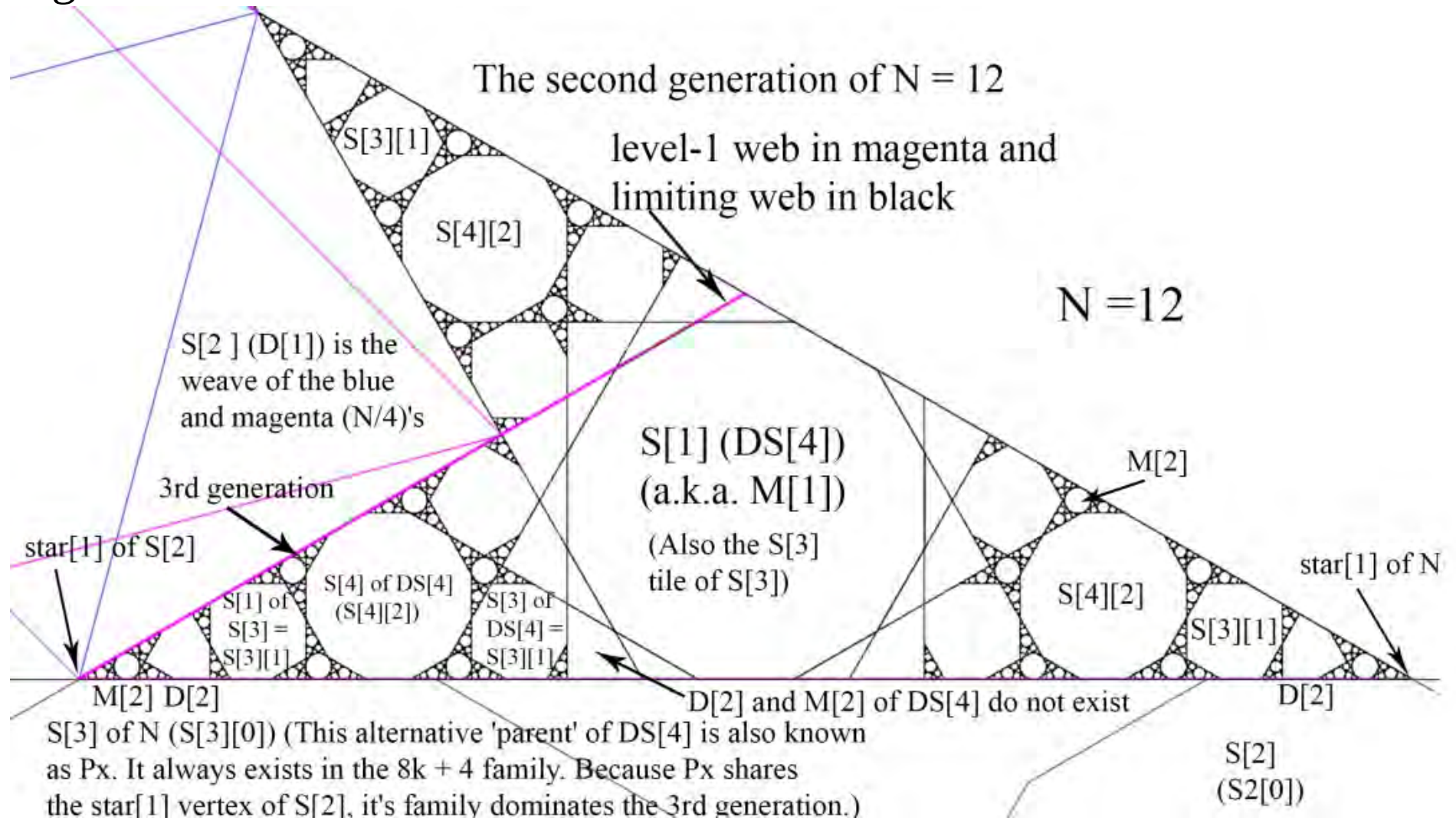

The First Family of S[3] does include DS[4] as well as a next-generation S[3][1] at the S[1] position. This S[3][1] survives the web along with DS[4]. (The S[2] of S[3] does not survive.) This S[3][1] and symmetric copies then generate M[2]s as their S[3] tiles. We will show that in the limit each S[3][k] will account for 3 M[k+1]s, and since there will 3 of these S[3][1]s in each 'dart' and there are 3 darts in each generation, the overall growth rate should be 27 as expected. Here the count of M[2]s is a little short with only 24 but that is because the mutation of S[2] destroys potential M[2]s (and D[2]s) which would exist on the edges of DS[4]. We will show why these mutations will play a diminishing role with each new generation.

Below is an enlargement of one dart. The only M[2] that is generated by a (mutated) D[2] is at star[1] of N. The remaining 5 M[2]s are generated as step-3 tiles of the S[3][1] so they have normal webs which will contain M[3]s at the 'next-generation' step-1 positions. Therefore of the 23 'darts' shown here only 3 are affected by the mutations - and this ratio will decrease with each new generation to yield a limiting count of 3 M[k]s for each S[3][k–1] and a limiting temporal scaling of 27 for the M[k]s (and D[k]s). See Example 5.2 for the resulting fractal dimension.

**Fig 12.3** Enlargement of one 'dart' of the 2nd generation

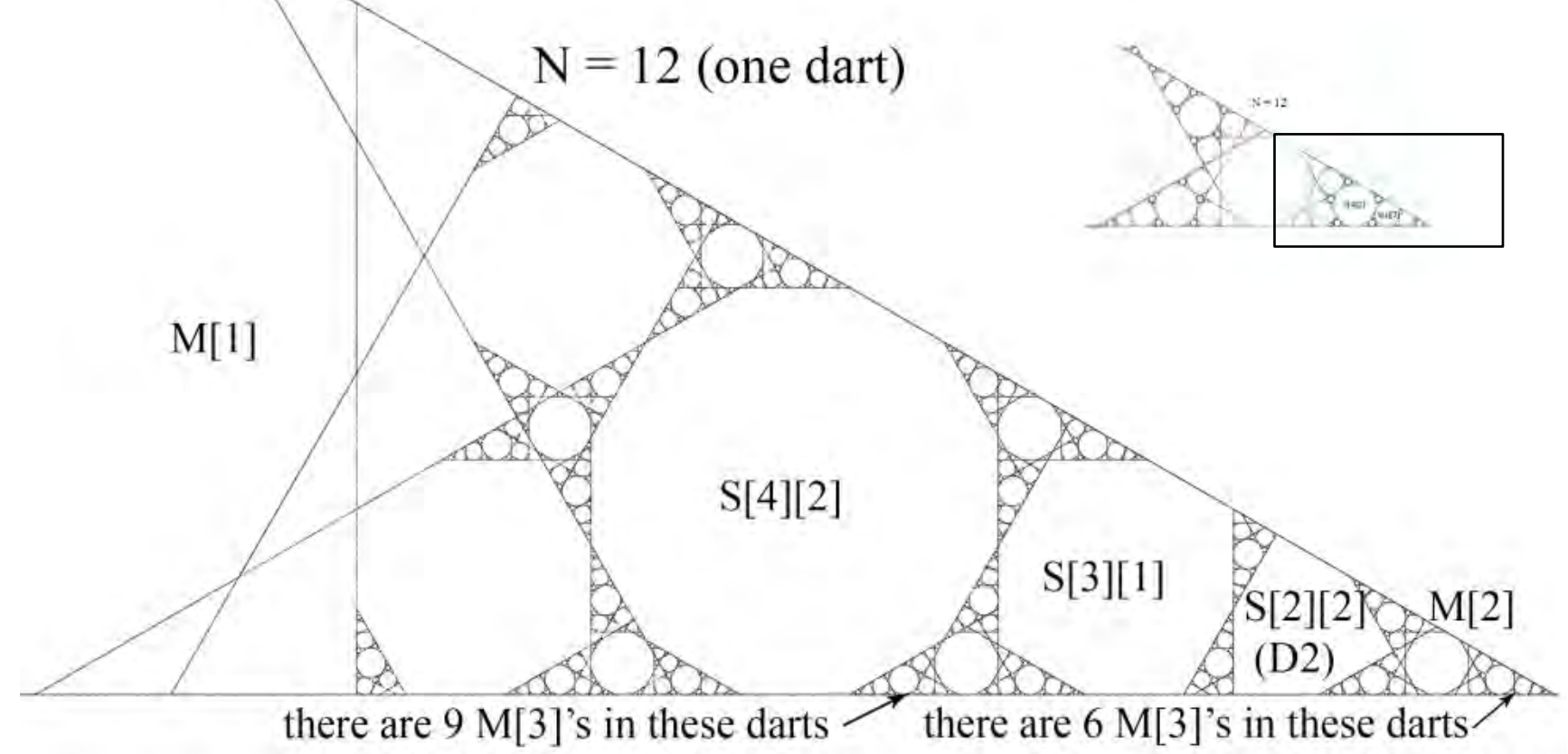

Another way to verify the temporal scaling for quadratic self-similar webs is to simply count the growth of the tiles using the $\tau$-periods. Here the $\tau$-periods of the M[k] in the canonical invariant 'star-region' are 60, 942, 28292, 775356, 21055308,.. with ratios of 15, 30, 27.41, 27.15. These are the combined periods of the M[k] at GenStar and their reflections about S[4]. Even though the local web has perfect reflective symmetry with respect to S[4], the dynamics are different and this combined count helps to minimize these differences. The dynamics of any composite N-gon allows for the possible 'decomposition' of expected orbits unto groups of orbits with smaller periods. This makes it difficult to match tile counts with periods, but for self-similar webs, the effect of these exceptions diminishes with each generation and in the limit the $\tau$-ratios will match the geometric ratios.

The **8k+4 Conjecture** states that volunteer Px 'parents' of DS[4] always exist and come in two versions. N = 12 is the charter member of a mod-16 family 12+16j where Px has DS[4] at S[3 +8j]. The next member is N = 28 where DS[4] is S[11] of Px and this in turn defines Px. In all cases the subsequent dynamics local to Px look promising. The matching mod-16 family is based on N = 20 and Px is simply a D tile of DS[4]. What these mod 16 families share is a secondary role of the predicted DS[4], with volunteer Px tiles now playing the major role in the next-generation web evolution local to S[2]. See N = 20 to follow.

**A comparison of the fractal dimension of the quadratic N-gons: N = 5, 8, 10 and 12**

N = 5, 8, 10 and 12 have φ(N)/2 = 2, so they have quadratic complexity- where the only non-trivial scale is GenScale[N]. Since the webs are naturally recursive, a single scale should yield a self-similar web and we will derive the similarity dimension of these webs below. Since the geometric scaling is known, the only issue is the 'temporal' scaling – which describes the limiting growth in the number of tiles. For self-similar webs this temporal scaling can be derived from a simple 'renormalization' process – where a representative portion of the web is scaled by GenScale[N] and mapped to itself under $\tau^k$ as shown by the magenta lines below. (The web for N = 10 is identical to the web for N = 5 and the N = 8 and 12 cases are closely related since their cyclotomic fields are generated by {$\sqrt{2}$,$i$} and {$\sqrt{3}$,$i$}.).

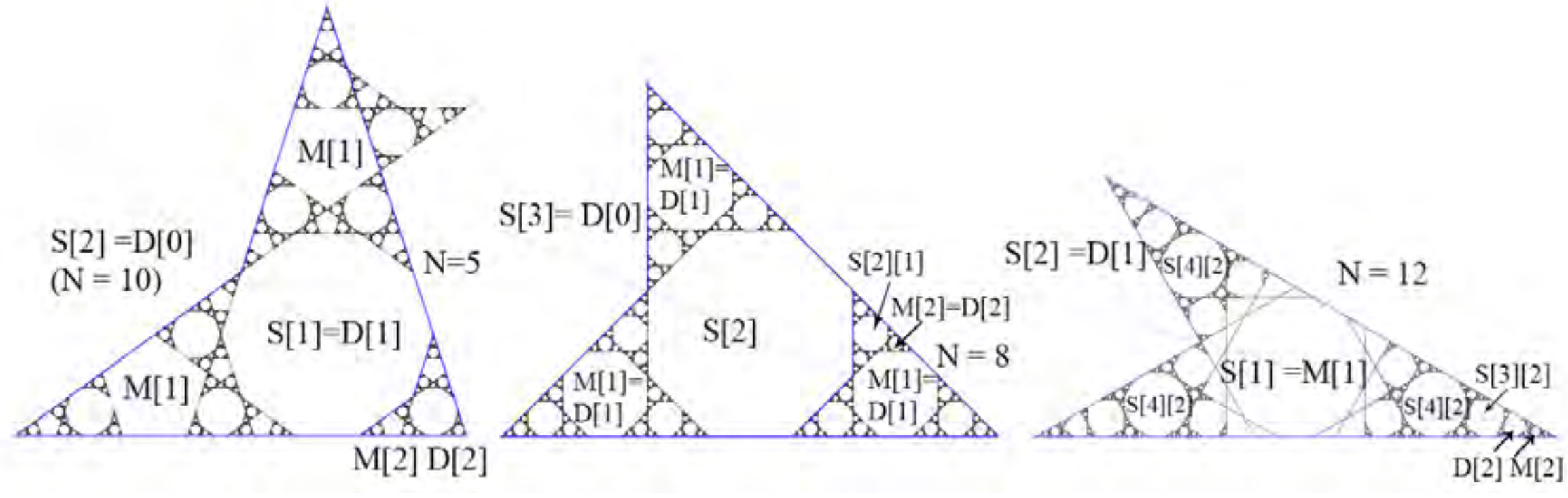

(i) As with N = 5 described earlier, each self-similar 'dart' is anchored by a D tile. There are 2 components to each dart and they are anchored by same-generation M tiles. Locally each of these M tiles is surrounded by 3 next-generation D tiles, so the D tiles will scale by 6 in the limit.

(ii) N = 8 is twice-even so the D[k] and M[k] are identical and we will just describe the D[k] converging to star[1] of N. The primary dart is anchored by S[2] and each of the 3 components are anchored by D[1]s. Each D[1] is surrounded by 3 S[2][1]s and each S[2][1] accounts for 3 D[2]s, so the Ds scale by 9 in the limit (with global periods 16 (D[1]), 96,1008, 8640,79056...)

(iii) N = 12 is also twice-even and each dart is anchored by an M tile and there are 3 components. Each component is anchored by a next-generation S[4] and each of S[4]s is surrounded by 3 same-generation S[3]s so there will be 9 S[3]s, and in the limit each S[3] will account for 3 M[2]s, so the Ms scale by 27 in the limit. (Here the combined S[3]s of each component only accounts for 6 M[2]s but that is because M[1] has a one-time mutated web).

Therefore the Hausdorff-Besicovitch fractal dimension of the three webs are:
(i) N =5: Log[6]/Log[1/GenScale[5]] ≈ 1.2411 where GenScale[5] = Tan[π/5]Tan[π/10]
(ii) N = 8: Log[9]/Log[1/GenScale[8]] ≈ 1.2465 where GenScale[8] = $\text{Tan}[\pi/8]^2$
(iii) N = 12: Log[27]/Log[1/GenScale[12]] ≈ 1.2513 where GenScale[12] = $\text{Tan}[\pi/12]^2$

It is no surprise that these are increasing, but this applies only the quadratic family. For the cubic family and beyond, the webs are probably multi-fractal, with a spectrum of dimensions. However it is possible that the maximal Hausdorff dimension will increase with the algebraic complexity of N, with limiting value of 2. (The M[k]-D[k] of the 8k+2 family cannot yield any fractal dimension beyond N = 10 because Log[N/2+1]/Log[1/GenScale[N/2] is decreasing.)

**•N = 13**

N = 13 has algebraic complexity 6 and is in the 8k+5 family. The salient feature of this family is the existence of a DS[1]. Because of the Rule of 8, this DS[1] is typically isolated and the question is whether it can support a next-generation family on the edges of S[2]. Since N is odd the DS[k] are typically not part of the First Family of S[2] but DS[1] is an exception because by definition it is the 'D' tile of the virtual S[1] of S[2] and therefore it is congruent to an S[2] of S[2] and a valid member of what we call the 3rd generation. The question is whether DS[1] has an edge geometry that can support a next-generation. We will discuss this issue below, but first we present an overview of this second generation on the edges of N.

**Figure 13.1** Overview showing lines of symmetry

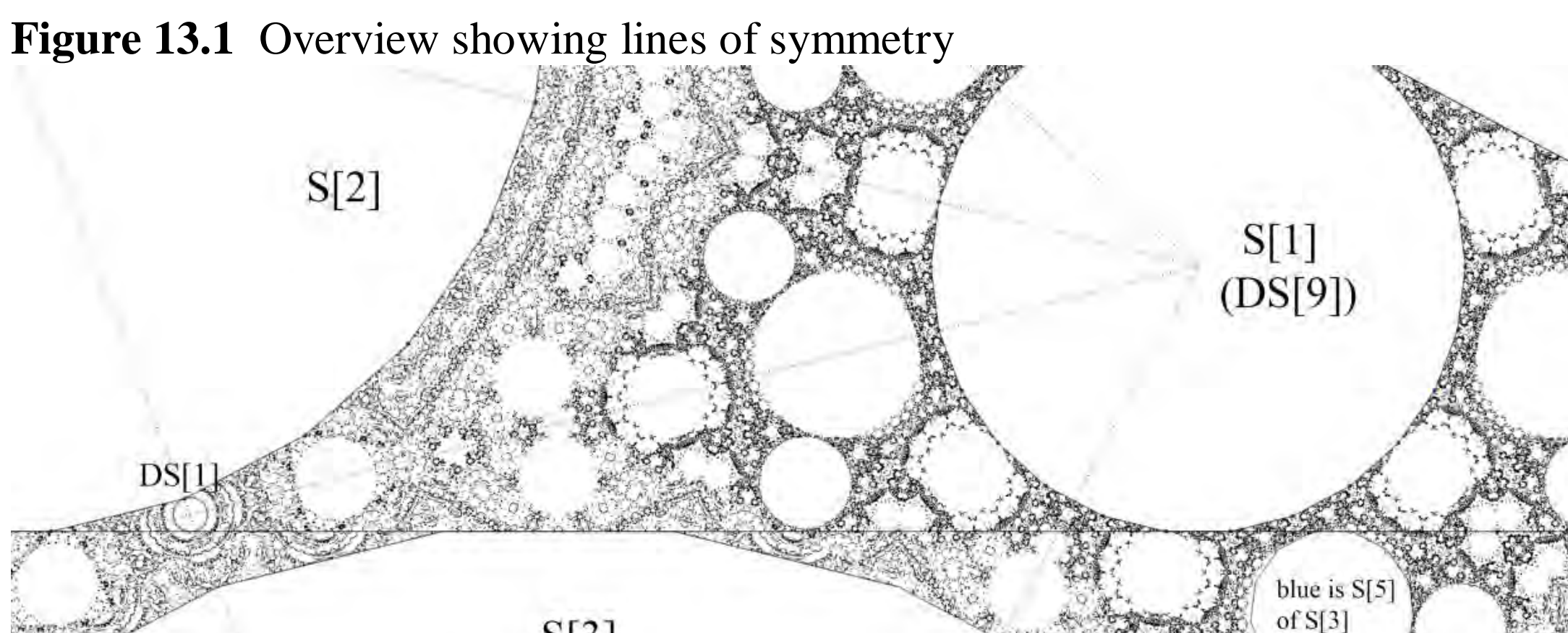

There is no doubt that the geometry here is strongly influenced by S[3], but there are no shared tiles except for a close match with S[5] of DS[3] shown here in blue. Of course this match would persist for the copies rotated around S[1], but this is largely independent of S[2]. Below is the geometry local to DS[1] where star[N-2] must match star[2] of S[2]. As expected the local web is step-2 and there are vertex-based copies of S[4] tiles at step-2 intervals. These displaced S[4]s have step-2 webs of their own and these webs have a promising tile structure.

**Figure 13.2** Detail of DS[1] which has a step-2 web

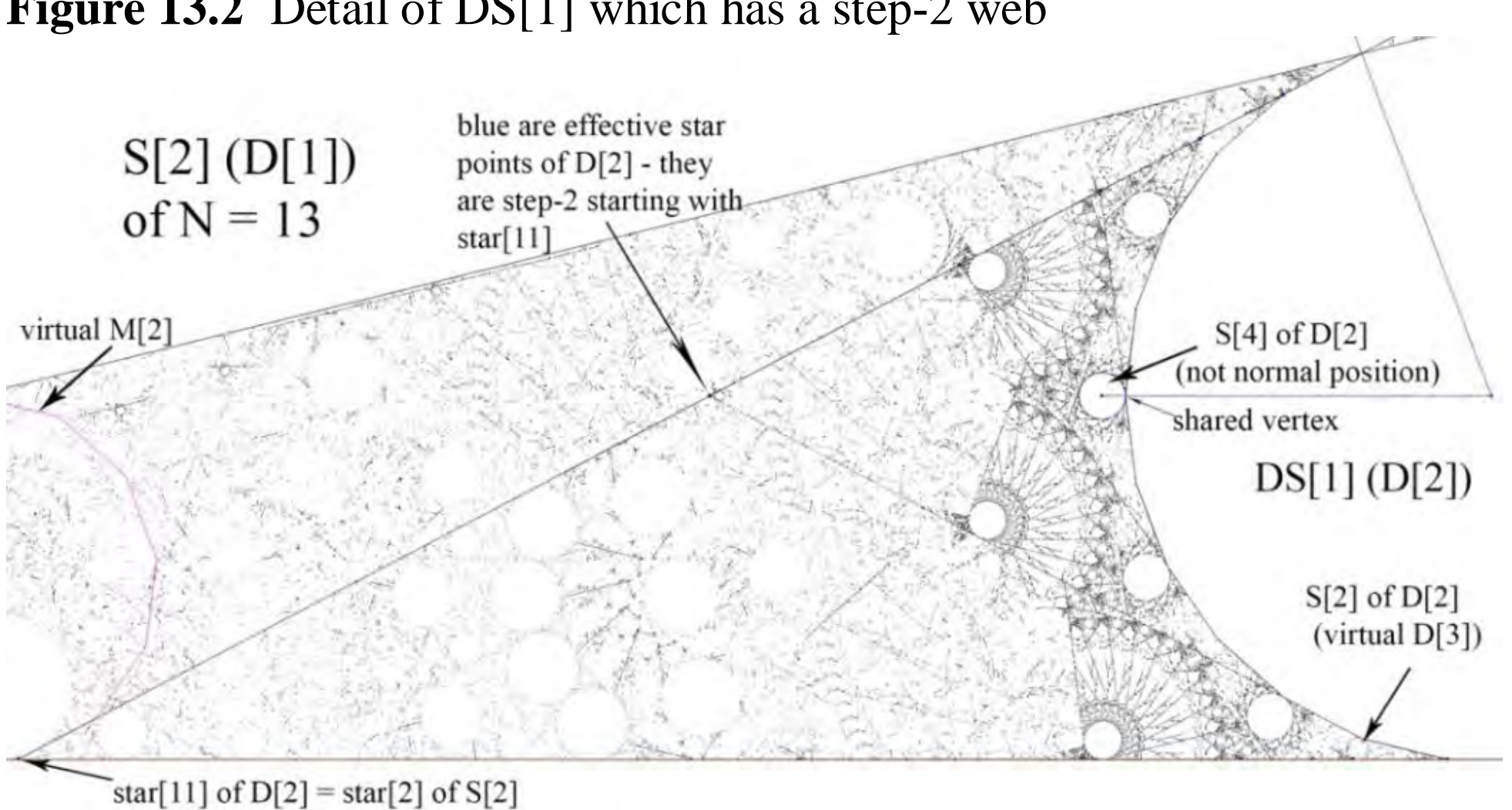

**•N = 14**

N = 14 has complexity 3 and is the first non-trivial member of the 8k+6 family. To trace the evolution here it is useful to take a step backwards and look at the First Family – which by definition is also the First Family of N = 7.

**Figure 14.1** The First Family of N = 14 and N = 7.

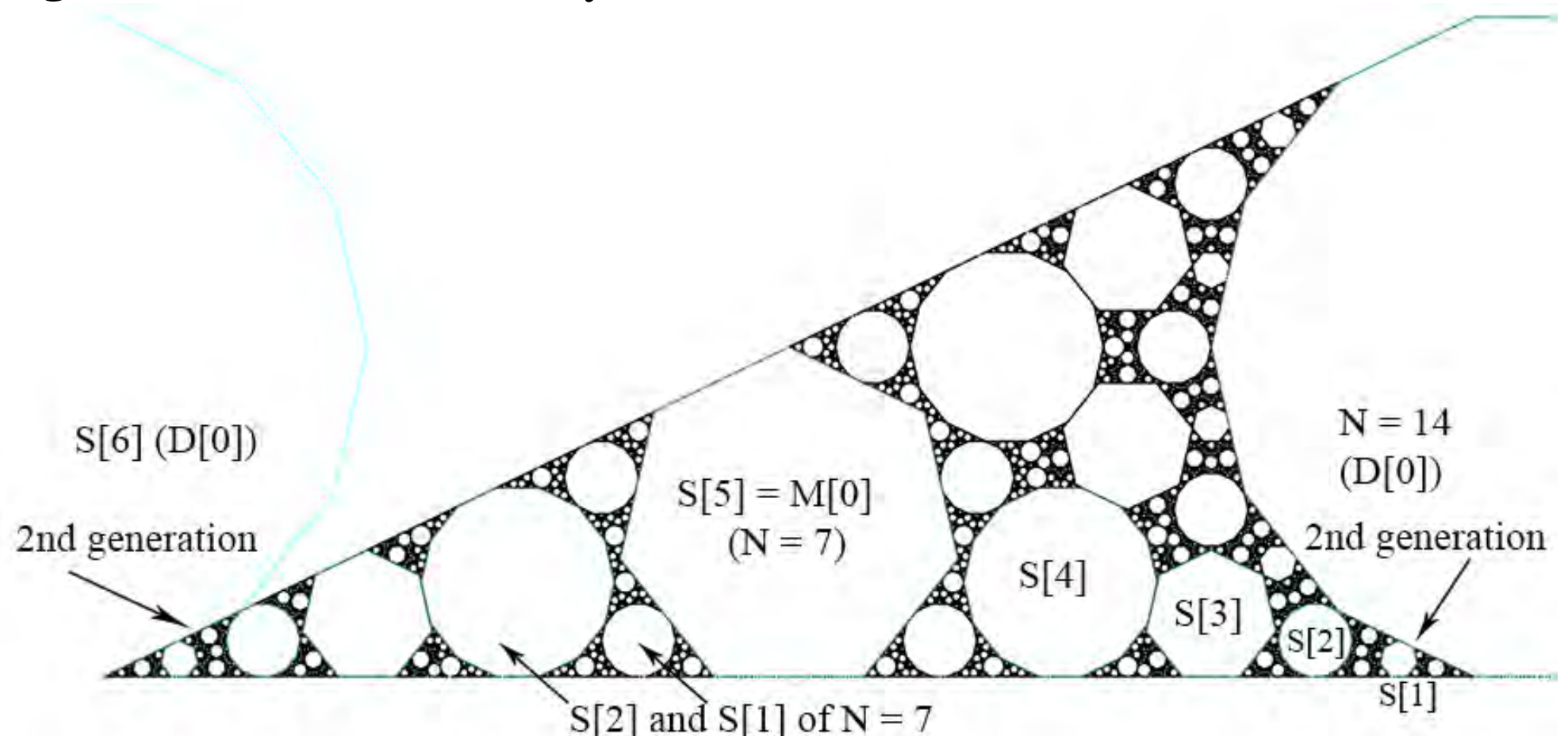

Even at this coarse resolution it is clear that the second generation presided over by S[2] is very different from the first. Of course the orientation is reversed, but the large S[4] tile which should be shared by S[2] and S[1] is missing. As shown below, it is replaced by a 'volunteer' which we call a Portal M tile because it is an N/2-gon with an invariant mantle. These PMs have a divided allegiance between S[1] and S[2]. Their characteristic equation is discussed below.

The Rule of 4 predicts a DS[1] along with S[1] at DS[5]. For N even, the DS[k] are simply the S[k] of S[2], so S[1] is always the penultimate M tile in the First Family of S[2]. Therefore S[2] and S[1] have the same D-M realationship as D and M[0] (N = 7) and it makes sense to call them D[1] and M[1] and acknowledge that it is possible to repeat this recursively with D[1] as the new N. In theory this could generate sequences of D[k[ and M[k] converging to star[1] of S[2] and N.

**Figure 14.2** The 2$^{nd}$ generation showing a blue 3$^{rd}$ generation at star[1] of S[2]. This 3$^{rd}$ generation will be an exact copy of the 1$^{st}$ generation, scaled by hM[1]/hM[0] =x =GenScale[7].

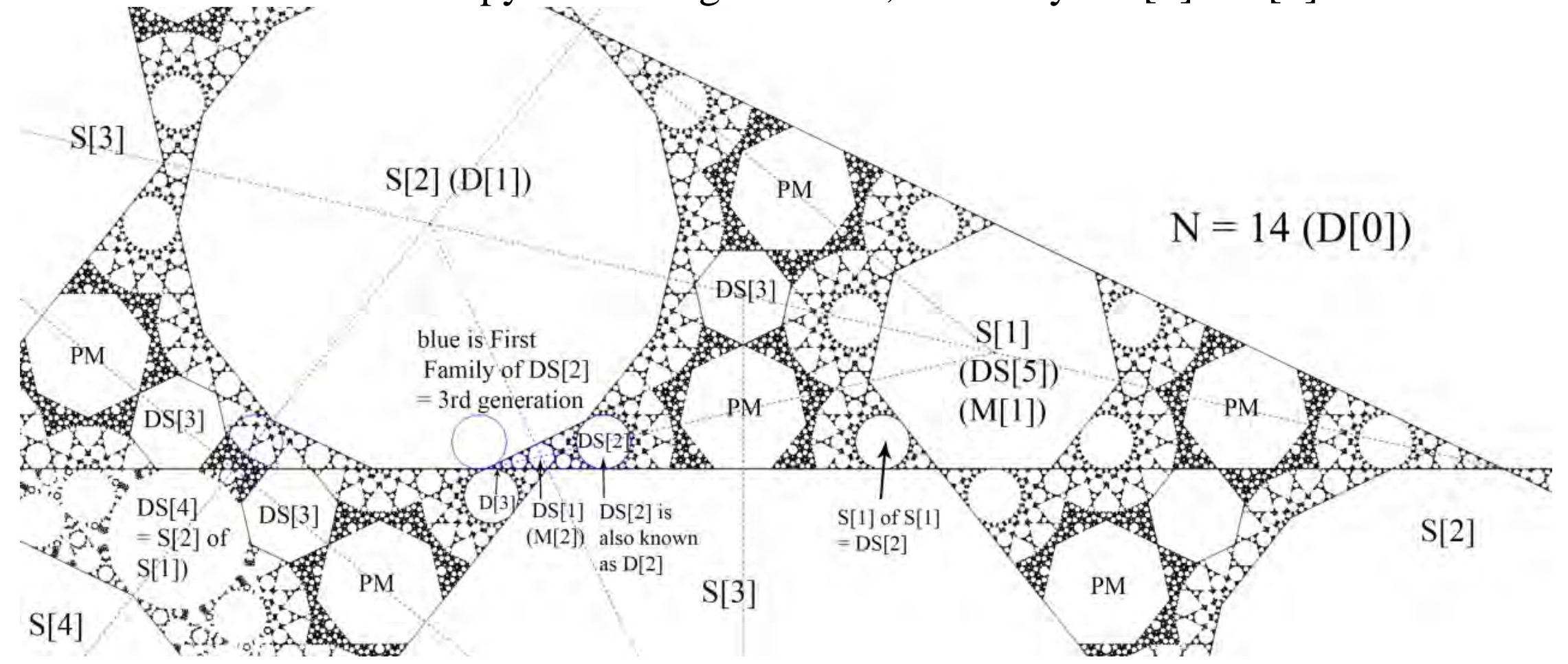

Therfore M[k] and D[k] could serve as patriarch and matriarch of generation[k+1]. For the 8k+6 family, the predicted DS[1] can serve as the M[2] in this sequence and the question is whether a volunteer D[2] will exist. Based on cases like N = 22 to follow, this D[2] may not exist. But for N = 14, D[2] tiles do exist and, like the 8k + 2 family, they develop along with DS[3]s. This implies that there may indeed be squences of M[k] and D[k] converging to star[1] of S[2],

This matches the convergence in the 8k+2 Conjecture, but N = 14 is the only known case where there is stict alternation of self-similar even and odd generations. This is no doubt a signature of the cubic nature of N = 14. These odd and even sequences have identical geometric scaling of x = GenScale[7] = Tan[Pi/7]·Tan[Pi/14], but the two sequences appear to have distinct temporal scaling of 8 (N/2+1) and 25 respectively so this appears to be a multifractal sequence. The blue 3rd generation shown above is scaled by hM[2]/hM[0] which is $x^2$ and M[2] has $\tau$-period 98. The first few M[k] periods are 14, 98, 2216, 17486, 433468, 3482794, 86639924, 396527902.

The DS[3] tiles here are volunteers which form in conjunction with the PM tiles. These PM tiles only arise in even-generations so these are also called PM[2]s. For N = 14, the S[2] tile of S[3] is D[1] but the S[1] tile of S[3] does not exist and in its place are the PM-DS[3] pairs. These PM sit symmetrically on the edges of S[3] between DS[2] and the conguent S[1] of S[1] so their parameters are easy to find. The close connection to S[1] is supported by the relative polynomial **AlgebraicNumberPolynomial[ToNumberField[hPM/hS[1],GenScale[7]],x]** = $(1+x^2)/2$.

There is a dual orthogonal convergence at left-side star[3] of D[1] with alternating PM and DS[3] tiles on the left and alternating virtual and real D[k] on the right. The first D[2] in this sequence is virtual as shown in Figure 14.1 above. This yields a four generation sequence with local geometry which is a mixture of single-scale self-similarity in blue and multi-scale dynamics in magenta as shown in the insert. The temporal scaling here appears to be 113 and this is consistent with the convergence at star[2] of N.

**Figure 14.3** The dual convergence at star[3] of S[2] – which is also star[2] of S[3]

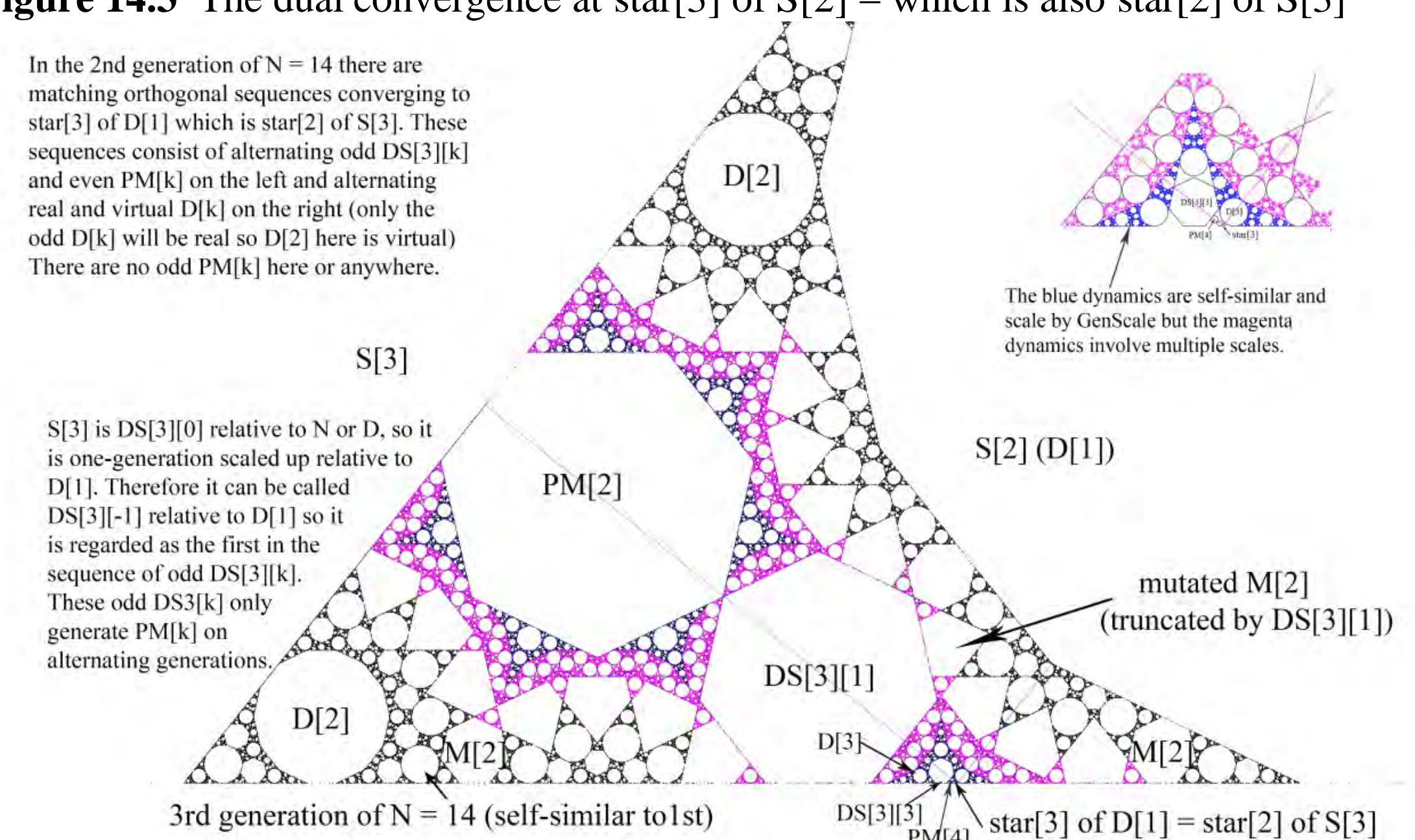

## • N = 15

N = 15 is really the 2nd member of the 8k+7 family because N = 7 can be regarded as the charter member. N = 15 has quartic complexity along with N = 16, 20, 24 and 30.

Below is an overview of the First Family of N = 15. As N grows S[3] and S[4] align with S[1] and S[2] respectively, but here the geometry is dominated by the S[3]-S[2] interaction, where star[5] of S[2] is equal to star[3] of S[3], and DS[3] which plays a major part in the local geometry of S[2], is also an S[4] tile of S[3]. This DS[3] will be mutated but first we look at mutations in S[5] and S[3] of N. For N-odd, the webs steps are doubled from N/2-k to N-2k, so S[5] will have $k' = 15\text{-}10 = 5$. The S[k] are 2N-gons, so the mutation should be the weave of two 30/gcd(30,5)-gons as shown below. S[6] is not shown here but it is also mutated since $k' = 15\text{-}12$ and 30/gcd(30,3) = 10. This is similar to S[3] with $k' = 9$ and 30/gcd(30,9) = 10, so both S[3] and S[6] will be the weave to two decagons in barely perceptible mutations.

**Figure 15.1** An overview of N = 15 showing mutations in S[3] and S[5] (S[6] is also mutated)

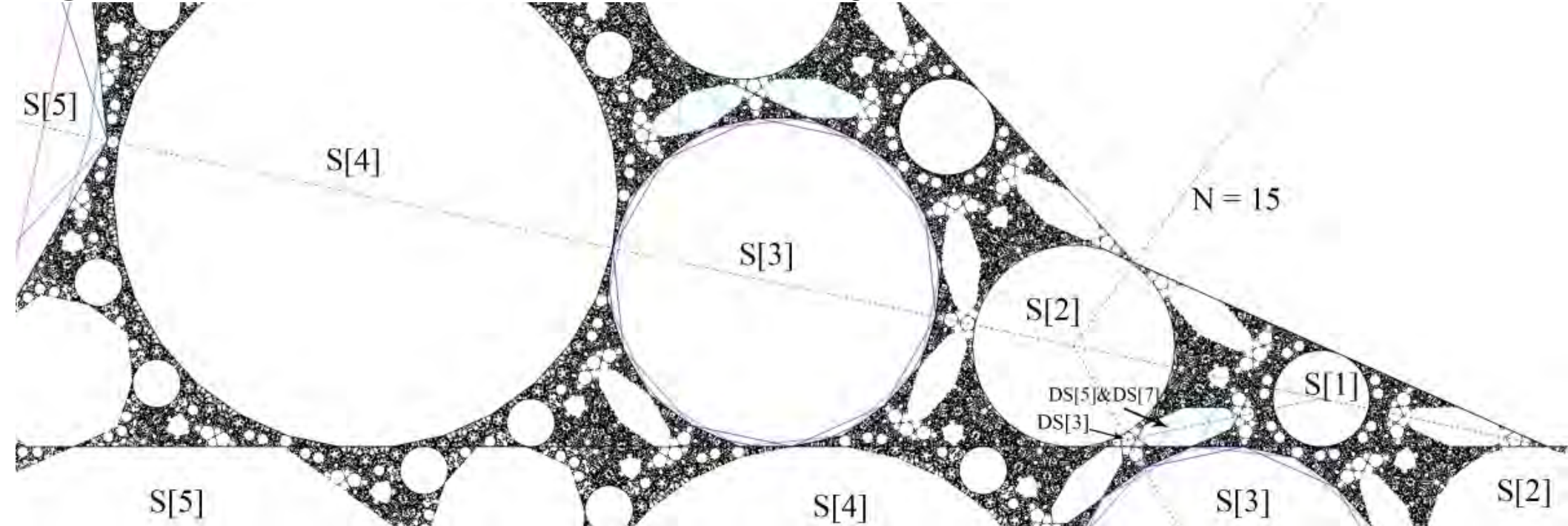

A salient feature of the 8k+7 family is the existence of a DS[3] which will co-exist with dual DS[1]s to form a step-4 cluster on the edges of S[2], anchored by DS[3]. This cluster is bounded on the right by star[5] of S[2] and as noted above this star point is shared with S[3]. DS[3] is mutated since S[2] is twice-odd and the web steps for the DS[k] are $k' = N\text{-}k$. Therefore S[3] has $k' = 12$ with 30/gcd(30,12) = 5 and DS[3] is the weave of two pentagons as shown below. DS[5] is virtual but if it existed it would be the weave of two triangles since $30/\gcd(30,k') = 3$. These adjacent mutations together with the web interaction of S[2] and S[3] are probably contributing factors for the dodecagon Dx tiles. N = 18 has a similar geometry. See Figure 18.2

**Figure 15.2** The DS[3] region is dominated by mutations with 3 and 5 as adjacent odd factors.

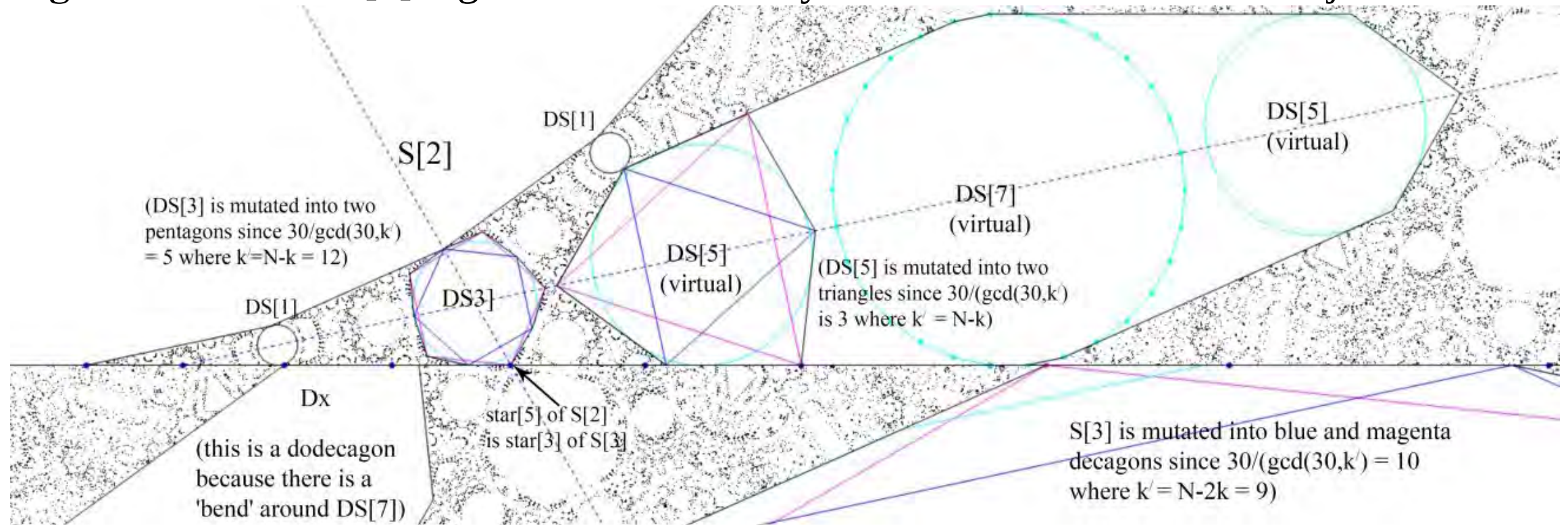

For the 8k+7 family the volunteer DS[1]s will always be S[N-3] tiles in the First Family of the underlying DS[3] and it is not unusual for DS[3] to foster local tiles as well. DS[3] will always have a step-6 limiting web since k′ = N-3 and Mod[2N,N-3] will be 6 for N >7. Here DS[3] has ‘normal’ S[1] and S[6] tiles and it seems that S[1] tile has its own S[6]. In a manner similar to the S[2] mutation in the 8k+4 family, the most interesting local geometry of DS[3] occurs at the inner blue vertices which are also vertices of the underlying DS[3]. At the magenta vertices DS[3] shares extended edges with the S[6]s which are similar to the S[4]s in the 8k+4 mutations.

**Figure 15.3** Detail of the web local to S[2]

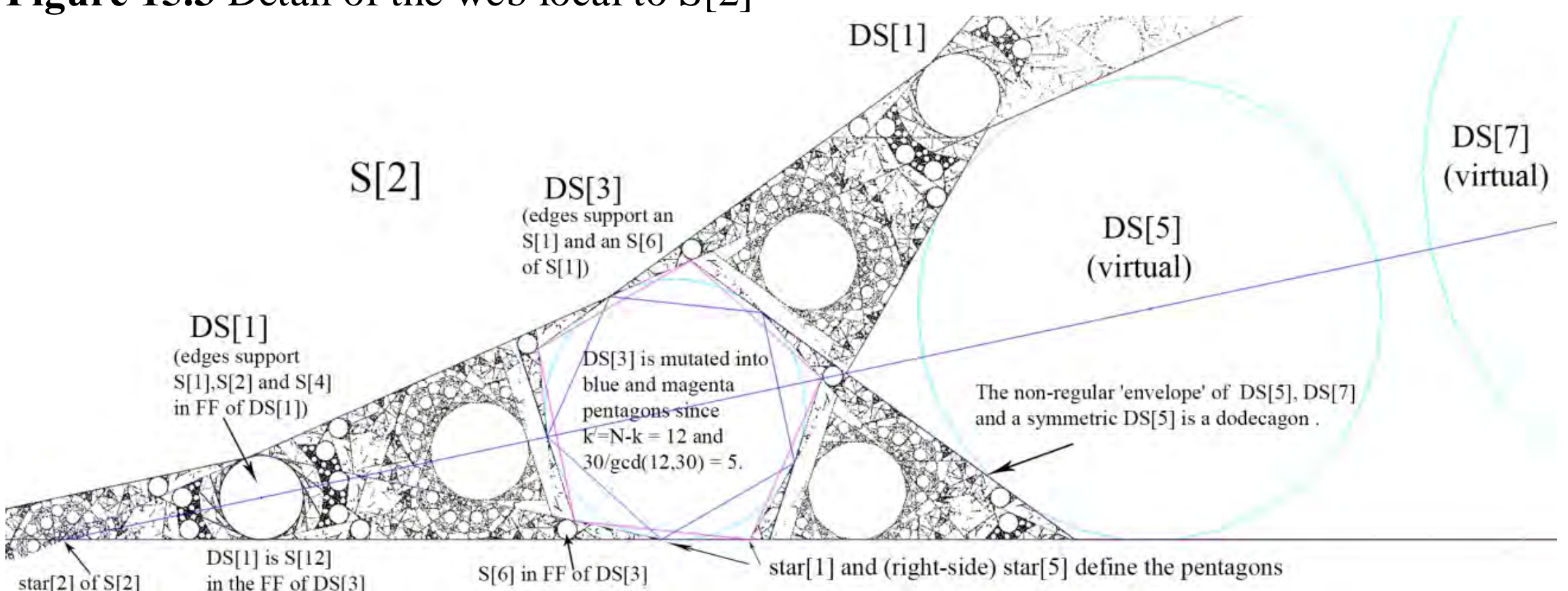

DS[1] has maximum k′ = 14 with gcd(30,14) = 2 but for N odd the mutation condition for the DS[k] is gcd(2N,k′) > 2, so DS[1] is not mutated. Since Mod[2N,N-1] is always 2, DS[1] will have a limiting step-2 web as predicted by the 8k+7 Conjecture. These step-2 webs typically have good support geometry and here it appears that the S[1], S[2] and S[4] survive the web. In a very gender specific environment, these S[1] are among the few N-gon survivors.

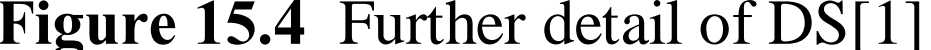

**Figure 15.4** Further detail of DS[1]

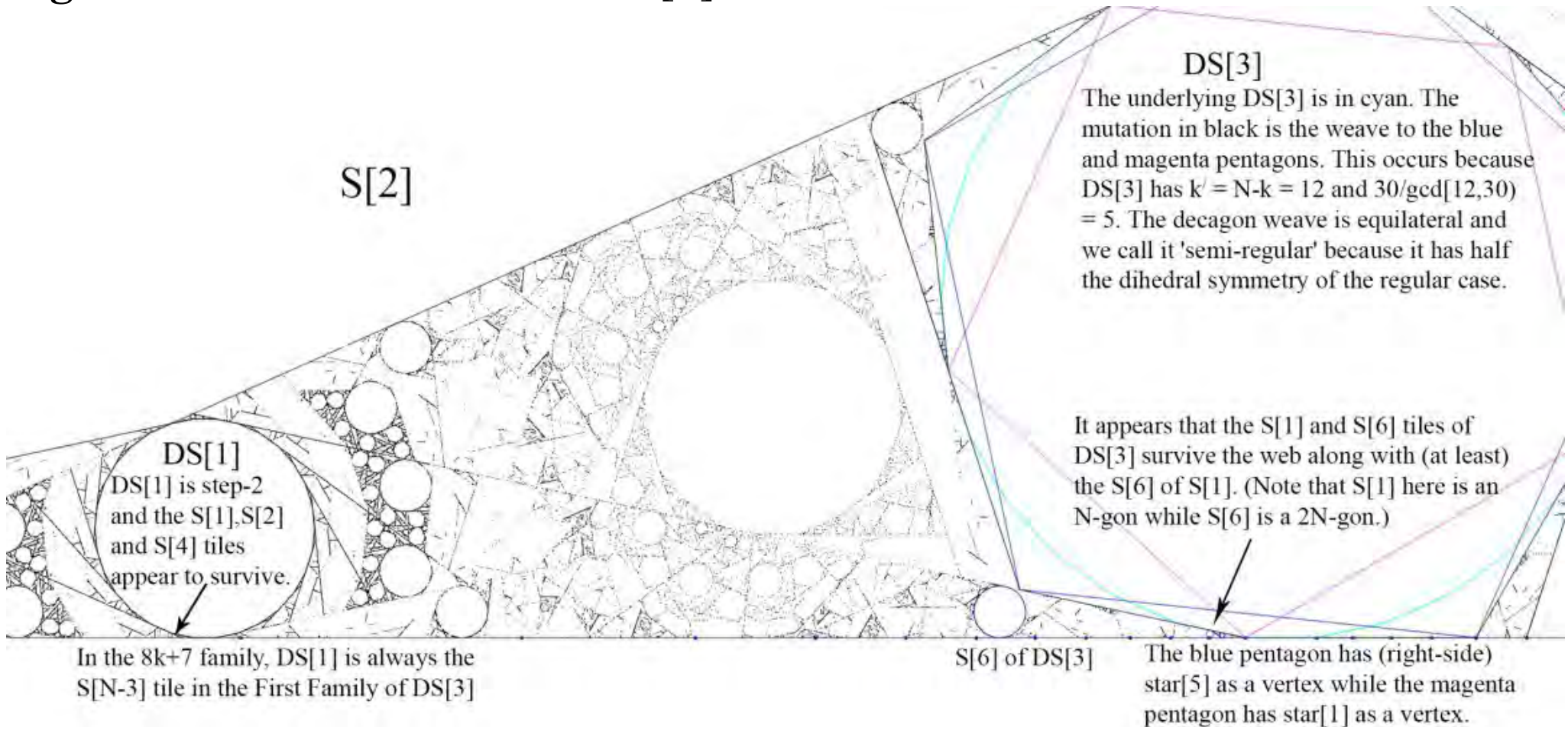

The DS[1] survivors seem to benefit from the step-2 web combined with star[2] being ‘effective’. This occurs because the count-down starts from star[N-2] of DS[1] which matches the star[2] reference point of S[2].

**• N = 16**

N = 16 has ‘quartic’ complexity along N = 15, 20, 24, and 30. Both N = 16 and N = 24 are in the 8k family so the Rule of 4 implies that after S[1] at DS[6] there will be an isolated DS[2] which is an S[2][2]. In the case of N = 16 this tile can function as a D[2] and even though there is no M[2], there will be a chain of D[k] converging to star[1] of S[2] and the matching star[1] of N. This convergence will be discussed below, but here we look at the overall geometry of this region and discuss the influence of S[3] and S[4]. Since S[2] and the rotated S[4] share a vertical ‘tower’, they will share their star[3] points when N is twice-even. This is still true here even though S[4] is mutated. Likewise S[1] shares its star[2] point with S[3].

**Figure 16.1** - Overview of the edge geometry showing the influence of S[3] and S[4]

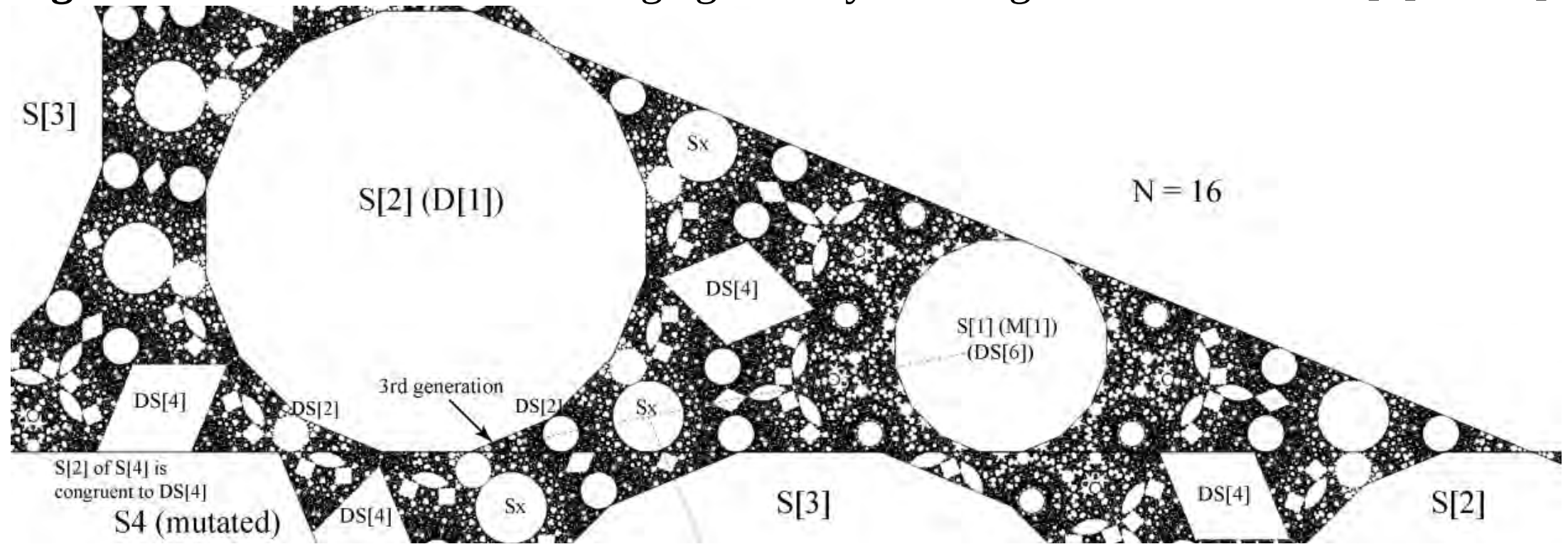

S[4] is mutated since $k' = N/2-k = 4$, and gcd(4,N) = 4. Therefore the mutated S[4] will consist of the weave of two squares $M_1$ and $M_2$. The Mutation Conjecture says that $M_1$ and $M_2$ will be based on star[1] and right-side star[3] of the underlying S[4]. This is illustrated on the left below using DS[4] as a stand-in for S[4]. MuDS4 in black is a scaled reflected copy of the mutation of S[4] shown above. The revised Mutation Conjecture predicts that the DS[4] of S[2] will share the same mutation as S[4] itself, but for N = 16 it seems that the interaction with S[4] turns this into a simple extended edge mutation. These ‘lazy’ mutations are actually very common as witnessed by N = 11,15 and 18. They are not ‘canonical’ riffles or weaves of two regular polygons but if 2-gons are allowed this DS[4] mutation can be regarded as a 2-step canonical mutation. Since MuDS4] will have step-8 symmetry it has an effective $k' = 8$ and N/gcd(8,N) = 2 so MuMuDS4 can be regarded as the rhombus weave of two ‘tuples’ shown in magenta and blue below.

**Figure 16.2** - The two-step mutation of the (virtual) right-side DS[4] tile of S[2]

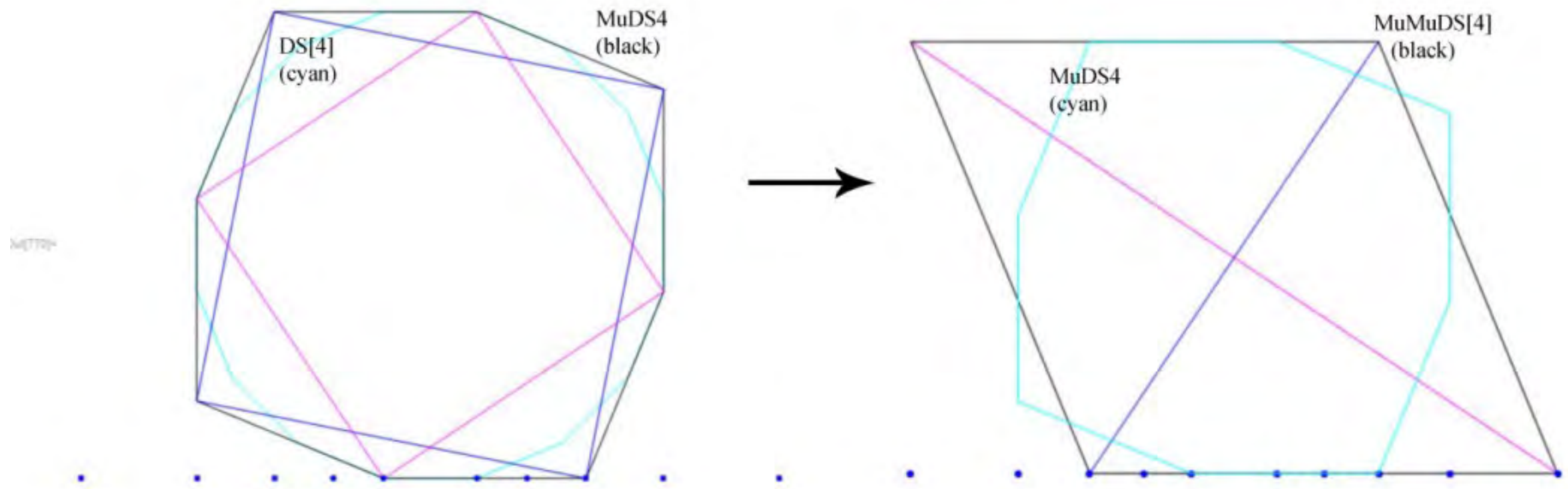

The remaining DS[4]s shown here are rotations of MuMuDS[4] about cS[2]. However the DS[4][2] and DS[4][3] in Figures 16.4 and 16.5 below have just the normal MuDS4 mutation.

The rotated S[3] clearly interacts with both S[2] and S[1]. The Sx 'volunteers' shown above are weakly conforming to S[2] but they are located symmetrically with respect to the edges of S[3] so their centers are aligned. This makes it easy to find their parameters. Sx may be responsible for the copies of DS[2] that surround it but it is also located symmetrically relative to two of these copies of DS[2] and this is reminiscent of N =14 where the PM tile is symmetric with respect to DS[2] and the congruent S[1] of S[1]. For N = 18 the DS[2]s are replaced with DS[4]s and this raises the question about the expected geometry between congruent tiles. For the canonical First Family spacing from D to D, the midpoint is an M tile so Sx and these PM could be regarded as 'local' M tiles. Indeed the characteristic polynomials for the PM are more closely aligned with S[1] than any other major tile. But the D-M symmetry breaks down when N is twice-even. D is no longer in the First Family of M, and S[1] is largely independent of S[2].

For Sx the characteristic polynomial relative to S[1] is slightly more manageable than the others

**AlgebraicNumberPolynomial[ToNumberField[hSx/hS[1],GenScale[16],x]** $= \frac{3}{4} + 28x - \frac{55x^2}{4} \pm \frac{x^3}{2}$

However Sx has a linear relationship with the neighboring DS[2] tiles since hSx/hDS[2 ] = 2 –x where x = GenScale[16] = hDS[1]/hS[2] = $\text{Tan[Pi/16]}^2$. Therefore with the height 1 convention for N, hSx = hDS[2]·(2- GenScale[16]) = $(\text{Tan[Pi/8]}\cdot\text{Tan[Pi/16]})^2\cdot(2 - \text{Tan[Pi/16]}^2) \approx .0133084$.

The Twice-even S[1] Conjecture says that the (possibly virtual) S3x tile in the step-2 family of S[1] will also be an S[2] relative to S[3]. S[2] itself is not in the First Family of S[1] but it is always in the step-2 family. Here it is S7x where S[7] is the (virtual) D tile of S[1].This occurs because the ratio of the edge lengths of S[1] and S[2] is scale[2].

**Figure 16.3** The web local to S[1] showing the early level-8 web and First Family of S[1] in magenta. The effective star points of S[1] are step-2 counting down from star[N/2-1] which is star[1] of S[2] – so the effective star points are just the odd integers bounded by [N/2-1].

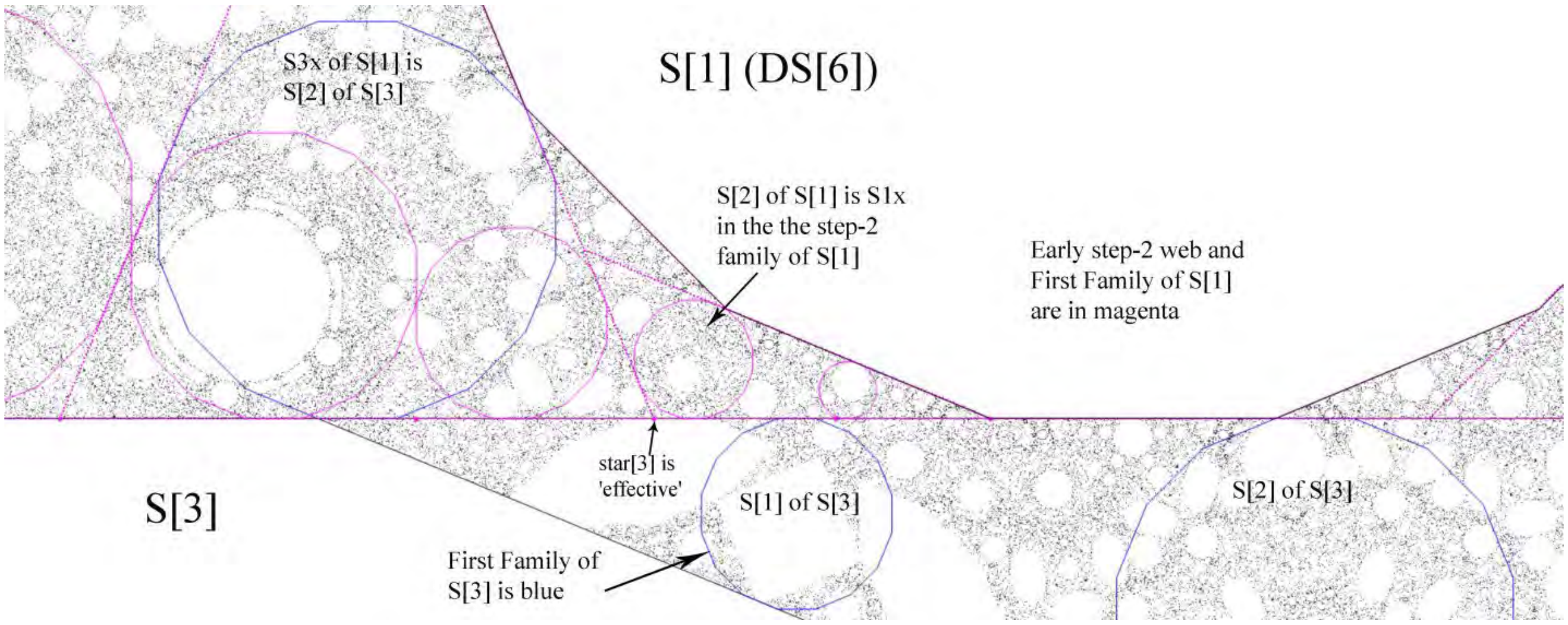

Typically these step-2 webs are not very supportive of First Family tiles and there are no obvious survivors here from the magenta First Family of S[1].

Returning to S[2], note that the DS[2]s evolve with two very different local geometries. The 'even-vertex' DS[2]s like the DS[2] shown on the right below, have minimal local structure at this scale, but the matching 'odd' cases which are reflections of the (virtual) D tiles, offer a haven for family tiles. In the blue First Family of D[2], M[2] is virtual but M[3] is real. This sequence of D[k] tiles and their virtual clones inside S[2] will converge to star[1] of S[2].

**Figure 16.4** - DS[2] and its clone at star[1] of S[2]

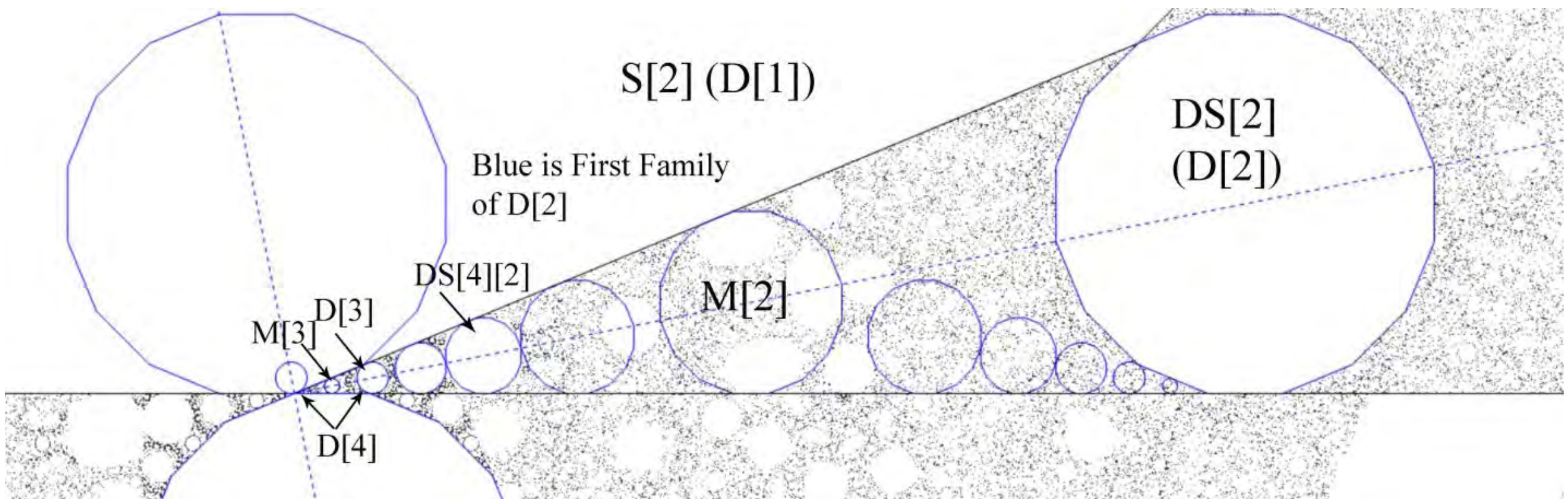

Even though M[2] is missing, the First Family of the surrogate D[3] contains a 'normal' M[3] and D[3] as well as matching S[3][2] and (mutated) S[4][2]. This is what we call the 3rd generation, presided over by D[2]. The geometry is very different from the first or second generations with just these four neighboring survivors. The detail below shows a unique 4th generation evolving local to D[3], but without M[4]s. Dual D[4]s survive and appears to foster D[5]s – which also exist in reflective form at star[1] of S[2]. The S[4][k] all appear to share the same mutation of S[4] and their local geometry appears to be consistent with D[3] and S[3][3] even without the 'towers' of the First Generation.

**Figure 16.5** Detail of the D[3] region showing the 4th generation which is not similar to any previous generation but it appears that sequences of D[k] will survive – probably retaining the even-odd distinction of the DS[2] with missing even M[k].

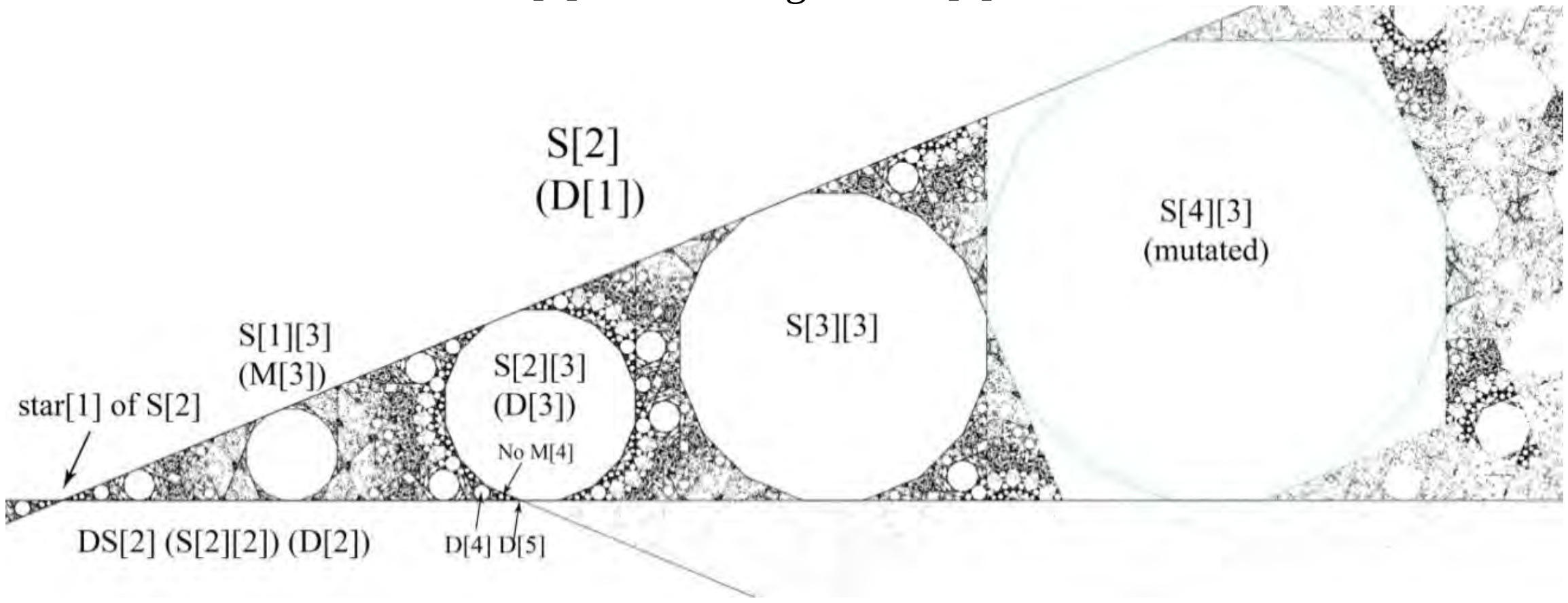

The periods of the first 9 D[k] are: 8, 16·2, 8·57, 16·154, 8·2609, 16·6898, 8·121863, 16·322826. This appears to be a multi-fractal sequence with distinct even and odd temporal scaling in keeping with the 'quartic' nature of N.

● **N = 17**

N = 17 is order 8 and the second member of the 8k + 1 family. Since N = 9 was strongly influenced by mutations, this is an important test case for the family. As predicted by the 8k+2 Conjecture there is an extended family structure at the foot of D and this will be discussed in the context of N = 34. The 8k+1 Conjecture predicts that there will be a volunteer DS[2] to go along with the predicted DS[5]. Below in Figure 17.3 we give a plausible explanation for this fact. Except for the special case of N = 9, this DS[2] will have a step-4 web because k′ = N- k and Mod[2N,N-2] will always be 4. This matches S[1] which is DS[N-4], so here S[1] matches up with DS[2] in the same way that S[1] and DS[1] have matching step-2 webs for N even. This is part of a surprisingly uniform transition from N even to N odd.

**Figure 17.1** The limiting web showing the influence of S[3] and S[4]

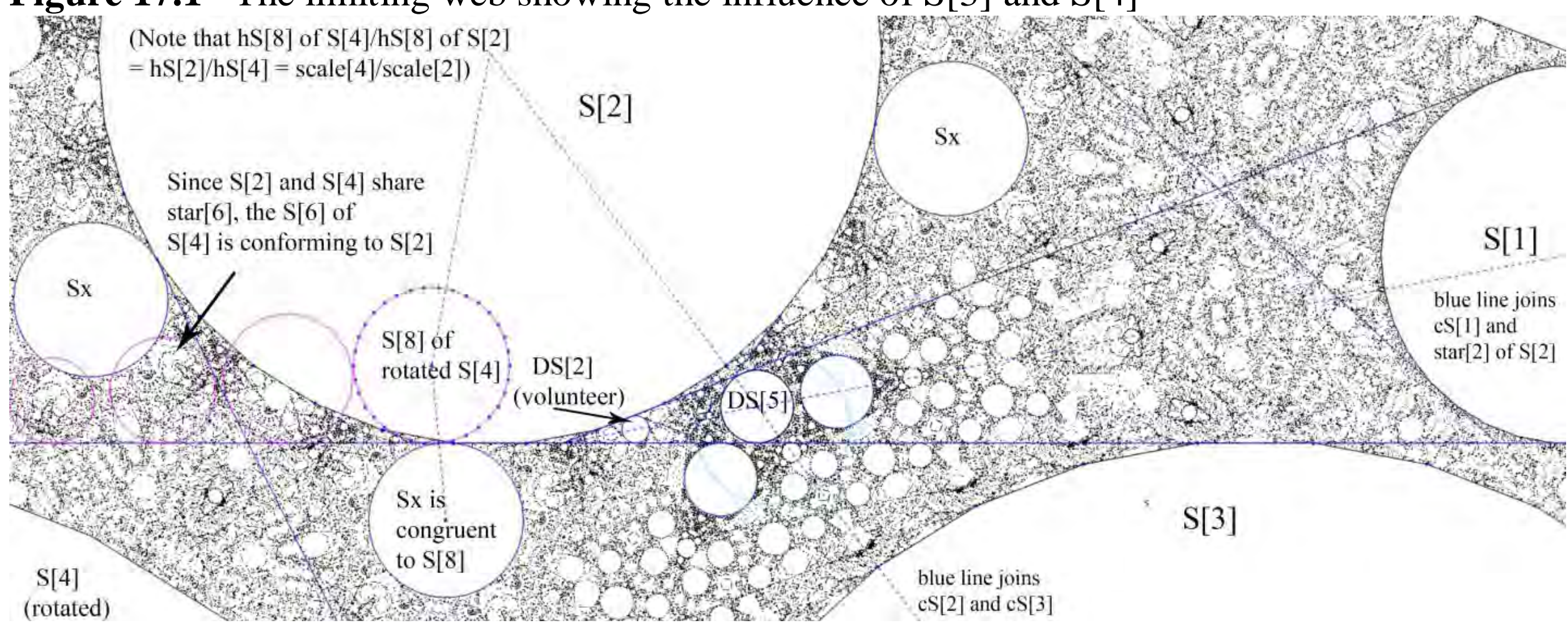

At first glance it seems that the DS[5] clusters might define the large Sx volunteers, but Sx is actually congruent to the S[8] tile of the rotated S[4]. In a manner similar to N = 11 and N = 16, S[3] also plays an important role here. DS[5] is aligned with the centers of S[2] and S[3] and the local geomery is oriented toward S[3].

For N odd, S[2] and the rotated S[4] share their star[6] points so the S[6] tile in the Fist Family of S[4] could play the part of a surrogate tile of S[2]. Here the S[8] of S[4] apparently shares a vertex with the Sx volunteer, so this virtual S[8] can be used to define Sx as shown above. This Sx has an interesting characteristic polynomial relative to S[2]

**AlgebraicNumberPolynomial[ToNumberField[hSx/hS[2],GenScale],x]** yields

$$\frac{3}{2}+\frac{219x}{2}-\frac{2675x^2}{4}+511x^3+446x^4-53x^5-\frac{195x^6}{4}-\frac{7x^7}{2}$$

The relationship between DS[5] and the rotated S[3] seems very robust. Algebraically there is no doubt that DS[5] is more closely aligned with S[3] than with either S[2] or S[1]:

**AlgebraicNumberPolynomial[ToNumberField[hDS[5]/hS[2],GenScale],x]** yields

$$-\frac{21}{64}+\frac{1575x}{64}-\frac{1045x^2}{64}-\frac{2137x^3}{64}-\frac{423x^4}{64}+\frac{237x^5}{64}+\frac{81x^6}{64}+\frac{5x^7}{64}$$ compared with hDS[5]/hS[3] below

$$\frac{1}{16}+\frac{3x}{32}-\frac{309x^2}{16}-\frac{557x^3}{32}-\frac{21x^4}{16}+\frac{73x^5}{32}+\frac{9x^6}{16}+\frac{x^7}{32}$$

Since DS[2] will be conforming to star[2] of S[2], it shares its left side star[N-2] with star[2] of S[2]. Because the count-down of effective star points is mod-4, this implies that the smallest effective star point of DS[2] will be star[3].This is often a good sign since it allows for the formation of a next-generation S[3] or S[2]. This is similar to D acting as N = 34 where this star[3] effective star point of S[2] fosters an S[3] tile that in turn implies the 8k+2 Conjecture.

**Figure 17.3** The star[2] region of S[2] showing the DS[2] volunteer

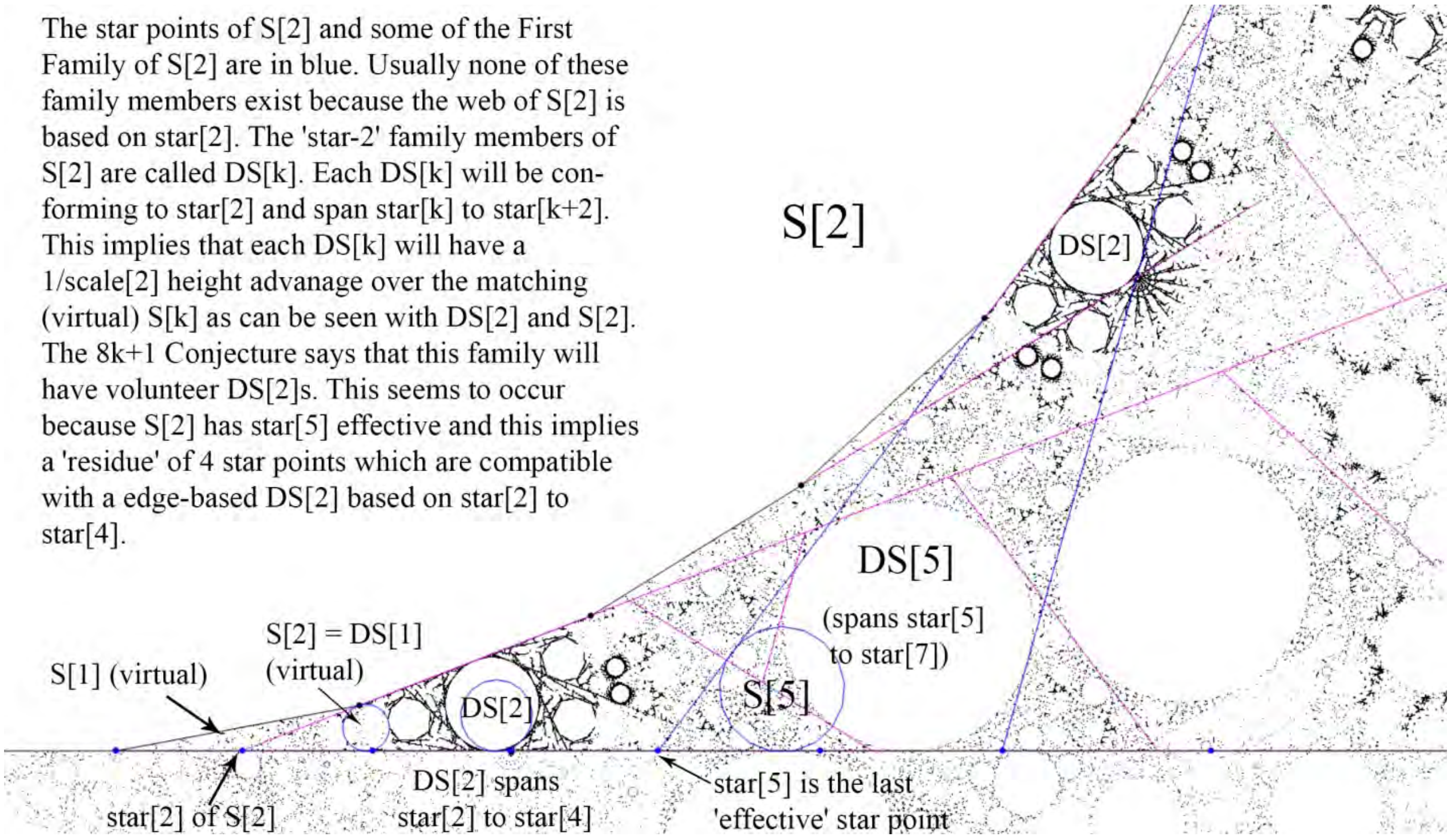

**Figure 17.4** Detail of DS[2] showing the early step-4 web

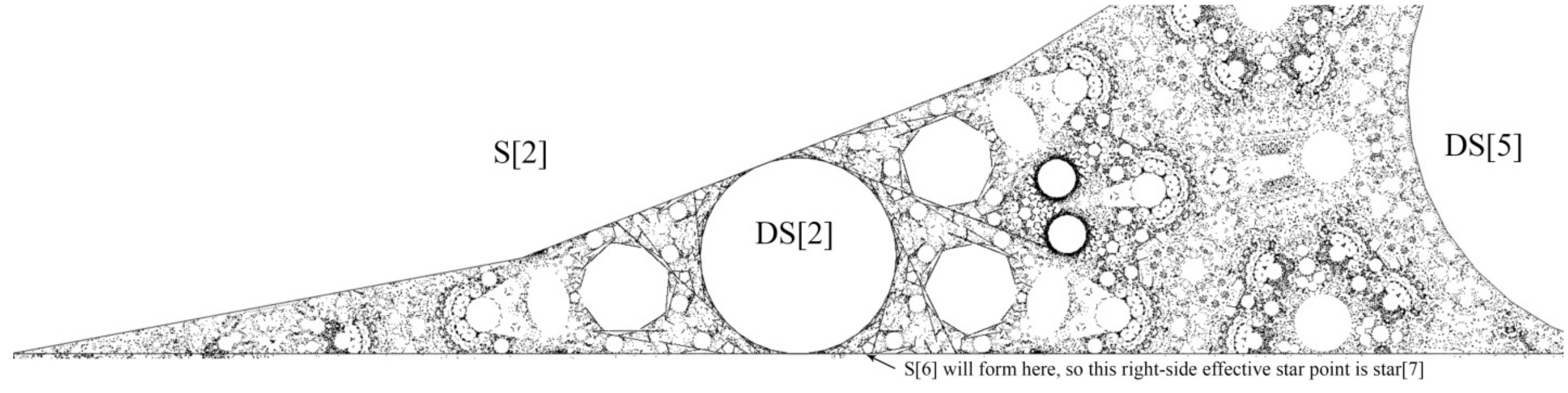

This web is highly fractured but clearly has overall step-4 symmetry as can be seen by the rotated copies of S[6] above and below. This implies that star[7] is effective and star[3] is also effective. This in turn appears to support an S[2] and matching S[1] pairs. This is very similar to N = 25 to follow.

**• N = 18**

N = 18 has cubic complexity and is the 2nd member of the 8k+2 family. This means there should be sequences of D[k] and M[k] tiles converging to star[1] of N and by reflective symmetry there will be an equivalent convergence at star[1] of S[2] acting as D[1]. This does not always imply self-similarity of generations, but for N = 18 it appears that the 3rd generation at the foot of S[2] is self-similar to the 2nd generation and this chain of scaled 2nd generations should continue. Even though the first generation is not be part of this sequence, this is still a strong version of convergence involving generations instead of just D[k] and M[k]. The revised 8k+2 Conjecture says that all the DS[k] predicted by the Edge Conjecture will survive and this leaves open the possibility of self-similar generations, but we believe that in general these generations will vary.

**Figure 18.1** – The edge geometry

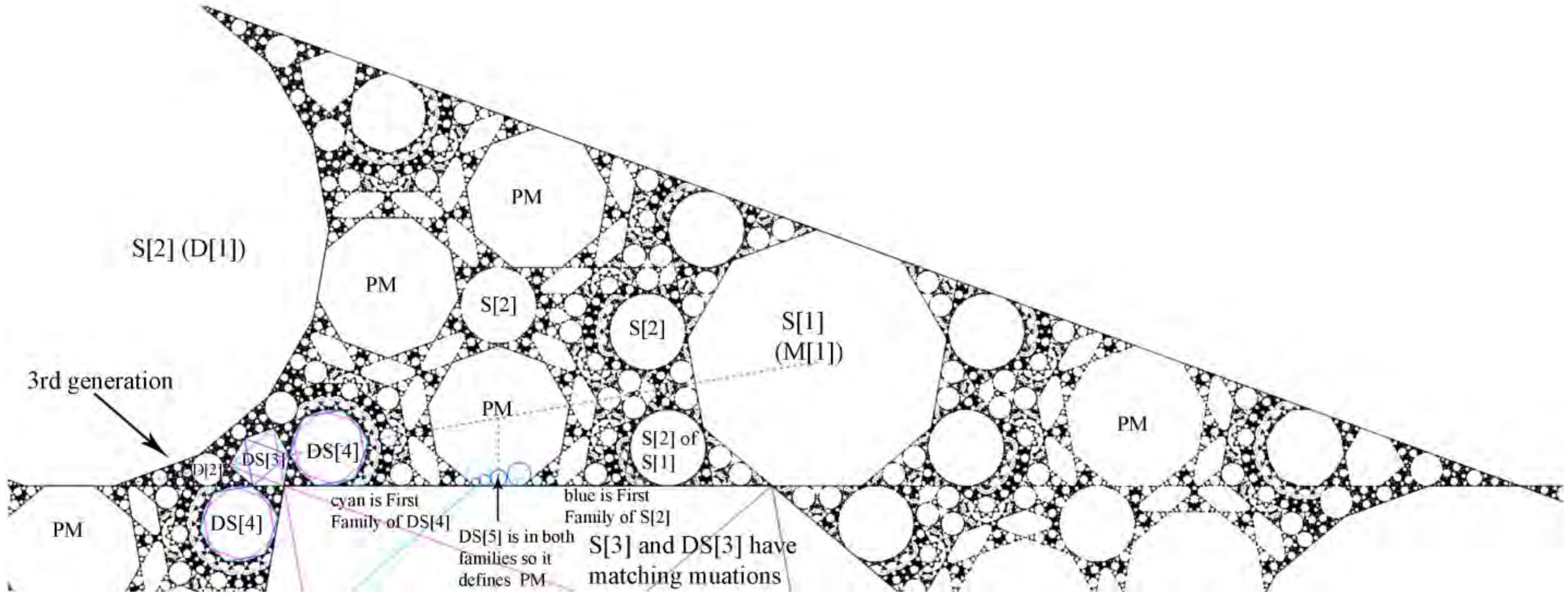

In the twice-odd family, the odd S[k] and DS[k] will be 'mutated' into N/2-gons and it is our convention to display them in this form. The revised Mutation Conjecture of Section 2, predicts that there will be further (matching) mutations of S[3] and DS[3] since both have k′ = 9-3 = 6 with 9/gcd(9,6) = 3. Therefore S[3] and DS[3] will be the weave of blue and magenta triangles as shown above for S[3] (but with reversed orientation because their webs are reversed).

The Mutation Conjecture says that since N/2-1-k′ = 2, the mutation base will be run from star[1] to opposite side star[2] of the cyan underlying S[3] as shown here. (S[4] and DS[4] have k′ = N/2 - 4 = 5 so they should not be mutated but DS[4] has a star[1[ to opposite-side star[2] mutation to match DS[3] but the triangles of the DS[3] mutation are replaced with hexagons. This may be a 'joint mutation', similar to DS[3] and DS[7] of N = 15. DS[4] is congruent to an S[2] of S[1] and these are not mutated.)

The PM tiles are conforming volunteers similar to the Portal M tiles of N = 14, In that case Figure 14.2 shows that PM is symmetric with respect to DS[2] and the matching S[1] of S[1]. The situation is the same here with DS[4] and S[2] of S[1]. Therefore the First Families are identical and the DS[5] tile shown here is in both families, so it defines the midpoint of PM, and the center must be on the blue line from cS[1] to star[1] of S[2]. The characteristic polynomial of PM relative to the 'usual suspects' can provide clues about its algebraic origins.

**AlgebraicNumberPolynomial[ToNumberField[hPM/hN,GenScale[9]],x]** = $\frac{11}{3}-\frac{169x}{3}-\frac{25x^2}{3}$

This improves relative to hS[3] at $\frac{5}{2}-36x-\frac{11x^2}{2}$ and relative to hS[1] it is simply $1-7x-x^2$

Since hS[1] is GenScale[18] =Tan[Pi/18]$^2$, hPM = GenScale[18]$\cdot(1-7x-x^2)$ where x is GenScale[9] = Tan[Pi/18]·Tan[Pi/9] is the default generator of the scaling field $S_9 = S_{18}$.

The elongated octagons which surround PM are collectively called Dx. Our Web Conjecture says that their edges (relative to the edge of N) , should be in the scaling field $S_9$ along with PM. One way to see this is to embed First Family members from DS[4] (or S[2]) as shown above.

**Figure 18.2** The non-regular Dx octagons at the vertices of the PM tiles

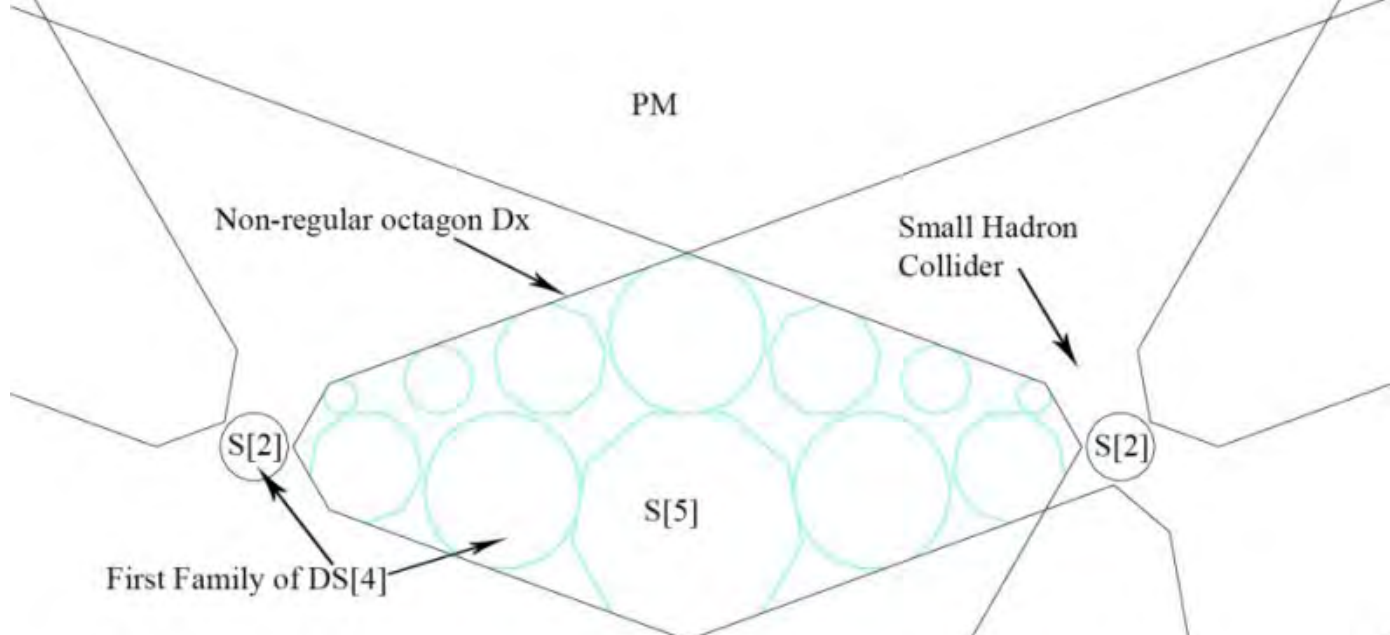

All the embedded tiles are formed from rational rotations (kπ/18) of the First Family tiles of DS4 so every vertex of the Dx octagon is either a vertex of a canonical tile or a star point of a canonical tile so the Dx are 'canonical' tiles. Just like N = 7 and N = 14 the small scale geometry here is a mixture of relatively predictable self-similar geometry and very unpredictable geometry. The region around the Small Hadron Collider is an example of the latter.

**Figure 18.3** The large scale geometry on the left and the small scale on the right

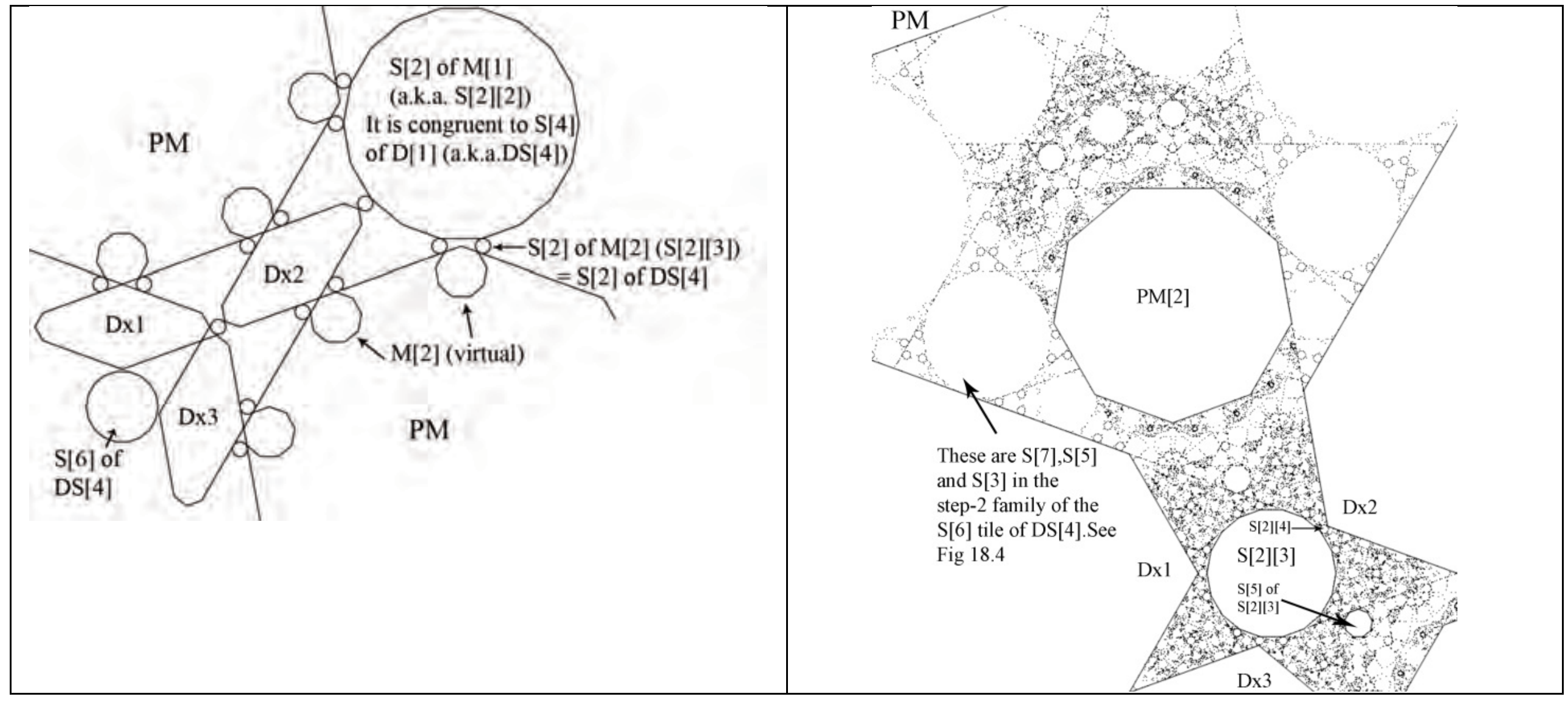

**Figure 18.4** The 3rd generation of N = 18 and its amazing environs

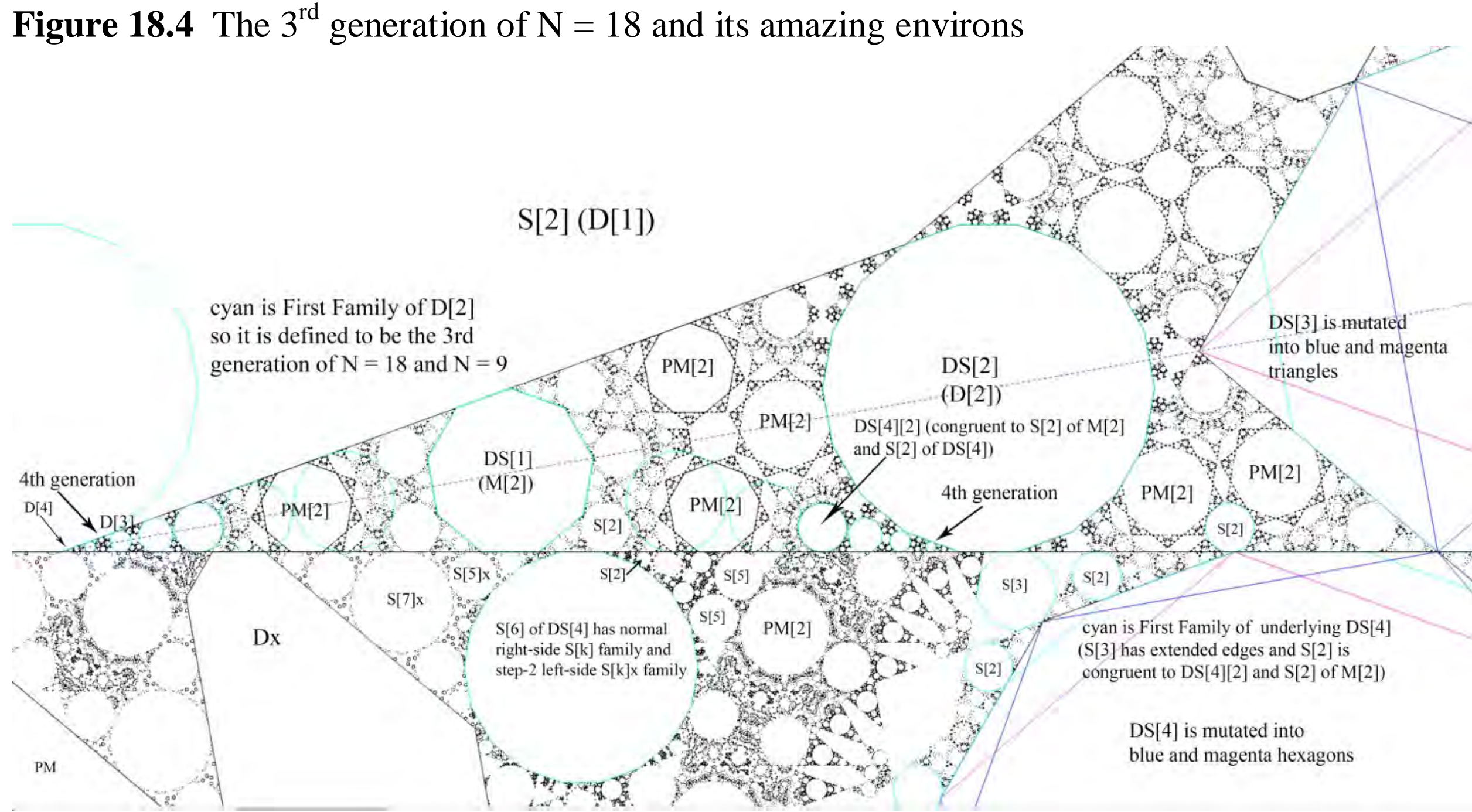

This graphic appears to support the hypothesis stated earlier that all subsequent generations in the convergence to star[1] of S[2] are self-similar to the 2nd generation. Here the 3rd generation scaled by GenScale[9], appears to mimic the 2nd, even though the geometry below the horizontal axis is quite different. Both above and below the horizontal axis the prominent S[6] tile of DS[4] clearly has traditional step-1 First Family edge geometry on one side and what we call 'step-2' geometry on the opposite side where the S[k]x are constructed from 3 adjacent star points instead of 2 , in a manner similar to the DS[k[ for N-odd. Typically step-2 tiles can occur with S[1] tiles for N twice-even, for example N = 24.

Taking a step back it appears that the above self-similarity extends out to the S[3] and DS[3]tiles.

**Figure 18.2** – A vector plot of this region on the edge of N = 18

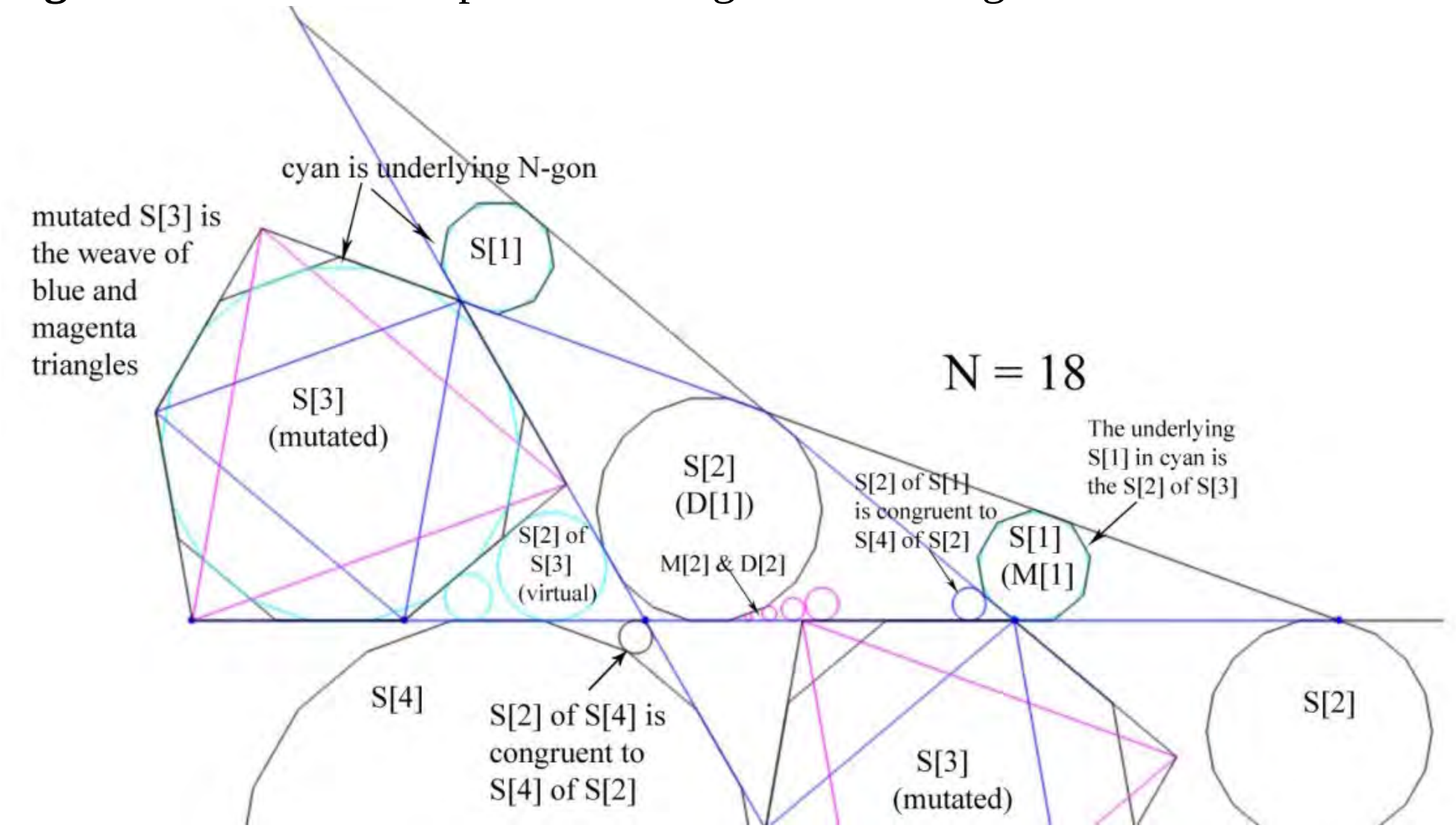

**•N = 19**

N = 19 is the second non-trivial member of the 8k+3 family so the Edge Conjecture implies the existence of a DS[7] along with S[1] at DS[N-4] = DS[15]. Algebraically N = 19 has complexity 9 along with N = 27 which is the next member of the 8k+3 family.

There may be an 8k+3 Conjecture which says that for 'most' family members there will be a conforming volunteer between pairs of DS[k]. The one known exception is N = 35. For N = 19 the first such volunteer between DS[7] and S[1] is known as $D_1$. It can serve as a surrogate 'M' tile for the DS[k] family of S[2] as it will be approximately half-way between S[2] and S[1]. For N = 11 earlier, $D_1$ is S[1], which is exactly the midpoint between S[2] and N. The $D_2$ tile shown here for N = 19 can be regarded as a conforming volunteer between DS[7] and S[2].

The existence of $D_1$ implies a rough form of self-similarity with respect to the first generation. As N increases and S[1] is further isolated from S[2] this pseudo symmetry will become more relevant. See N = 27, 35 and 43 below. This isolation is accelerated in the 8k+3 family since DS[7] is the first DS[k].

**Figure 19.1** The $\tau$ web of S[2]

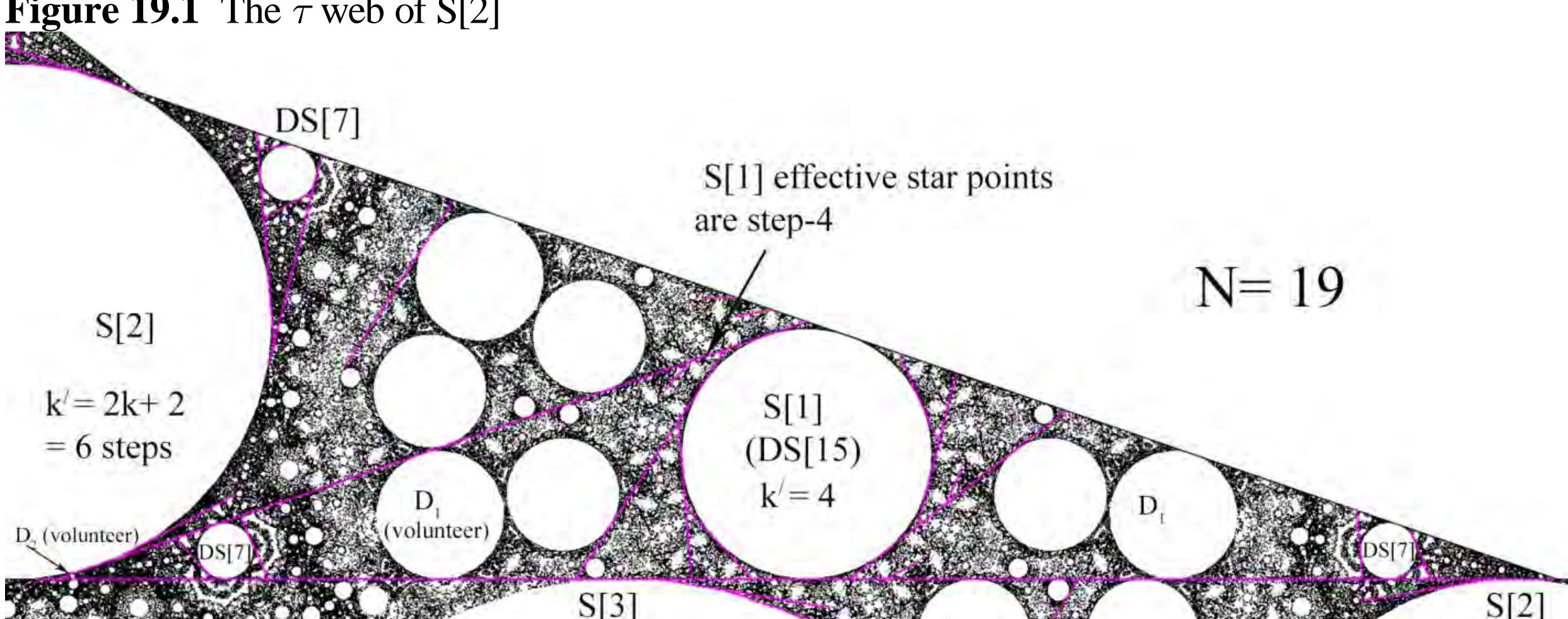

**Figure 19.2** The star[2] region of S[2]

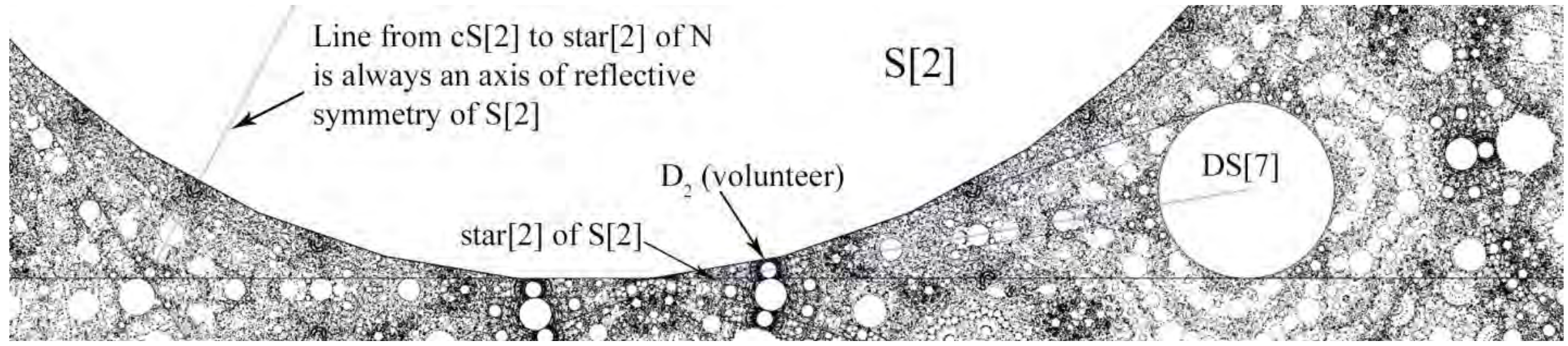

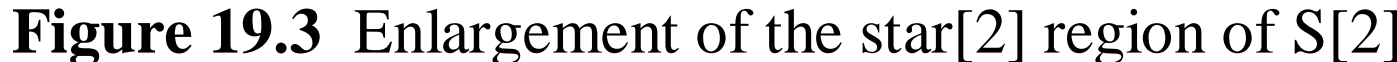

**Figure 19.3** Enlargement of the star[2] region of S[2]

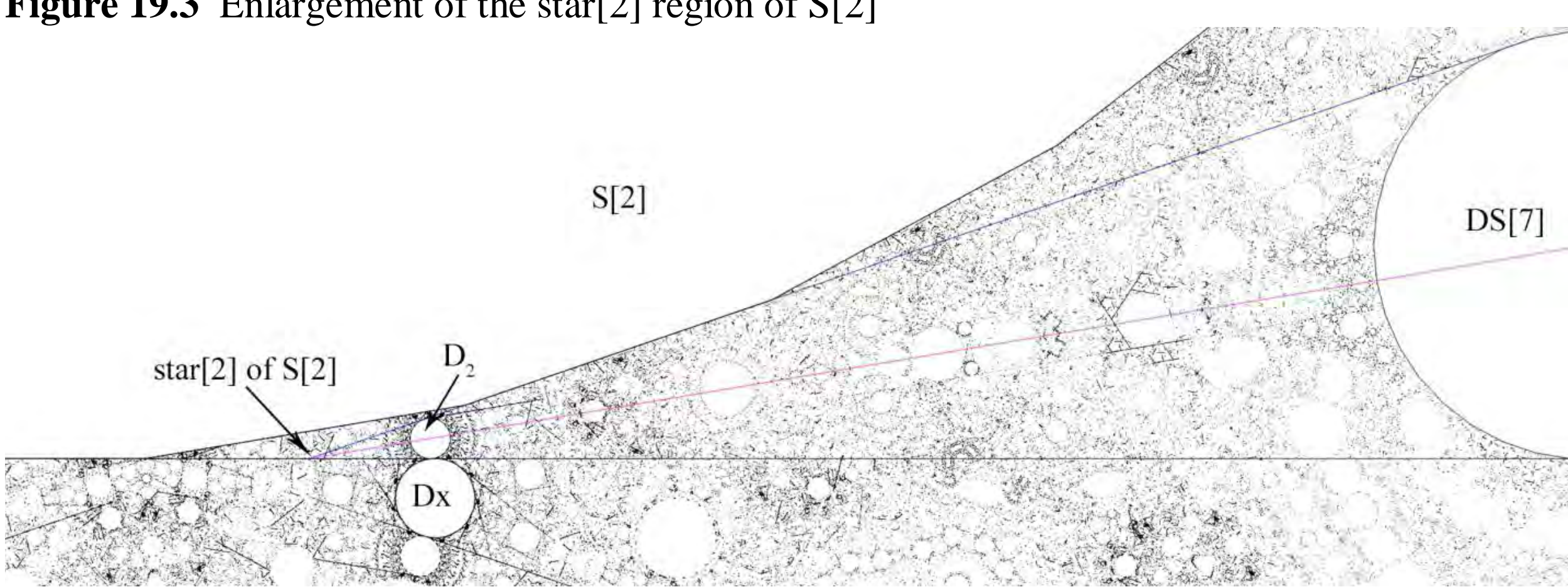

This $D_2$ tile has no obvious relationship with DS[7] or S[1]. For all odd N-gons, the virtual matching M tiles of the DS[k] may help to show hidden relationships – but that is not true here. The (known) families of these M tiles show no clear relationship with $D_2$. This $D_2$ is clearly a 2N-gon and appears to share a vertex with Dx below. This Dx has period 2546 and period doubling so the center (and height) can be found to arbitrary precision. This means that $D_2$ is also known to arbitrary precision, but we do not have the exact value of a second star point.

It is easy to find the parameters of 'star[2]' family tiles like DS[7] because they are strongly conforming to S[2]. For example DS[7] shares star[2] and star[7] with S[2] (while the virtual S[7] of S[2] shares star[1] and star[7]). Therefore if d is the distance between star[2] and star[7] of S[2] the Two Star Lemma says that:

hDS[7] = d/(Tan[17*Pi/38]-Tan[12*Pi/38]) ≈ .0187136 (inside DKHO where sN = 1).

(The two star points of DS[7] are star[38/2-2] and star[38/2-7] to match star[2] and star[7] of S[2]). Note that hDS[7] is also hS[7] /scale[2] of S[2]. (This calculation has to be done with care because it involves two different versions of scale[2]. First to find hS[7] of S[2], this is an N-gon so hS[7] = hDS[7]·GenScale/scale[2] where DS[7] is a (known) tile in the FF of N and scale[2] is relative to N. But now it has to be scaled up by dividing by scale[2] of 2N which by definition is $s_1/s_2$ = Tan[Pi/38]/Tan[2*Pi/38].)

The characteristic equation for hDS[7]/hN is valid for any N. In $\tau$-space by default hN = 1 so AlgebraicNumberPolynomial[ToNumberField[hDS7,GenScale[19]],x] is

$$\frac{55}{128}-\frac{2009x}{64}+\frac{3493x^2}{64}+\frac{4127x^3}{64}-\frac{89x^4}{8}-\frac{1539x^5}{64}-\frac{381x^6}{64}-\frac{3x^7}{64}+\frac{9x^8}{128}$$

where GenScale[19] = Tan[Pi/19]·Tan[Pi/38]. The step sequence for DS[7] is {2,2,1,1} which places it half-way between S[2] at {2} and S[1] at {1}.

**Figure 19.5** A ‘Deep-Field’ image of the region local to S[2]

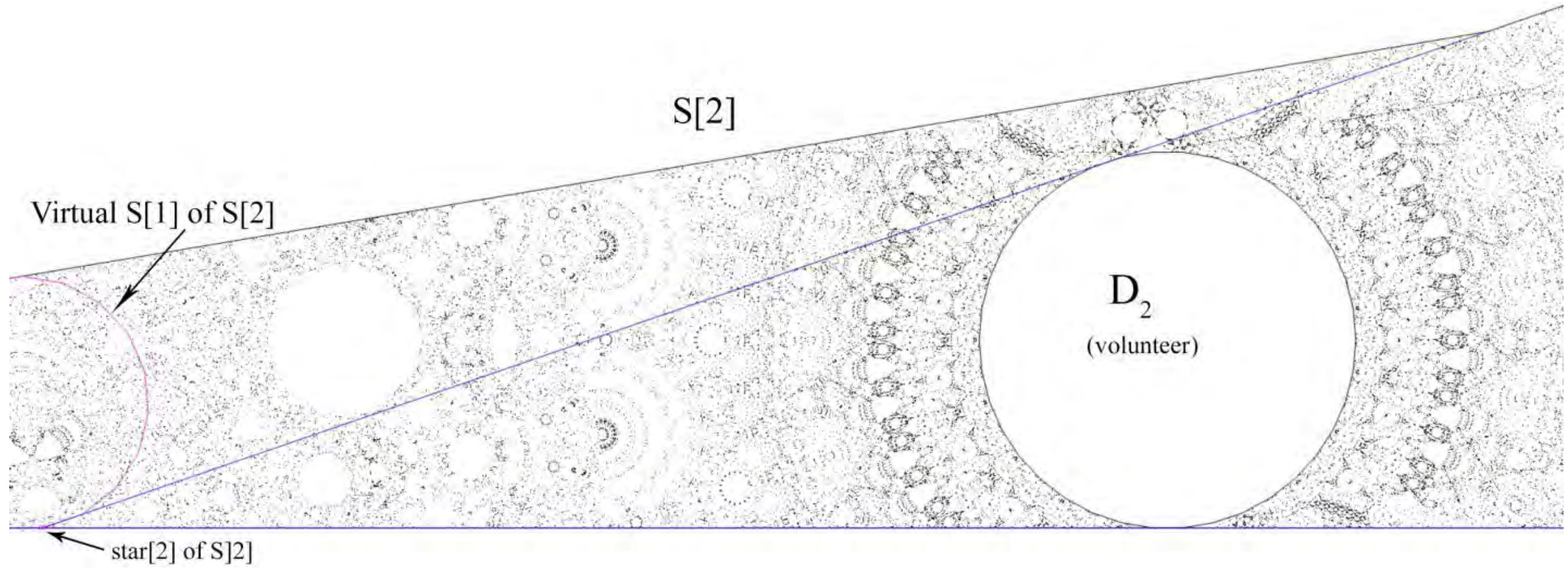

We have found no relationship between the weakly conforming volunteer $D_2$ and the magenta virtual S[1] of S[2]. It is clear that $D_2$ supports a whole ‘sub-culture’ of tiles and this geometry is repeated at the star[1] point of S[2] and the matching star point above $D_2$. For N odd, the geometry of the star[2] point is identical to the star[1] point of N and the star[1] point of S[2] is ‘out-of bounds’ when it comes to the star[2] family of S[2].

**Figure 19.6** The $\tau$-web local to S[1] has effective star points which are step-4. In the 8k+ 3 family this implies that star[1] of S[1] will be effective as show below in the early web. This star point countdown starts with star[N-2] of S[1] which is always star[2] of S[2].

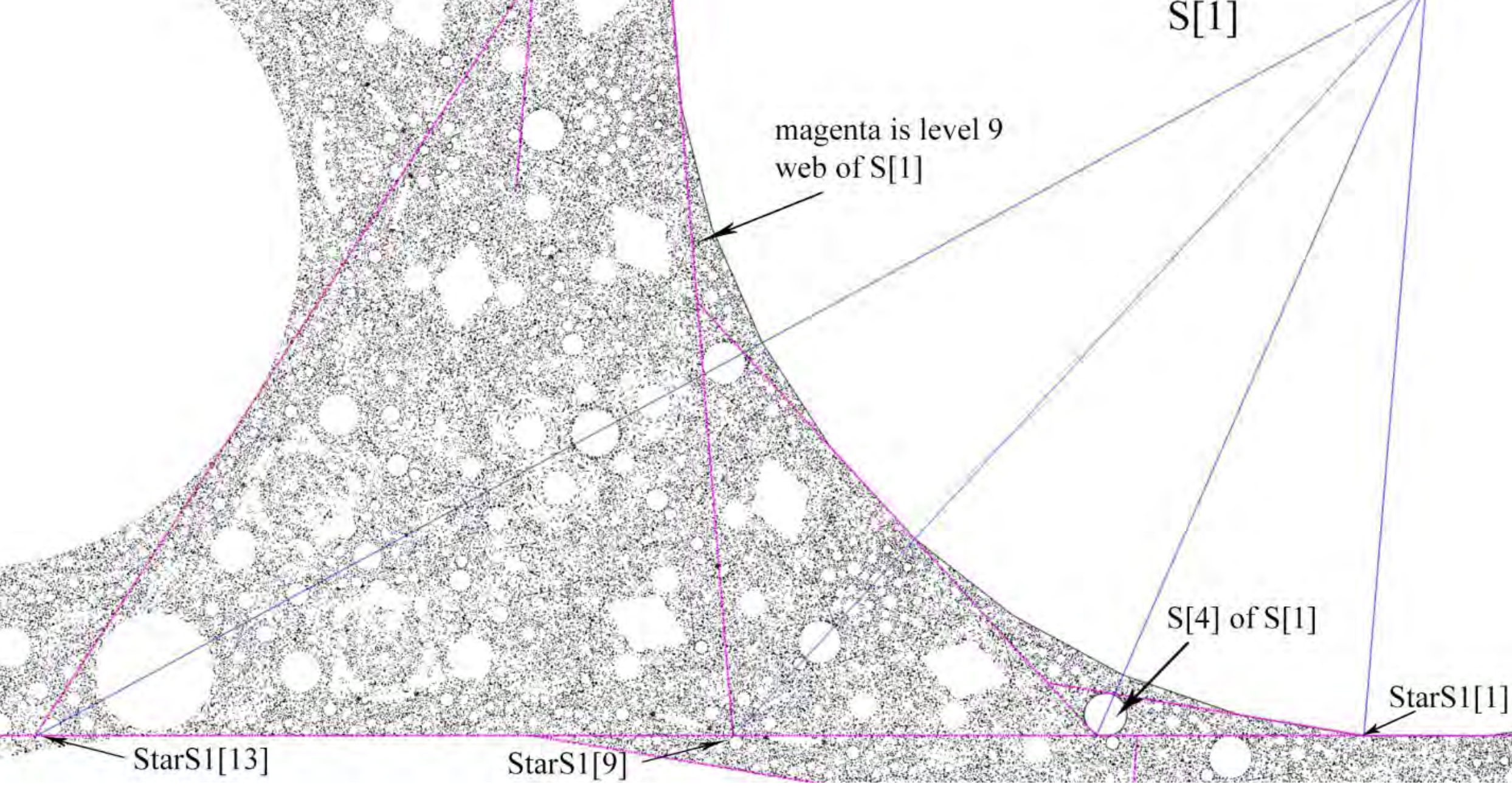

Note that S[1] supports an S[4] tile just like N = 11. Every effective star point has the potential to support families of tiles because these star points form early in the web. Any sequence that is conforming to the star point will have centers on the blue lines of symmetry.

**N = 20**

N = 20 has complexity 6 and is the second non-trivial member of the 8k+4 family. The Edge Conjecture predicts the survival of a DS[4] as well as S[1] at DS[8]. The Mutation Conjecture of [H5] says that S[6], S[5] and S[2] will be mutated as shown below. The 8k+4 family is the only one where S[2] is mutated and this combined with DS[4]s, will drive the local geometry.

**Figure 20.1.** Towers of S[k] tiles showing S[1] and S[2] sharing star points with S[3] and S[4]

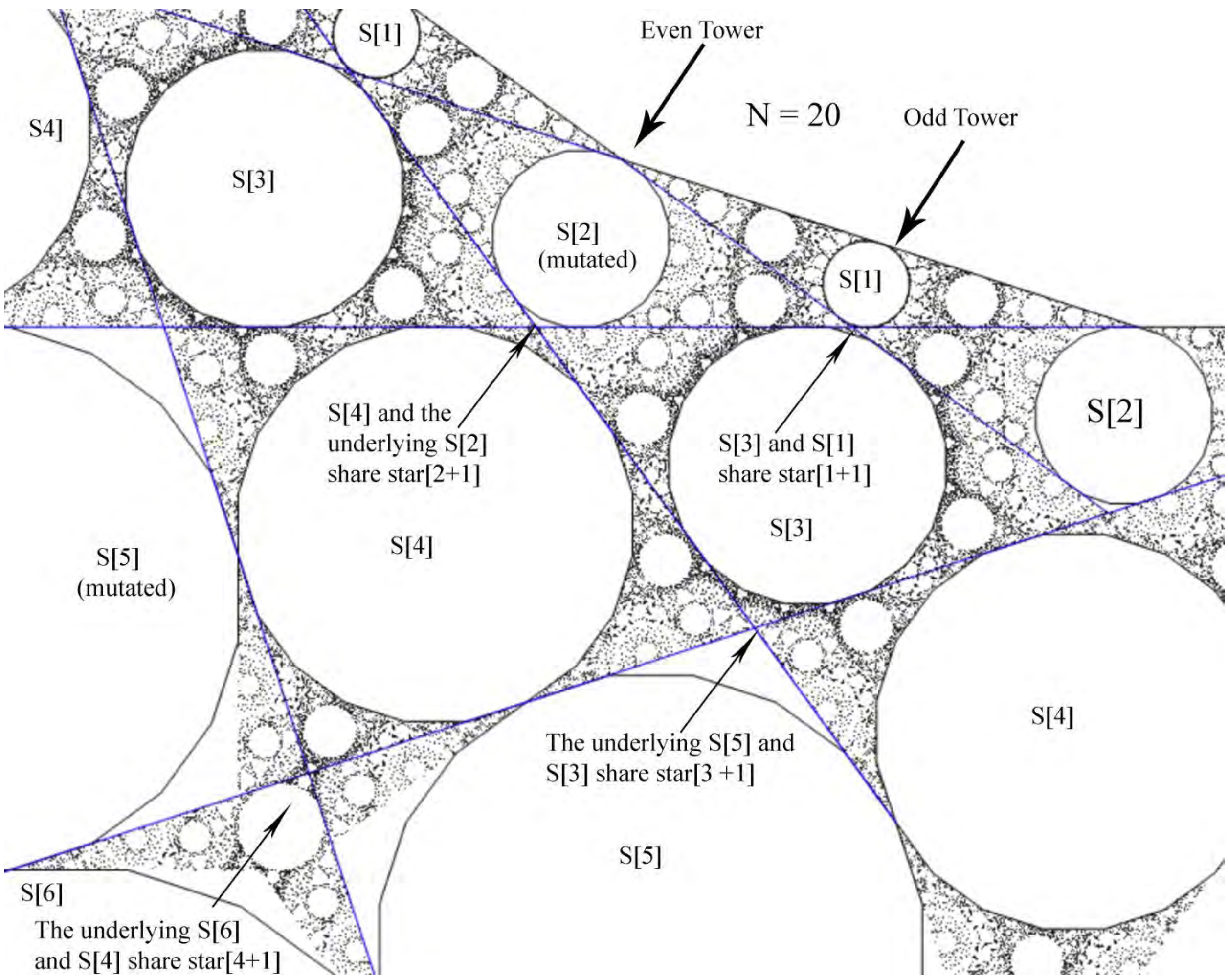

Since S[2] has k′ = N/2 – 2 = 8 and gcd(20,) = 4, S[2] will consist of two 'woven' N/4-gons, so here it will be an equilateral decagon which is the 'Riffle' or weave of two regular pentagons with slightly different radii. The 8k+4 Conjecture gives the expected implications of this mutation. This is a 'mod-16' conjecture and the two branches are anchored by N = 12 and N = 20. In all cases the predicted DS[4]s will be 'almost-vertex' tiles of the larger N/4-gons and there will be real or virtual Px tiles that will be actual vertex tiles of the smaller N/4-gons. These Px tiles will always be 'parent' tiles of the DS[4] but that will occur in two different ways.

For the N = 12 branch of the form 12 + 16j, DS[4] is the S[3+8j] of Px so for N = 12 it is simply S[3]. For the branch anchored by N = 20 the Px 'parents' will be D tiles relative to DS[4] but these Px may be virtual after N = 36. Every S[k] of S[2] has such a virtual counterpart by the FFT.

**Figure 20.2** S[2] is the weave of two N/4-gons for all members of the 8k+4 family.

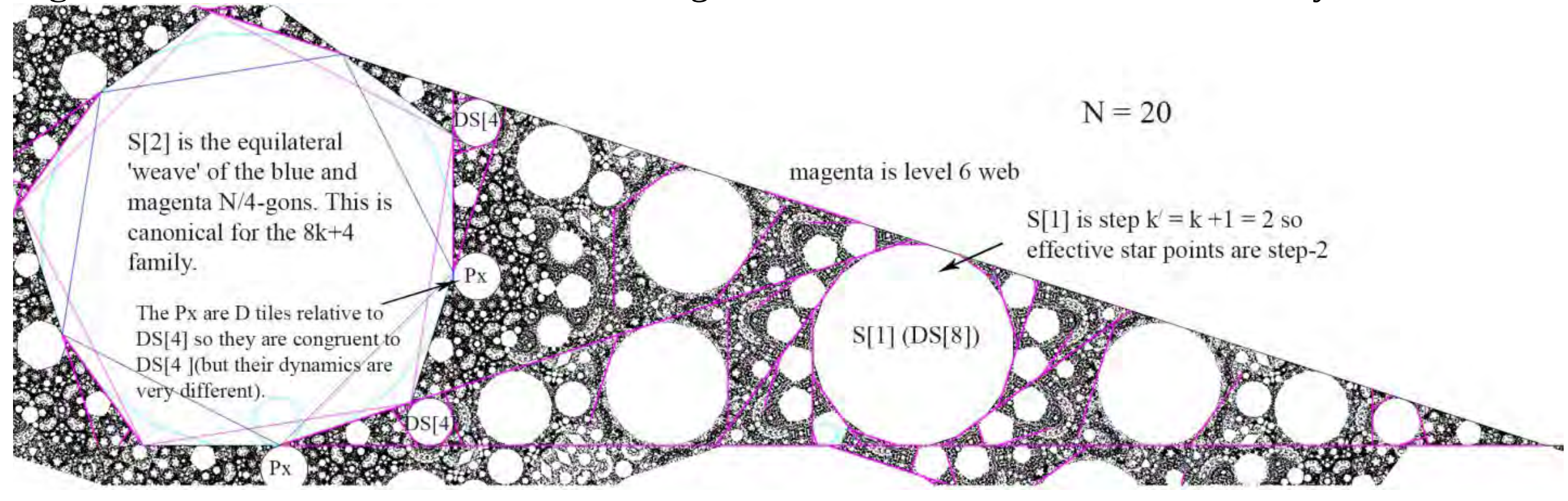

.

Here the Px have S[3] tiles which survive the web and they show promise of self-similar geometry. These S[3] tiles of Px will have height hS[3]·hDS[4] = hN·hS[3]·hS[2]·hS[4]= $\text{GenScale}^3$ / (scale[3]·scale[2]·scale[4]) ≈.000477, so this S[3] is a (modified) $3^{rd}$ generation tile that appears to foster a $4^{th}$ generation at star[1] of S[2]. This may be unique to N = 20 and in general we have no evidence that the (real or virtual) Px at star[1] of S[2] foster a next generation.

gure **20.3** The geometry local to the DS[4] tiles.

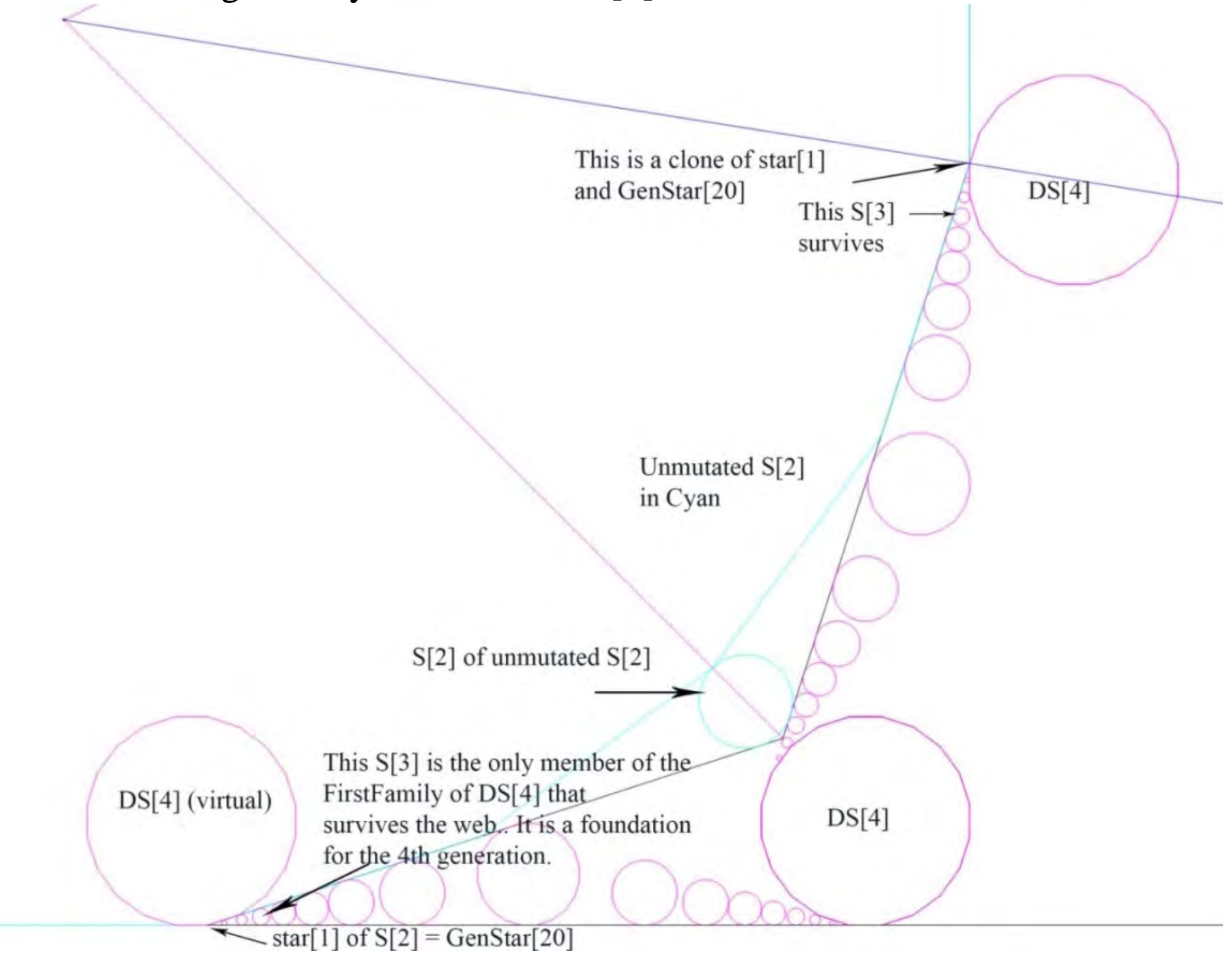

**Figure 20.4** The 3rd and 4th generations at star[1] of S[2]. The large conforming 'volunteer' appears to be mutated in the same fashion as S[2].

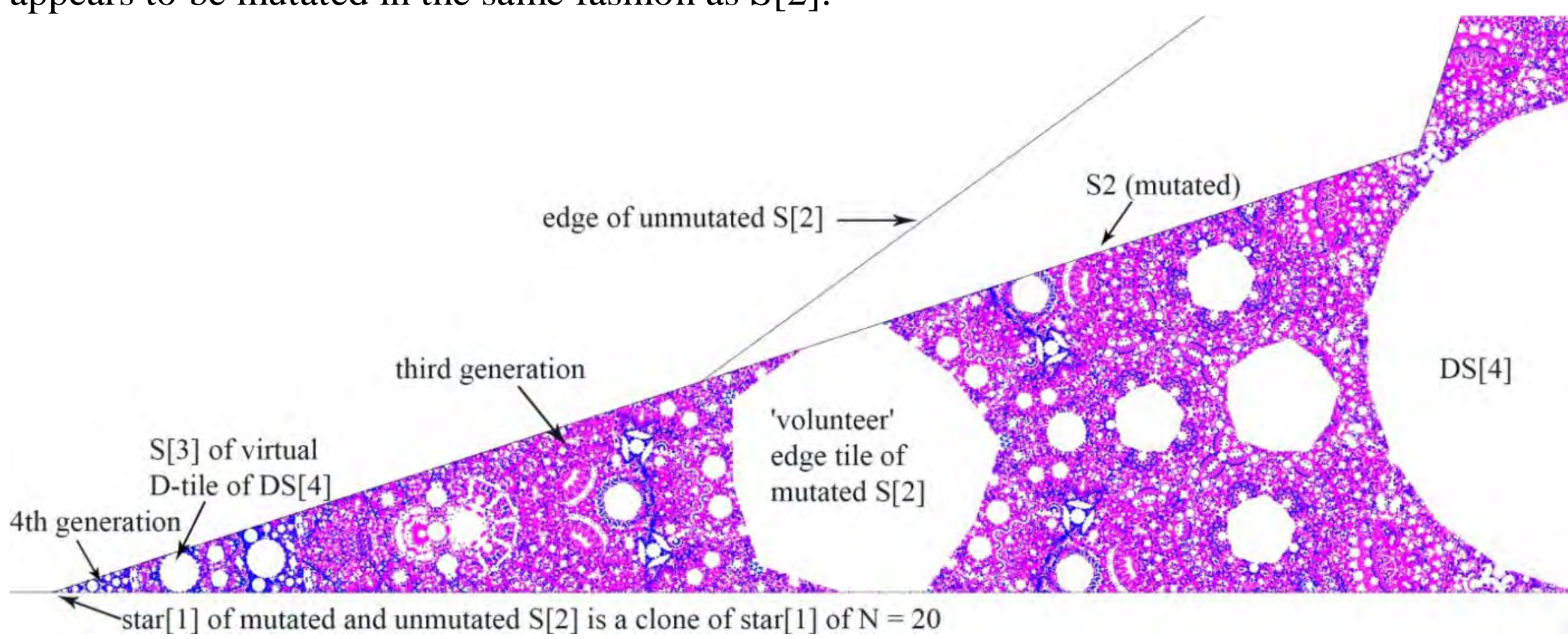

When N is even, S[1] will be in the First Family of all the S[k] but as N increases the major influence on the local geometry of S[1] will be the rotated S[3] as shown in the 'tower' plot above. S[1] will share its star[2] point with S[3] so the S[2] of S[3] will be adjacent to S[1]. This tile will have a scale[3] height relationship with the S[2] of S[1] and it will be the S2x tile in the 'step-2' family of S[1].

**Figure 20.5** The geometry local to N when N is twice even

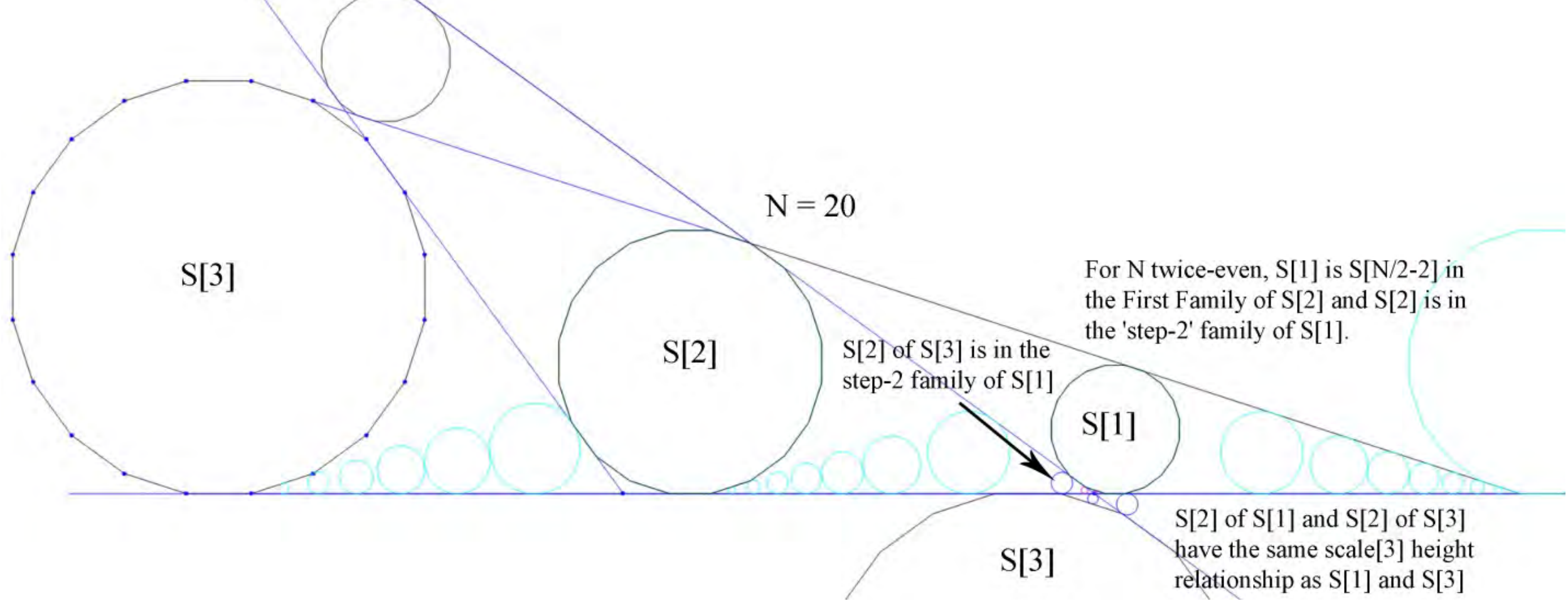

The tiny magenta S[2] of S[1] is also in this step-2 family of S[1] because it has the canonical scale[2] relationship with S[1] and it will have a scale[3] relationship with the larger S[2] of S[3] which is similar to the scale[2] relationship between S[1] and S[2].The only difference with N twice-odd is that S[2] is now 'mutated' and an N/2 gon. This actually fits better with S[3] and they share star[1] points which means they share a common vertex. See N = 18 for a vector plot like the one above. As expected the N odd case is double the twice- even case with star[4] shared between S[1] and the rotated S[3]. See N = 11.

Because S[2] of N is mutated it should be no surprise that this S[2] of S[3] is also mutated as shown below. The Twice-even S[1] Conjecture states that S[1] has the potential to support a Step-2 (or ‘star[2]’) family. Since the S[1] web is step-2, the odd star points are ‘effective’ and each one can support an Skx tile which is a D tile relative to S[k] for k odd.

**Figure 20.5** The mutated Sx tile shown here is in the Step-2 Family of S[1]

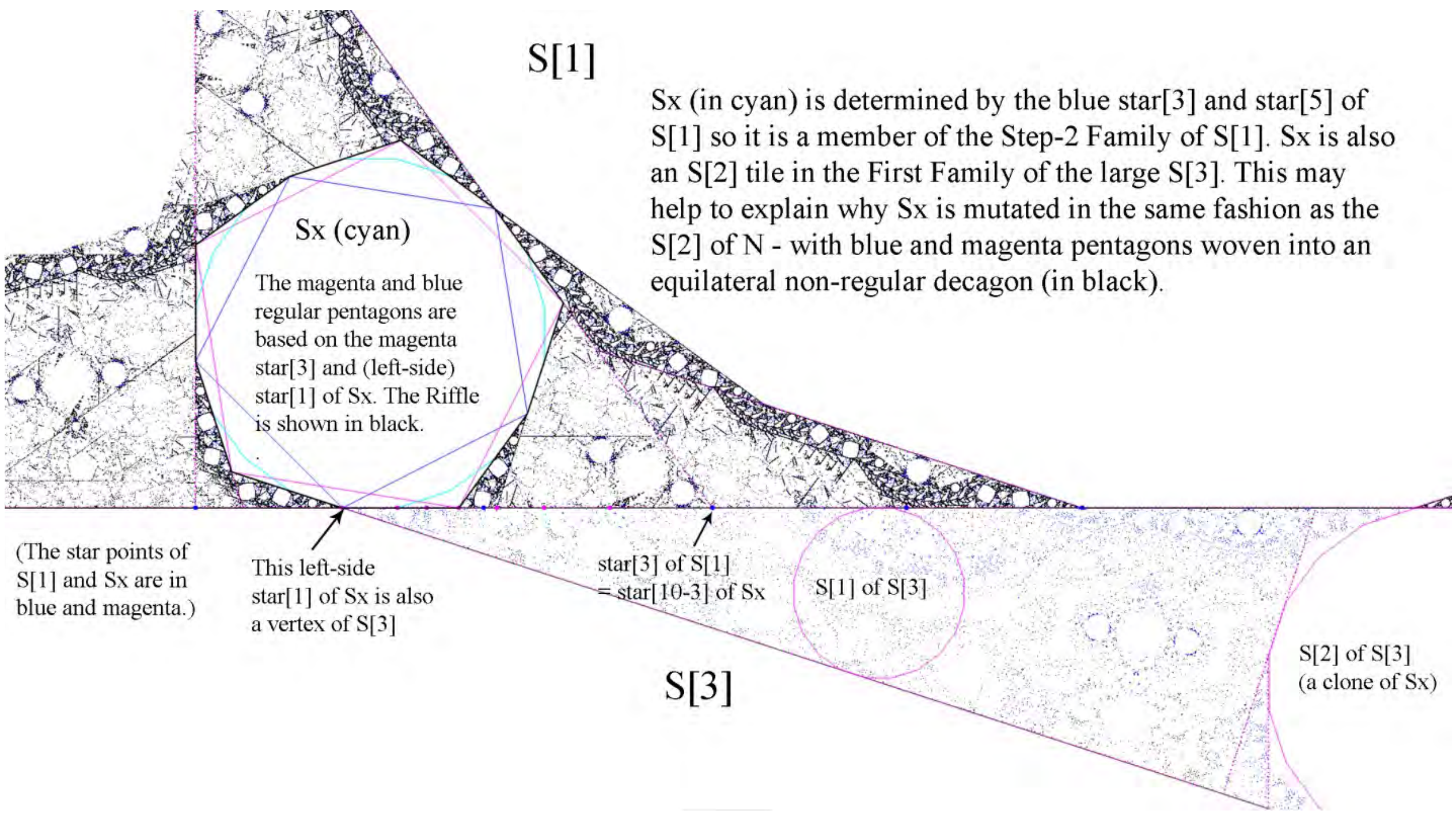

Since the ideal Sx is known it is easy to obtain the parameters of the mutated version. Here the mutated Sx is the Riffle or weave of two pentagons. This makes the mutation similar to S[2] with k′ = N/2-k = 8 and S[6] with k′ = 4 . Both have N/gcd(k′ , 20) = 5.

**• N = 21**

N = 21 has algebraic order 6 along with N = 13, 26, 28, 36 and 42. Both N = 21 and N = 13 are 8k+5 so there is a well-defined 8k+2 generation structure on the edges of D acting as 2N. We would like to know what effect this has on N. There is no explicit 8k+5 Conjecture but every member of this family has a DS[1] and such tiles always have potential for extended tile structure. Since these are star[2] families, DS[1] is the 'D' tile relative to an S[1] which is usually virtual. But for N = 5 and N = 21 this S[1] actually exists. Therefore N = 21 has the foundation for a self-similar 3rd generation presided over by M[2] and D[2]. For N = 5 these D[k], M[k] sequences survives at both D and N, but typically in the 8k+5 family these is no M[2] to match the existing D[2]. However it seems that D[2] always has surviving First Family members with potential for further generation structure.

**Figure 21.1** The early web of S[1] and S[2]

This black web was generated by an 'egalitarian' scan of the interval from star[1] of S[2] to star[1] of N at level 800, so the dark semi-invariant region around S[2] is an artifact of N = 21 and not an artifact of the scan depth. For all N, the S[2] region bounded by S[2] and the penultimate DS[k] shows signs of invariance, so for N-odd, DS[N-12] may serve a 'shepherd' satellite of S[2] as shown here.

The 'effective' star points of S[1] and S[2] are determined by the magenta early web. For N odd these retrograde webs are always 2k + 2 steps, so the effective star points of S[2] are step-6 while the effective star points of S[1] are step-4 starting with star[19] (N-2] at star[2] of S[2]. By convention the $\tau$-web uses an initial interval that spans S[1] and S[2] so the 'effective' star points of S[2] are determined by the combined web which is step-8 to match the Rule of 8 – and the step-8 rotational symmetry of S[2].

The S[1] effective star points are not modified by S[2], so they are simply step-4. Since the right-side web of S[1] is clock-wise, these step-4 star points can be tracked on the horizontal axis between S[1] and S[2]. These effective star points of S[1] are 19,15,11,7 and 3 but just the first four of these are shown here. It is no surprise that DS[9] is mutated and in general there will be mutations of the DS[k] when gcd(k,N) >1.

As indicated above, N = 21 is apparently the lone 8k+5 non-trivial case where the existence of DS[1] implies a matching M[2]. This is a unique because the ‘star[2] family’ of S[2] is usually not backwards compatible with the ordinary First Family of S[2]. In all cases these two families are related by the fact that S[k] is congruent to the ‘M’ tile of DS[k] and in the case of DS[1] here, this congruence is an equivalence.

The First Family S[k] of S[2] are almost always virtual, so the DS[9] here would be a regular 2N-gon with an N-gon M tile congruent to a virtual S[9] of S[2]. It is easy to construct these DS[k] with the Two Star Lemma because they always share star[2] and star[k] with S[2]. This means that DS[1] must share star[2] and star[1] of S[2] as shown in the enlargement below. Therefore the ‘M’ tile of DS[1] will always be the M[2] tile of S[2] by the First Family Theorem, and they will form a D-M pair. It appears that N = 5 and N = 21 are the only cases where this M[2] tile actually exists

**Figure 21.2** The web local to DS[1] showing the (ideal) First Family of DS[1] in blue

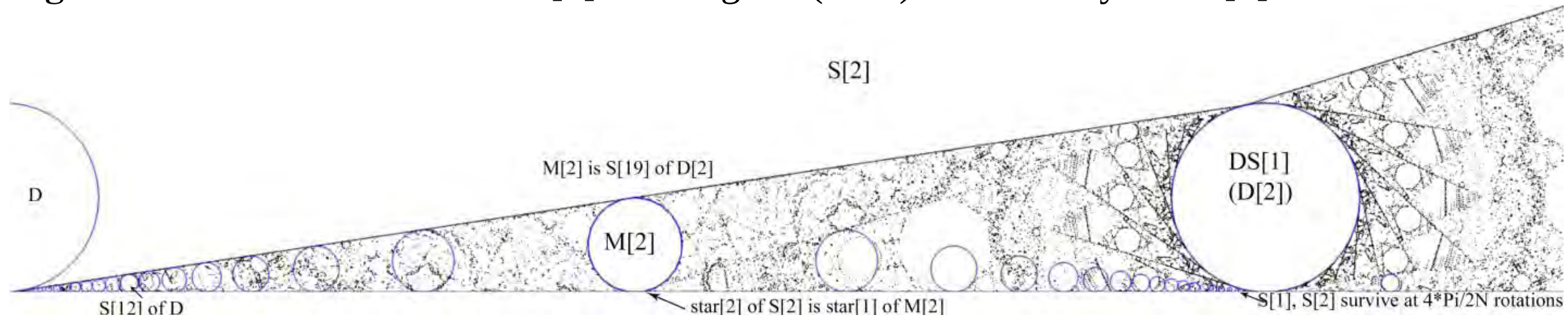

In all known cases of the 8k+5 family, DS[1] has a step-2 web as shown here. This usually implies a slight asymmetry in the web geometry. The left-side web appears to have survivors at S[1] and S[2[ as well as M[2] at S[19], but the right-side web may have more survivors. An M[2] seems to exist on both sides of DS[1]. Note that D at star[1] also has and S[12] survivor.

**Figure 21.3** The right-side web local to DS[1]

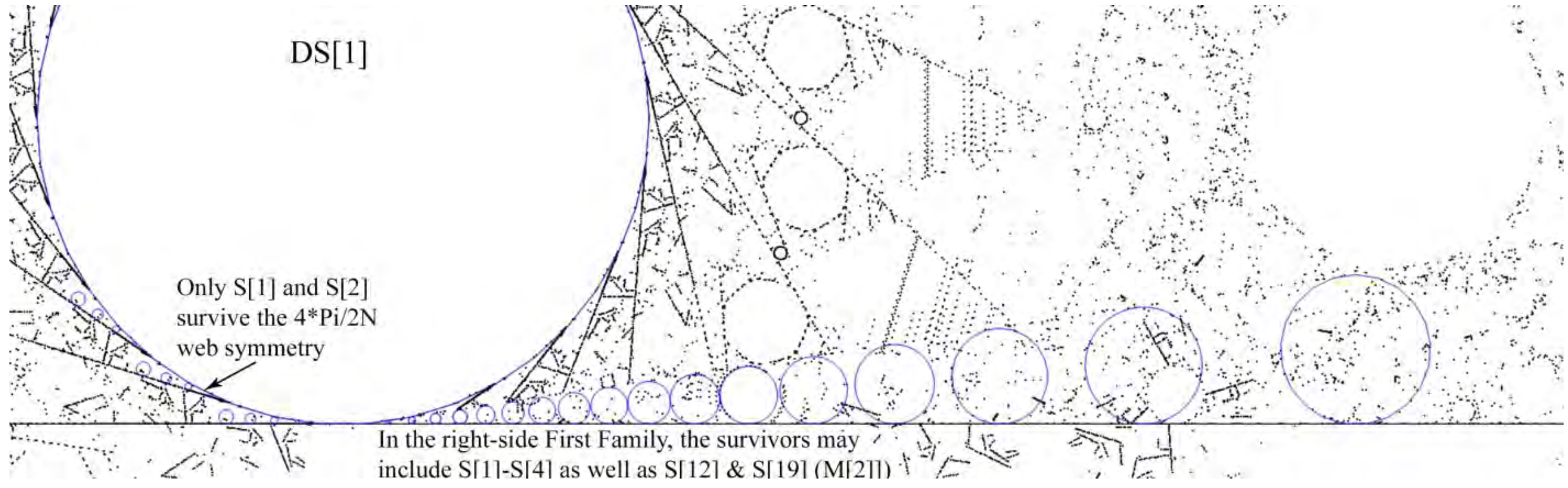

**• N = 22**

N = 22 and the matching N = 11 are the only quintic regular polygons so the web can offer some insight into the geometry of a degree 5 number field. As shown in Figure 2.1, the edge geometry is shared by the S[1]-S[5] tiles of N = 22. As in all twice-odd cases, the edge geometry of N = 22 is what we call the $2^{nd}$ generation of both N = 11 and N = 22.

N = 22 is in the 8k+6 family so there is a DS[1] serving as an M[2], but unlike N =14, there is no volunteer DS[3] whose edges can help to form a matching D[2]. It appears that the weakly conforming Mx tile shown below is the remnants of a DS[3] that never formed. Therefore there is no obvious path toward self-similarity, but on smaller scales there are many colonies of tiles which may be invariant. Some of these colonies are local to Mx and the small Sx which is a satellite of both Px and DS[5]. In [H5] we show how the Two-Star Lemma can be used to construct Sx because it shares a star point with both larger tiles.

**Figure 22.1** The web showing the shared geometry across the horizontal extended edge of N.

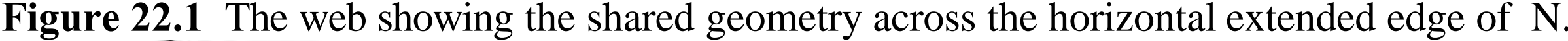

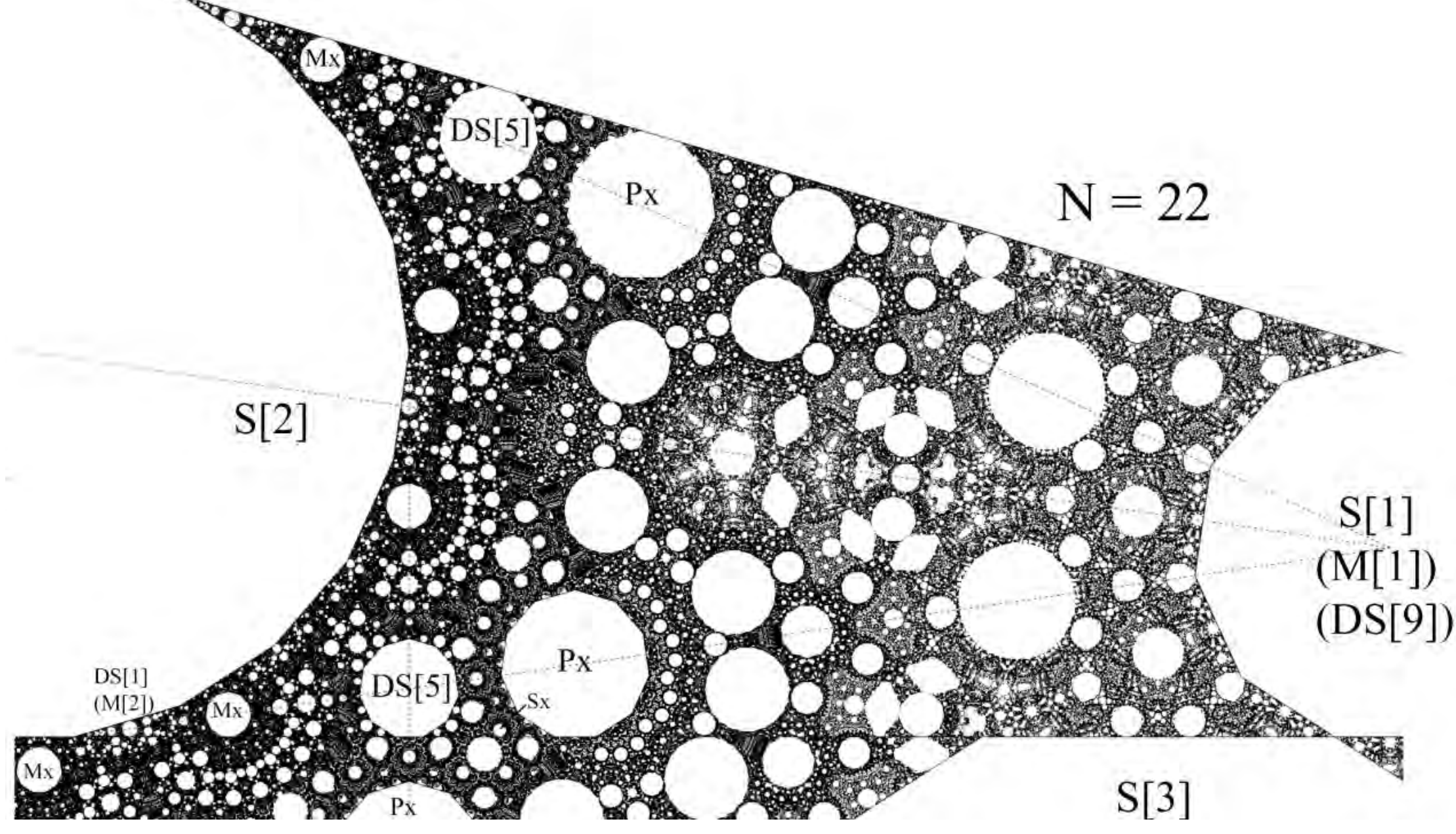

**Figure 22.1** For N even, the First Family of S[2] are what we call the DS[k]. By the Rule of 4 the survivors include DS[1], DS[5] and S[1] at DS[9]. There are no other survivors.

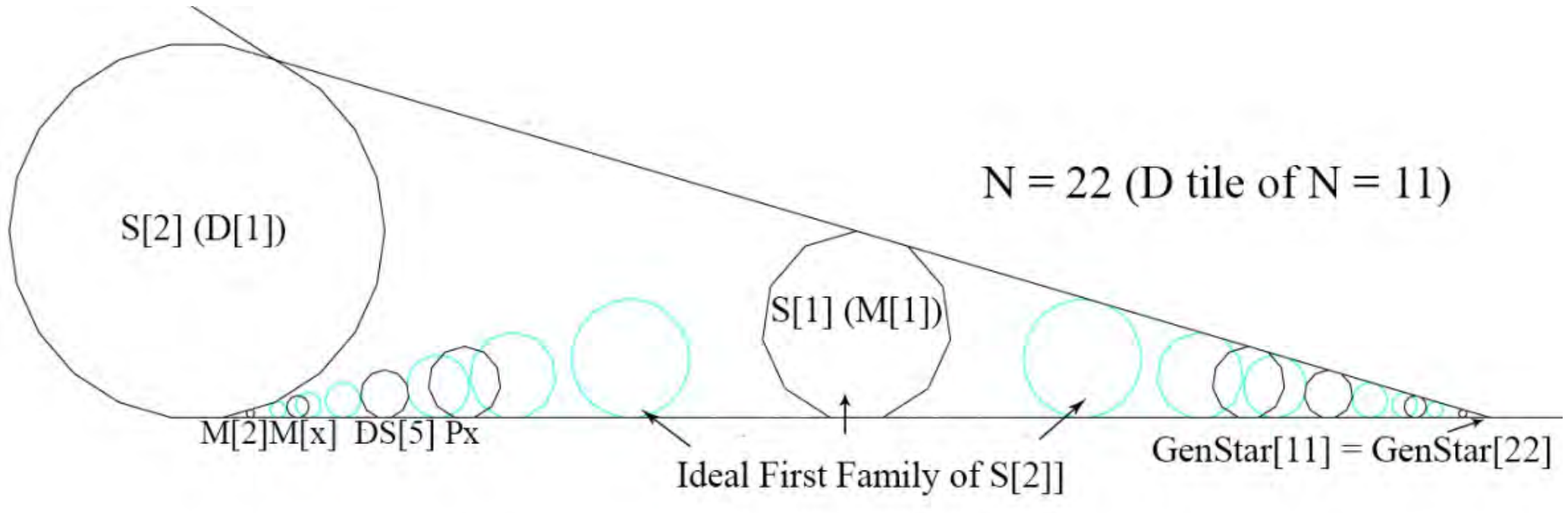

The Mx and Px tiles are only weakly conforming to S[2] because they share star[1] of S[2] but no other star point. By the Two-Star Lemma their parameters can be found if another star point can be identified, and in the appendix of [H5] we show how it is possible to use the web evolution of $\tau$ or the Digital Filter map or the Dual Center map, to find a second star point. The case of Px was easier and the corresponding polynomials are shown in Table 22.1 below.

The Px tile can be regarded as a volunteer 'M' tile shared between S[1] and S[2] in a manner similar to the PM tiles of N = 14 and N = 18. Studies of combined S[1]-S[2] webs show clearly how S[1] and S[2] share in the early evolution of the DS[k] and the volunteers, but S[1] itself evolves in a largely independent fashion.

The graphic below compares the influence of S[1] and S[2] using blue for the right-side web of S[1] and magenta for the right-side web of S[2]. It is clear that Px is like a 'swing-state' between these two influences. This will have a lasting effect of the limiting geometry as can be seen in the plot above.

**Figure 22.2** The level-4000 webs of S[1] (blue) and S[2] (magenta)

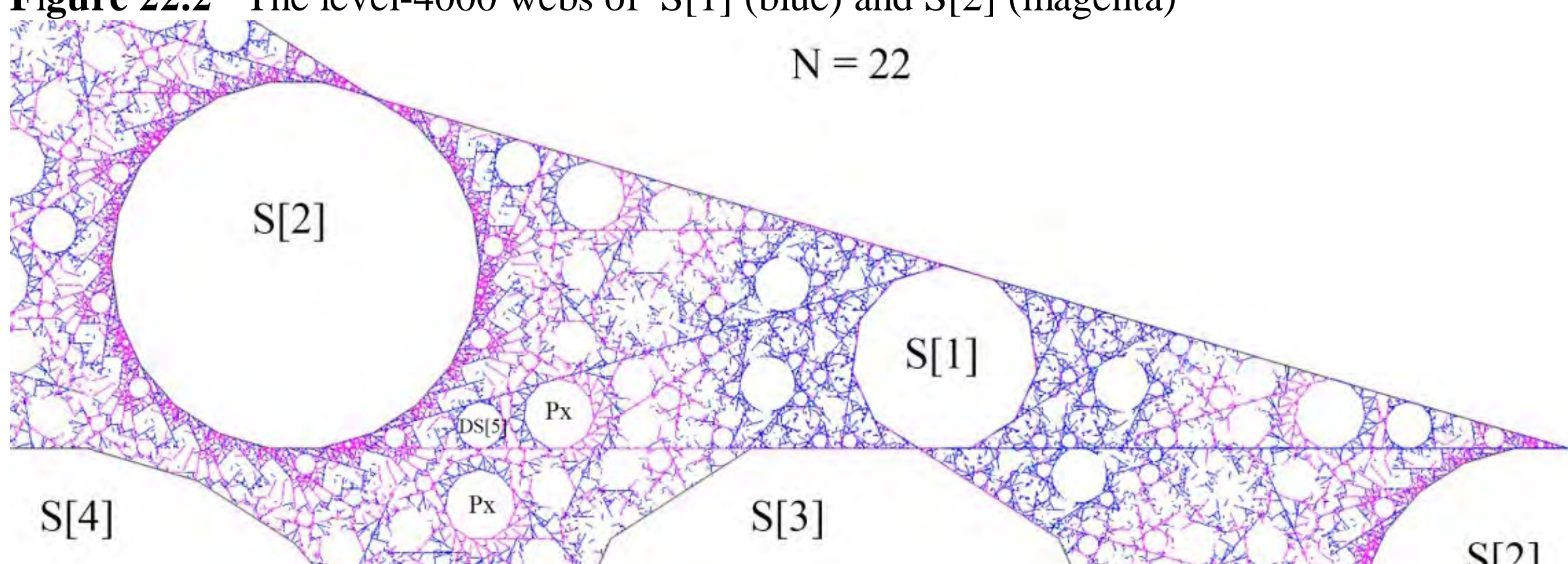

In the same fashion it is possible to trace the early influence of S[3] by generating its local (right-side) web and comparing it with the joint S[1]-S[2] web. As expected the influence is much greater for S[1] where it shares a vertex.

**Figure 22.3** The level-300 magenta web of S[3] showing its influence on S[1] and S[2]

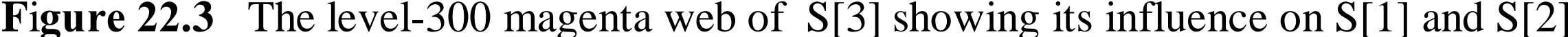

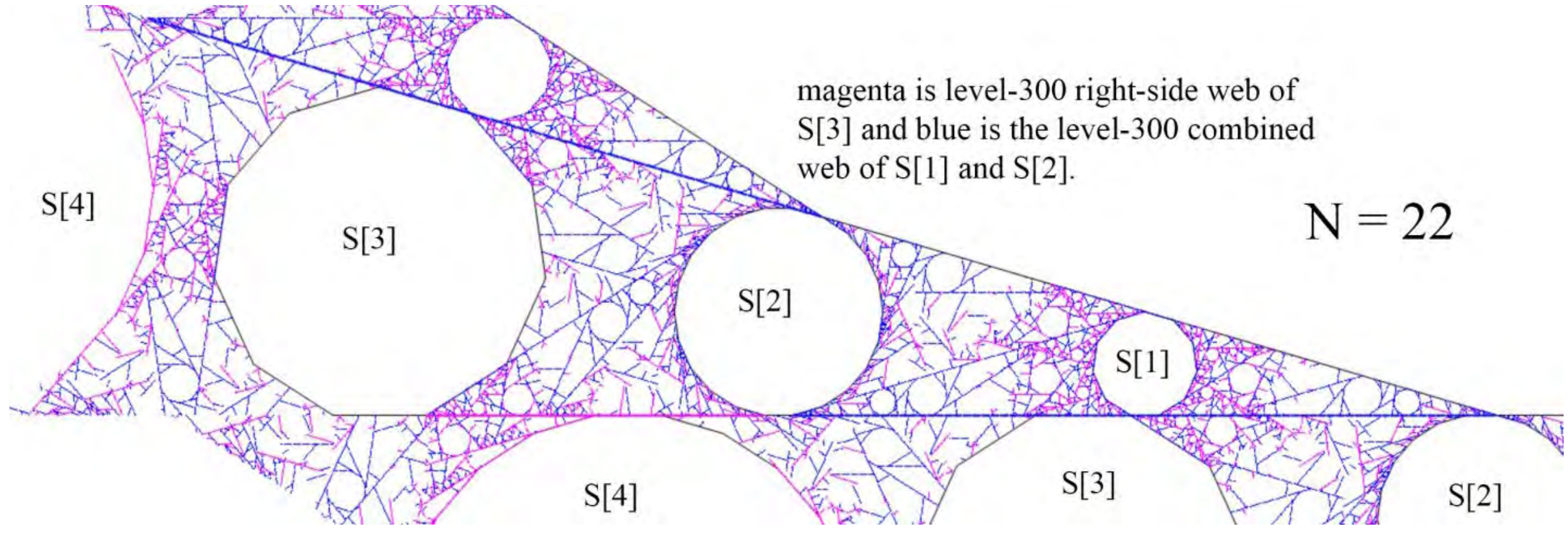

Returning to the Px region, it is clear that the Sx tile is more closely aligned with Px than DS[5], but they both have copies of Sx as satellites.

**Figure 22.4** The region local to Sx showing colonies of tiles.

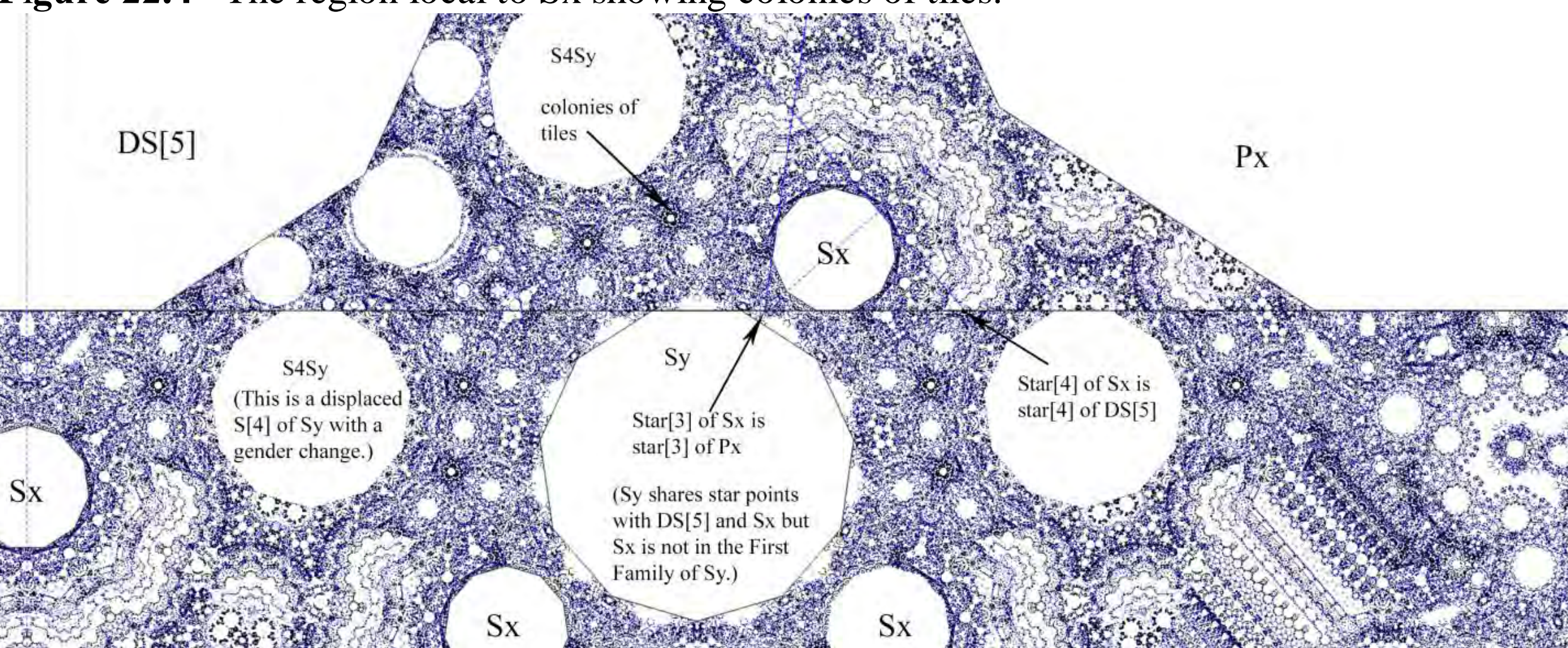

At this scale it is hard to see, but the (right-side) star[4] of Sx supports a small tile which is weakly conforming to this mutual star point. It is highly likely that there is a convergent sequence of such tiles. Throughout this region there are clusters of tiles which appear to be invariant. These colonies have the potential to foster endless chains of future generations with their own unique geometry The edges of Px support similar colonies. It is likely that almost all N-gons will have tiles on all scales and it is equally likely that these tiles will exist in diverse environments. Here this future geometry would be expected to retain some quartic influence.

The geometry of the star[1] region of S[2] is just a reflection about S[1] of the geometry at star[1] of N – which we call GenStar[22] or GenStar[11]. Every star point has a corresponding scale and here it is called GenScale[11]. This is our traditional generator of the scaling field $S_{11} = S_{22}$. When N is 8k+2, there should be a matching convergent sequence of tiles, but here the sequence is almost entirely virtual. However the local geometry of this GenStar point is still replicated everywhere, so at all scales there should be copies of the origin. This means there should be colonies like those seen above throughout the geometry. Below we will track some of these invariant colonies in the vicinity of Mx.

**Figure 22.5** The star[1] region of S[2]

One way to probe a star point is to generate the parameters of ideal tiles that would exist and then iterate the centers under $\tau$ to look for temporal scaling. It is not surprising that there are no obvious coherent sequences here, but it is easy to find colonies similar to those in the vicinity of Sx and the $S_k$ tiles of N = 11, See Figure 11.7. Here we will examine one of these 'island' colonies in the vicinity of Mx.

**Figure 22.6** The local geometry of Mx

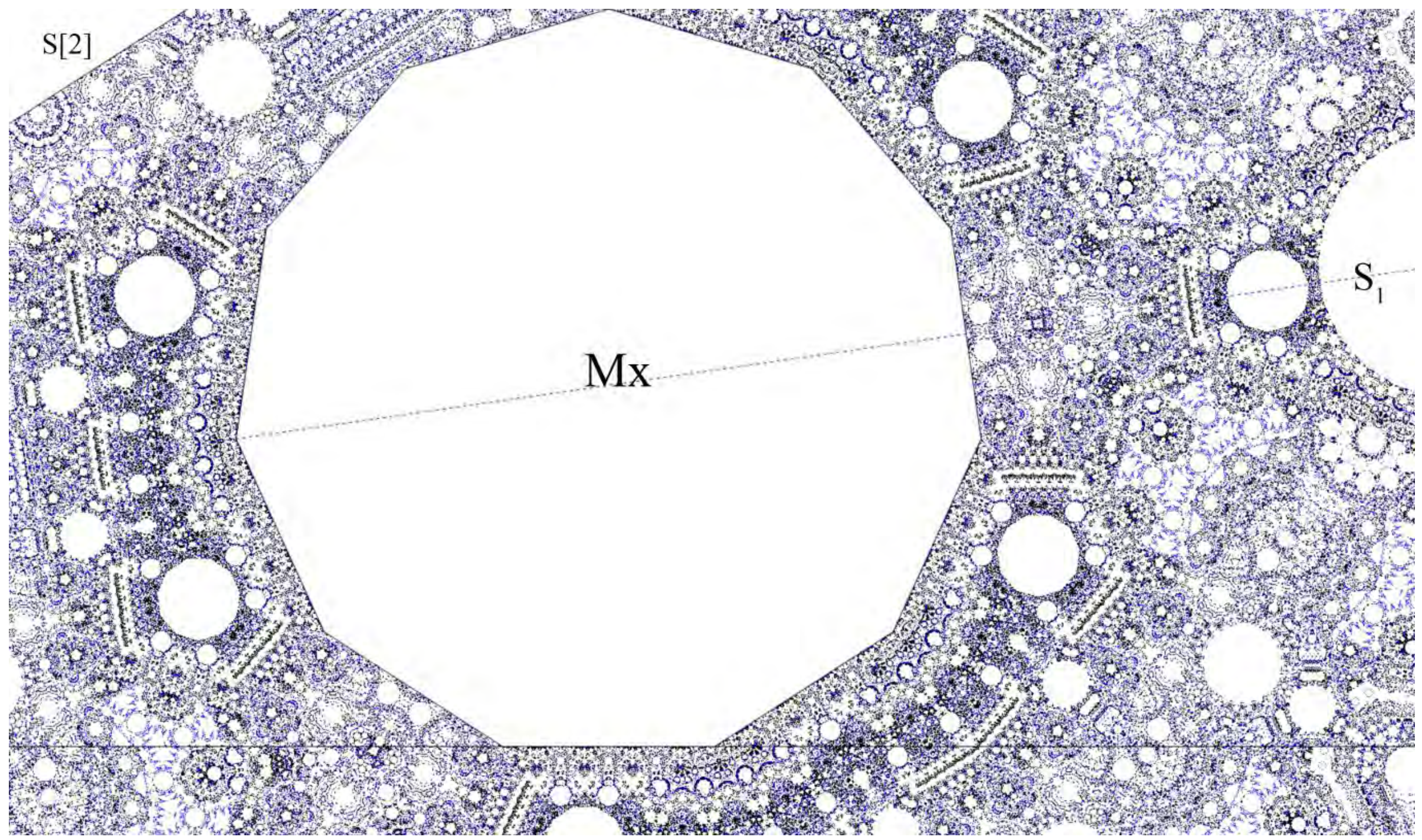

Over a period of years we have studied these invariant islands local to Mx. One of them is called N40 and it is shown at the right below. In the left panel N40 is shown with its reflection near the $S_2$ tile which is just off the screen above. See the overview in Figure 22.5. It was no surprise to discover that the geometry of N40 is congruent to a region close to star[1] of S[2]. This region is between the virtual M[3] and M[4] tiles so the scale is about GenScale[11]$^3 \approx .000075$. The N40 pair shown here 'point' to a virtual clone of star[1] at the center of the small edge tile of $S_2$. There is a similar (but not identical) pair that point to the center of the satellite of $S_1$.

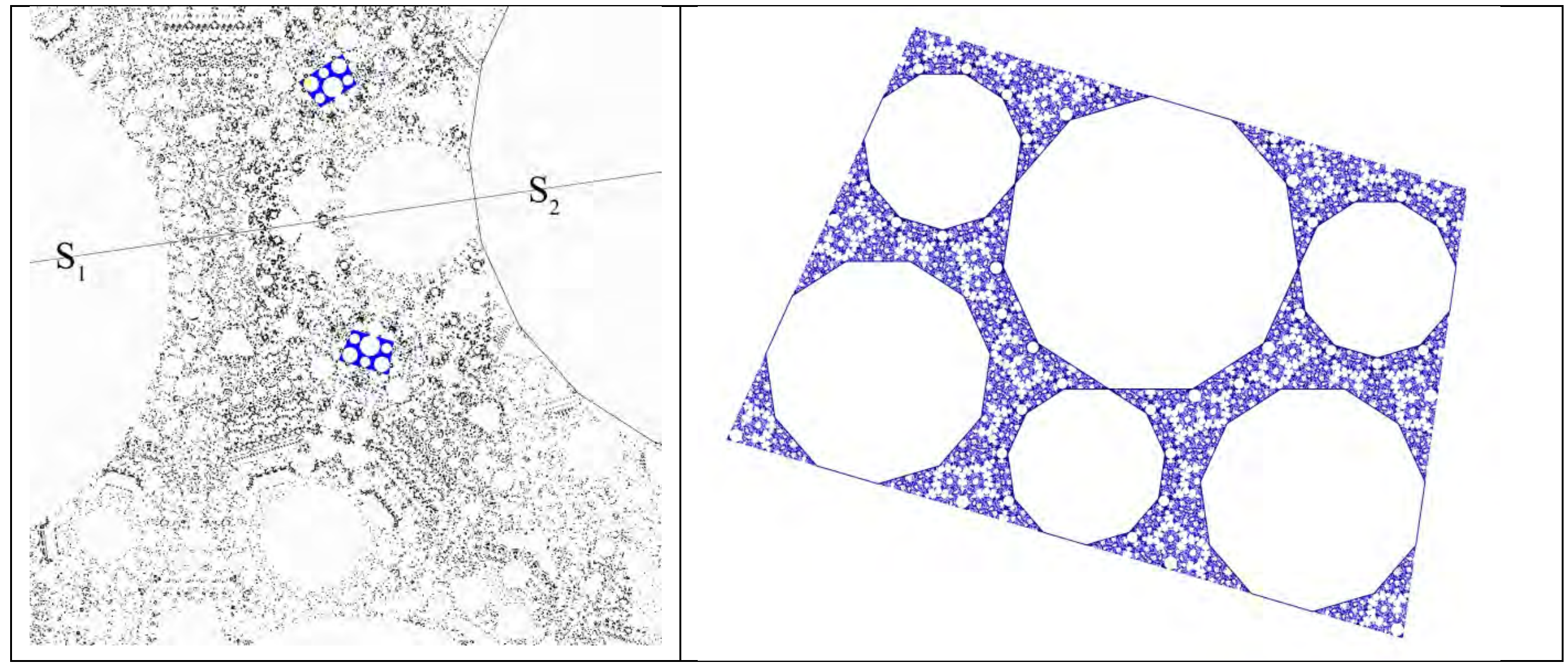

These islands have their own scaling and geometry which appears to be more uniform than the surrounding geometry. This may be because they all have a similar origin in the vicinity of the $4^{th}$ or $5^{th}$ (virtual) generations at GenStar. In this sense they could provide some insight into 'the geometry of 'future' generations. One clear indication of its origin is that it is symmetric with respect to its center line – which is just a continuation of the blue line of symmetry of running from the center of S[1] to star[1]..

. We have seen with N = 14 and N = 18 that in these twice-odd cases, the S[1] tile may have 'hidden' structure based on its dual role as N/2-gon and N-gon. This is also true for N = 22 as shown below. The First Family of S[1] has no local web survivors, but the First Family of the underlying N-gon has S[1] and S]3] as survivors of the web. Note that the virtual S[1] in this family sits on the line of symmetry joining the center of S[3], the center of S[1] and the origin.

**Figure 22.9** The step-2 web local to S[1] showing the dual nature of S[1]

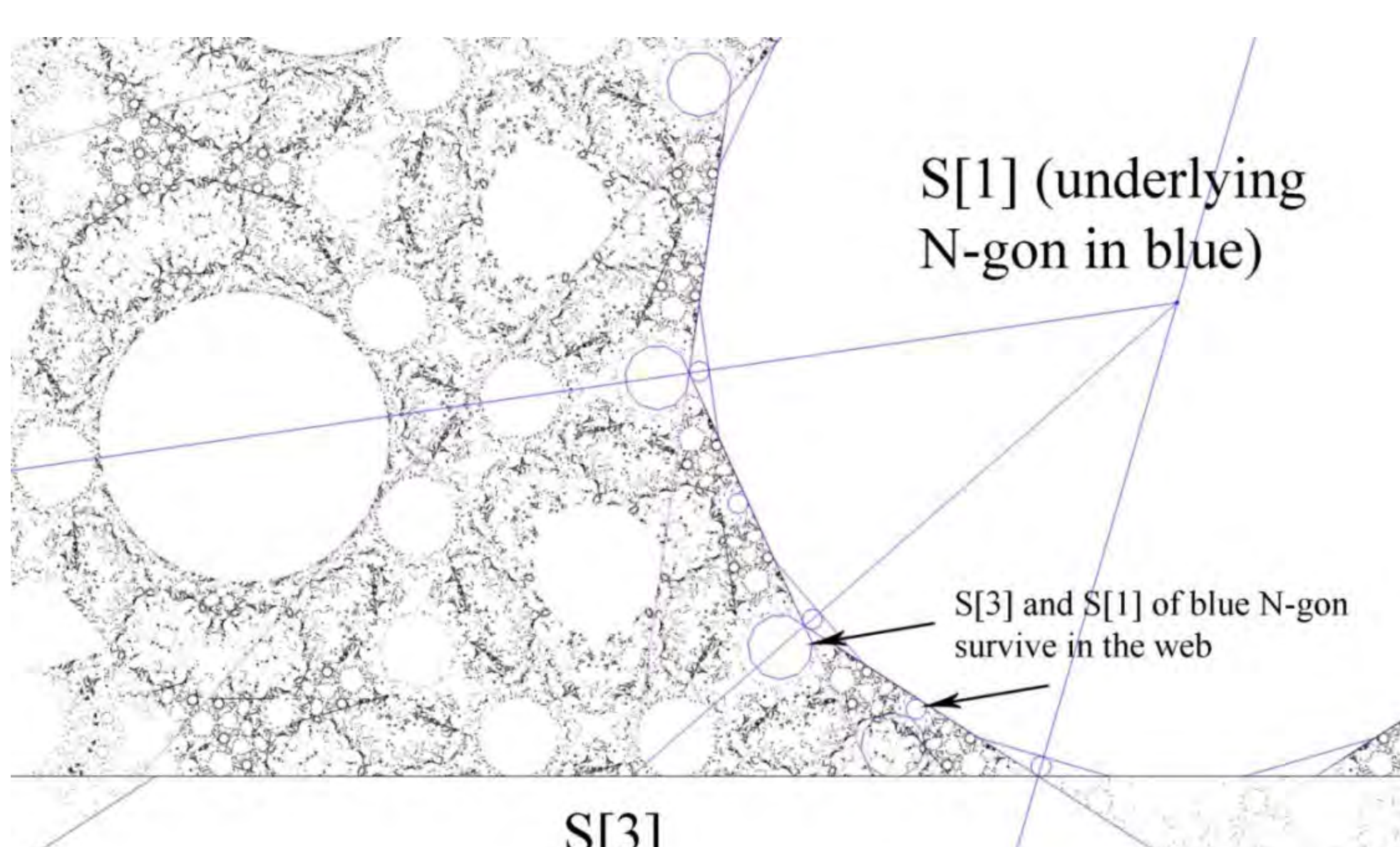

The table below gives what we call the 'characteristic polynomials' for height and midpoint of the major tiles: S[1], S[2], DS[5], Px , Mx and Sx. These are polynomials in x = GenScale[11] that yield the relative height and horizontal displacement of these tiles. In general, these midpoint polynomials are redundant and can be derived from the height polynomials, but they can be useful to simplify the construction of the corresponding polygons. For conforming polygons there is a linear relationship between height and displacement from the star point.

Except for S[2] at D[1], all of the tiles in the tables are regular 11-gons so by convention their heights will be given relative to the S[9] tile of N = 22 (a.k.a M) . This tile is the surrogate N = 11 within the First Family of N = 22 and hM/hN = Tan[π/22]/Tan[π/11] ≈ 0.489664. When N is twice-odd there can be issues with using GenScale[N] to generate the scaling field, so we will use x = GenScale[N/2] as generator, so for N = 22, x = GenScale[11] = Tan[Pi/11]·Tan[Pi/22].

In the table below we have chosen to measure the displacement of a tile relative to star[1] of N instead of the origin. This is more meaningful than the displacement from the origin, and in the Dc map this star[1] point is the origin. Therefore we will use the sN = 1 convention for the displacements, with the origin at star[1] of N. This means that for all N-gons, the critical S[1] tile will have displacement -1/2. When N is even star[1] of S[2] can serve as the local GenStar[N[ point and it will have displacement -1. When N is odd, the star[2] point of S[2] can serve as the local GenStar[N] and it has displacement -1 with the Dc convention.

**Table 22.1** The characteristic polynomials for height and midpoint of the major tiles

| **Polygon P** | **polynomial for hP using hN = 1 so hM =Tan[π/22]/Tan[π/11]** | **polynomial for Midpoint of P using sN = 1 and origin at star[1] of N** |
|---|---|---|
| S[1] (M[1]) | hS[1]/hM = hM[1]/hM = x | hS[1] = x*hM*hN so MidS[1] = -x*hM*hN*Tan[5π/11] = -1/2 (for all N) |
| S[2] (D[1]) | hS[2]/hN = hD[1]/hD = x | MidS[2]/sN = $-\frac{3}{2}-\frac{x}{2}$ |
| DS[5] (S[5] of D[1]) | hDS[5]/hM = $\frac{1}{8}[1-22x+8x^2+6x^3-x^4]$ | MidDS[5]/sN = $-\frac{3}{16}+\frac{11x}{8}+\frac{5x^2}{4}+\frac{x^3}{8}-\frac{x^4}{16}$ |
| Px | hPx/hM] = $\frac{1}{4}[-1+25x+5x^2-x^3]$ | MidPx/sN = $-\frac{35x}{8}-\frac{25x^2}{8}-\frac{x^3}{8}+\frac{x^4}{8}$ |
| Mx | hMx/hM = $1-23x-\frac{27x^2}{2}+\frac{x^4}{2}$ | MidMx/sN = $-1+\frac{87x}{4}+13x^2+\frac{x^3}{4}-\frac{x^4}{2}$ |
| Sx | hSx/hM] = $-\frac{9}{4}+52x+\frac{63x^2}{2}+x^3-\frac{5x^4}{4}$ | MidSx/sN = $-\frac{1}{2}+\frac{31x}{4}+\frac{25x^2}{4}+\frac{x^3}{4}-\frac{x^4}{4}$ |

When N is twice-odd, sM = sN so either one can be used for reference for midpoints. For the Dc map these are both 1 and this will simplify the calculations for midpoints in the same way that the height 1 convention simplifies the height calculations. Because of the reflective symmetry of the edge geometry any existing polygon P has a left and right-side version, and we will us the right-side version here so that the (negative) horizontal displacement of the midpoint of P can be measured relative to the origin at star[1] of N.

These height polynomials are the standard versions that we have derived in [H5] and they apply for any hN. This applies to the midpoint polynomials also but inside DKHO they will be shifted by -½ to match the shift in centers, so for the S[5] tile of N, the traditional MidS[5] is star[5][[1]]-sN/2 with polynomial MidS[5]/sN = $-\frac{61}{16}+\frac{55x}{8}+5x^2+\frac{x^3}{8}-\frac{3x^4}{16}$ and inside DKHO this is still valid with sN = 1 so MidS[5] = $-\frac{53}{16}+\frac{55x}{8}+5x^2+\frac{x^3}{8}-\frac{3x^4}{16}$ where the constant terms is off by -½.
In the table above we have included the polynomial for hS[5] of D[1] even though this is easily derivable from the polynomial for hS[5] of N – because these two just differ by x. To find the polynomial for S[5] of N, the First Family Scaling Lemma says that hS[1] /hS[5] = scale[5] = Tan[π/22]Tan[5π/22] and hS[1]/h[N] = Tan[π/22]$^2$. Graphically it is easy to import DS[5] and DS[1] (M[2]) as part of the First Family of S[2] using :**FFS2 =TranslationTransform[cS[2]] /@(FirstFamily*GenScale);**

● **N = 23**

N = 23 has complexity 11 along with the matching N = 46. N = 23 is in the 8k+7 family so it will have early DS[3]s and the 8k+7 Conjecture makes a number of predictions about the efficacy of these DS[3]s. The most important prediction is that each DS[3] will evolve in tandem with dual DS[1]s which will be S[N-3] tiles in the (2N) First Family of DS[3]. In the case of N = 7, DS[3] is the S[1] tile of N and indeed it supports DS[1]s as S[4] tiles. The evolution of DS[1] is discussed as part of the analysis of N = 7 because this tile plays an important role as an S[2][2] or D[2] candidate for ‘next-generation’ S[2]. We believe that this true throughout the 8k+7 family, but N = 15 is somewhat special because of the mutation of DS[3]. That makes N = 23 an important test case for the 8k+7 Conjecture. The symmetry diagrams for N = 39 in Figure 2.7 shows the symmetry that is shared by all members of this family. Here with N = 23 we trace the first 7 interations of the S[1]-S[2] web and see that it is remarkably similar to the case of N = 7. This helps to explain the shared symmetry.

**Figure 23.1** The level-7 web in magenta

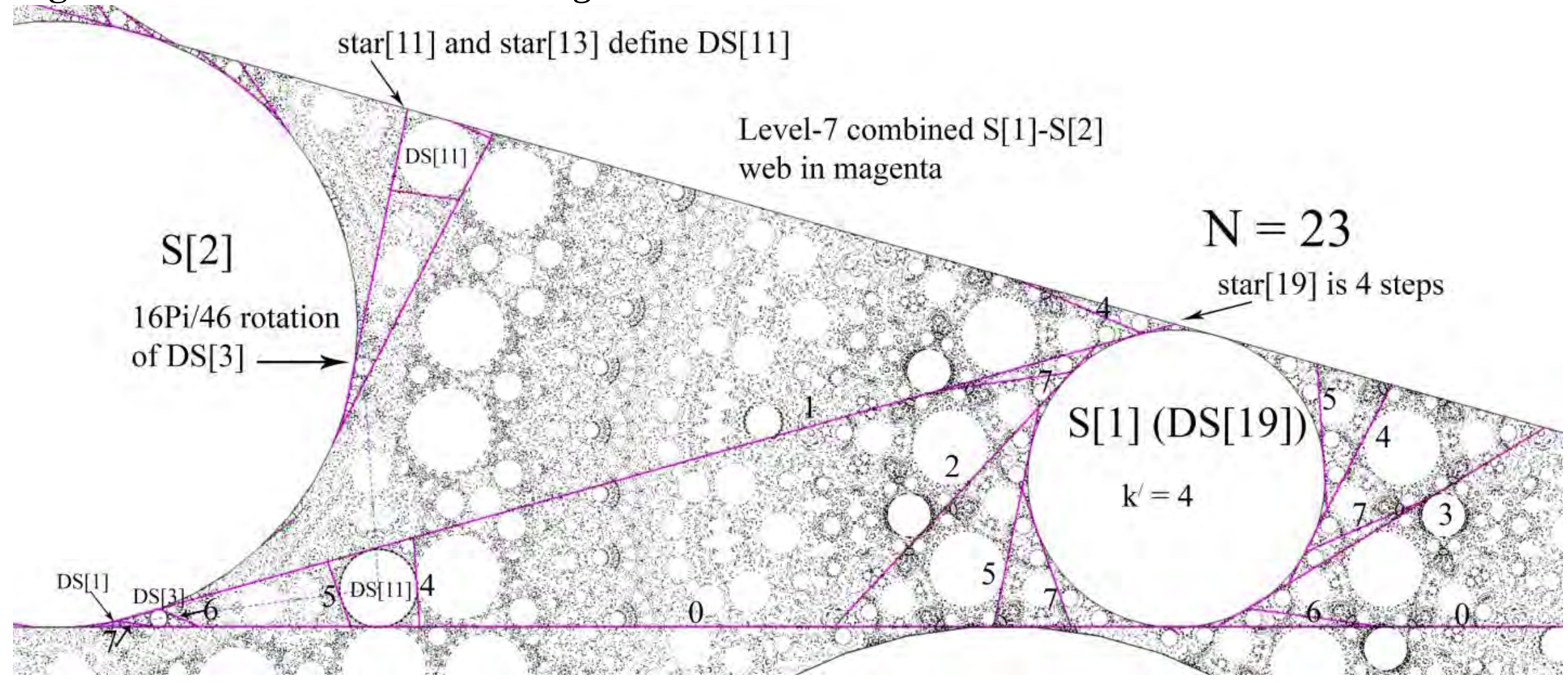

The predicted DS[k] count down mod-8 from S[1] at DS[N-4] so here they are DS[11] and DS[3]. Since N is odd the web is based on star[2] of S[2] and this extra displacement from star[1] to star[2] means that DS[3] and DS[11] will be displaced and upscaled relative to the (virtual) S[3] and S[11] in the First Family of S[2].

Since S[2] is a 2N-gon, the DS[k] will evolve with k′ = 2N/2-k = N –k as in the N-even case. For S[1] at DS[19] this matches the 2k+2 steps for N-odd. By symmetry these 4 steps can be partitioned into 2 + 2 to match the star[2]s at top and bottom of S[2] as in Figure 2.8. The 12 steps of DS[11] are also partitioned into 6 + 6 as shown by the blue dotted lines above. In the 8k+7 family the web of DS[3] always seems to evolve in a step-6 fashion since Mod[2N,N-3] is 6 for N > 7. Even though DS[1] is a volunteer it should have k′ = N-1 with Mod[2N,N-1] = 2 and indeed DS[1] has a step-2 web as shown below. Therefore we conjecture that DS[1] and DS[3] will have a step-2 and step-6 (cw) webs for the 8k+7 family beyond N = 7 where Mod[14,4] = 2 does give the correct symmetry with respect to the underlying S[1] in the N = 14 family.

Here DS[1] has S[2]s at step-1 and a trio of S[1], S[2] and S[4] survivors at step-2, while DS[3] has at least S[1] and S[2] at step-2 intervals (and overall step-6 symmetry).

**Figure 23.3** The web local to DS[1] showing the step-2 and step-6 webs of DS[1] and DS[3]

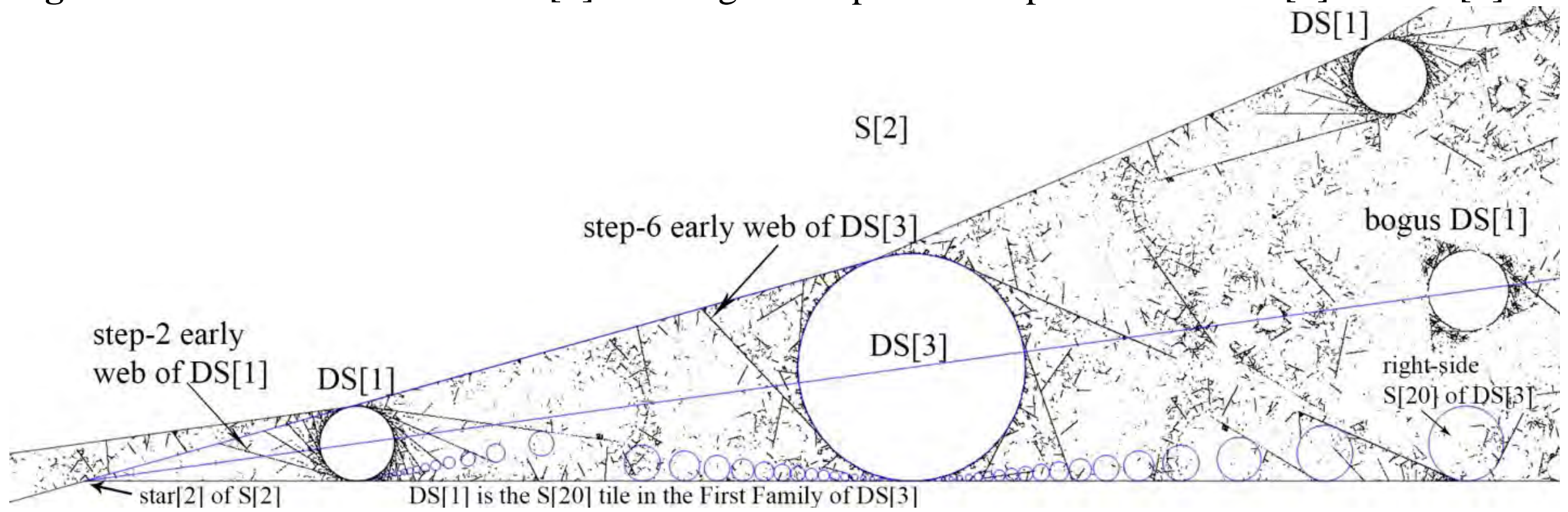

**Figure 23.3** Detail of the web local to DS[1] showing S[2]s at step-1 and a mix of S[1]s, S[4]s and S[6]s at step-2.

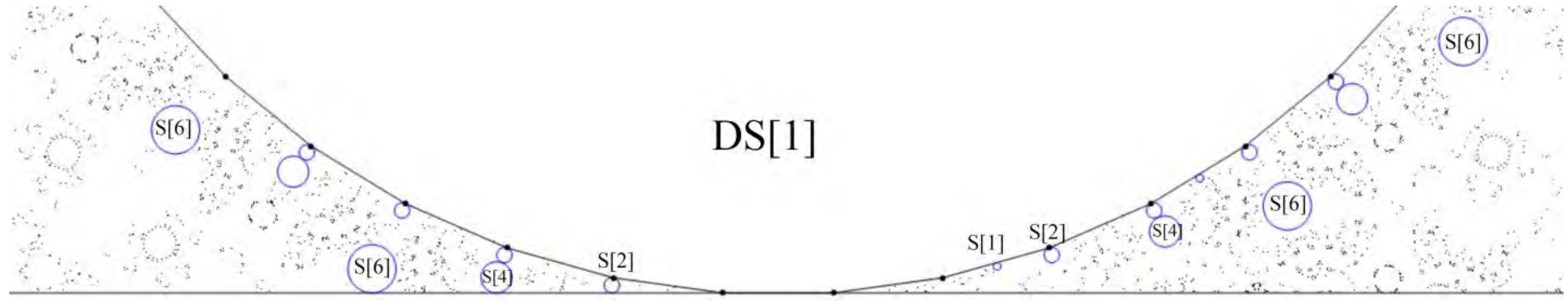

**Figure 23.4** Detail of the step-6 web local to DS[3]

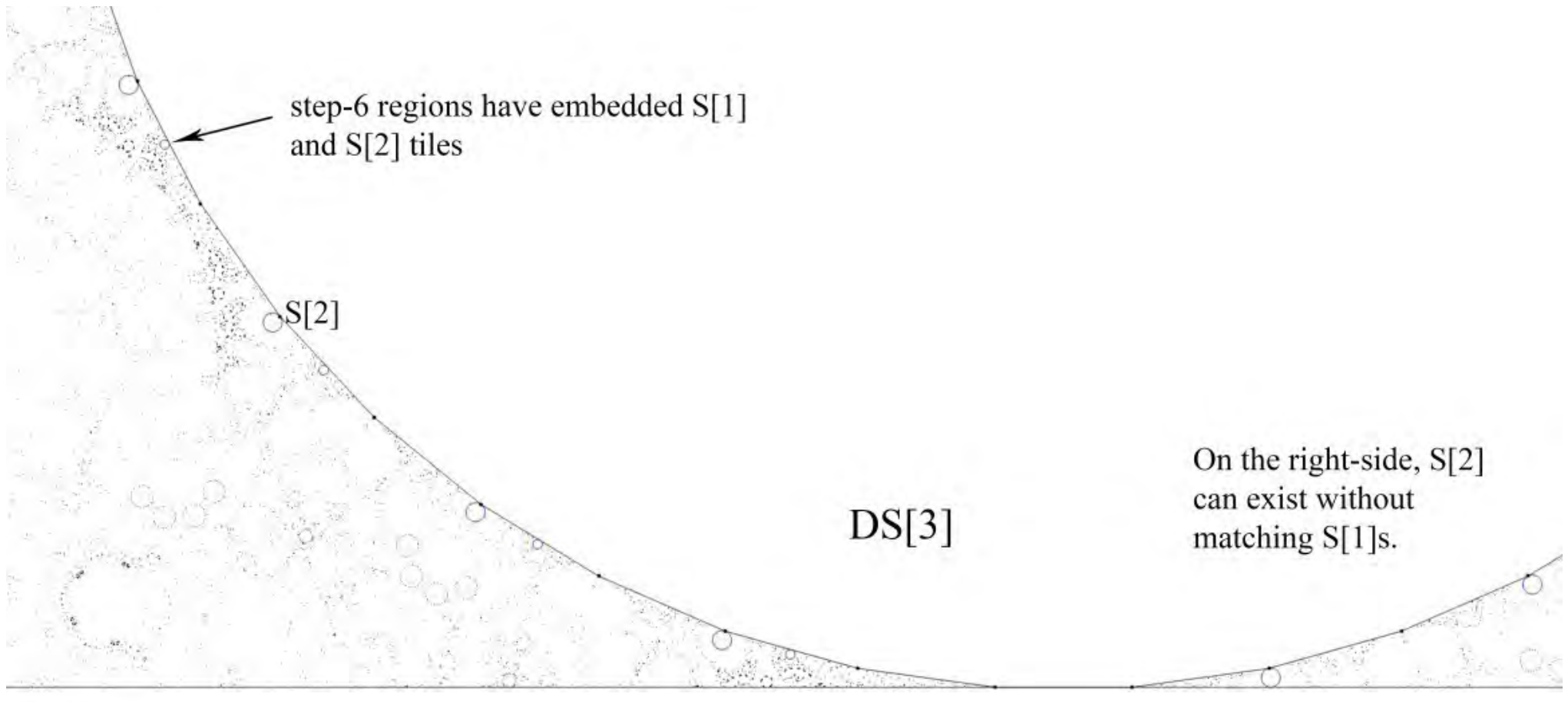

Any DS[k] web can be regarded as ccw or cw but based on our ccw convention for S[1] and S[2], it is natural to adopt the opposite cw orientation for the DS[k] – excluding S[1]. For tiles like DS[1] and DS[3] which are vertex tiles of S[2], this opposite orientation makes the most sense. Each new generation should have web orientation reversed from the previous. Here S[2] can be regarded as a D[1] patriarch of the 2nd generation and DS[1] would be expected to play the role of a D[2] which is patriarch of a third generation This means that the 4th generation, presided over by a D[3], should have a ccw web in a manner similar to S[2]. This should imply some form of shared dynamics and shared geometry with the 2nd generation. But as N grows, resolving 4th generations becomes a non-trivial issue.

We will use N = 23 as a ‘token’ N-odd example to compare the S[1] and S[2] web development with the case of N = 22 above. From the plot below it is clear that they are very similar, except that the steps are doubled and based on star[2] of S[2]

**Figure 23.5** The S[2] (right-side) web in magenta) and S[1] (right-side) web in blue.

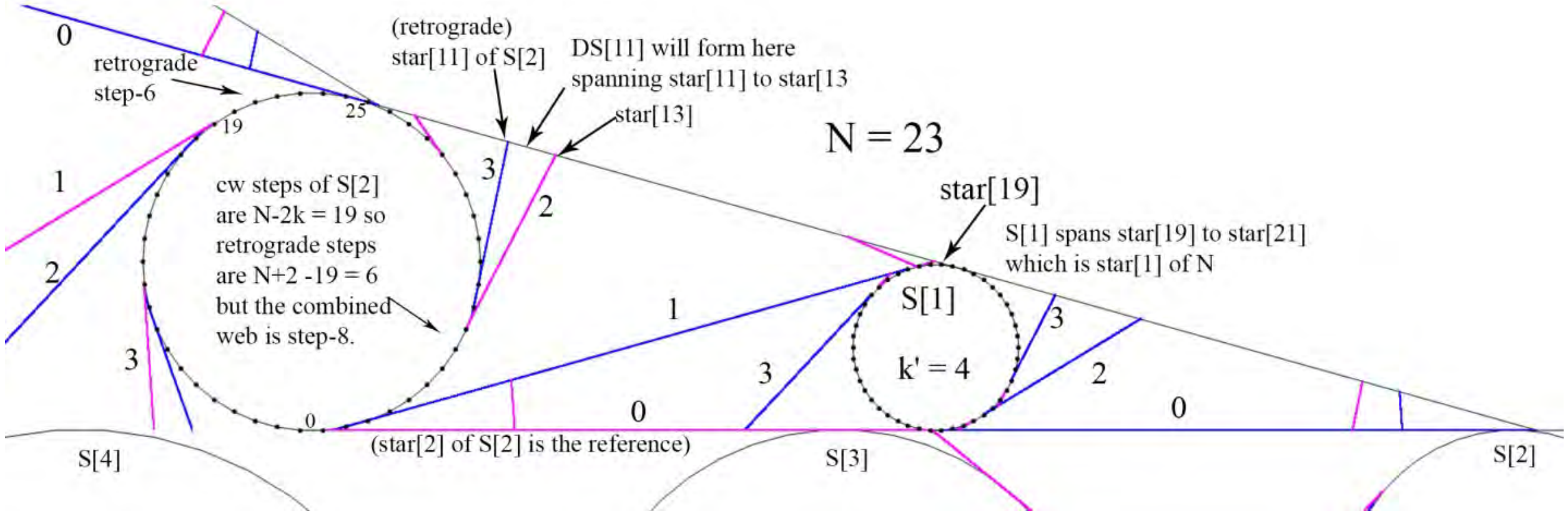

**Figure 23.7** The combined web local to S[1] has monochrome islands that may be invariant

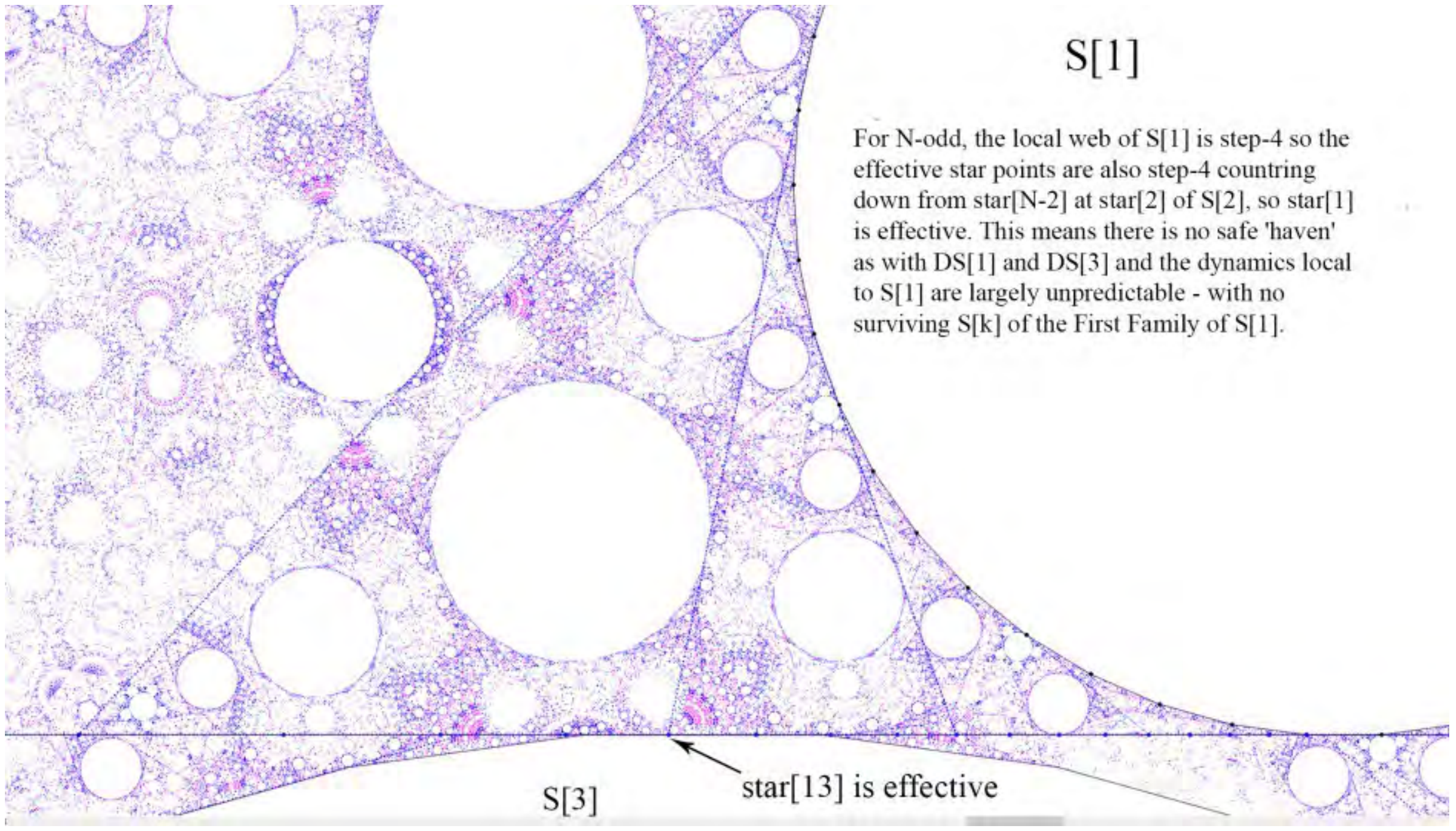

## • N = 24

N = 24 is in the 8k family. It is a quartic polygon along with N = 16 and N = 30. Algebraically it is also related to N = 8 and N = 12 since n|m implies that for cyclotomic fields $\mathbb{Q}_n \subseteq \mathbb{Q}_m$. $\mathbb{Q}_{24}$ is generated by $\sqrt{2}$ $(2Cos[2\pi/8]), \sqrt{3}$ $(2Cos[2\pi/12])$ and $i$, while $\mathbb{Q}_8$ and $\mathbb{Q}_{12}$ are generated by $\{\sqrt{2}, i\}$ and $\{\sqrt{3}, i\}$ respectively.

**Figure 24.1** The local (right-side) families of the S[k] as predicted by the GFFT of [H5]

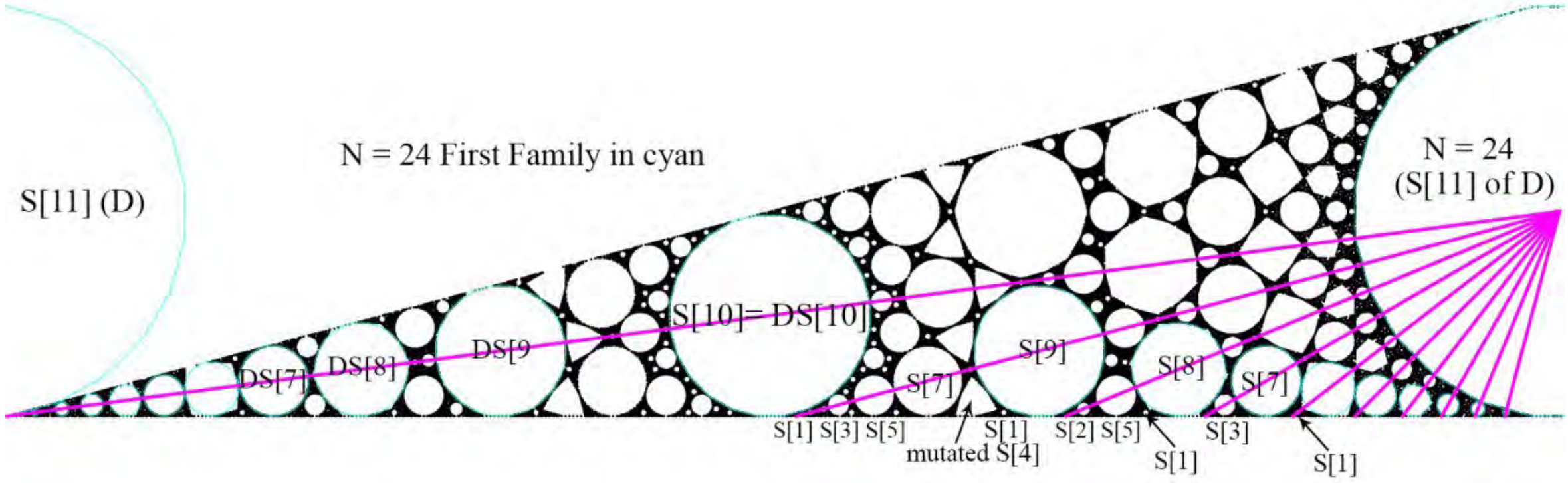

The actual S[k] are shown in cyan. As expected there are many mutations. For N twice-even, the M tile is S[N/2-2] which is S[10] here. The step sizes are k′ = N/2-k so M is step-2, but it is not mutated because the mutation condition for an S[k] is gcd(k,N) > 2 . This twice-even case is more forgiving than the twice-odd case because the two cycles that form each S[k] are N/2-1 steps apart and now this is odd, so the cycles are non-redundant when gcd(k,N) = 2. Therefore the mutated S[k] of N = 24 are S[3], S[4], S[6], S[8] and S[9]. It is not hard to predict the form of these mutations and it appears that most local families inherit the mutations of the S[k].

**Figure 24.2** The mutations of the First Family for N = 24. The S[1] tile is never mutated for N twice even, but S[2] is mutated when N is 8k+4.

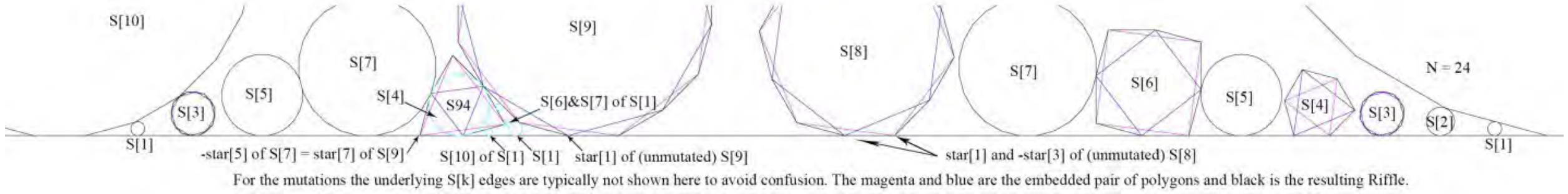

**Figure 24.3** The web local to S[1], S[2] and S[3]. As shown by the insert, every external 'factor polygon' of N = 24 shares star points with N.

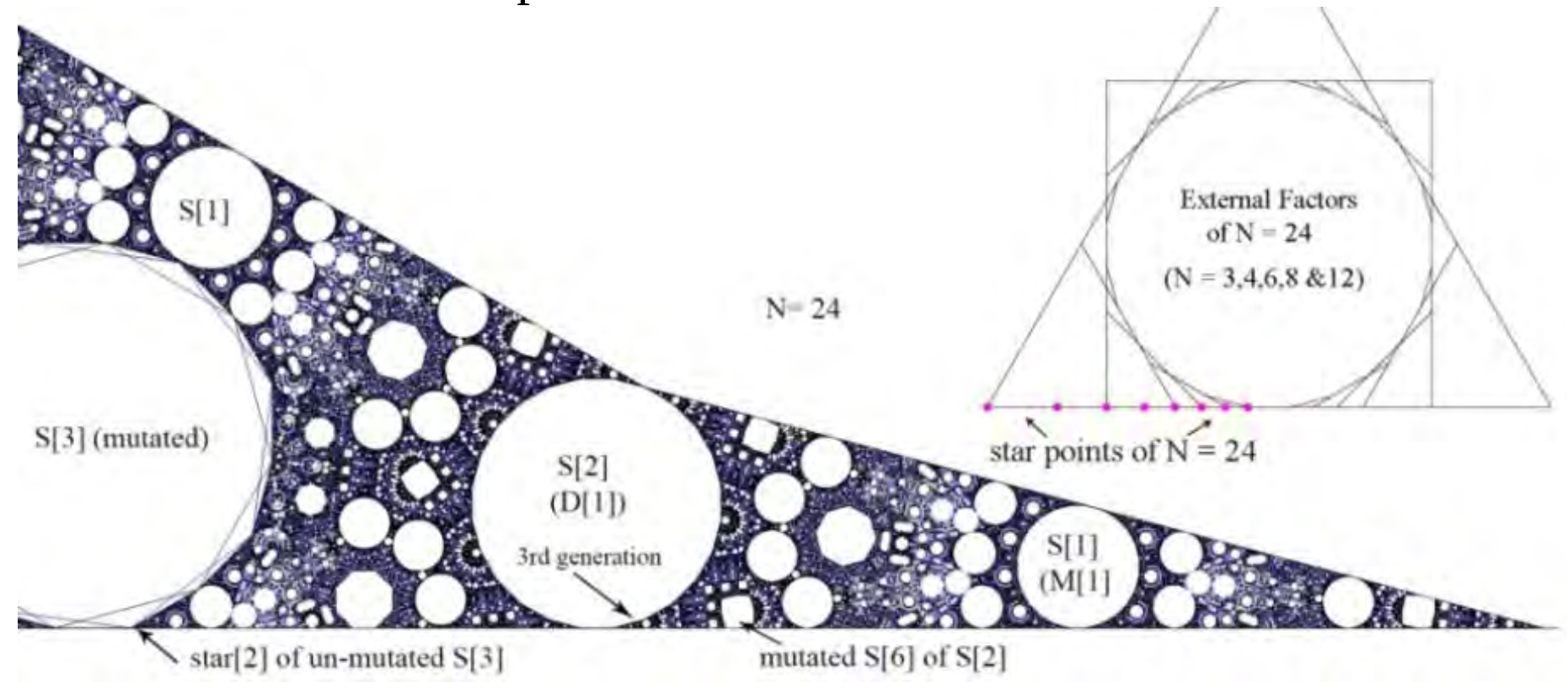

Since N = 24 is a member of the 8k family it has an isolated D[2]. But unlike N = 16, the local web of D[2] does not generate an M[3] or D[3] to make up for the missing M[2] so there is little hope of a 'normal' generation evolution at the foot of S[2].

The lack of M[3] and D[3] can be seen in the 3rd generation enlargement below. This failure may be due to the interaction of the mutations for N = 24. From experience with N= 9 , 12 and 16, it is clear that individual mutations can evolve in a predictable fashion, but there is no theory that attempts to explain how distinct mutations interact. Here it is clear that D[2] is not formed in a 'normal' step-2 fashion and the only candidate for an M[3] is highly mutated.

The S[6] tile of D[1] is also mutated but it is possible to construct the resulting tile based on the star points of the unmutated S[6]. However each generation becomes more difficult to track and the 4th generation shown below is almost unrecognizable.

**Figure 24.4** The 3rd generation

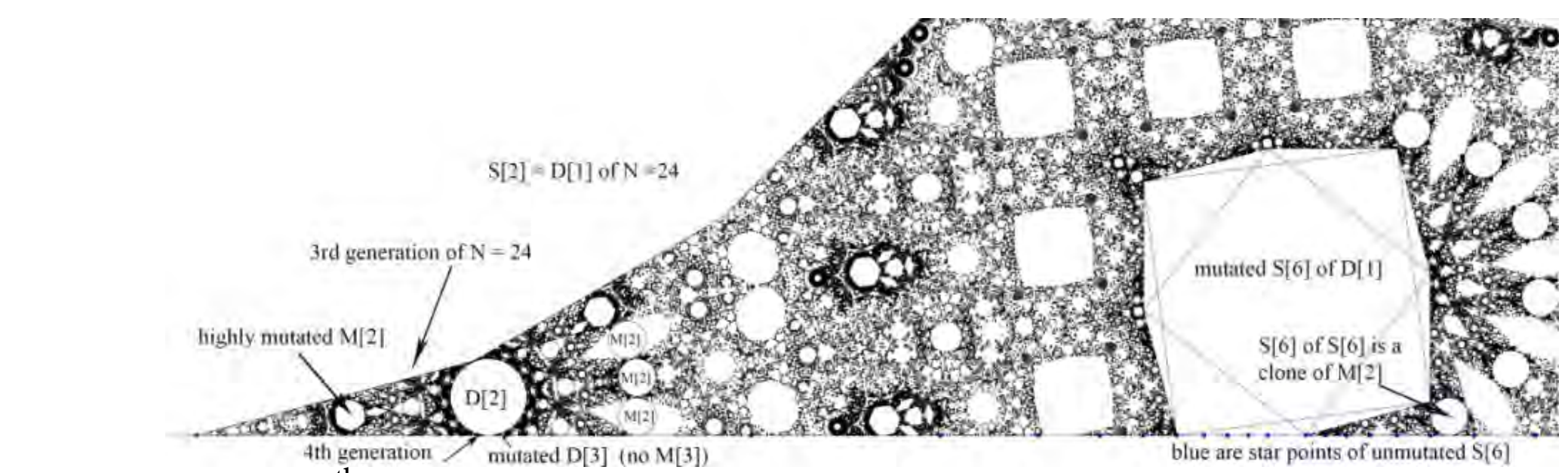

**Figure 24.5** The 4th generation is highly mutated

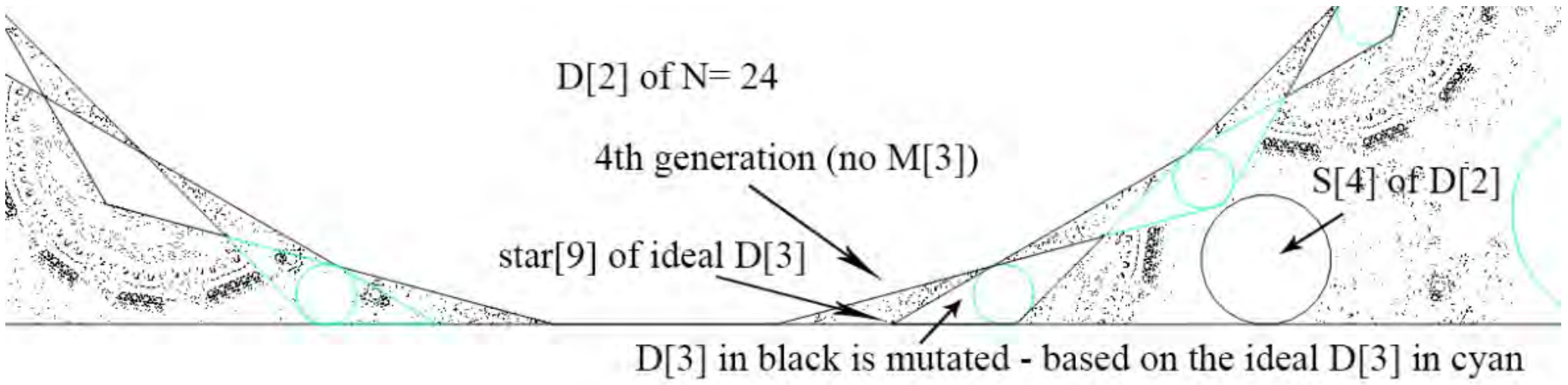

The canonical generation scaling is hS[1]/hN = $\mathrm{Tan}[\pi/24]^2$ = GenScale[24]. The '3rd generation' shown above is not a true 3rd generation because D[2] exists but there is no matching M[2] – so there is no sign of chains converging to GenStar or star[1] of D[1] – but such chains could exist elsewhere The mutated octagonal S[6] of S[2] is a scaled version of the original mutated S[6] – composed of two squares of different radii.

The S[1] region may be more promising as shown below. The Twice-even S[1] Conjecture says that since S[1] is step-2, the local web could support tiles in a 'Step-2' family where each tile is based on consecutive odd star points. The graphic below shows the ideal Step-2 family tiles of S[1] in blue. All of these blue tiles must be conforming to star[2] of S[1], so they are similar to a star[2] family of S[2] when N is odd. In both cases the tiles span 3 star points.

**Figure 24.6** The early step-2 web of S[1] is in magenta and the First Family and 'step-2' family are in cyan and blue. Their relationship is the same as S[1] and S[2] – namely scale[2].

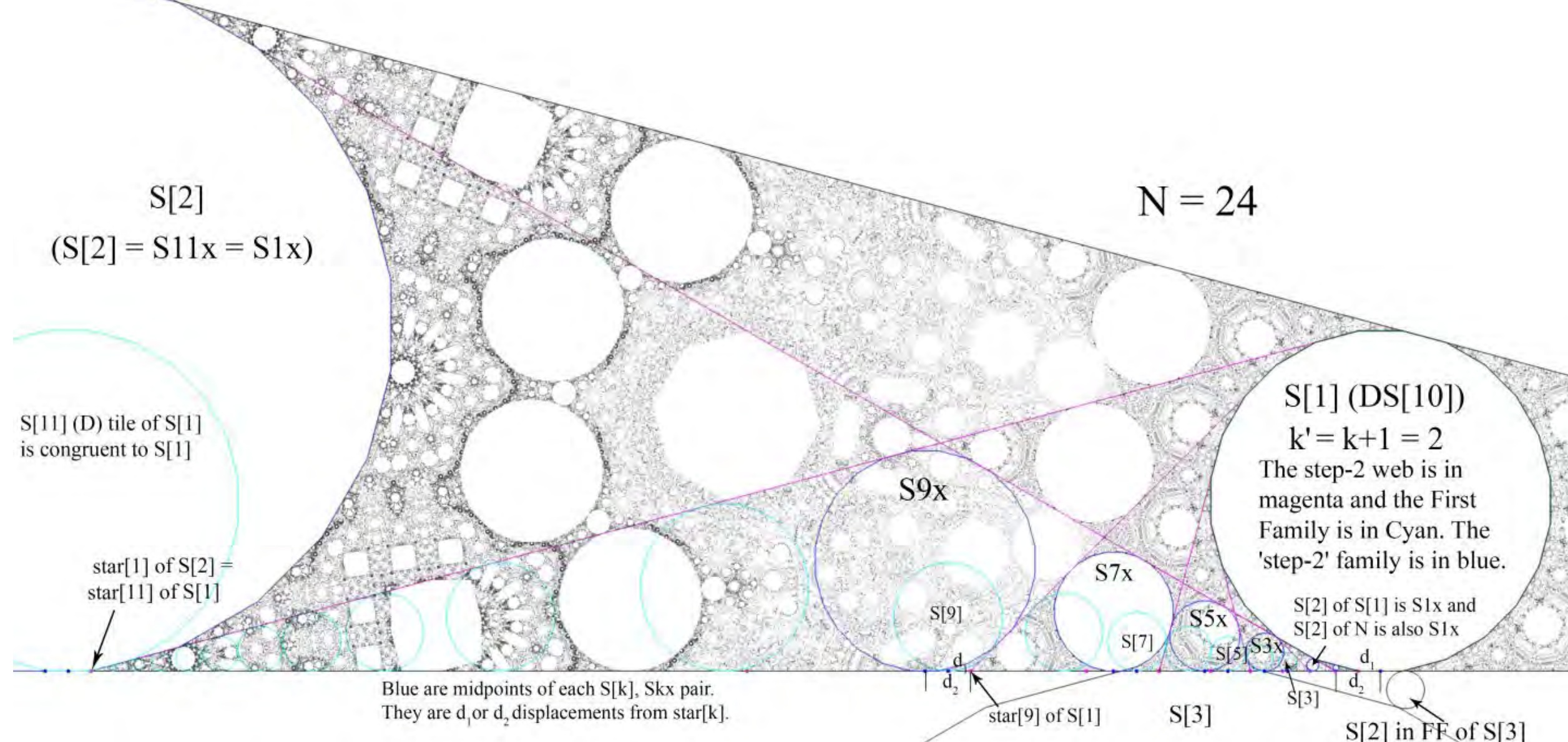

In a manner similar to the N-odd case, each blue Skx tile is a D tile relative to the matching S[k]. This occurs because their displacements from star[k] are related in the same fashion as N-odd where the traditional $d_1$ = sS1/2 displacement of S[k] is expanded out to $d_2$ at star[2] of S[1], so $d_1/d_2$ is scale[2] of N and heights scale the same as sides.

Therefore all the Skx are conforming to star[2] of S[1] and this includes S[2]. Since S[1] is DS[10] this means that star[10] of S[2] is star[2] of S[1] (and star[2] of N). By symmetry S[1] must be able to construct S[2] with the opposite displacement using star[11] as shown here. (Note that for N even the edges of S[2] and S[1] are related by scale[2], so for N twice-odd they are the same length because of the mutation of S[1].)

As N increases, S[3] is on a collision course with S[1], but the spacing for the 8k family seems to be determined by the fact that the S[2] tiles in the First Family of S[3] is the same as the S3x tile of S[1]. These are both vertex tiles and the fit is perfect. Note that S7x also shares a vertex (and hence dynamics) with S[3].

All of the Skx of S[1] are formed in the same fashion and here are the steps for S3x. (As above, the displacement $d_2$ used here is the displacement of star[2] of S[1] from the midpoint of S[1].)

```
hS3 (of S[1) = hS[3]*hS[1]; hS3x = hS3/scale[2]; MidS3x=starS1[[3]]-{d2,0};
cS3x =MidS3x +{0,hS3x}; starS3x=Table[MidS3x+{hS3x*Tan[k*Pi/24],0},{k,12-1}];
v1=starS3x[[1]]; S3x=RotateCorner[v1,40,cS3x];
```

Using star[1] of S3x to construct S3x is a way to avoid using the radius. In the twice-even family, based on our default orientation, there are no vertices in 'cardinal' positions.

**•N = 25**

N = 25 has algebraic complexity 10 along with N = 33, 44 and 50. It is a member of the 8k+1 family so it has a DS[5] as well as a volunteer DS[2].

**Figure 25.1** The level-7 web in magenta.

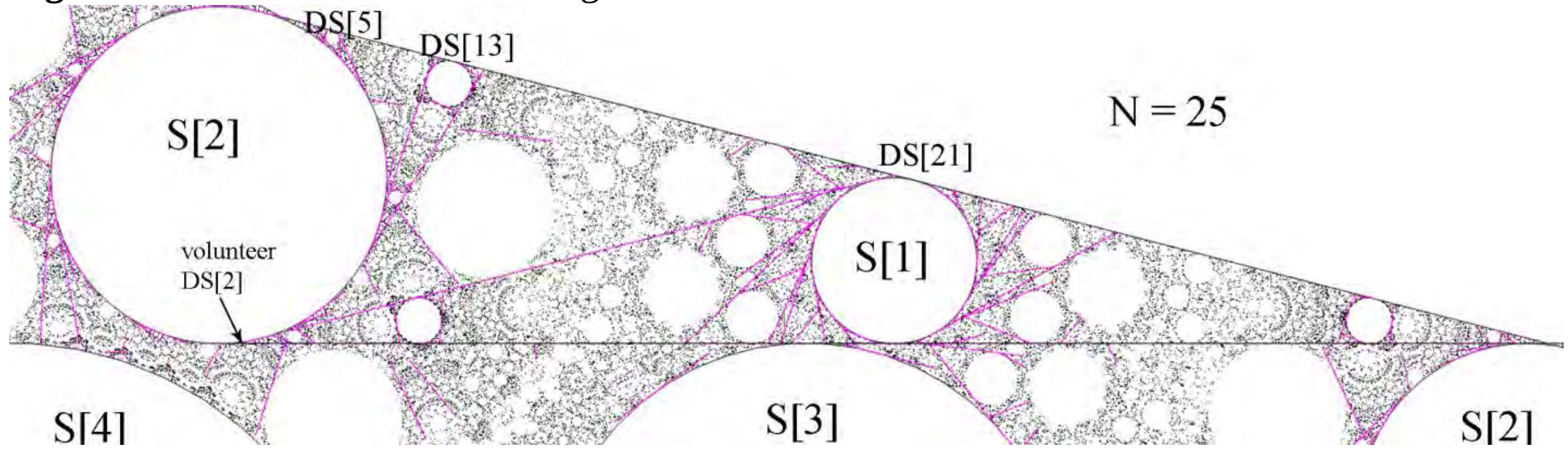

**Figure 25.2** The web geometry local to S[2] showing lines of symmetry

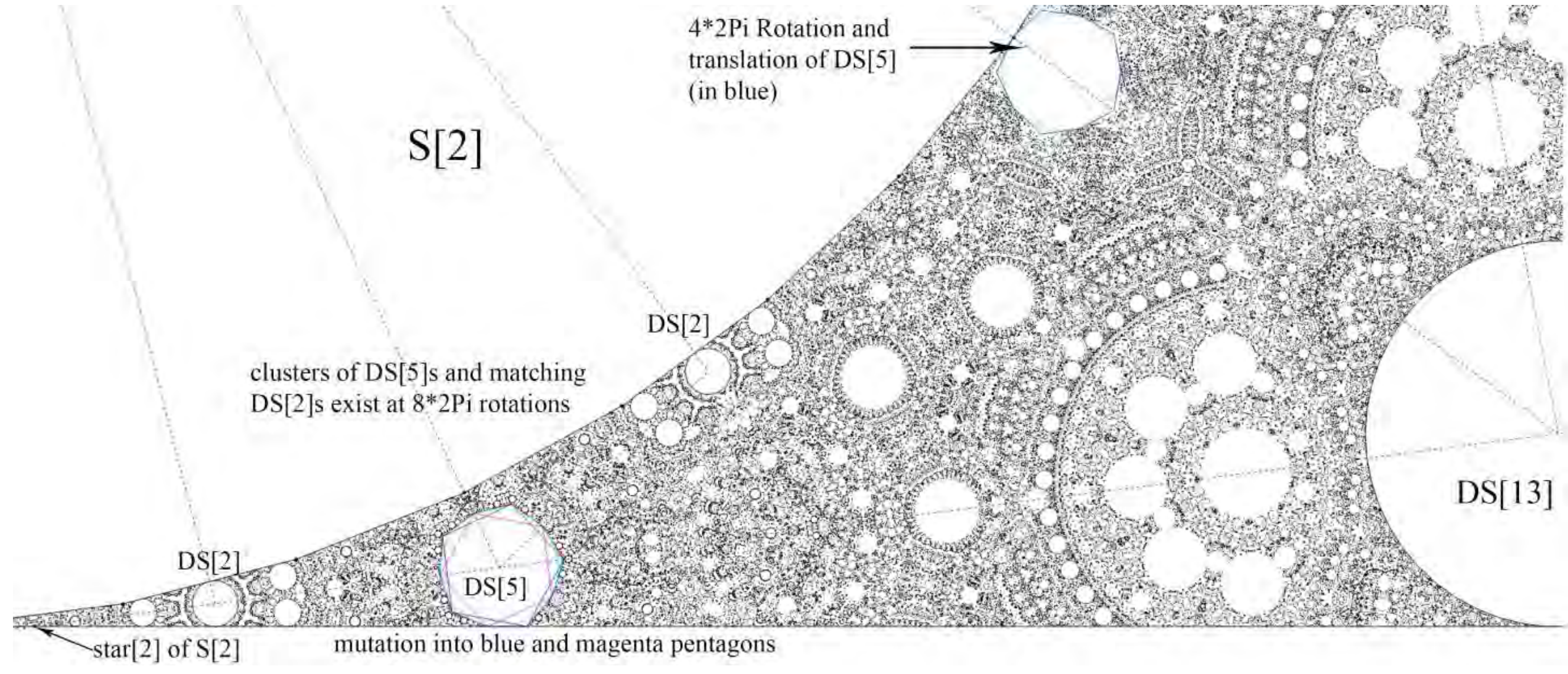

For DS[5], $k' = N\text{-}k = 20$ and $2N/\gcd(2N,20) = 5$ so the mutation consists of 2 pentagons with a horizontal span of 10 star points. The Mutation Conjecture says that the smallest star point should be the minimum of $N\text{-}2 \text{ -}jk'$ which is 23-20, so the magenta pentagon is anchored at star[3] of the underlying S[5] and the blue pentagon is at (right-side) star[7]. As always for N-odd, S[2] has step-8 rotational symmetry because this is the combined S[1]-S[2] web. Any step-8 web would be expected to have some degree of step-4 symmetry and here the step-4 version of DS[5] is a very close match for the web, but this volunteer has very different local geometry compared with DS[5].

**Figure 25.3** The region local to DS[2] is semi-variant and DS[2] has a step-4 web so it has 8*Pi/2N rotational symmetry with star[3] effective. There are clearly survivors from the First Family of DS[2].

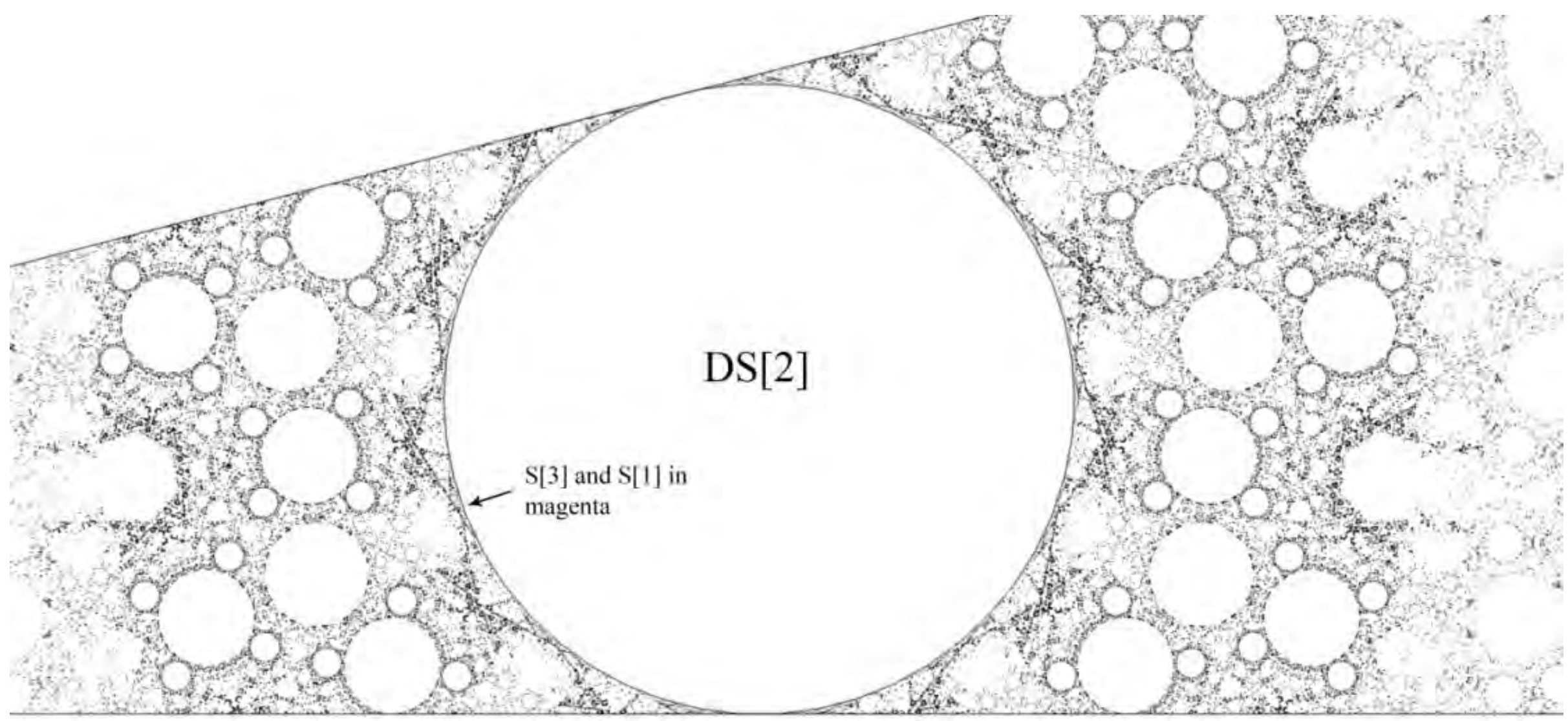

**Figure 25.4** The early web local to S[1] is always step-4 for N-odd since k′ = 2k+2.

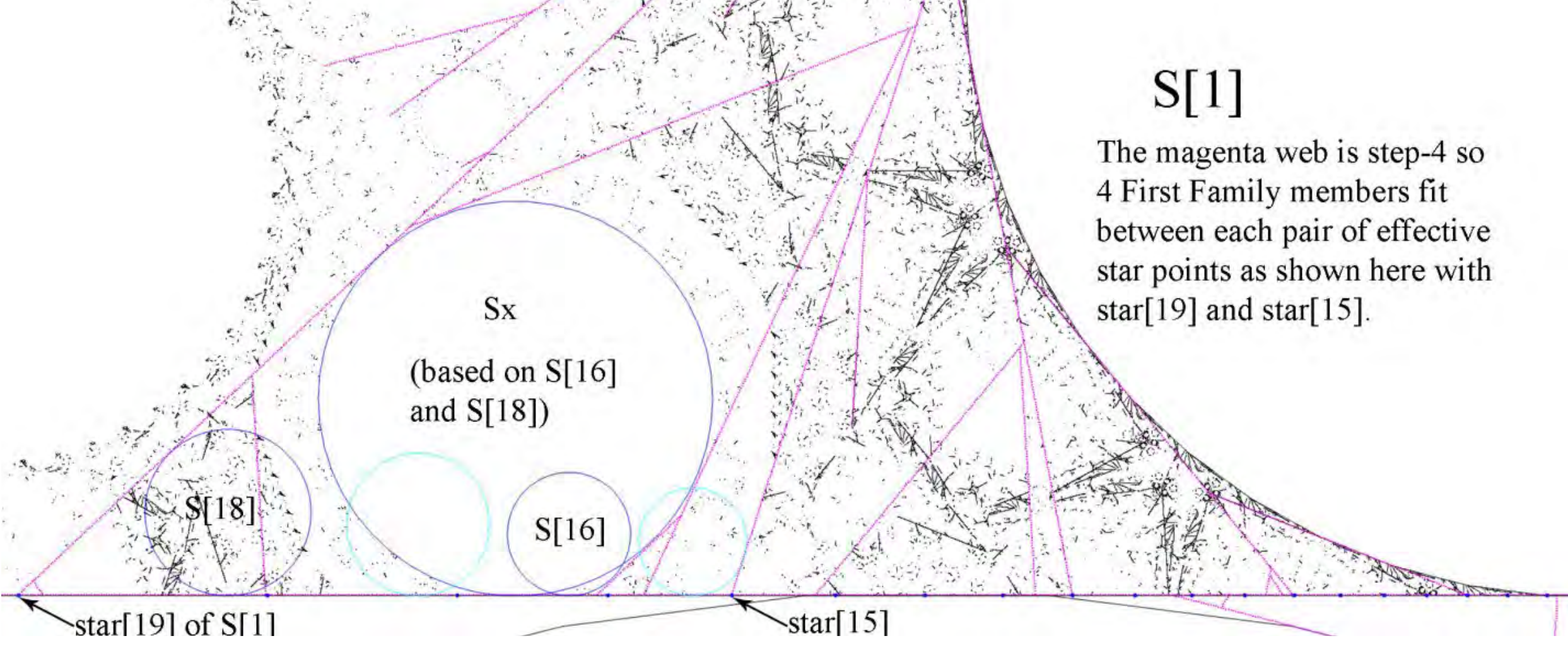

S[1] here is a 2N-gon which is congruent to the S[2] tile of D. This S[2] tile is in the twice-odd family of N so it would be expected to have a hybrid step3-4 web shared with the local S[1] but this S[1] is an N-gon and there seems to be little common ground with S[1] here. This step-4 web fractures early and this may be due to the influence of S[3] at the bottom. There is little sign of cohesive tile structure except for the Sx tile whose parameters appear to be based on the virtual S[16] and S[18] in the First Family of S[1].

There is no theory for ‘2 out of 4’ constructions as there is for ‘2 out of 2’ with N twice-even, but this 3 out of 4 web splitting does occur naturally as the web evolves. In general these step-4 webs are not well-behaved. Typically they have step-4 rotational symmetry with large-scale structures like those shown here.

**Appendix I :'Deep Field' maps of the edge geometry for N = 19 and N = 200**

In his Wikipedia article on the outer billiards map, (https://en.wikipedia.org/wiki/Outer_billiard) Richard Schwartz listed the foremost unanswered questions. One of these was: *Show that outer billiards relative to a regular polygon has almost every orbit periodic.*

This says that the points with non-periodic orbits should have Lebesgue measure 0. This Hamiltonian 'phase-space' conjecture has a long history beyond outer billiards and has never been proven, but we believe it is true. From a practical standpoint it means that every web W should be dominated by periodic tiles and this 'white matter' should have 'full measure' leaving only measure 0 for 'dark matter'. For the quadratic cases of N = 5,8,10 and 12, this is clearly true because the web W has a simple fractal structure where the non-periodic exceptional points are at most countable. For regular N-gons, these non-periodic orbits cannot originate inside tiles because every point in a tile has the same period and the 'inner star' region around N is invariant and bounded. Therefore no tile with non-zero measure can have a non-periodic orbit.

Any non-periodic point can potentially be of value because it may be possible to use the orbit to illuminate the tile borders. Indeed some of these orbits appear to be locally 'dense' in the limiting web. The star points of N are technically non-periodic, because they have no image under $\tau$ or $\tau^{-1}$, but the neighborhoods of these 'saddle-points' have the potential to act as 'candles' to illuminate the web structure. No point on an extended edge of N can be periodic because these points have no inverse image, but such points could be non-periodic and never quite map to an extended trailing edge. N = 5 has such orbits. Any orbit with very long period can possibly be used to trace the details of the web.

Below is the region adjacent to DS[7] of N = 19. This web is generated by choosing initial points on the x-axis within .001 of star[7] of S[2] and iterating each point at depth 30 million with the Dc map. (The orbits of star points of tiles like S[2] could be mapped under $\tau$ but such orbits would tend to terminate quickly at an extended trailing edge of N. )

**Figure A1**: A 'deep-field' image for N = 19 showing pockets of 'dark matter' around DS[7]

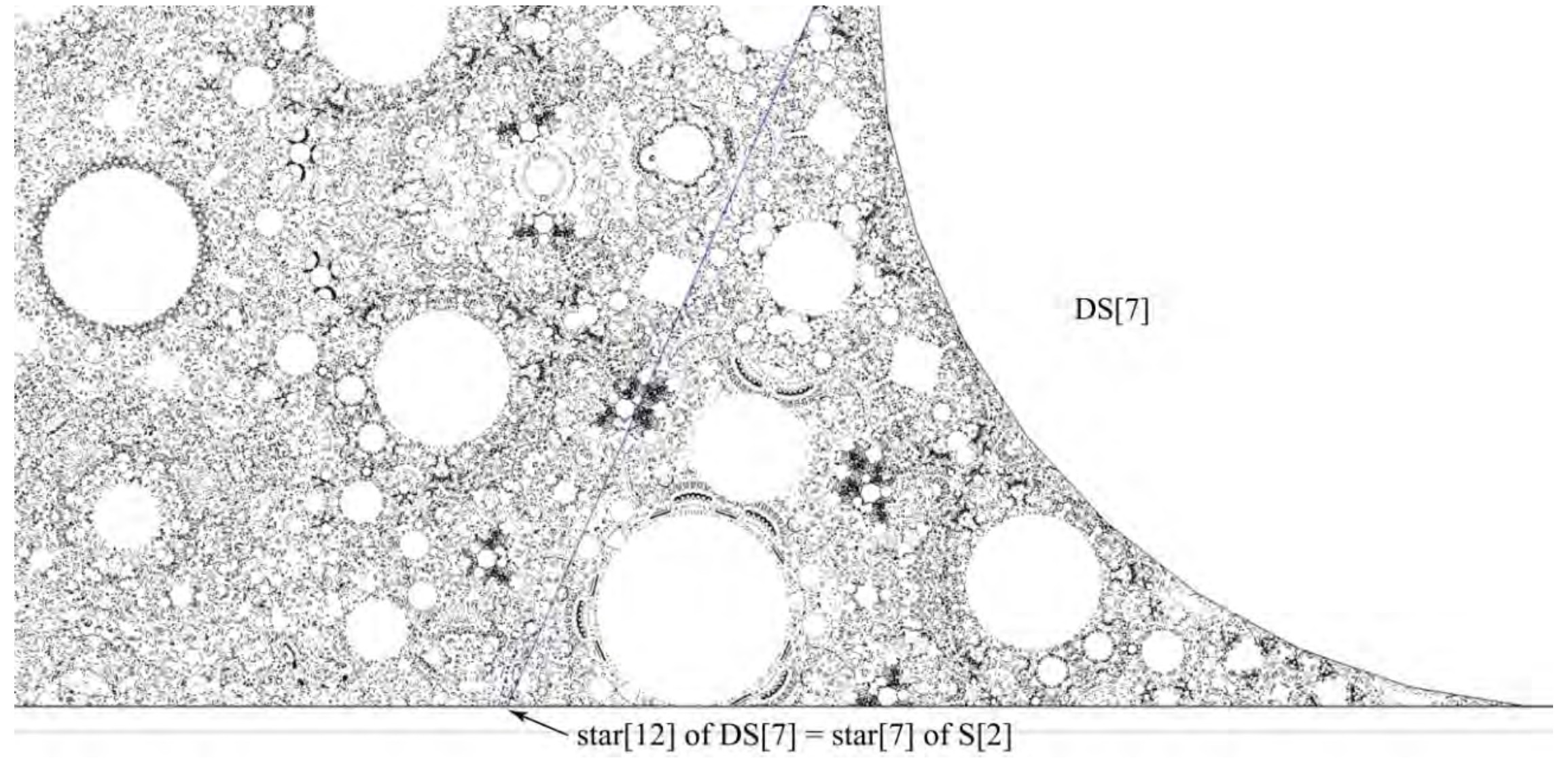

We conjecture that any ‘dark matter’ always dissolves into tile structure on closer examination so ‘most’ N-gons will have structure on all scales. By default, the images generated by Mathematica are vector based Postscript files and we typically use 35-decimal place accuracy for each point, so there is virtually no loss of detail on enlargement. But from a practical standpoint it is necessary to convert images from vector form to ‘raster’ pixel form for display or printing. This is usually done with a program like Photoshop or Adobe Illustrator. To keep the file size down, most of the images in this paper use a modest 200 dpi which would enable one or two levels of screen enlargement. These sample raw images were scanned by Photoshop at 600 dpi to give 7200 by 3500 pixels which is about 25 Mb raw and 4Mb compressed. The original Postscript file from Mathematica was 400 Mb. Even on a fast computer it can be a time-consuming process for Mathematica to generate these files. This N = 19 data set had about 1 million points and took more than 20 minutes to generate the Postscript file.

**Figure A2**: A ‘deep-field’ enlargement of the volunteer $D_2$ tile for N = 19

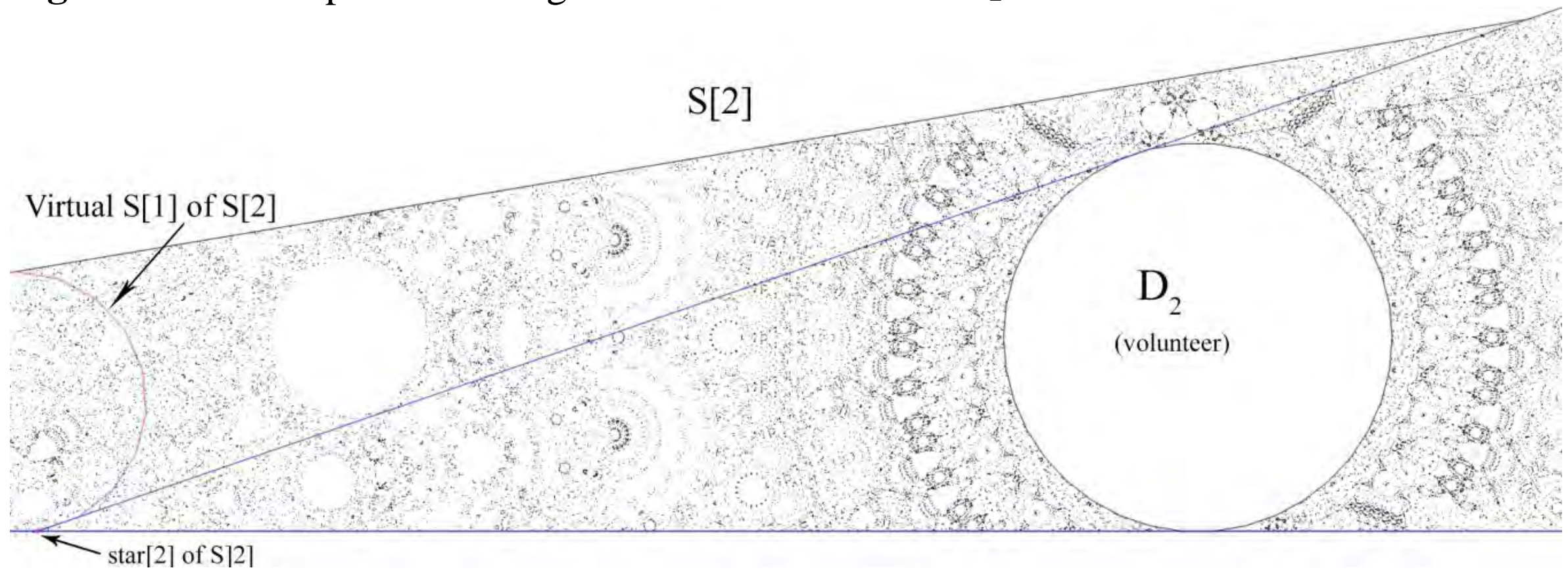

It is not clear whether large N values with higher algebraic complexity may yield denser webs with more potential for non-zero Lebesgue measure. With the Dc map the edge length is fixed at 1 and there is a nominal price to pay for larger N values – primarily the smaller rotation angle w = 2Pi/N and inherent loss of accuracy.

Below are some ‘deep-field’ scans of N = 200 with algebraic complexity 40. Mathematica only takes a few seconds to generate the First Family. D is now S[99] and with the height 1 convention for N, the center of D is at {-(cot[Pi/N] + tan[Pi/N),0} $\approx$ {-63.6725,0}. N = 200 is the 25$^{th}$ member of the 8k family so there is a tiny DS[2] web survivor to begin a mod-4 chain of survivors ending at S[1] which is DS[98]. Since N is even, these DS[k] are the same as the S[k] of S[2] and the only mutations will be k = 10,30,50,70, 90 all with gcd(200,k) = 10 except for k =50. The natural reference here for points on the x-axis are the 99 star points of S[2] which are

**StarS2 =Table[MidpointS[2] + {hS[2]*Tan[k*Pi/200],0} , {k, 1, 200/2-1}];**

**Figure A3** A 600 dpi deep-field image of the S[1]-S[2] web for N = 200

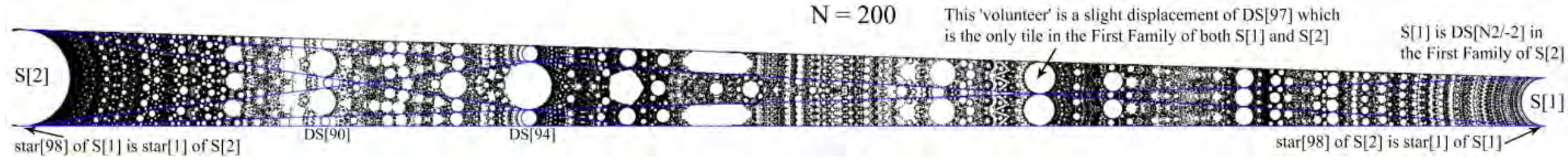

A typical collection of scan points for Dc might extend from star[40] to star[70] of S[2]
**H =Table[x, {x,Tw[StarS2[[40]]][[1]],Tw[StarS2[[70]]][[1]], .000131}]**
Where Tw converts between normal tau-space and Dc space. These 150 points can then be iterated 'overnight' to depth of about 40 million in a modest computer.

For N even the retrograde steps of S[2] are $k' = k+1$, but the 'effective' steps are 4 since the web of S[2] is linked with the web of S[1]. This means that the surviving DS[k] of S[2] will be mod-4 counting backwards from DS[N/2-2] as S[1]. It also means that S[2] will have step-4 symmetry and rotation by 8*Pi/200 will yield a locally equivalent web. One salient feature of this combined web is that the blue DS[k] survivors will have local geometry which is linked to both S[2] and S[1]. The left-side geometry will tend to be more closely tied to S[1] than the right-side. This can be observed with the penultimate DS[94] below which lives in both worlds. Beyond DS[94], the geometry may change dramatically as the largely autonomous web of S[1] dominates.

**Figure A4** The DS[94] region of N = 200. DS[94] is the last predicted survivor before S[1]

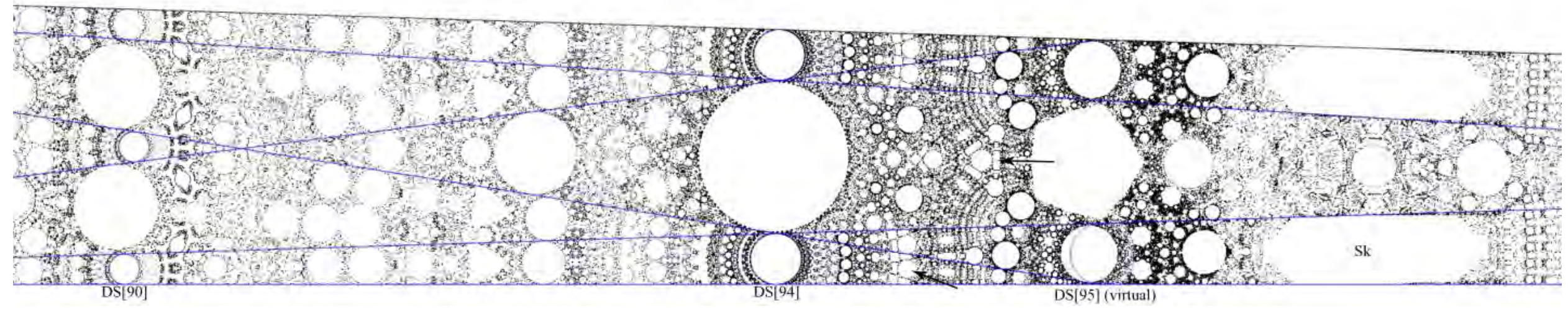

The Sk tiles are apparently octagons similar to those observed with N = 15 and 18. They are weakly conforming to S[2] and S[1] since they share blue edges. They may be tied to the breakdown of the virtual DS[95] shown here in magenta. Smaller Sk are shown at the arrows. Typically they are composed of linked regular tiles with shared edges.

**Figure A5** The S[2] region shows the strong step-4 symmetry expected of the 8k family.

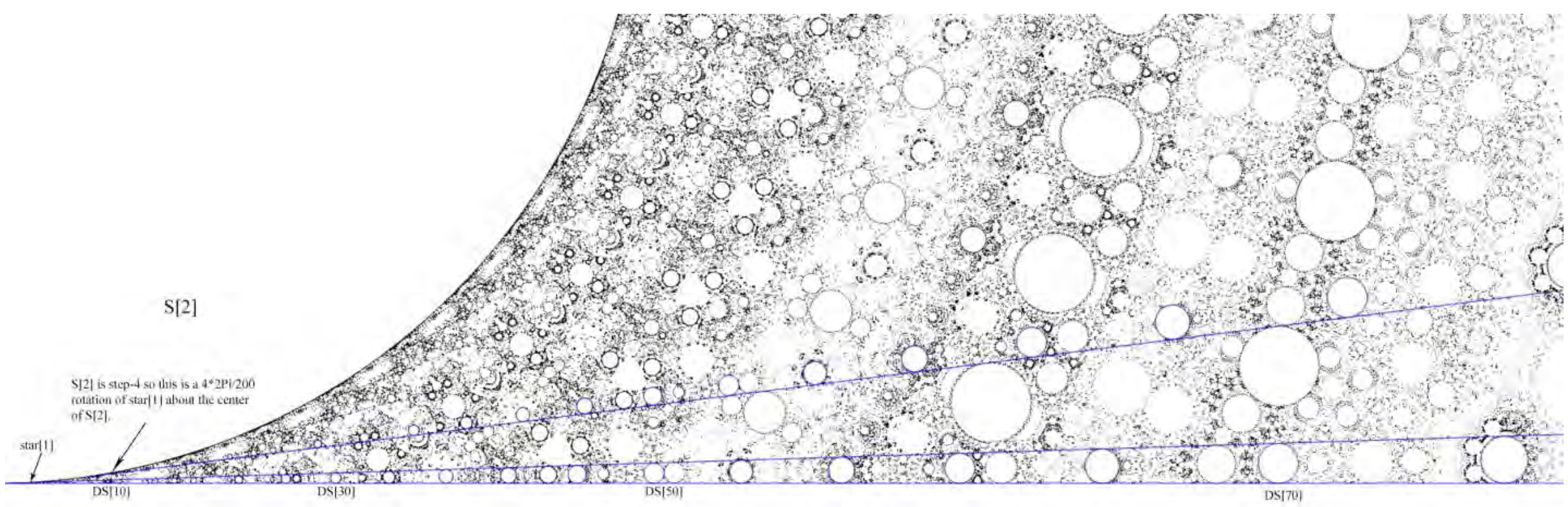

The smallest mod-4 survivor is DS[2] which is the S[2] of S[2] so it is no surprise that it has web features that resemble S[2]. GenScale[200] = $(1-\cos[\frac{\pi}{200}]^2)\sec\left[\frac{\pi}{200}\right]^2$ $\approx .00024\ 67$ so it will be a challenge to explore future generations. These images of N = 19 and N = 200 each took about 5 billion iterations to generate a paltry 1 million points. In the limit with the convention of edge length 1, the rotation angle of N would vanish and the N gon would become an unbounded ray.

## Appendix II: The invariant regions of N = 202 and the matching 'projections'

N =202 is the 25th member of the 8k+2 family so it can provide some insight into what to expect for large N when the S[2] tile supports a well-defined family of D[k] and M[k] tiles with known geometric and temporal scaling. Since N is twice-odd the last surviving S[k] in the inner invariant region will be S[Floor[N/4]] = S[50. As noted in the Introduction this will imply that S[50] will share an edge with S[51] and this will define the bounds of the black region.

**Figure 1** The shared edge of S[50] and S[51] creates a dynamical barrier in the early web. The S[k] of N evolve as k′ =N/2-k so S[50] has k′ = 51 and S[50] and S[51] are dual.

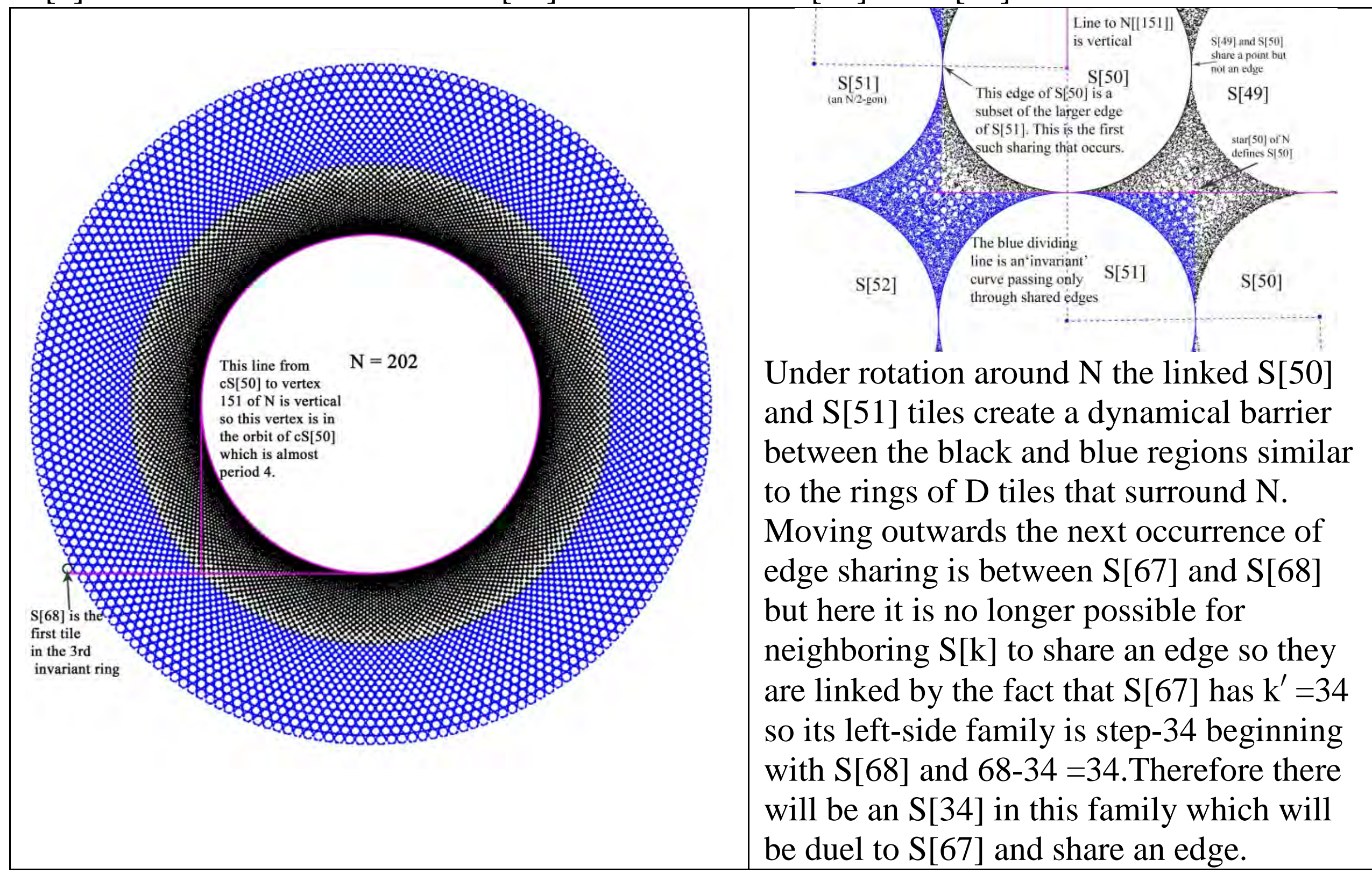

Under rotation around N the linked S[50] and S[51] tiles create a dynamical barrier between the black and blue regions similar to the rings of D tiles that surround N. Moving outwards the next occurrence of edge sharing is between S[67] and S[68] but here it is no longer possible for neighboring S[k] to share an edge so they are linked by the fact that S[67] has k′ =34 so its left-side family is step-34 beginning with S[68] and 68-34 =34.Therefore there will be an S[34] in this family which will be duel to S[67] and share an edge.

Appendix II of [H5] is devoted to the large-scale geometry of N which is dominated by concentric rings of maximal 'D tiles. The first such ring always guarantees that the region local to N is invariant, but invariance occurs at all scales and there will typically be well-defined 'secondary' invariant regions such as the black and blue regions above. For these rings of D tiles the S[1] tiles of each D will exist and share an edge with D, so these D tiles will have a strong link consisting of both vertex sharing and edge sharing. As noted in the Introduction, this sharing has its origin with S[1] and N. Because D is formed with web steps k′ = N/2-(N/2-1) = 1, it can share an edge with its dual S[1].

For S[50] and S[51] it is convenient when your 'hugging partner' lives next door but the twice-odd family is the only one where this is possible and it cannot be repeated. The next best is finding a pair like S[67] and S[68] which share through a smaller intermediary. Further out the S[k] will be too far apart and they can only link up via their shared local families.

In terms of edge geometry the S[2] tile of N is the usually the natural choice of 'next-generation' N because S[1] is an N/2-gon . The local web shown here is the joint S[1]-S[2] web but both of these tiles have their own locally invariant webs which are strongly influenced by the S[4] and S[3] tiles below.

**Figure 2** The second generation presided over by S[2] and S[1] acting as D[1] and M[1]

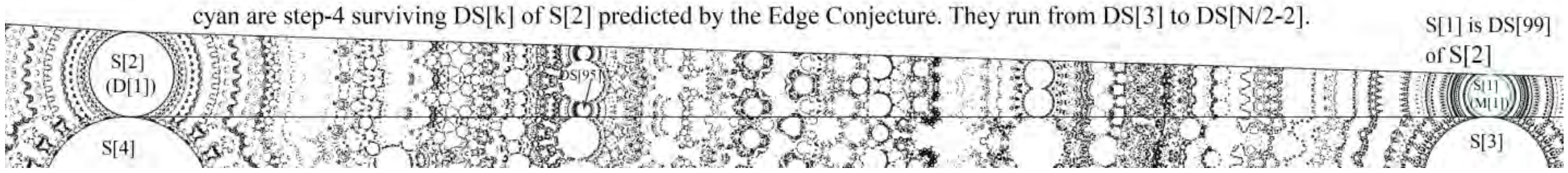

The Edge Conjecture of [H5] predicts that the surviving First Family tiles of S[2] will include at least the step-4 S[k] of S[2] which we call DS[k] to avoid confusion. Here in the 8k+2 family, the smallest survivor will be DS[3] which is simply the S[3] of S[2]. This tile will play a very important role as web 'influencer'.

As shown below this DS[3] tile will in turn foster S[1] and S[2] tiles of S[2] and these will be preserved at step 8 rotations around S[2]. The 8k+2 Conjecture predicts that these DS[k] will exist at all subsequent generations and the periods of the DS[k] will satisfy known difference equations. We call the generations D[k] with N as D[0] and S[2] as D[1].

Resolving the 3rd generation here is a challenge because the scale is small. For N twice-odd the generations scale by GenScale[N/2] = Tan[Pi/N]·Tan[Pi/2N] which here is about 0.0004839. This is the scale between S[2] and N so it is the height of S[2] a.k.a. D[1]. Then recursively the S[2] tile of S[2] is D[2] which has height $\text{GenScale}^2 \approx 2.34209 \times 10^{-7}$. This geometric scaling applies to each new generation and here we also know the matching temporal scaling.

**Figure 3** The 3rd generation of N on the edges of S[2] acting as D[1]. As in the plot above, the cyan tiles are surviving members of the First Family of S[2] so they are remnants of the 2nd generation. As noted above, the Edge Conjecture says that for the 8k+2 family these survivors will be step-4 starting with the S[3] of S[2] which we call DS[3]. On the right below are the next 2 predicted survivors: DS[7] and DS[11]. The last one is DS[N/2-2] which is S[1] (a.k.a. M[1])

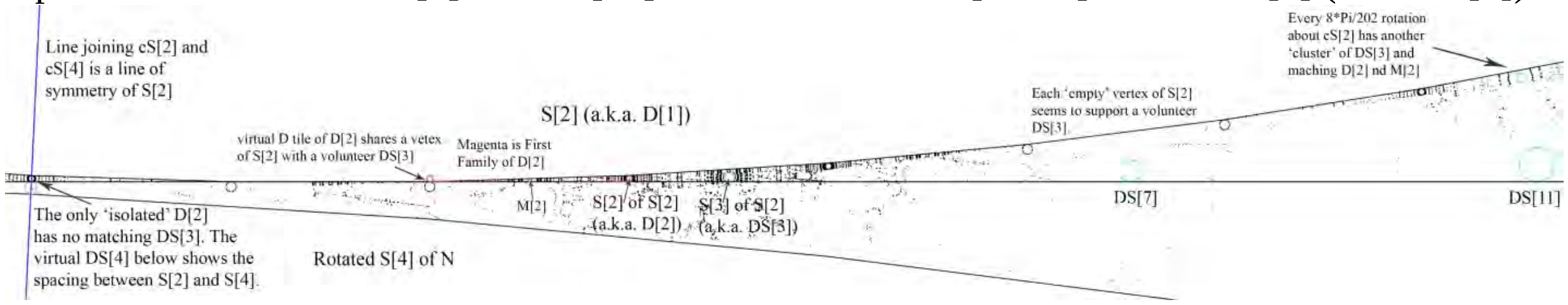

In the web evolution, this DS[3] tile will have extended edges which support 'volunteer' D[2] and M[2] tiles symmetrically on both sides and the 8k+2 Conjecture states that there will be matching 'clusters' at $8\pi/N$ (ccw) rotations about S[2] so here there will be 25 such clusters and one lone DS[2] is shown on the left above. Recursively this implies the period $D_k$ of cD[k] will satisfy $D_k = nD_{k-1} + (n+1)D_{k-2}$ where n = N/2, so $D_{[k]} = -(n(-1)^k - (n+1)^k)/(n+2)$ and Table[D[k],{k,1,6}] ={101, 10201, 1040603, 106141405, 10826423411, 1104295187821}

Finding the cD[k] is easy because their heights scale by GenScale and cD[7] has a period so high that it may as well be non-periodic, but its orbit will be very local to S[2] and will be part of a locally invariant region on the edges of D[1]. S[1] has its own locally invariant web.
To illuminate this web structure we often appeal to the complex-valued Dc map, but it has dynamics very different from $\tau$, so it can be useful sometimes to go back to first principles and look at the symbolic dynamics of $\tau$-orbits. Every orbit defines a unique sequence of vertices of N visited and the 'first derivative' of this sequence will yield the 'step-sequence' which records the number of N-vertices skipped.

Using the Mathematica First Family notebook at dynamicsofpolygons define N by setting npoints = 202 and then Evaluate the notebook. All the First Family S[k] and their parameters are then available and Vp[cS[50],5] will generate the first 5 'step sequences' in the $\tau$-orbit of cS[50] which are simply {50,50,50,50,50}. Here we want to be more literal and get the actual vertices visited in the orbit: IND[cS[50],5] will use the Vp routine to find the first 5 'indices' which are the actual vertices visited: {151, 201,49, 99, 149}. Recall from Figure 202.1 that cS[50] maps vertically to N[[151]] so this is the first vertex in the orbit.

**Figure 4** The first 4 points in the orbit of cS[50] are V = V[cS[5],4] as shown here in blue and compared with the $\tau^2$ orbit in black. The right-side is the full orbit which is
**Graphics[{Black, Line[PIM[cS[50],50,1], Magenta, poly[N],Blue, Arrow[V]}]**

| Blue $\tau$-orbit of cS[50] and black $\tau^2$ orbit | The period 50 $\tau^2$ orbit of cS[50] |
|---|---|

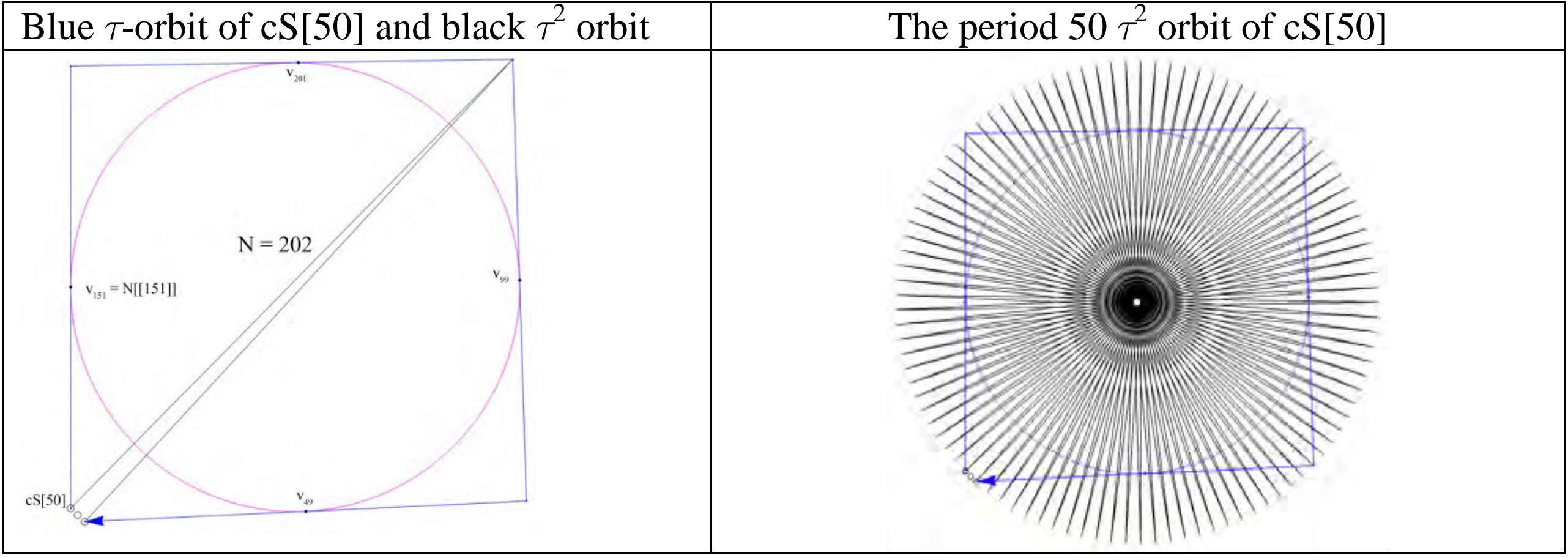

The black orbit on the left is an 'accelerated' version using the 'return map' $\tau^2(p) = p + 2(v_k - v_j)$ where $v_j$ is the initial vertex and $v_k$ is the return vertex. Here we set p1 = cS[50] then the top corner is p2 = p1 + 2($v_{201}$-$v_{151}$) and the return is p3 = p2 + 2($v_{99}$–$v_{49}$), so p3 = p1 + 2Q where Q = {$v_{201}$ –$v_{151}$ + $v_{99}$ – $v_{49}$}. All return orbits will have an alternating sum of vertices in Q.

This is easy to implement in Mathematica. For a given point p the routine PIM[p,12,1] will take the Ind indices generated by p in pairs and calculate the first 12 points in the return orbit . Here PIM[cS[50],3,1] = {$p_1$,$p_2$,$p_3$} and this can be plotted as points or lines. For a periodic orbit the return orbit will also be periodic but the indices have to be twice the period, so the plot on the right above needs IND[cS[50],100]

**Projections.** The ‘1’ in the PIM command seems to serve no purpose but in fact any of the first 100 integers would work here and the orbits would typically be quite different. P3 = PIM[cS[50, 100, 3] would use the same indices of cS[50] but the matching vertices would be those of a ‘mod-3’ version of by simply rotating each vertex to the right by 3. It should be clear that there are only $\phi$(N)/2 independent ways to remap the vertices of N in this fashion because this is the rank of the vector space determined by these vertices. For N = 202 this gives odds up to 99.

The orbits generated by these re-mappings are what Richard Schwartz calls ‘projections’ or ‘arithmetic graphs’. Since the vertices of N are ‘algebraic numbers’ these projections are exact and can be analyzed algebraically. In S1(2007] Richard used the return map to help prove that a certain ‘kite’ had an unbounded orbit and later in his study of the octagon in [S2] he used projections to aid in the analysis of orbits. This new tool has given us a new perspective on the geometry of complex obits. Here we will use them to track orbits in the invariant regions surrounding N. Since the projections are based on the return map they will preserve symmetry and periods. Here is a simple example with cS[50].

**Figure 5** On the right is the period 50 P3 projection of cS[50] in green and the normal P1 return orbit in black. On the left is the first iteration comparing the black return map $\tau^2$ and the green P3 projection which involves first remapping the vertices of N mod 3. **Graphics[{Green, Line[PIM[cS[50],50,3], Black, Line[PIM[cS[50],50,1], Magenta, poly[N]}]**

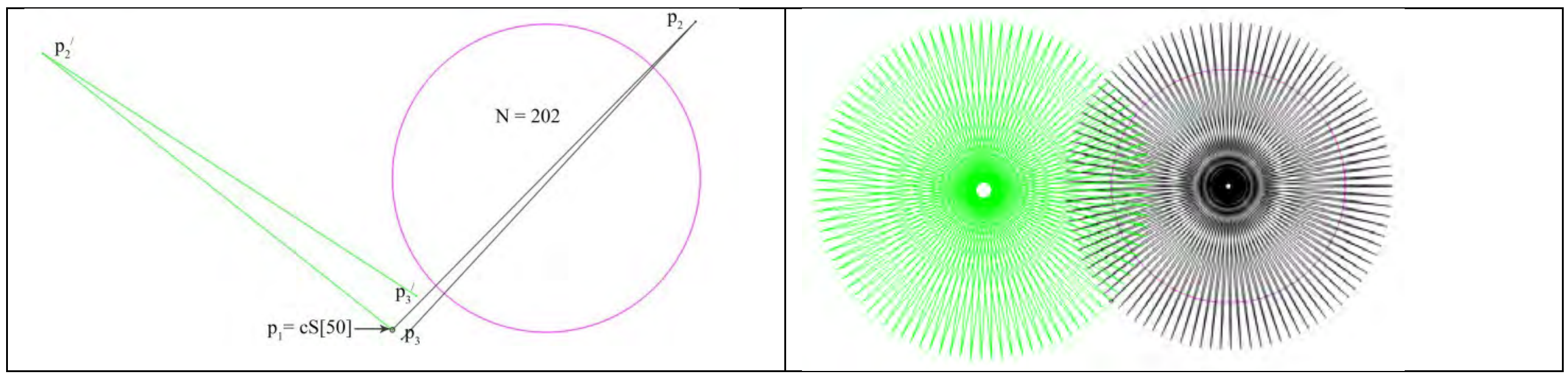

The ‘3’ in the PIM routine generates the mod-3 versions of N called Wc[[3]]. Set M3 = Wc[[3]] to be the new N and now $v_{151}$ is M3[[151]] inside PIM , so $p_2'$ above is cS[50] + 2($v_{201}$-$v_{151}$) just like $p_2$ but these two vertices are now evaluated with respect to M3 – which looks exactly like N.

Every projection of a point p is based on the same indices but there are usually some that provide a better perspective of the orbit. Here with N = 7 we might prefer P3. Using p = cS[1] the orbit is period 7 so **IND[p,14**] ={6,7,1,2,3,4,5,6,7,1,2,3,4,5} (vertex 1 at top)

**Figure 6** The 3 Projections for the period 7 orbit of S[1] for N = 7

| P1 = PIM[p,7,1] | P2 = PIM[p,7, 2] | P3 = PIM[p,7,3] | The actual orbit |
| --- | --- | --- | --- |

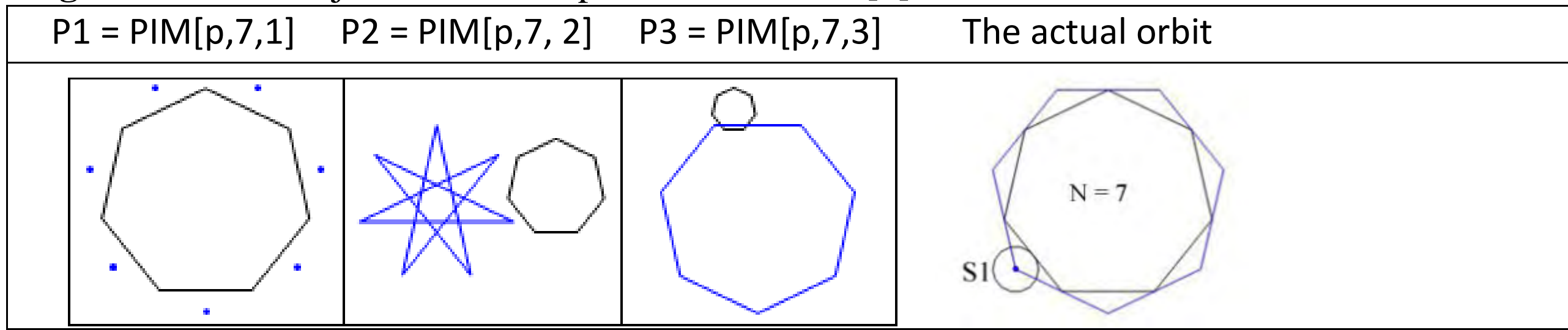

Projections can uncover structure and self-similarity in an orbit . For N = 8 there is a sequence of D[k] tiles converging to star[3] of N. The P2 projection reduces N to a rotated square as shown by the 'sponge' plot. The P3 projection is a self-similar 'Koch' snowflake. For the next cD[5], this image will be just one 'bud' of the larger snowflake. See N= 8 Summary.

**Figure 7** The P1,P2 and P3 projections for N = 8 based on cD[4] with period 8640. By convention the P1 projection is plotted as points instead of lines.

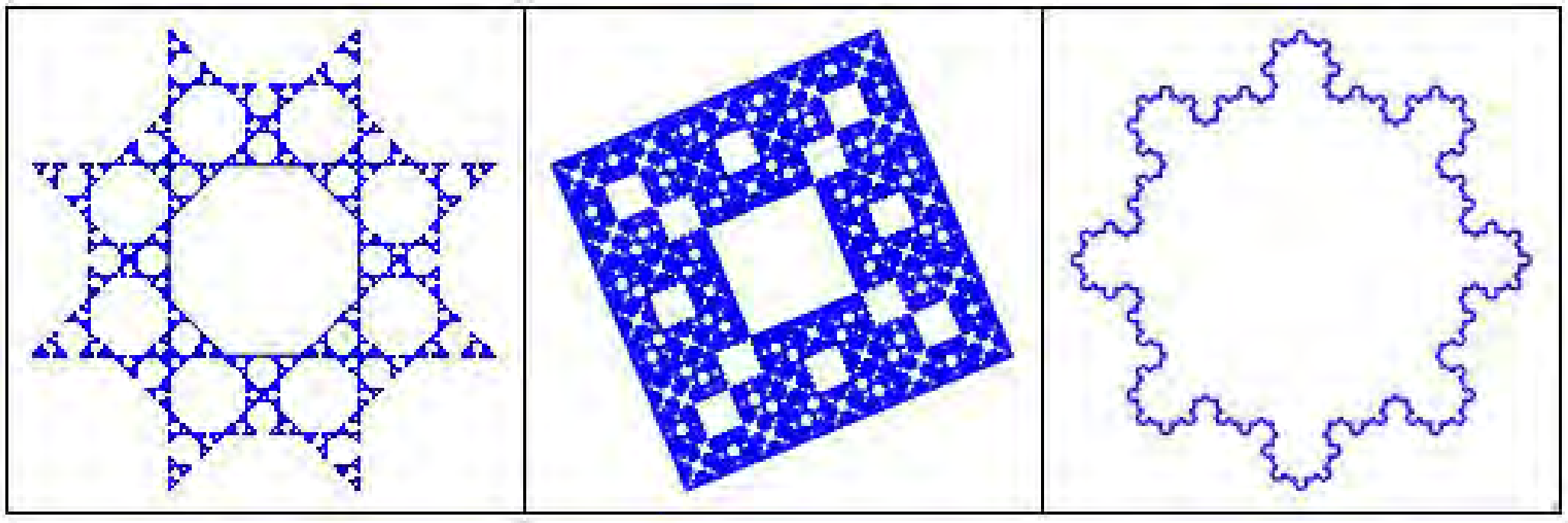

For N = 9 below the P3 projection reduces to N = 3 and this 'resonance' makes for interesting geometry. In the web S[3] will be mutated

**Figure 7** P1,P2 ,P3 and P4 projections for N = 9 based on M[2] with period 648

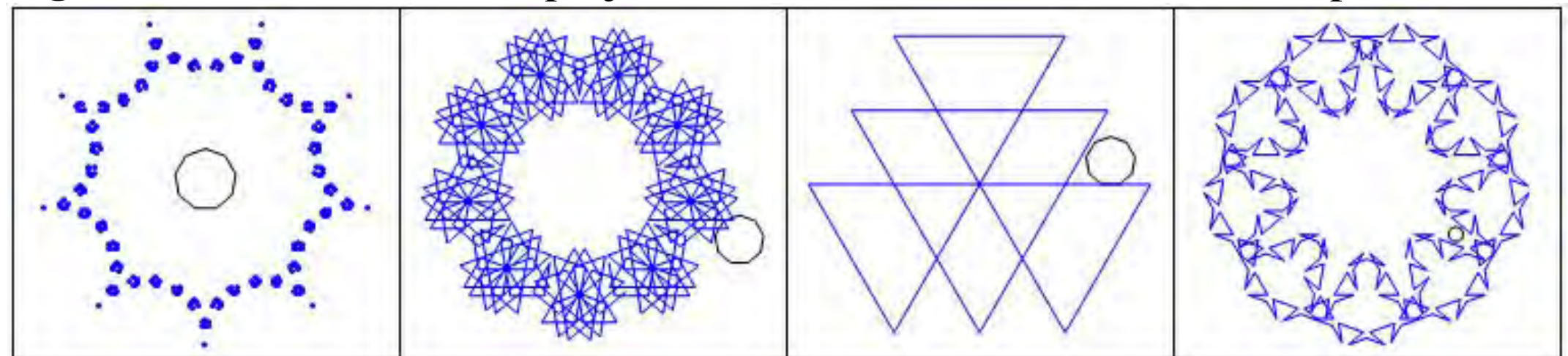

N = 26 below is in he 8k+2 family so there is a sequence of D[k] (and M[k]) tiles converging to star[1] of S[2] acting as D[1]. The ratio of periods approaches N/2+1 so P7 is the prefect projection where each pair of vertices of N fosters a new generation. The first three projections are plotted here, showing self-similarity. The S[6] and S[7] tiles lie on the boundary of the inner invariant region.The P5 projection here is based on a slight displacement of vertex 6 of S[7].

**Figure 8** Projections for N = 26 in the 8k+2 family

| The star[1] region of S[2] | P7 projections of cD[2],cD[3] & cD[4] | P5 projection of a point adjecent to S[7] |
|---|---|---|
| S[2] (D[1])<br>DS[2] (D[2])<br>star[1] of S[2] | These are P7 projections of cD[2], cD[3] and cD[4] for N = 26. The periods are 69(magenta), 2379(black) and 33293(blue). The scale is very large and N = 26 is just a tiny dot at the origin.<br>N = 26 | N = 26. 1000K lines in the P5 projection of p ... This point has period that is at least 500 million.<br>The initial interlude shows a semi-periodic beginning which may reflect the proximity of p in a neighbor which has period about 650K. |

**Figure 9** N = 32 has a strong resonance at P8 and secondary resonance at P12. Below are projections using p1 = S[1][[5]] +{.00000001,0} and p10 = S[10][[2]] +{.00000001,0}. The p10 point is just inside the initial invariant region surrounding N.

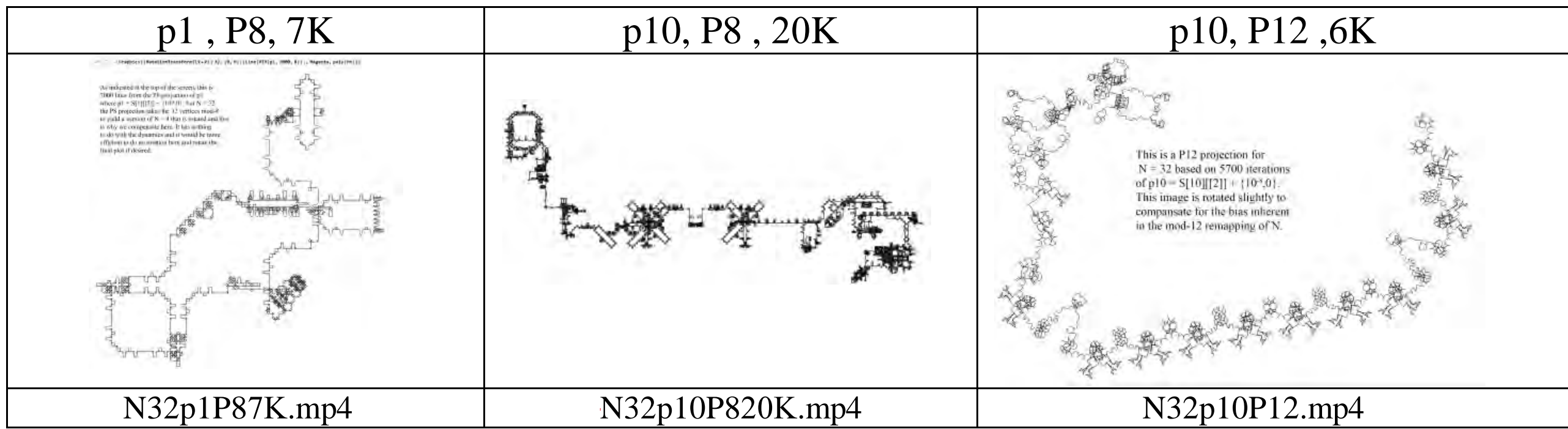

| p1 , P8, 7K | p10, P8 , 20K | p10, P12 ,6K |
| --- | --- | --- |
| N32p1P87K.mp4 | N32p10P820K.mp4 | N32p10P12.mp4 |

N = 202 is in the 8k+2 family so the growth rate of periods of D[k] tiles is N/2+1 =102. This implies that the P50 or P[51] projections may be able to track the natural self-similarity of D[k] orbits in the same fashion as P7 for N = 26. This growth rate makes it hard to plot 'generations' of D[k] . The first P50 plot below has D[1] (period 101) in blue and a portion of D[2] (period 10,202) in magenta and D[3] (period 1,040,603).in black. This black orbit will fold over and return showing that the entire D[2] orbit is just one half-circle of the D[3] orbit. But the D[4] projection at period 106,141,405 can no longer track the self-similarity which is only approximated by the P50 remapping. The third plot below shows that the P8 projection of cS[7] also diverges from the expected regular behavior after only about 40 K iterations. Since S[50] is the 'shepherd' of the inner region, the p49 = star[49] + {$10^{-8}$, $10^{-8}$} point can serve as a representative of the dynamics of the inner region.

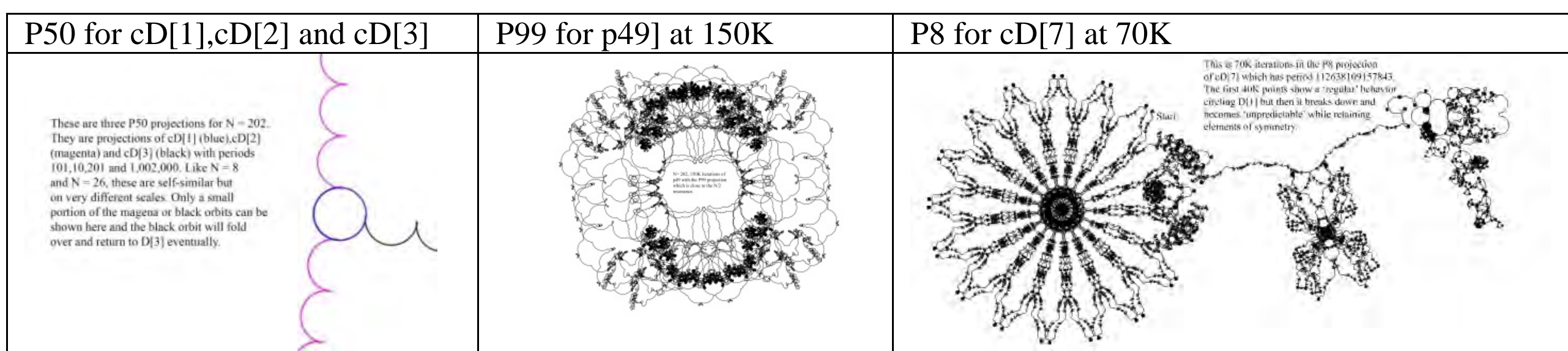

| P50 for cD[1],cD[2] and cD[3] | P99 for p49] at 150K | P8 for cD[7] at 70K |
| --- | --- | --- |

Moving outwards from N, orbits will have additional 'noise' due to extra rotations around N and some projections can act as 'filters'. Below is 500K of a P7 projection in the 3rd invariant region.

**Figure 11** N = 202 , P7 projection for S[68[[[30]]+ {$10^{-8}$,0}

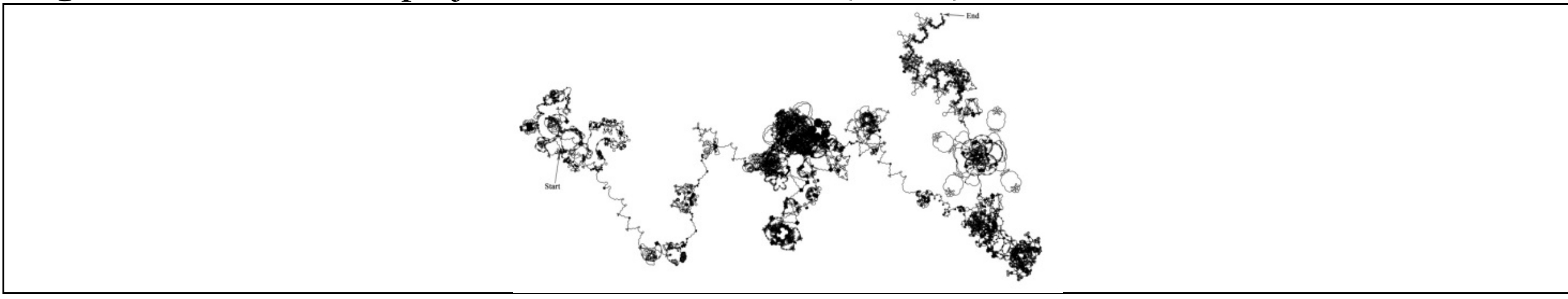

Be sure to check out the animations of projections at dynamicsofpolygons They are guaranteed to amuse your cat and also cure insomnia.

**Links**

The author's web site at DynamicsOfPolygons.org is devoted to the outer billiards map and related maps from the perspective of a non-professional. This is a very safe site which has been on-line for more than 12 years. It has full SST encryption and no commercial content. The main menu has links for accessing PDFs, Animations, Software, Images, and further Links. Under Software there is a Mathematica notebook called FirstFamily.nb which can be downloaded. It will generate the First Family for any regular N-gon and allow explorations using $\tau$, Df or Dc.The notebook called NonRegular.nb will work for any N-gon and the InnerBilliards.nb notebook introduces the case of orbits inside an N-gon.

For someone willing to download the free Mathematica CDF Reader there are many 'manipulates' that are available at the Wolfram Demonstrations site - including an outer billiards manipulate by the author and two other manipulates based on the author's results in [H2]. At the DynamicsOfPolygons site there are more cdf manipulates and there are also animations of orbits in the form of Projections. Using Adobe Reader or Adobe Acrobat to read this file, the image downloads can be opened as a temporary 'tab' in the current window for convenience but it may be necessary to go to Edit then Preferences and check the box that allows this.